\documentclass[11pt,reqno, twoside]{amsart}
\usepackage{amscd}
\usepackage{amsfonts}
\usepackage{amsmath}
\usepackage{amssymb}
\usepackage{amsthm}
\usepackage{fancyhdr}
\usepackage{latexsym}
\usepackage[colorlinks=true, pdfstartview=FitV, linkcolor=blue, citecolor=blue, urlcolor=blue]{hyperref}
\usepackage{enumitem}      
\usepackage{mathtools}            
\usepackage{upgreek}
\usepackage{indentfirst} 
\usepackage{blindtext, rotating}   
\usepackage{soul} 
\usepackage{ulem} 
\setstcolor{red}
\usepackage{color}
\usepackage{tikz}
\newcommand{\eee}[1]{\begin{equation}#1 \end{equation}}
\newcommand{\sss}[1]{\begin{subequations}#1\end{subequations}}
\newcommand{\ddd}[1]{\begin{alignat}{2}#1\end{alignat}}
\newcommand{\nn}{\nonumber}
\newcommand{\p}{\partial}

\renewcommand{\b}{\mathcolor{blue}}

\newcommand{\ve}{\varepsilon}
\def\re{{\textnormal{re}}}
\def\im{{\textnormal{im}}}
\def\err{{\textnormal{err}}}
\def\asymp{{\textnormal{asymp}}}
\def\dom{{\textnormal{dom}}}
\def\Vo{v_o}
\def\Vs{v_s}
\def\Vd{\widetilde v_s}
\def\Vw{v_w}
\def\qs{q_s}
\def\qd{q_p}
\def\Re{\mathop{\rm Re}}
\def\Im{\mathop{\rm Im}}
\makeatletter
\def\overl@ss#1#2{\vcenter{\offinterlineskip
        \ialign{$\m@th#1\hfil##\hfil$\crcr#2\crcr<\crcr } }}
\def\gl{\mathrel{\mathpalette\overl@ss>}}
\makeatother
\def\@#1{{\mathbf{#1}}}
\usepackage{scalerel,stackengine}
\stackMath
\newcommand\wwhat[1]{%
\savestack{\tmpbox}{\stretchto{%
  \scaleto{%
    \scalerel*[\widthof{\ensuremath{#1}}]{\kern-.6pt\bigwedge\kern-.6pt}%
    {\rule[-\textheight/2]{1ex}{\textheight}}
  }{\textheight}%
}{0.5ex}}%
\stackon[1pt]{#1}{\tmpbox}%
}
\newcommand\wc[1]{%
\savestack{\tmpbox}{\stretchto{%
  \scaleto{%
    \scalerel*[\widthof{\ensuremath{#1}}]{\kern-.6pt\bigwedge\kern-.6pt}%
    {\rule[-\textheight/2]{1ex}{\textheight}}
  }{\textheight}%
}{0.5ex}}%
\stackon[1pt]{#1}{\scalebox{-1}{\tmpbox}}%
}
\makeatletter
\renewcommand\subsection{\@startsection{subsection}{2}%
  \z@{-0.8\linespacing\@plus-0.7\linespacing}{0.7\linespacing}%
  {\normalfont\bfseries}}
\renewcommand\subsubsection{\@startsection{subsubsection}{3}%
  \z@{-0.8\linespacing\@plus-0.7\linespacing}{0.7\linespacing}%
  {\normalfont\itshape}}
\makeatother
\makeatletter
\def\mathcolor#1#{\@mathcolor{#1}}
\def\@mathcolor#1#2#3{%
\protect\leavevmode
\begingroup
\color#1{#2}#3%
\endgroup
}
\makeatother
\theoremstyle{plain}  
\newtheorem{theorem}{Theorem}[]

\theoremstyle{definition}

\newtheorem{remark}{Remark}[section]
\newenvironment{Proof}[1][\proofname]
{\proof[\textnormal{\textbf{#1.}}]}{\endproof}
\newcommand{\bp}{\begin{Proof}}
\newcommand{\ep}{\end{Proof}}

\DeclareMathSizes{12}{12}{7}{5}
\numberwithin{figure}{section}
\numberwithin{equation}{section}
\usepackage{geometry}
\allowdisplaybreaks
\numberwithin{theorem}{section}

\begin{document}
\title{Long-time asymptotics for the focusing nonlinear Schr\"odinger equation with nonzero boundary conditions in the presence of a discrete spectrum}
\author{Gino Biondini, Sitai Li \& Dionyssios Mantzavinos$^*$}
\date{September 18, 2019. \textit{Revised}: April 27, 2020. \mbox{}$^*$\!\textit{Corresponding author}: mantzavinos@ku.edu}
\begin{abstract}
The long-time asymptotic behavior of solutions to the focusing nonlinear Schr\"odinger (NLS) equation on the line with symmetric, nonzero boundary conditions at infinity is studied in the case of initial conditions that allow for the presence of discrete spectrum.
The results of the analysis provide the first rigorous characterization of the nonlinear interactions between solitons and the coherent oscillating structures produced by localized perturbations in a modulationally unstable medium.
The study makes crucial use of the inverse scattering transform  for the focusing NLS equation with nonzero boundary conditions, as well as of the nonlinear steepest descent method of Deift and Zhou for oscillatory Riemann-Hilbert problems.
Previously, it was shown that in the absence of discrete spectrum the $xt$-plane decomposes asymptotically in time into two types of regions: a left far-field region and a right far-field region, where to leading order the solution equals the condition at infinity up to a phase shift, and a central region where the asymptotic behavior is described by slowly modulated periodic oscillations.
Here, it is shown that in the presence of a conjugate pair of discrete eigenvalues in the spectrum a similar coherent oscillatory structure emerges but, in addition, three different interaction outcomes can arise depending on the precise location of the eigenvalues: 
(i) soliton transmission, 
(ii) soliton trapping, and 
(iii) a mixed regime in which the soliton transmission or trapping is accompanied by the formation of an additional, nondispersive localized structure akin to a soliton-generated wake.
The soliton-induced position and phase shifts of the oscillatory structure are computed, and the analytical results are validated by a set of accurate numerical simulations.
\end{abstract}


\maketitle

\markboth
{Long-time asymptotics for focusing NLS with nonzero boundary conditions and discrete spectrum}
{G. Biondini, S. Li \& D. Mantzavinos }

%
%

\section{Introduction}

In this work, we characterize the long-time asymptotic behavior of solutions to the focusing nonlinear Schr\"odinger (NLS) equation 
formulated on the line with symmetric, nonzero boundary conditions at infinity
and initial conditions that allow for the presence of discrete spectrum. 
Specifically, we consider the initial value problem (IVP) 
\sss{\label{ds-fnls-ivp-u}
\ddd{
&iq_t + q_{xx} + 2 |q|^2 q=0,  &&x\in\mathbb R, \ t>0,
\label{ds-fnls-ivp-u-eq}
\\
&q(x,0) = {f}(x),  && x\in \mathbb R,
\label{ds-fnls-ivp-ic-u}
\\[-0.2em]
&\kern-0em{\lim_{x\to\pm\infty} q(x,t) = {} } q_\pm e^{2iq_o^2 t}, \quad && t\geqslant 0,
\label{ds-fnls-ivp-bc-u}
}
}
where $q_\pm$ are complex constants such that
\eee{
|q_\pm|=q_o>0,
}
and the initial datum  $f(x)$ generates a conjugate pair of discrete eigenvalues in the spectrum (as discussed in detail in Section~\ref{ds-rhp-s}). 
The nonzero boundary conditions \eqref{ds-fnls-ivp-bc-u} are referred to as \textit{symmetric} and imply that the initial datum also tends to nonzero values at infinity:
${\lim_{x\to\pm\infty} f(x)} = q_\pm$. 
In particular, throughout this work we assume that  
\eee{\label{ds-fnls-ic-space}
e^{\pm q_o x} \left({f}-q_\pm\right) \in L^1(\mathbb R_\pm)
}
with $L^1(\mathbb R_\pm)$ denoting the spaces of Lebesgue integrable functions over $\mathbb R_\pm$. 
This is a standard assumption when the long-time asymptotic analysis is performed via inverse scattering transform techniques. Well-posedness results for IVP \eqref{ds-fnls-ivp-u} with rough initial data  are available via harmonic analysis techniques, e.g. see the recent work \cite{m2017} by Mu\~noz where local well-posedness is shown in Sobolev spaces $H^s$ with $s>\frac 12$.

The boundary conditions \eqref{ds-fnls-ivp-bc-u} motivate the transformation
\eee{
{q(x, t) \longmapsto q(x, t) e^{2iq_o^2t},}
}
which turns IVP \eqref{ds-fnls-ivp-u} into the  convenient form
\sss{\label{ds-fnls-ivp}
\ddd{
&i{q}_t + {q}_{xx} + 2\left(|{q}|^2-q_o^2\right){q}=0,  \quad  &&x\in\mathbb R, \ t>0,
\label{ds-fnls-ivp-fnls-eq}
\\
&{q}(x, 0) = {f}(x),  && x\in \mathbb R,
\\
&{\kern-0.6em\lim_{x\to\pm\infty}q}(x, t) = q_\pm,  && t\geqslant 0, 
\label{ds-fnls-ivp-bc}
}
}
where, importantly, the boundary conditions at infinity are now {independent of time.}

The focusing NLS equation \eqref{ds-fnls-ivp-fnls-eq} is a prime example of a  completely integrable system \cite{zs1972,AS1981}. 
As such, it can be written in the form of the compatibility condition $X_t - T_x + [X,T] = 0$ of the Lax pair
\eee{\label{ds-fnls-lp}
\Psi_x = X\Psi, 
\quad
\Psi_t = T\Psi, 
}
where $\Psi = \Psi(x, t, k)$ is a  $2\times 2$ matrix-valued function  and
\eee{
X=ik\sigma_3+Q,
\quad
T = -2ik^2\sigma_3 + i\sigma_3\left(Q_x-Q^2-q_o^2 I \right)-2kQ
}
with  $k\in \mathbb C$ and
\eee{
\sigma_3
=
\begin{pmatrix}
1 &0
\\
0 &-1
\end{pmatrix},
\quad
Q
=
\begin{pmatrix}
0 &q
\\
-\bar q &0
\end{pmatrix}.
}
The Lax pair \eqref{ds-fnls-lp} can be used to analyze IVP \eqref{ds-fnls-ivp}  by means of the celebrated inverse scattering transform. 
For rapidly vanishing initial conditions, in which case $q_o=0$, this task was accomplished  by Zakharov and Shabat in 1972  \cite{zs1972}.  
For nonvanishing initial conditions, however, which is the case relevant to the problem considered here, only partial results  were available (e.g.,  \cite{m1979}) until the recent work by Kova\v{c}i\v{c} and the first author \cite{bk2014}. 
There, the authors were able to develop the complete inverse scattering transform formalism for IVP \eqref{ds-fnls-ivp} and, in particular, to associate its solution to that of a  matrix Riemann-Hilbert problem. The work was then extended to asymmetric and one-sided boundary conditions in \cite{JMP55p101505} and \cite{PV2015}, respectively.

The results of \cite{bk2014}  provide a starting point for the \textit{rigorous} analysis of the long-time asymptotic behavior of the solution of IVP \eqref{ds-fnls-ivp}. This task is far from trivial due to the fact that, in the case of nonvanishing initial conditions, the focusing NLS equation exhibits \textit{modulational instability} (also known as Benjamin-Feir instability \cite{bf1967}), namely, 
the instability of a constant background with respect to long-wavelength perturbations \cite{ZO2009}. 

For example, in the special case of constant initial data $f(x) = q_o$ it is straightforward to verify that problem \eqref{ds-fnls-ivp} admits the constant solution $q(x, t) = q_o$. Seeking a solution of \eqref{ds-fnls-ivp} in the form of the localized perturbation $q(x, t) =q_o \left[1 + \ve \nu(x, t)\right]$ with $\nu = O(1)$ and $\ve\ll 1$ yields to $O(\ve)$ a linear equation with \textit{zero} conditions at infinity, which can therefore be solved explicitly via Fourier transform. The associated dispersion relation is $\omega = k \sqrt{k^2-4q_o^2}$, which becomes purely imaginary for small wavenumbers (i.e. long wavelengths) characterized by $|k|<2q_o$. Hence, $\nu$ grows exponentially as $t\to\infty$, indicating instability. But, of course, the linearization becomes invalid once $\nu$ grows to $O(\ve^{-1})$. 
The question of what happens  to the solution of the focusing NLS equation beyond this point  is referred to as \textit{the nonlinear stage of modulational instability}. 

Despite interesting results concerning the behavior of solutions with periodic boundary conditions \cite{akhmedievkorneev, forestlee, trillowabnitz}, the nonlinear stage of modulational instability for the focusing NLS equation on the infinite line remained essentially open for more than fifty years. Recently, it was conjectured in \cite{zakharovgelash,gelash} that the nonlinear stage of modulational instability is governed by the formation of certain
breather pairs termed ``super-regular solitons''. However, this conjecture was disproved in \cite{SIAP75p136}, where it was shown that solitons are not generically the main vehicle for the modulational instability; instead, the signature of the instability in the inverse scattering transform lies in the portion of the continuous spectrum associated with the nonlinearization of the unstable Fourier modes and manifests itself via exponentially growing jumps in the Riemann-Hilbert problem. The problem was then settled in \cite{bm2016,bm2017}.
First, the inverse scattering transform formalism of \cite{bk2014} was suitably modified to yield   a  Riemann-Hilbert problem  convenient for carrying out a long-time asymptotic analysis. 
The asymptotic behavior of the solutions of this Riemann-Hilbert problem was then studied using the Deift-Zhou nonlinear steepest descent method \cite{dz1993, dz1995} 
and borrowing ideas from \cite{boutet, bv2007, jm2013}. Eventually, it was shown  in \cite{bm2017} that the solution of IVP \eqref{ds-fnls-ivp} \textit{remains bounded at all times} and, more specifically,  at leading order it takes on the following asymptotic forms (see Figure \ref{ds-bifurc-f}): 
\vskip 3mm
\begin{enumerate}[label=(\roman*), leftmargin=10mm]
\item For $|x|>4\sqrt 2 q_o t$, the solution is described by two plane waves, one for $x<0$ and one for $x>0$, whose amplitudes are equal to the ``boundary data'' $q_-$ and $q_+$ respectively; 
\vskip 2mm
\item For $|x|<4\sqrt 2 q_o t$, the solution is described by slowly modulated periodic oscillations whose amplitude is given in terms of the well-known Jacobi elliptic  snoidal solution of focusing NLS. 
\end{enumerate}

Importantly, in both of the above regions the spatial structure of the leading-order asymptotics is independent of the initial datum $f$. That is,  within the class of initial data \eqref{ds-fnls-ic-space}, \textit{generic} localized perturbations of a constant background display the \textit{same} long-time behavior in \textit{all} modulationally unstable media governed by the focusing NLS equation on the infinite line. 
In this sense, the results of \cite{bm2017} demonstrate that the asymptotic state of the nonlinear stage of modulational instability is \textit{universal}. 
These analytical predictions were recently confirmed, and the resulting behavior was observed, in optical fiber experiments \cite{kser2019}. 
Moreover, it was shown in \cite{blmt2018} that this behavior is not limited to the focusing NLS equation, but instead it is a common feature of more general NLS-type systems. 
In this regard, we note that the focusing semilinear Schr\"odinger equation with power nonlinearity (which is not integrable besides the cubic case)  and nonzero boundary conditions at infinity with perturbations in Sobolev spaces was recently studied via harmonic analysis techniques~\cite{m2017}.

However, the analysis of \cite{bm2017} was carried out for initial data~\eqref{ds-fnls-ic-space} such that \textit{no discrete spectrum is present} in the Riemann-Hilbert problem emerging from the inverse scattering transform. This is a major assumption at the technical level (as will become evident while the analysis unfolds in the forthcoming sections) but, more importantly, a significant restriction from a physical point of view since, as is well-known, discrete spectrum is the mechanism generating \textit{solitons}. 
Hence, in the case of IVP \eqref{ds-fnls-ivp}, an empty discrete spectrum excludes the possibility of describing solutions that contain solitons.

In this work, we perform the long-time asymptotic analysis of the focusing NLS IVP \eqref{ds-fnls-ivp}  without the assumption of an  empty discrete spectrum  that was used in \cite{bm2017}. Specifically, we consider initial data $f$ satisfying \eqref{ds-fnls-ic-space} such that the analytic scattering coefficients arising in the inverse scattering transform have a single pair of conjugate simple poles in the complex spectral plane. This is clearly the simplest scenario that allows for the presence of solitons. As in the case of zero boundary conditions at infinity, 
each conjugate pair of discrete eigenvalues contributes a soliton to the solution of NLS. Hence,  in the case considered here there is exactly one soliton present.

The simultaneous presence of a discrete spectrum and a nonvanishing reflection coefficient allows one to study the interactions 
between solitons and radiation (i.e. the components of the solution of the NLS equation arising from the reflection coefficient).
In the case of zero boundary conditions at infinity,  problems of this kind were first studied in the 1970s \cite{JMP17p710, JMP17p714, JETP44p106}.
Those studies, however, employed formal methods. 
Moreover, and most importantly for our purposes, they were limited to the case of a zero background (i.e. $q_o=0$).
In the context of the focusing NLS IVP \eqref{ds-fnls-ivp}, the presence of a discrete spectrum affords us the ability to rigorously study --- for the first time --- 
the interaction between solitons and radiation on a modulationally unstable  background.

%
%

\section{Overview of Results}
\label{overview-s}

\noindent
\textbf{Definitions and notation.}
Before we can state our results precisely,  we need to introduce some notation and provide definitions of various quantities that will appear throughout this work. 
\vskip 1mm
\begin{enumerate}[label=$\bullet$, leftmargin=4mm, rightmargin=0mm]
\advance\itemsep 1mm
\item For any complex-valued function $f$, we denote $f_\re := \text{Re}(f)$ and $f_\im := \text{Im}(f)$. Complex conjugation is denoted by an overbar. 

\item The complex square root $\left(k^2+q_o^2\right)^{\frac 12}$, with $k\in\mathbb C$ being the spectral variable introduced through the Lax pair \eqref{ds-fnls-lp},  
is expressed in terms of a single-valued function $\lambda(k)$, which is uniquely defined by taking the branch cut along the segment 
\eee{
B := i[-q_o, q_o]
}
of the complex $k$-plane and defining 
\eee{\label{ds-lambda-def}
\lambda(k) 
= 
\left\{
\def\arraystretch{1.2}
\begin{array}{ll} 
\sqrt{k^2+q_o^2}, &k\in\mathbb R^+\cup B, \\ -\sqrt{k^2+q_o^2}, &k\in\mathbb R^-, 
\end{array}
\right.
}
so that $\lambda(k)\sim k$ as $k\to \infty$. 

\item 
The phase function $\theta(\xi, k)$ is defined by 
\eee{\label{ds-theta-def}
\theta(\xi, k) =  \lambda(k)  \left(\xi-2k \right),
}
where $\xi$ is the similarity variable
\eee{\label{ds-xi-def}
\xi = \frac xt
}
which, as usual, is the key independent parameter in the calculation of the long-time asymptotics. Importantly, $\theta$ is Schwarz-symmetric, i.e. $\theta(\xi,\bar k) = \overline{\theta(\xi,k)}$.
\item 
As in \cite{bm2017}, a key role in the analysis will be played by the function
\eee{\label{ds-h-abel-i}
h(\xi, k) = \frac 12 \left(\int_{iq_o}^k + \int_{-iq_o}^k \right) dh(\xi, z)
}
defined via the Abelian differential 
\eee{\label{ds-dh-abel-i}
dh(\xi, k) = -4\, \frac{\left[k-k_o(\xi)\right]\left[k-\alpha(\xi)\right]\left[k-\bar \alpha(\xi)\right]}{\gamma(\xi, k)} \, dk,
}
with $\alpha$ and $\gamma$ defined below and
\eee{\label{ds-ko-def-i}
k_o = -\alpha_\re + \frac \xi 4.
}
Note that $h$ is also Schwarz-symmetric, i.e. $h(\xi,\bar k) = \overline{h(\xi,k)}$.

\item
The complex quantity $\alpha$ and the elliptic parameter $m$ of the slowly modulated genus-1 oscillations are uniquely determined by the solution of the modulation equations \cite{egkk1993,kamchatnov}
\sss{\label{ds-mod-eqs-i}
\ddd{
&\frac \xi 2 = 2\alpha_\re+\frac{q_o^2-\alpha_\im^2}{\alpha_\re},
\quad
m^2=\frac{4q_o\alpha_\im}{\alpha_\re^2+\left(q_o+\alpha_\im\right)^2},
\\
&\left[\alpha_\re^2+\left(q_o-\alpha_\im\right)^2\right]K(m) 
= 
\left(\alpha_\re^2-\alpha_\im^2+q_o^2\right)E(m),
\label{ds-mod-eq-bm-i}
}
}
with $K(m)$, $E(m)$ being  the complete elliptic integrals of the first and second kind respectively. 
\item
The function 
\eee{\label{ds-gamma-def-i}
\gamma(\xi, k)
:=
\left[\left(k^2+q_o^2\right)\left(k-\alpha\right)\left(k-\bar \alpha\right)\right]^{\frac 12}
}
is  uniquely defined by taking branch cuts along $B$ as well as an appropriate contour $\widetilde B$ connecting the points $\alpha$, $\bar \alpha$ and $k_0$. The  dependence on the similarity variable $\xi$ will often be suppressed from the arguments of $\alpha$, $\bar \alpha$, $k_o$, $\gamma$ and other quantities for  brevity.
\item We denote by $p, \bar p$ the single pair of conjugate simple poles that form the discrete spectrum of the Riemann-Hilbert problem associated with the focusing NLS IVP \eqref{ds-fnls-ivp} (see Section \ref{ds-rhp-s} for more details). 
The location of $p$ will play a crucial role in the analysis.
Thanks to the reflection invariance of the NLS equation (i.e. the fact that if $q(x,t)$ is a solution then so is $q(-x,t)$) and the symmetry
$k\mapsto \bar k$ of the spectrum of the scattering problem
(see Section~\ref{ds-rhp-s} for details), without loss of generality we may take $p$ to lie in the third quadrant of the complex $k$-plane.
\item 
Our analysis and  the corresponding results are intimately related to the value of $\xi$ relative to the following special values:
\eee{\label{ds-xistar-xip-def}
\Vo = -4\sqrt 2 q_o, 
\quad 
\Vs = 2\left[p_\re+\frac{\lambda_\re(p)}{\lambda_\im(p)} \, p_\im \right].
} 
The velocity $\Vo$ defines the edge of the modulated elliptic wave region, whereas $\Vs$ is the unperturbed velocity of a soliton produced by a discrete eigenvalue located at $k=p$ (see Figure \ref{ds-bifurc-f}). Note that $\Vs$ is the value of $\xi$ such that
\eee{\label{vs_def}
\Im\left[\theta(\xi, p)\right]=0
}
and that $\Im\left[\theta(\xi, p)\right]=0$ if and only if $\Im \left[\theta(\Vs, \bar p)\right]=0$.

\item
Besides $\Vo$ and $\Vs$, a key role will also be played by the solutions $\Vd$ and $\Vw$ of the equation
\eee{\label{ds-xih-def}
\Im\left[h(\xi, p)\right] = 0, \quad \xi\in (\Vo, 0).
}
Note that $\Im\left[h(\xi, p)\right] = 0$ if and only if $\Im\left[h(\xi, \bar p)\right] = 0$. 
The difference between  $\Vd$ and $\Vw$ is explained below (see also Sections~\ref{ds-esc-s} and~\ref{ds-trap-s} for more details).  

\begin{figure}[t!]
\begin{center}
\includegraphics[scale=0.5]{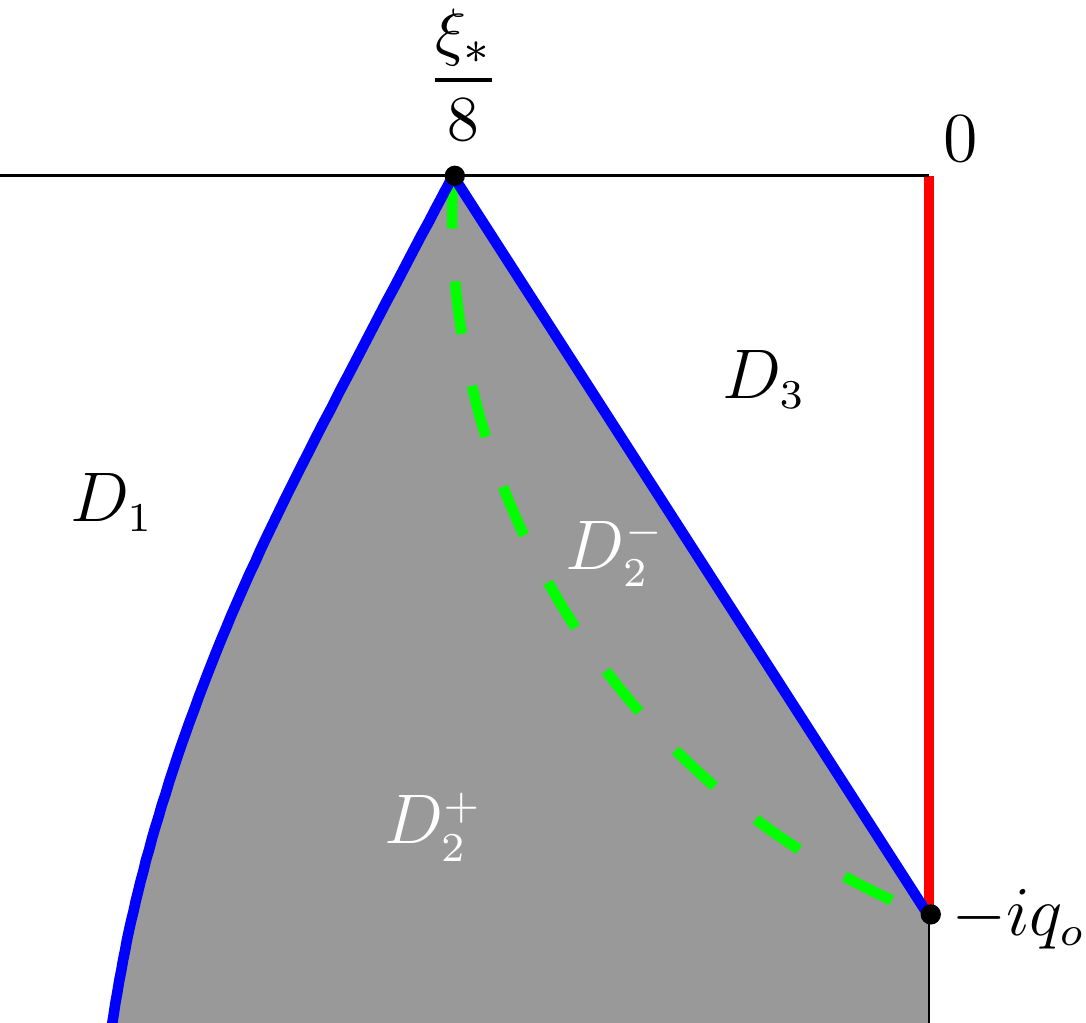}
\hspace*{0.3cm}
\includegraphics[scale=0.3]{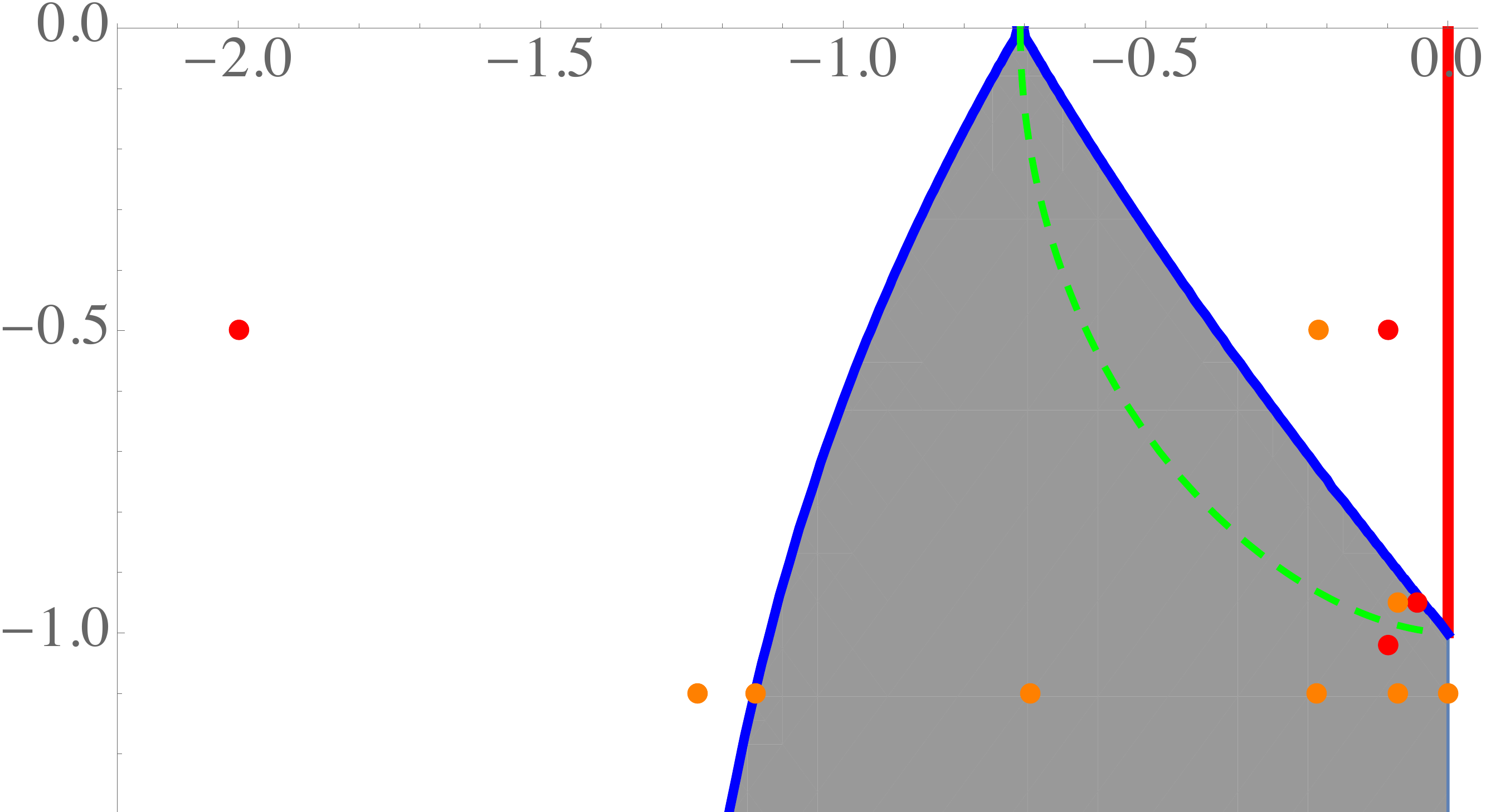}
\caption{
\textit{Left}: 
The third quadrant of the complex $k$-plane and the four regions $D_1$, $D_2^+$, $D_2^-$, $D_3$. 
Solid blue curve: $\Im[\theta(\Vo, k)]=0$; 
dashed green curve: the trace of the point $\bar \alpha$ (defined by \eqref{ds-mod-eqs-i}) 
as $\xi$ increases from $\Vo$ to $0$.
The case $p\in D_1$ corresponds to the \textit{transmission regime}, 
the case $p\in D_2^+$ to the \textit{trap regime}, 
the case $p\in D_2^-$ to the \textit{trap/wake regime}, 
and the case $p\in D_3$ to the \textit{transmission/wake regime}. 
\textit{Right}: The different choices of the pole $p$  used in the numerical simulations of Figures \ref{numerics1}, \ref{numerics2} (red dots) and \ref{h_root} (orange dots).
}
\label{ds-regions-f}
\end{center}
\end{figure}

\item
We will show in Sections~\ref{ds-esc-s} and~\ref{ds-trap-s} that
the third quadrant  $\mathbb C_\mathrm{iii}$ of the complex $k$-plane is divided into the four regions $D_1$, $D_2^+$, $D_2^-$, $D_3$  defined as follows.
Recall that, for a discrete eigenvalue at $k=p$, $\Vs$ is uniquely defined as the value of $\xi$ such that $\Im[\theta(\xi,p)]=0$. Then, $\mathbb C_\mathrm{iii}$ can be decomposed into
\sss{
\ddd{
&D_1\cup D_3 = \left\{ k\in\mathbb{C}_\mathrm{iii}: \Vs < \Vo < 0\right\},
\nn\\
&D_2 = \left\{ k\in\mathbb{C}_\mathrm{iii}: 0 > \Vs > \Vo \right\}.
\nn
}
}
These regions are shown in Figure~\ref{ds-regions-f} with $D_1\cup D_3$ in white and $D_2$  in gray.
The solid blue curve separating them corresponds to the values of $k$ for which  $\Vs=\Vo$ or, equivalently, to $\Im[\theta(\Vo, k)] = 0$. 
The dashed green curve corresponds to the trace of the point $\bar \alpha$ as $\xi$ increases from $\Vo$ to $0$. 

Note that:
\vspace*{1mm}
\begin{itemize}[label=$\circ$,leftmargin=4mm, rightmargin=0mm]
\advance\itemsep 1mm
\item 
The region where $\Vs<\Vo<0$ is divided by the blue curve $\Im[\theta(V_o,k)]=0$ into two disjoint domains.  
Among them, we take $D_1$ to be the infinite domain and $D_3$ the one adjacent to the imaginary axis.
\item
Similarly, the dashed green curve separates $D_2$ into two subdomains, $D_2^+$ and $D_2^-$, which we take as the portions of $D_2$ adjacent to $D_1$ and  $D_3$, respectively.
\item 
We will show that $D_1$ and $D_3$ differ with respect to the number of solutions of equation \eqref{ds-xih-def} that arise in the interval $(\Vo, 0)$. In particular, if $p\in D_1$ then  \eqref{ds-xih-def} does not have a solution in $(\Vo, 0)$, while if $p\in D_3$ then \eqref{ds-xih-def} possesses a unique solution $\Vw\in (\Vo, 0)$. 
\item 
Similarly,  we will show that if $p\in D_2^+$ then equation \eqref{ds-xih-def} possesses a unique solution $\Vd\in (\Vo, 0)$ while if $p\in D_2^-$ then \eqref{ds-xih-def} has two solutions $\Vd, \Vw\in (\Vo, 0)$ with $\Vd<\Vw$. 
\end{itemize}
\end{enumerate}

\vspace*{3mm}
\noindent
\textbf{The four interaction outcomes.}
Placing $p$ in each of the four regions $D_1$, $D_2^+$, $D_2^-$, $D_3$ gives rise to different,  inequivalent asymptotic regimes,  which we label as the \textit{transmission regime}, the \textit{trap regime}, the \textit{trap/wake regime}, and the \textit{transmission/wake regime} respectively. 

Specifically, in Sections~\ref{ds-esc-s} and~\ref{ds-trap-s}  we show that, depending 
on its location in the complex $k$-plane (see Figure \ref{ds-regions-f}), 
the presence of a discrete eigenvalue at $k=p$ gives rise to the following leading-order contributions 
in addition to the portion of the solution generated by the continuous spectrum:
\begin{enumerate}[label=(\roman*),leftmargin=12mm, rightmargin=0mm]
\advance\itemsep 1mm
\item 
In the transmission regime, i.e. when $p\in D_1$, a soliton along the ray $x=\Vs t$;

\item 
In the trap regime, i.e. when $p\in D_2^+$, a soliton along the ray $x=\Vd t$;

\item 
In the trap/wake regime, i.e. when $p\in D_2^-$, a soliton along  $x=\Vd t$ and a soliton wake along  $x=\Vw t$;

\item  
In the transmission/wake regime, i.e. when $p\in D_3$, a soliton along  $x=\Vs t$ and a soliton wake along  $x=\Vw t$.
\end{enumerate}
In particular, we will see that the above outcomes are determined by whether there exist solutions of equation~\eqref{vs_def} for $\xi\in(-\infty,v_o)$ and of   equation~\eqref{ds-xih-def} for $\xi\in(v_o,0)$.

\begin{figure}[t!]
\begin{center}
\includegraphics[scale=0.725]{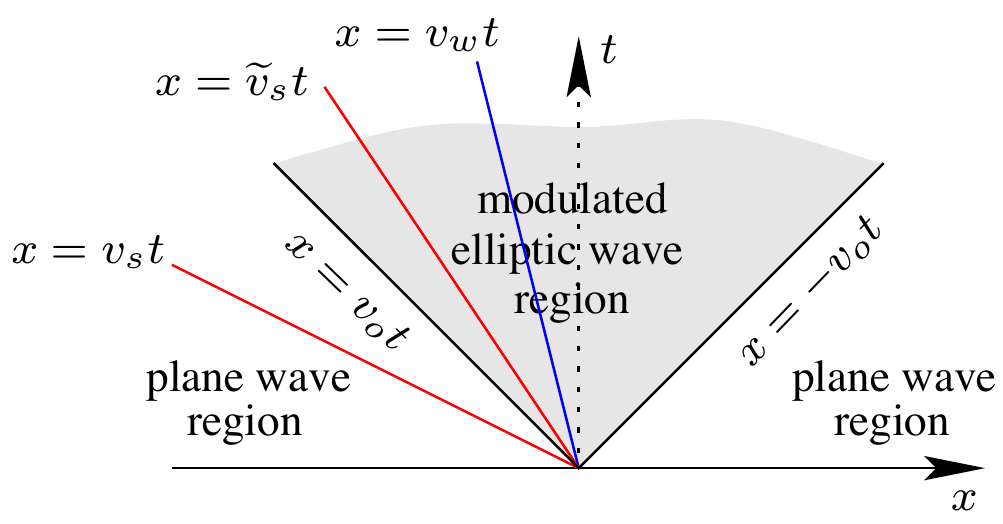}
\caption{Asymptotically in time, the $xt$-plane is divided into the plane wave regions $|x|>|\Vo|t$ and the modulated elliptic wave region $|x|<|\Vo|t$. Also shown are the locations of the $O(1)$ contributions generated by a discrete eigenvalue at $k=p$ in the four inequivalent cases corresponding to the regions of Figure \ref{ds-regions-f}.
}
\label{ds-bifurc-f}
\end{center}
\end{figure}

\begin{figure}[ht!]
\vspace*{3mm}
    \begin{center}
        \includegraphics[scale=0.25]{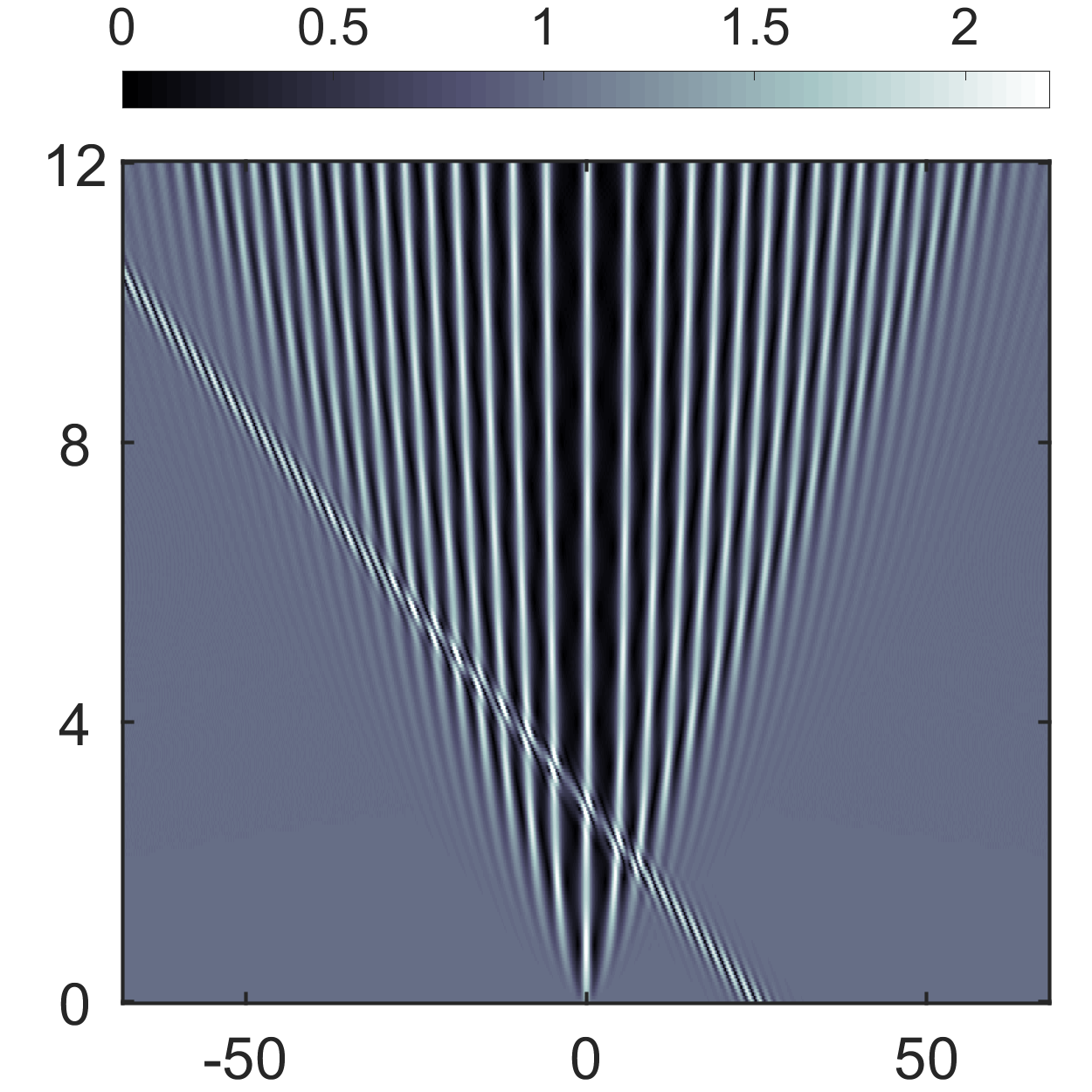}%
        \includegraphics[scale=0.25]{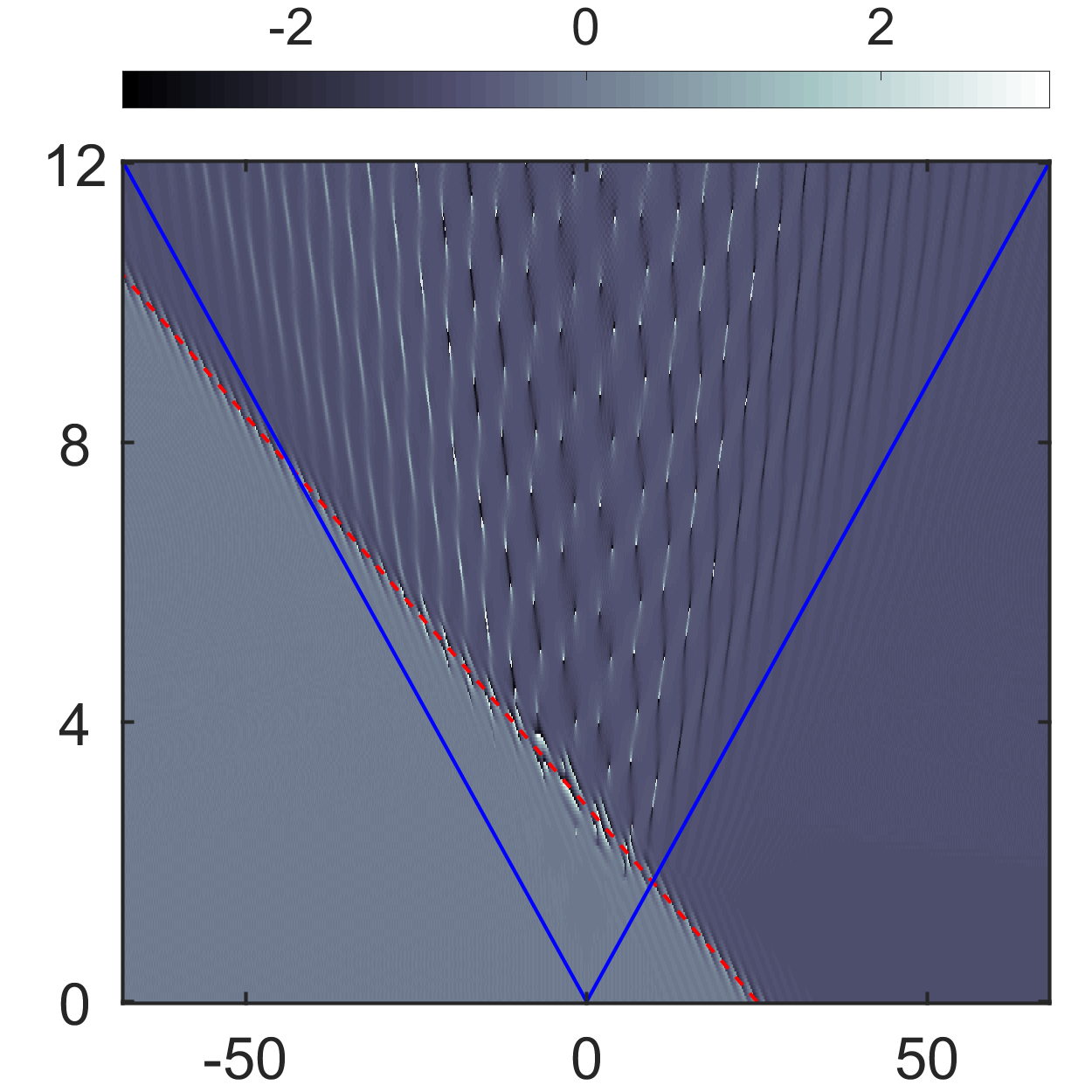}%
        \includegraphics[scale=0.25]{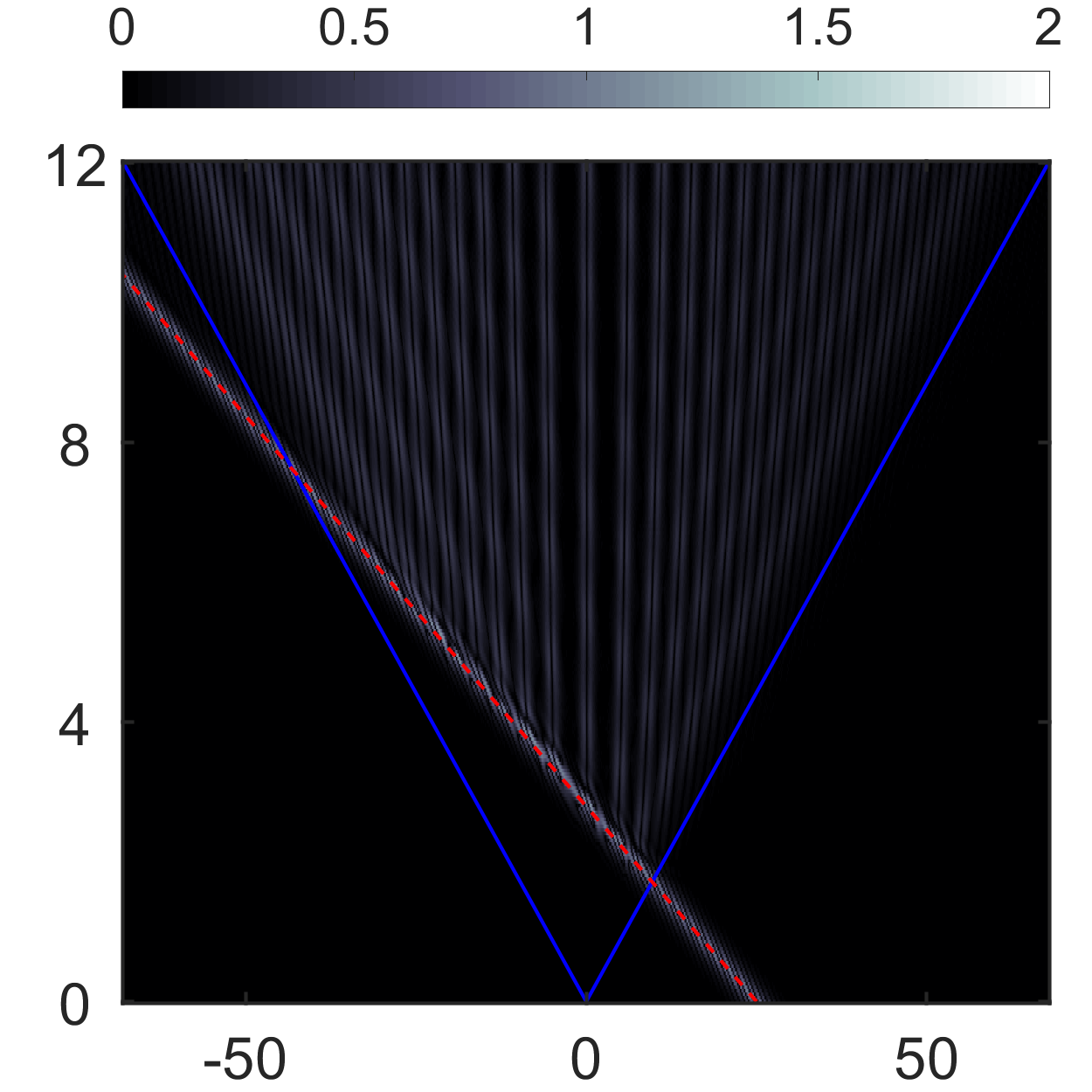}
        \\[2ex]
        \includegraphics[scale=0.25]{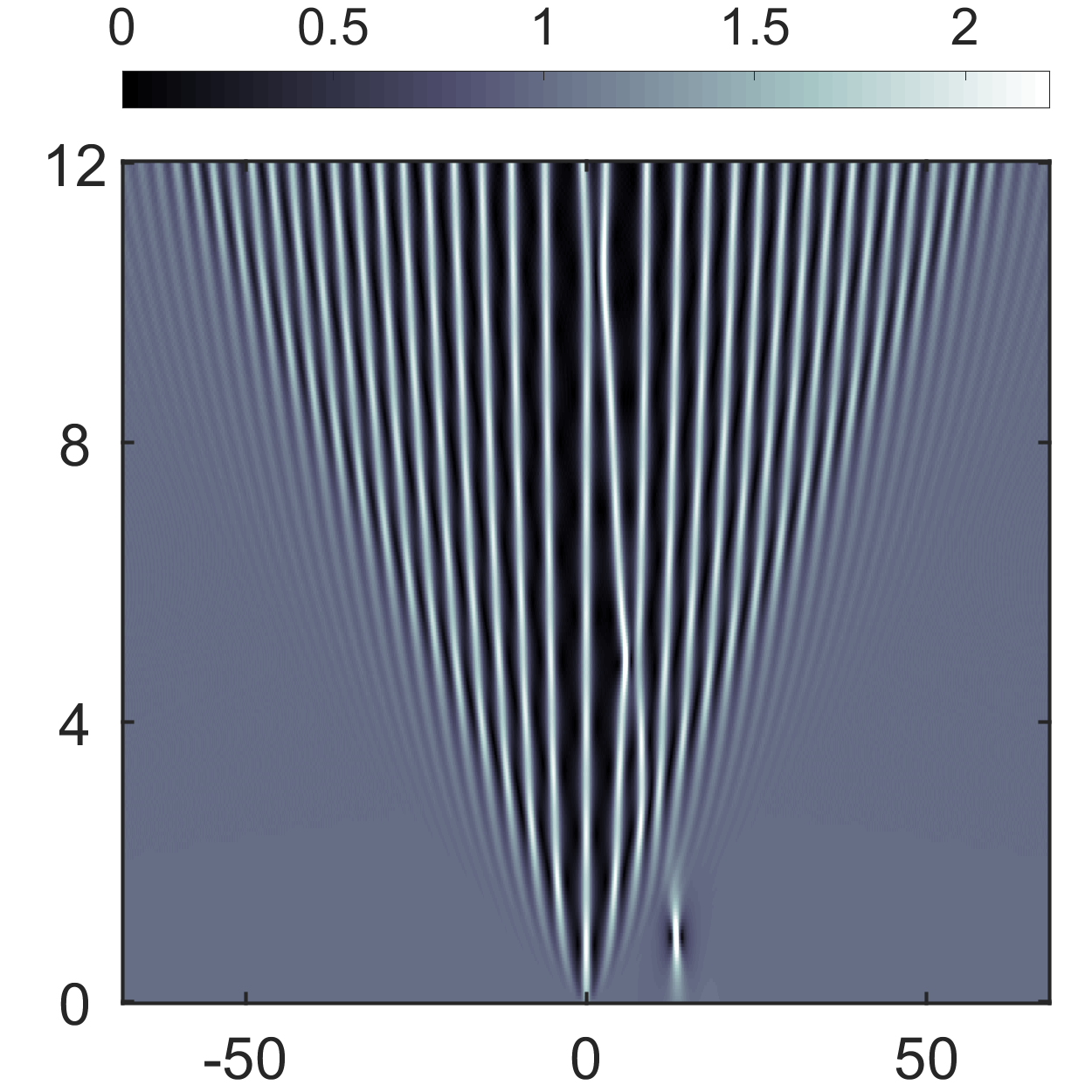}%
        \includegraphics[scale=0.25]{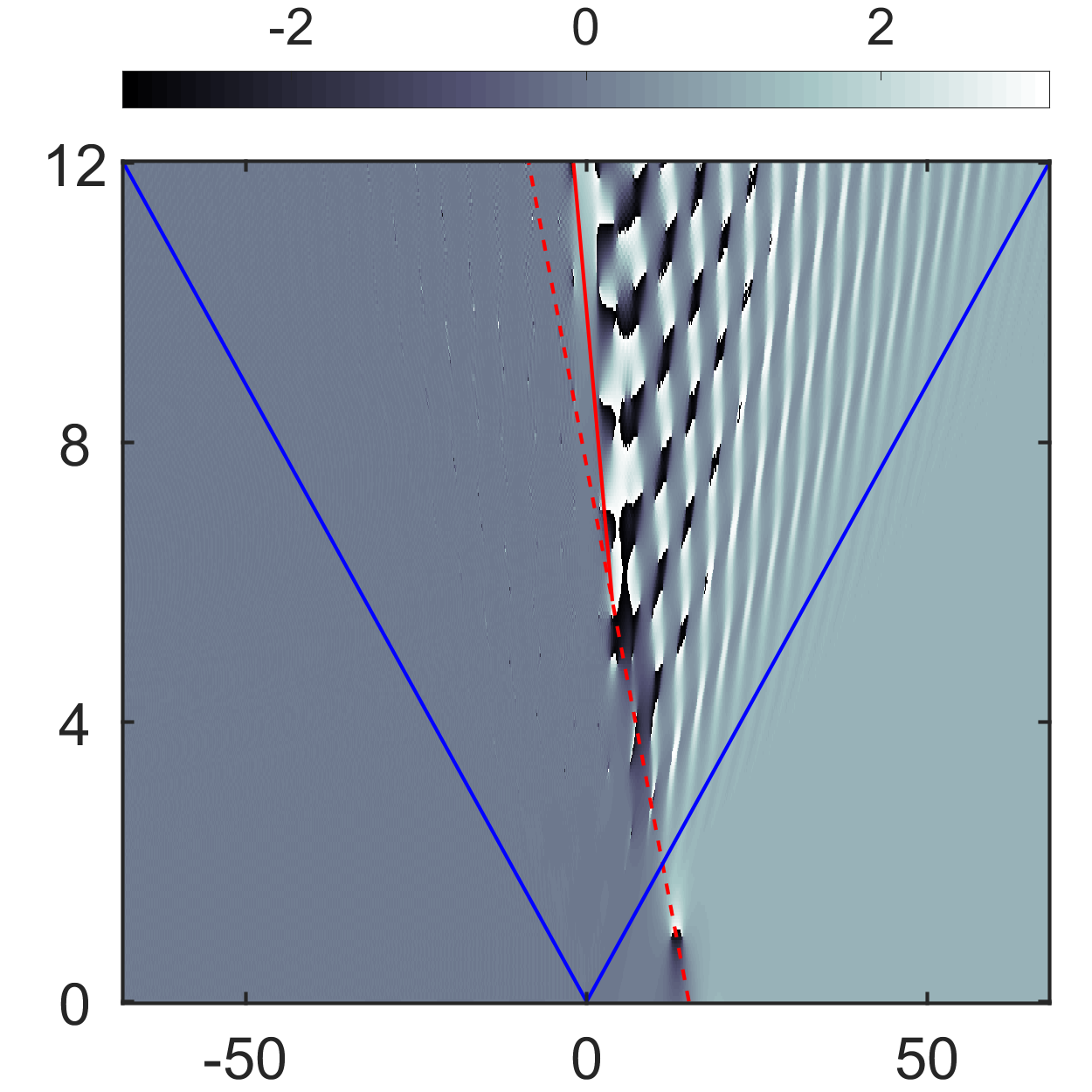}%
        \includegraphics[scale=0.25]{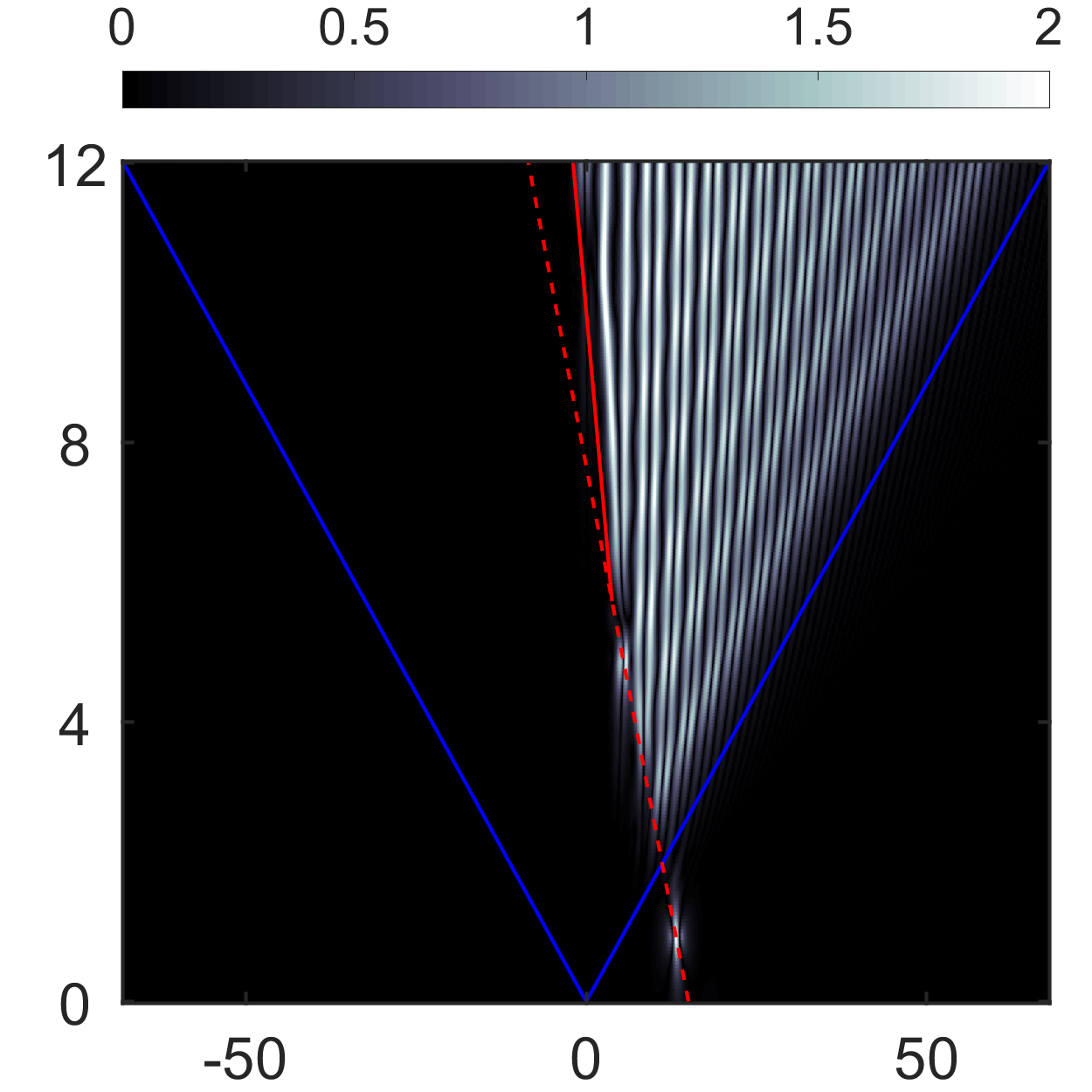}
\vspace*{2mm}
\caption{Numerical  solutions of the IVP \eqref{ds-fnls-ivp} with $q_o=1$ and different choices of $p$ for the two pure asymptotic regimes (red dots in Figure~\ref{ds-regions-f}).
\textit{Top row}: Transmission regime (Theorem \ref{ds-per-t}) with $p = -2 - 0.5i\in D_1$.
\textit{Bottom row}: Trap regime (Theorem \ref{ds-ptr-t}) with $p = -0.1 - 1.02i\in D_2^+$.
In all cases, the horizontal axis corresponds to $x$ and the vertical axis to $t$.
The grayscale used in  the plots is shown in the top inserts.
\textit{Left column}: Density plots of the solution amplitude $|q(x,t)|$.
\textit{Center column}: The complex phase difference between the solution in the left column and the solution  $q_\mathrm{wedge}(x,t)$   generated just by the localized disturbance (i.e. without the soliton),  illustrating  the asymptotic phase shift as $x\to\pm\infty$ induced by the soliton.
\textit{Right column}: The amplitude difference between the solution in the left column and $q_\mathrm{wedge}(x,t)$, illustrating the asymptotic position shift introduced by the soliton. 
\textit{Blue lines}: The boundary $x = \pm \Vo t$ between the wedge of modulated periodic oscillations from the left and right plane wave regions.
\textit{Dashed red lines}: The  original  soliton trajectory (velocity $\Vs$). 
\textit{Solid red lines}:  The  final soliton trajectory (which is either $\Vs$ or $\Vd$ depending on the regime). 
The position shift is very small in the transmission regime but becomes more noticeable in the trap regime. 
In both cases, the position shift is confined to the portion of the wedge lying above the soliton, in agreement with Theorems~\ref{ds-per-t} and~\ref{ds-ptr-t}.         
}
\label{numerics1}
\end{center}
\end{figure}

\begin{figure}[t!]
\vspace*{3mm}
    \begin{center}
        \includegraphics[scale=0.25]{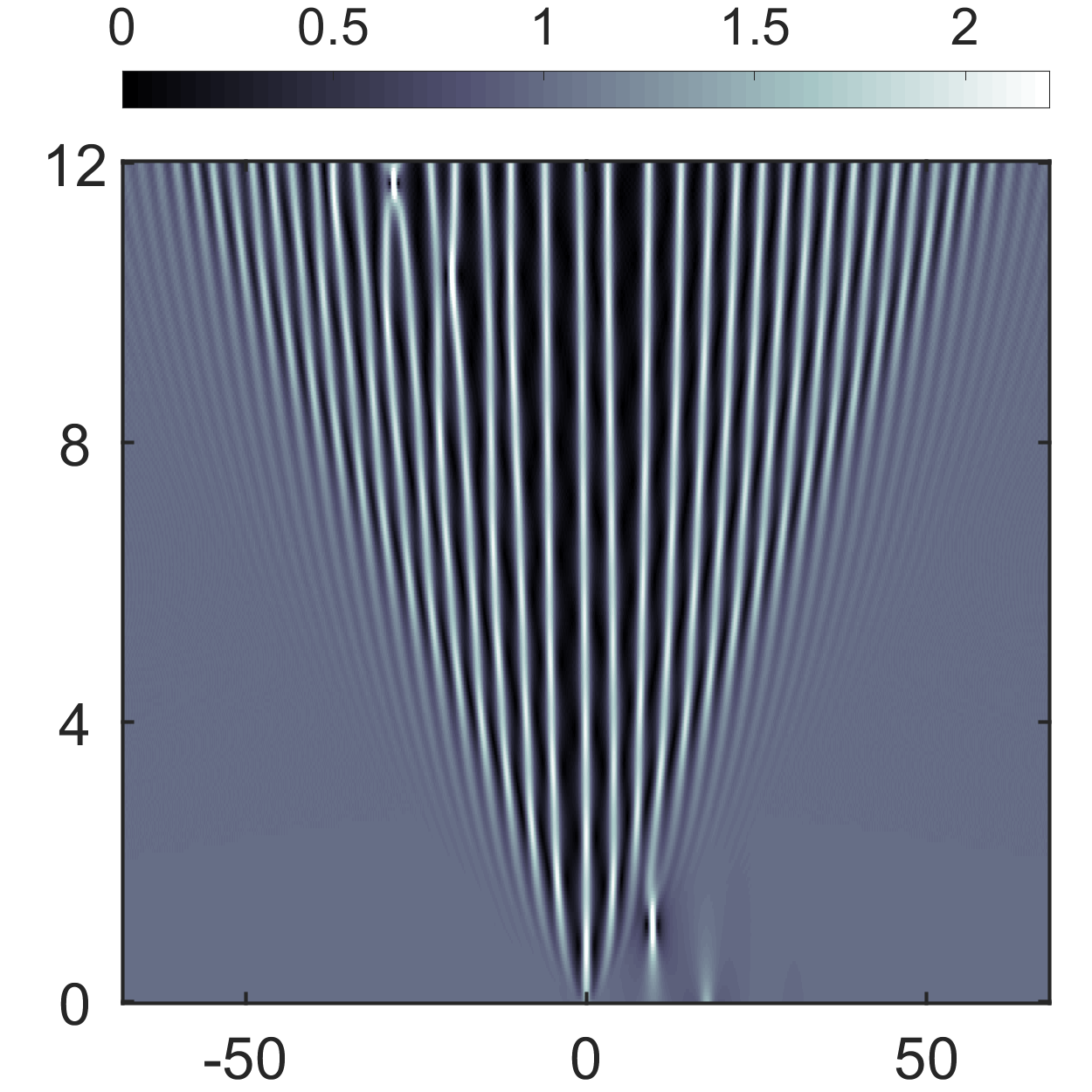}%
        \includegraphics[scale=0.25]{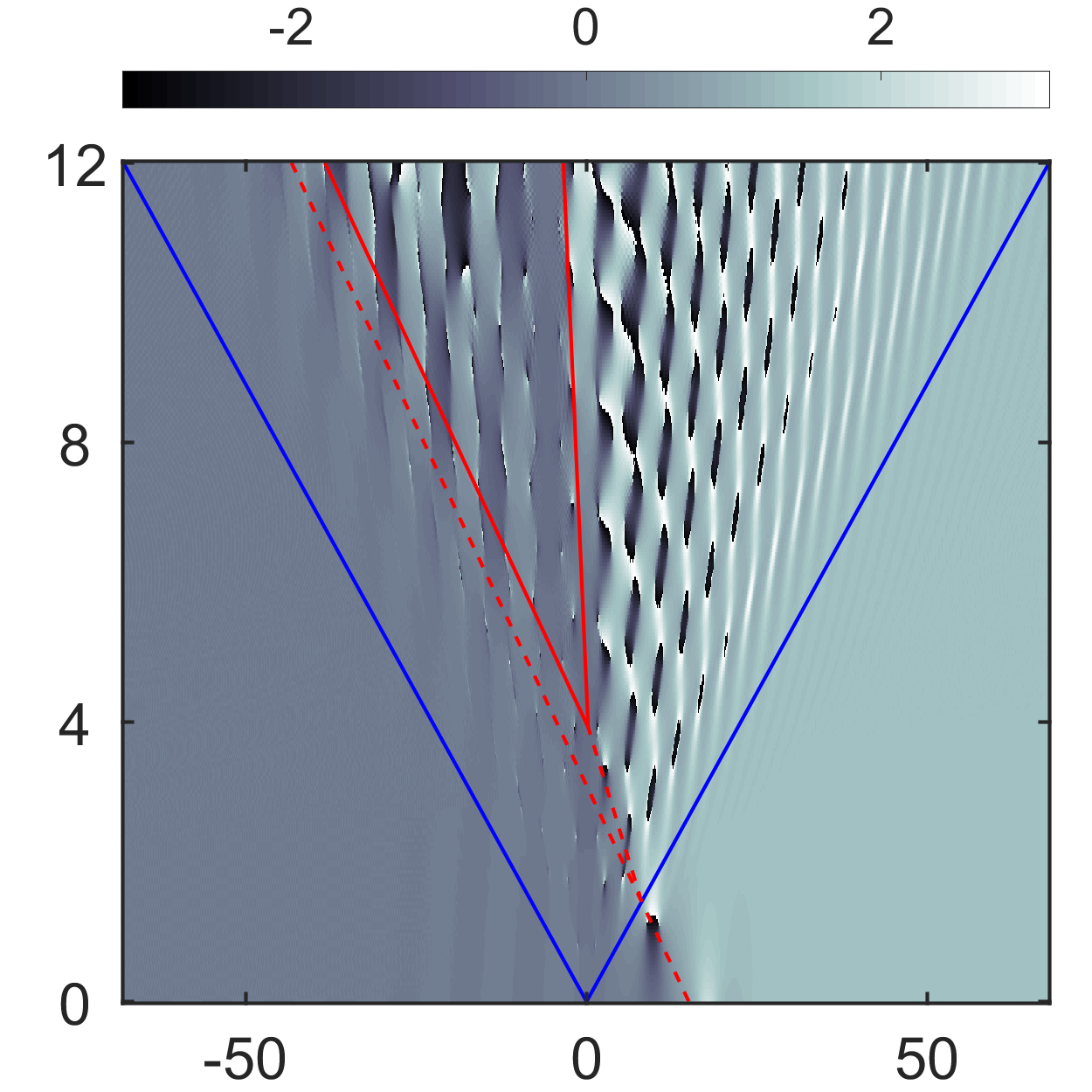}%
        \includegraphics[scale=0.25]{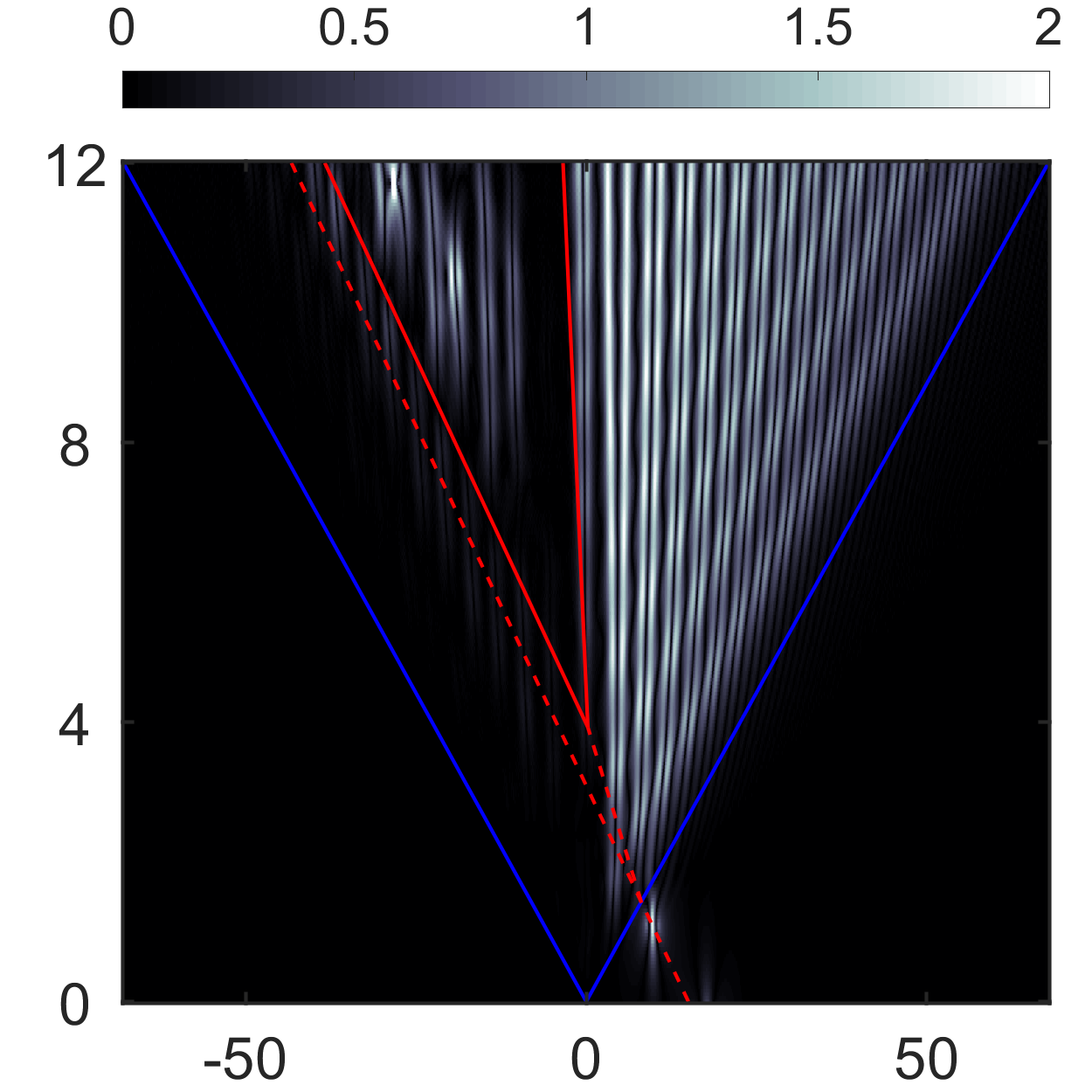}
        \\[2ex]
        \includegraphics[scale=0.25]{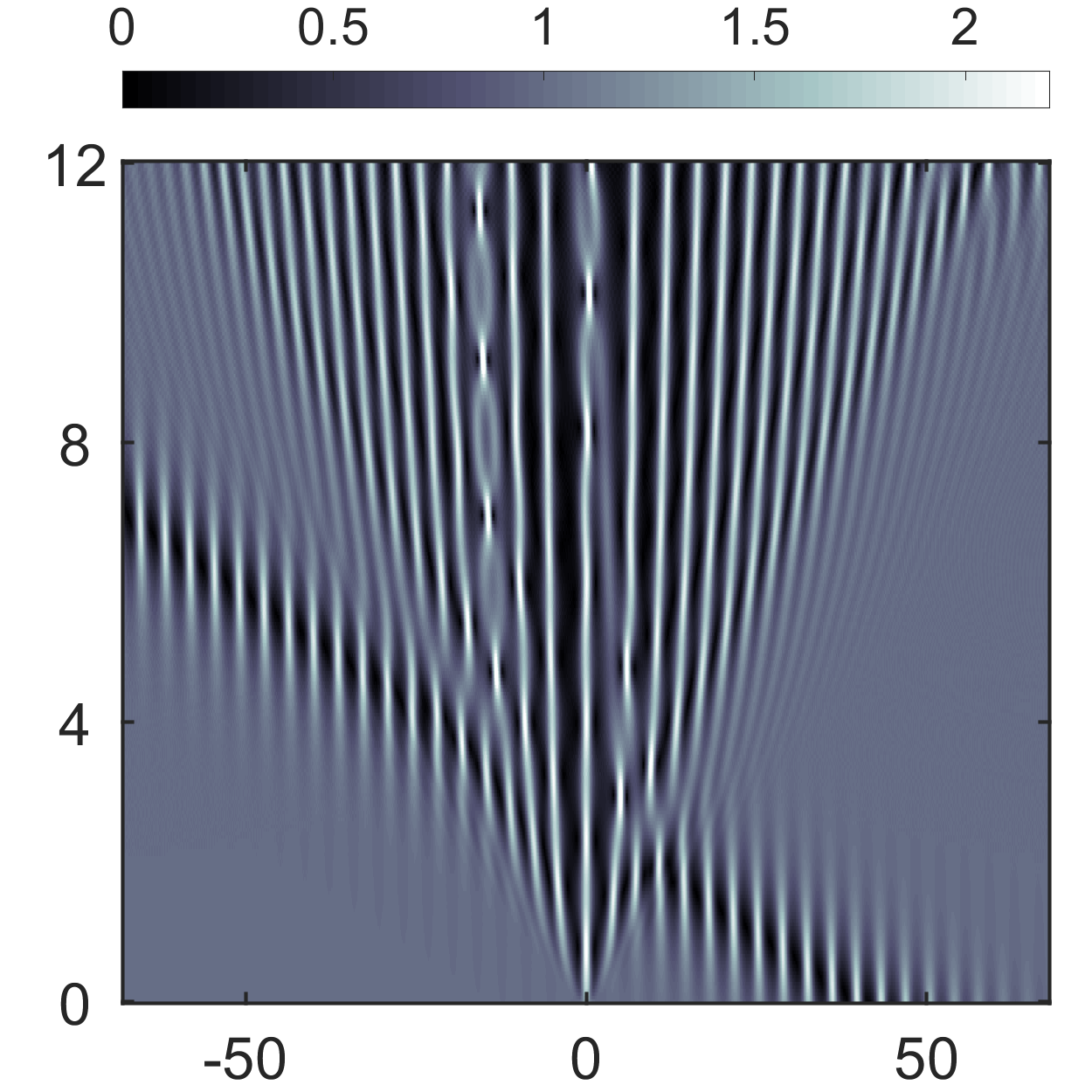}%
        \includegraphics[scale=0.25]{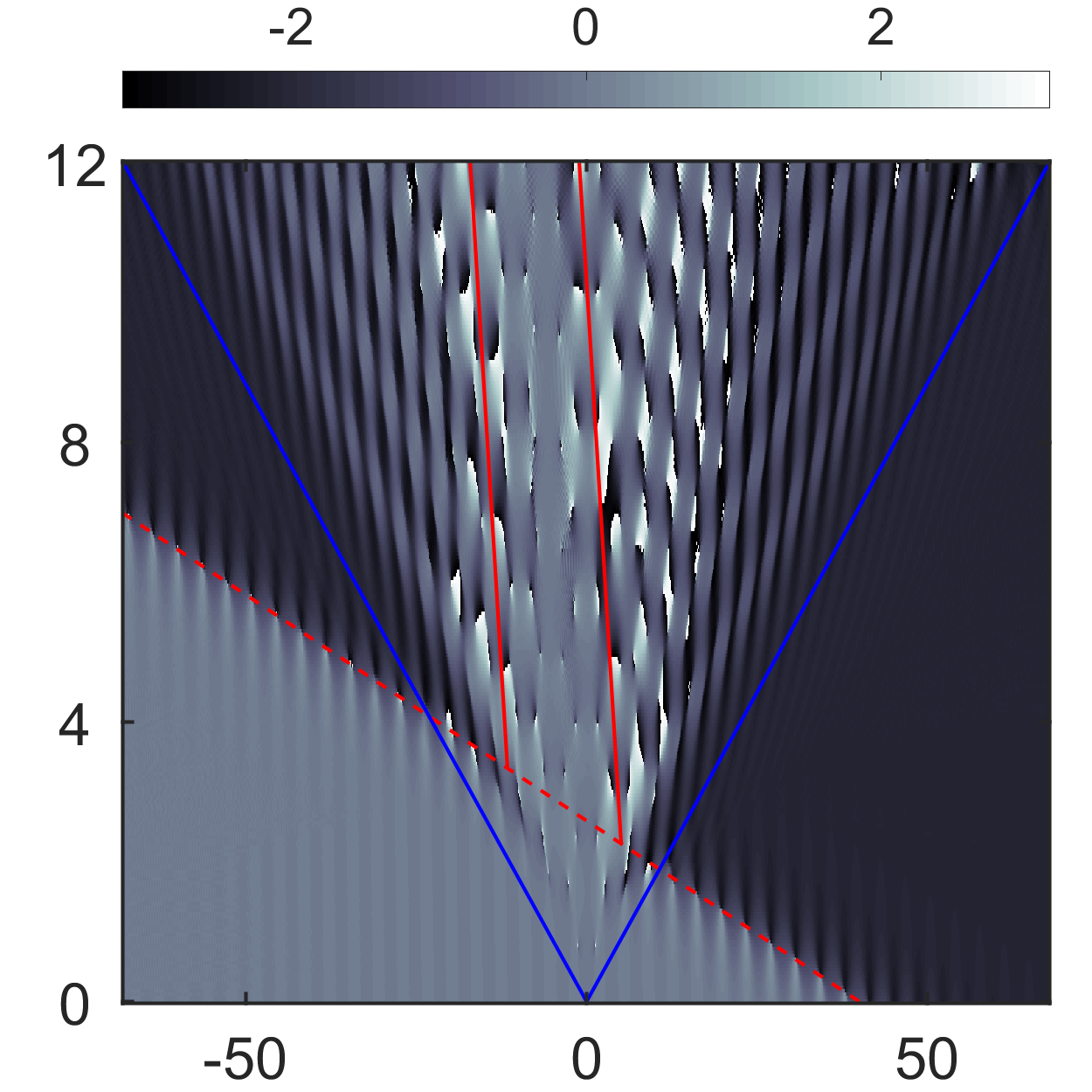}%
        \includegraphics[scale=0.25]{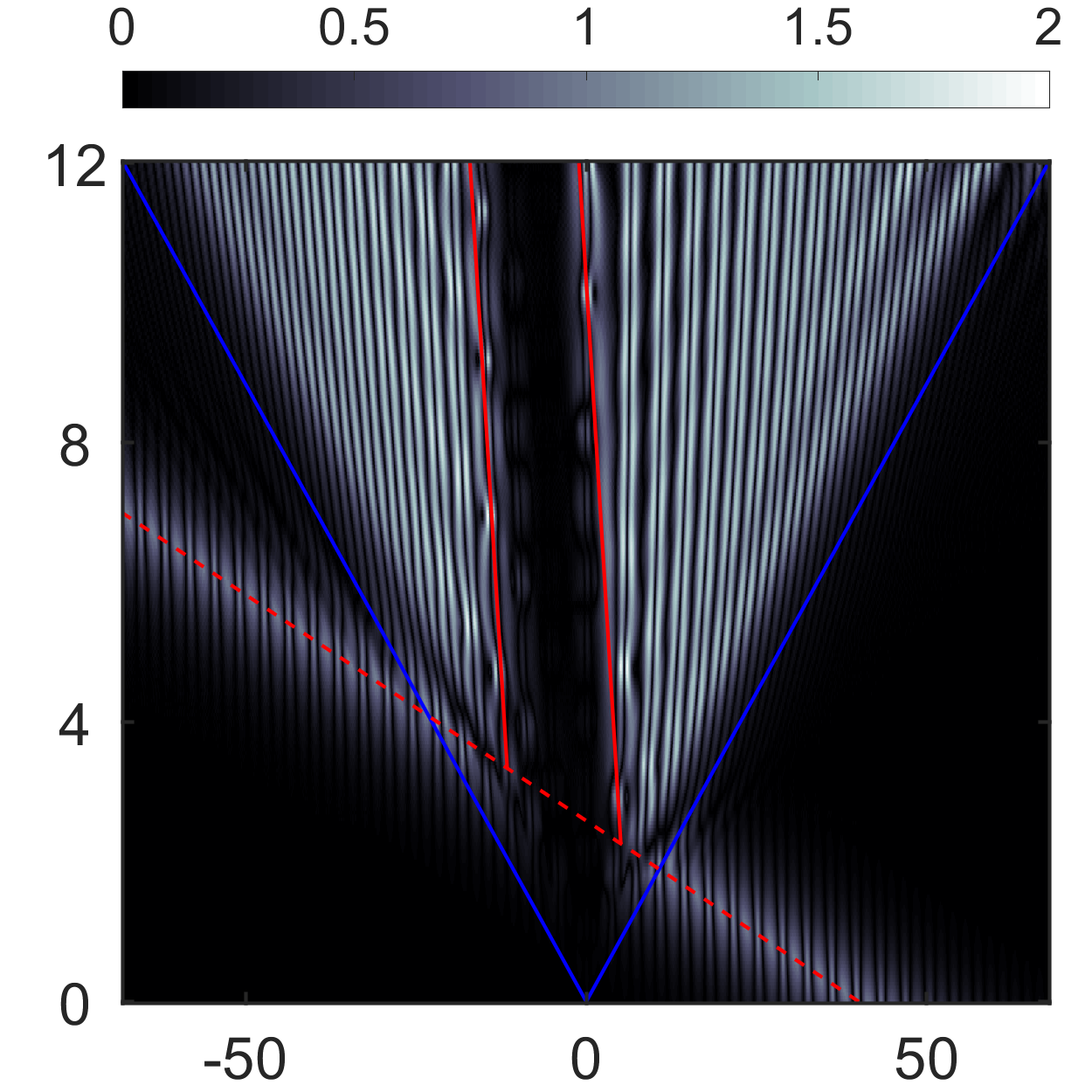}
        \vspace*{2mm}
\caption{
The analogue of Figure~\ref{numerics1} for $p$  corresponding to the two mixed asymptotic regimes (orange dots in Figure~\ref{ds-regions-f}).
Here, the solid red lines also identify the soliton wakes, propagating with velocity $\Vw$, where $|\Vw|<|\Vd|<|\Vs|$.
\textit{Top row}: Trap/wake regime (Theorem \ref{ds-twr-t}) with $p = -0.05 - 0.95i\in D_2^-$.
\textit{Bottom row}: Transmission/wake regime (Theorem \ref{ds-ewr-t}) with $p = -0.1 - 0.5i\in D_3$.
Note that the choice of $p$ in the trap/wake regime is very close to the boundary between $D_2^-$ and $D_3$ and, as a result, the dashed and solid red lines  corresponding respectively to the original and modified soliton velocities  are almost identical.  
Moreover, note that the two wakes seen in the transmission/wake regime have the same speed, and hence are detected as a single wake in the asymptotic analysis, which is why no wake-induced phase or position shift is observed in the asymptotics for $\xi>\Vw$.
}
\label{numerics2}
\end{center}
\end{figure}

\vspace*{3mm}
\noindent
\textbf{Long-time asymptotic results.}
We are now ready to give the precise form of the leading-order long-time asymptotics of the solution of the focusing NLS IVP \eqref{ds-fnls-ivp} 
in each of the four inequivalent regimes described above.
Numerical simulations with  the discrete eigenvalue chosen in each of the four regions of Figure~\ref{ds-regions-f},
illustrating the asymptotic results, 
are shown in Figures~\ref{numerics1} and \ref{numerics2}. 
For comparison purposes, Figures~\ref{numerics1} and \ref{numerics2}
also show the difference between $q(x,t)$ 
and the solution $q_\mathrm{wedge}(x,t)$ produced by an initial condition that generates the same reflection coefficient as $f(x)$ but no discrete spectrum. 
The numerical methods used in the  numerical  simulations were described in~\cite{blm2018}.
Recall that, since we are taking $\Re(p)<0$, all relevant velocities and all values of $\xi$ considered in 
Theorems~\ref{ds-per-t}--\ref{ds-ewr-t} below are negative.

\vspace*{1mm}

\begin{theorem}[\b{Transmission regime}]
\label{ds-per-t}
Suppose $p\in D_1$ and let $\Vs<\Vo<0$ be defined by \eqref{ds-xistar-xip-def}. 
Then  the solution $q(x, t)$ of the focusing NLS IVP \eqref{ds-fnls-ivp} exhibits the following asymptotic behavior as $t\to\infty$.
\vskip 2mm
\noindent 
\textnormal{(i)} 
If $\xi<\Vs$, then the leading-order asymptotics is described by the plane wave
\eee{\label{ds-qsol-pw-t}
q(x, t) = q_{\textnormal{pw}}(\xi) + O\big(t^{-\frac 12}\big), \quad t\to\infty,
}
where  
\eee{\label{ds-qpw-def}
q_{\textnormal{pw}}(\xi) := q_-\, e^{2ig_\infty(\xi)} 
}
and the real, constant phase $g_\infty(\xi)$ is given by \eqref{ds-esc-ginf-def-pw1}.
\vskip 2mm
\noindent 
\textnormal{(ii)} 
If $\xi=\Vs$, then the leading-order asymptotics is equal to a soliton on top of a nonzero plane-wave background, i.e.
\eee{\label{ds-esc-q-sol-lim-pw1-t}
q(x, t) 
=
q_{\textnormal{pw}}(\Vs) 
+
\qs(t)\,e^{2ig_\infty(\Vs)} 
+ O\big(t^{-\frac 12}\big), \quad t\to \infty,
}
with $q_{\textnormal{pw}}$ given by \eqref{ds-qpw-def}, $g_\infty(\Vs)$ defined by \eqref{ds-esc-ginf-def-pw1}, 
and the soliton $\qs$ given by
\eee{\label{ds-qstheta-def}
\qs(t)
=
 \frac{\left|R_p\right|\left(\bar{\mathcal A} \Lambda_1^2\bar q_- + \mathcal A \Lambda_2^2q_-  - 2 \mathcal B\Lambda_1 \Lambda_2 q_o\right) +  e^{i\left[2 \theta(\Vs,p) t + \arg(R_p)\right]} \Lambda_1^2\bar q_- +  e^{-i\left[2 \theta(\Vs,p) t + \arg(R_p)\right]} \Lambda_2^2q_-}
{4i\bar q_-\left\{\sqrt{|\mathcal A|^2-\mathcal B^2}\cosh\big[\ln\big(|R_p| \sqrt{|\mathcal A|^2-\mathcal B^2}\,\big)\big] + \Re \big(\mathcal A e^{i\left[2 \theta(\Vs,p) t + \arg(R_p)\right]}\big)\right\}}
}
with the constants $R_p$, $(\mathcal A, \mathcal B)$ and $(\Lambda_1, \Lambda_2)$ given by \eqref{ds-R-def},  \eqref{ds-ab-sol-def} and \eqref{ds-DE-def} respectively.
\vskip 2mm
\noindent 
\textnormal{(iii)} 
If $\Vs<\xi<\Vo$, then the leading-order asymptotics is given by the plane wave \eqref{ds-qsol-pw-t} up to a constant phase shift, namely
\eee{\label{ds-esc-qasym-pw2-t}
q(x, t) =  q_{\textnormal{pw}}(\xi)\,e^{4i\, \textnormal{arg}\left[p+\lambda(p)\right]} + O\big(t^{-\frac 12}\big),\quad t\to\infty.
} 
%
%
\noindent 
\textnormal{(iv)} 
Finally, if $\Vo<\xi<0$, then the asymptotic behavior of the solution is described at leading order by the phase-shifted modulated elliptic wave
\eee{\label{ds-esc-qasym-mew-t}
q(x, t) = \widetilde q_{\textnormal{mew}}(x, t)\,e^{4i\, \textnormal{arg}\left[p+\lambda(p)\right]} + O\big(t^{-\frac 12}\big),\quad t\to\infty,
} 
where 
\eee{\label{ds-qsol-mew-t}
\widetilde q_{\textnormal{mew}}(x, t)
=
\frac{q_o \left(q_o +\alpha_\im\right) }{\bar q_-}
\
\frac{
\Theta\Big(
\frac{\sqrt{q_o\alpha_\im}}{mK(m)}
\left(x-2\alpha_\re t\right)
-X_o  +2\upnu_\infty -\frac 12 - \frac{\widetilde \omega}{2\pi} \Big)
\Theta\left(\frac 12 \right)
}
{
\Theta\Big(
\frac{\sqrt{q_o\alpha_\im}}{mK(m)}
\left(x-2\alpha_\re t\right)
-X_o -\frac 12  - \frac{\widetilde \omega}{2\pi} \Big)
\Theta\left(2\upnu_\infty -\frac 12\right)
}
\,
e^{2i\left[g_\infty(\xi) -G_\infty(\xi)  t\right]}
}
with the Jacobi function $\Theta$ defined by \eqref{ds-Theta-def},  the complex quantity $\upnu_\infty$ given by \eqref{ds-vdefr}, and the real quantities $G_\infty$, $g_\infty$, $X_o$ and  $\widetilde \omega$  defined by equations \eqref{ds-Ginf-mew}, \eqref{ds-esc-ginf-def-mew}, \eqref{ds-Xdef} and \eqref{ds-esc-omegat-def-mew} respectively.
Importantly, all of these quantities depend on $x$ and $t$ only through the similarity variable $\xi$.
\end{theorem}

\vspace*{0mm}

\begin{theorem}[\b{Trap regime}]
\label{ds-ptr-t}
Suppose $p\in D_2^+$ and let $\Vd$ be the unique solution of equation \eqref{ds-xih-def} in the interval $(\Vo, 0)$. 
Then the solution $q(x, t)$ of the focusing NLS IVP \eqref{ds-fnls-ivp} exhibits the following asymptotic behavior as $t\to \infty$.
\vskip 2mm
\noindent 
\textnormal{(i)} 
If $\xi<\Vo$, then the leading-order asymptotics is given by the plane wave \eqref{ds-qsol-pw-t}.
\vskip 2mm
\noindent 
\textnormal{(ii)} 
If $\Vo<\xi<\Vd$, then the leading-order asymptotics is described by the modulated elliptic wave
\eee{\label{ds-trap-qasym-mew-t}
q(x, t) =  q_{\textnormal{mew}}(x, t) + O\big(t^{-\frac 12}\big),\quad t\to\infty,
} 
where $q_{\textnormal{mew}}$ is obtained from  \eqref{ds-qsol-mew-t} after setting $\widetilde \omega=0$, i.e.
\eee{\label{ds-qsol-mew-cpam-t}
q_{\textnormal{mew}}(x, t)
=
\frac{q_o \left(q_o +\alpha_\im\right) }{\bar q_-}
\
\frac{
\Theta\Big(
\frac{\sqrt{q_o\alpha_\im}}{mK(m)}
\left(x-2\alpha_\re t\right)
-X_o  +2\upnu_\infty -\frac 12  \Big)
\Theta\left(\frac 12 \right)
}
{
\Theta\Big(
\frac{\sqrt{q_o\alpha_\im}}{mK(m)}
\left(x-2\alpha_\re t\right)
-X_o -\frac 12  \Big)
\Theta\left(2\upnu_\infty -\frac 12\right)
}
\,
e^{2i\left[g_\infty(\xi) -G_\infty(\xi)  t\right]}.
}
\vskip 2mm
\noindent 
\textnormal{(iii)} 
If $\xi=\Vd$, then at leading order the asymptotics is equal to a soliton on top of a nonzero modulated-elliptic-wave background, i.e.
\eee{\label{ds-trap-q-sol-lim-pw1-t}
q(x, t) 
=
q_{\textnormal{mew}}(\Vd t, t) + \qd(t)  + O\big(t^{-\frac 12}\big), \quad t\to \infty,
}
where the modulated elliptic wave $q_{\textnormal{mew}}$ is defined by \eqref{ds-qsol-mew-cpam-t} and the soliton $\qd$ is given by
\vspace*{2mm}
\eee{\label{ds-qsh-def}
\qd(t)
=
2i\, 
\frac{2\mathcal B\rho_p\rho_{\bar p} W_{11}(p) W_{12}(\bar p)
-\left(1+\mathcal C\rho_{\bar p}\right) \rho_p W_{11}(p)^2
+\left(1+\mathcal A\rho_p\right)  \rho_{\bar p} W_{12}(\bar p)^2
}{\mathcal B^2 \rho_p\rho_{\bar p} + \left(1 + \mathcal C\rho_{\bar p}\right)\left(1+\mathcal A\rho_p\right)}
}
\vspace*{0mm}

\noindent
with $(\rho_p, \rho_{\bar p})$, $W$ and $(\mathcal A, \mathcal B, \mathcal C)$ given by \eqref{ds-rho-def}, \eqref{ds-trap-W-mew1} and \eqref{ds-trap-abc-def} respectively.
\vskip 2mm
\noindent 
\textnormal{(iv)} 
Finally, if $\Vd<\xi<0$, then the asymptotics is given by the phase-shifted modulated elliptic wave \eqref{ds-esc-qasym-mew-t}.
\end{theorem}

\vspace*{1mm}

\begin{theorem}[\b{Trap/wake regime}]
\label{ds-twr-t}
Suppose $p\in D_2^-$ and let $\Vd<\Vw$ be the two solutions of equation \eqref{ds-xih-def} in the interval $(\Vo, 0)$. Then, the solution $q(x, t)$ of the focusing NLS IVP \eqref{ds-fnls-ivp} exhibits the following asymptotic behavior as $t\to \infty$.
\vskip 2mm
\noindent \textnormal{(i)} If $\xi<\Vo$, then the leading-order asymptotics  is described by the plane wave \eqref{ds-qsol-pw-t}.
\vskip 2mm
\noindent \textnormal{(ii)} If $\Vo<\xi<\Vd$, then the leading-order asymptotics   is given by the modulated elliptic wave \eqref{ds-trap-qasym-mew-t}.
\vskip 2mm
\noindent \textnormal{(iii)} If $\xi=\Vd$, then the asymptotics  is characterized by  \eqref{ds-trap-q-sol-lim-pw1-t}, namely at leading order it is equal to the sum of the modulated elliptic wave \eqref{ds-qsol-mew-cpam-t} evaluated at $x=\Vd t$ and the soliton  \eqref{ds-qsh-def}. 
\vskip 2mm
\noindent \textnormal{(iv)} If $\Vd<\xi<\Vw$, then the leading-order asymptotics  is given by the phase-shifted  modulated elliptic wave \eqref{ds-esc-qasym-mew-t}.
\vskip 2mm
\noindent \textnormal{(v)} If $\xi=\Vw$, then at leading order the asymptotics is equal to a soliton wake on top of a nonzero modulated-elliptic-wave background, i.e.
\eee{\label{ds-trapwake-t}
q(x, t) 
=
q_{\textnormal{mew}, w}(t) + q_w(t)  + O\big(t^{-\frac 12}\big), \quad t\to \infty,
}
where the modulated elliptic wave $q_{\textnormal{mew}, w}(t)$ is given by \eqref{ds-qsol-mew-cpam-t} evaluated at $x=\Vw t$ but with $\omega$ in $X_o$ replaced by $\omega_w$ of \eqref{ds-wake-omega'-def}  and with $g_\infty$ replaced by $g_{w, \infty}$ of  \eqref{ds-ginf'-def},
%
and the soliton wake $q_w$ is defined by
\eee{\label{ds-trapwake-qw-def}
q_w(t)
:=
2i\, 
\frac{2  \mathcal B_w {\rho_p}_w {\rho_{\bar p}}_w  {W_w}_{11}(\bar p)  {W_w}_{12}(p)
-\big(1+   \mathcal C_w  {\rho_{p}}_w\big) {\rho_{\bar p}}_w  {W_w}_{11}(\bar p)^2
+\big(1+\mathcal A_w {\rho_{\bar p}}_w\big)  {\rho_{p}}_w  {W_w}_{12}(p)^2
}{\mathcal B_w^2  { \rho_p}_w  {\rho_{\bar p}}_w + \big(1 +  \mathcal C_w {\rho_{p}}_w\big)\big(1+ \mathcal A_w  {\rho_{\bar p}}_w\big)},
}
with $({\rho_p}_w, {\rho_{\bar p}}_w)$, $W_w$ and $(\mathcal A_w,  \mathcal B_w, \mathcal C_w)$ given by \eqref{ds-rrb-til}, \eqref{ds-wake-ww-sol} and \eqref{ds-abc-w-def} respectively.
\vskip 2mm
\noindent \textnormal{(vi)} Finally, if $\Vw<\xi<0$, then the leading-order asymptotics is the same with the one in the range $\Vd<\xi<\Vw$, namely it is given by the phase-shifted modulated elliptic wave \eqref{ds-esc-qasym-mew-t}.
\end{theorem}

\vspace*{0mm}

\begin{theorem}[\b{Transmission/wake regime}]
\label{ds-ewr-t}
Suppose $p\in D_3$, let $\Vs<\Vo<0$ be defined by \eqref{ds-xistar-xip-def}, and let $\Vw$ be the unique solution of equation \eqref{ds-xih-def} in the interval $(\Vo, 0)$. Then, the solution $q(x, t)$ of the focusing NLS IVP \eqref{ds-fnls-ivp} exhibits the following asymptotic behavior as $t\to \infty$.
\vskip 2mm
\noindent \textnormal{(i)} If $\xi<\Vs$, then the leading-order asymptotics  is given by the plane wave \eqref{ds-qsol-pw-t}.
\vskip 2mm
\noindent \textnormal{(ii)} If $\xi=\Vs$, then the asymptotics  is characterized by   \eqref{ds-esc-q-sol-lim-pw1-t}, namely at leading order it is given by the superposition of the plane wave \eqref{ds-qpw-def} and the soliton \eqref{ds-qstheta-def}.
\vskip 2mm
\noindent \textnormal{(iii)} If $\Vs<\xi<\Vo$, then the leading-order asymptotics   is described by the phase-shifted plane wave~\eqref{ds-esc-qasym-pw2-t}.
\vskip 2mm
\noindent \textnormal{(iv)} If $\Vo<\xi<\Vw$, then the leading-order asymptotics  is given by the phase-shifted modulated elliptic wave \eqref{ds-esc-qasym-mew-t}.
\vskip 2mm
\noindent \textnormal{(v)} If $\xi=\Vw$, then the  asymptotics is characterized by \eqref{ds-trapwake-t},  i.e. at leading order it is equal to the sum of the modulated elliptic wave $q_{\text{mew}, w}(t)$ and the soliton wake  \eqref{ds-trapwake-qw-def}.
\vskip 2mm
\noindent \textnormal{(vi)} Finally, if $\Vw<\xi<0$, then the leading-order asymptotics is the same with the one in the range $\Vo<\xi<\Vw$, namely it is given  by the phase-shifted modulated elliptic wave \eqref{ds-esc-qasym-mew-t}.
\end{theorem}

\begin{remark}[\b{Leading-order asymptotics for $x>0$}]
Since the pole $p$ lies in the third quadrant of the complex $k$-plane, it has no effect on the asymptotics for $x>0$ (equivalently, $\xi>0$; see Figure \ref{sign-structure-f}). In particular, for $x>0$ the leading-order asymptotics of IVP \eqref{ds-fnls-ivp} is described by  Theorems~1.1 and 1.2 of \cite{bm2017}, the only difference being that now one must also include the constant phase shift $4 \textnormal{arg}\left[p+\lambda(p)\right]$ and the position shift induced by the soliton arising for $x<0$.
\end{remark}

\begin{remark}
In  the appendix, we  explicitly verify that the expression \eqref{ds-qstheta-def} for the soliton obtained in the long-time asysmptotics agrees with the long-time asymptotics of the standard soliton solution of the focusing NLS with nonzero background.
\end{remark}

\begin{remark}[\b{Soliton versus soliton wake}]
\label{ds-shifts-r}
The soliton arising either at $\xi=\Vs$ or at $\xi=\Vd$ induces a constant phase shift (equal to $4\textnormal{arg}\left[p+\lambda(p)\right]$) as well as  a position shift (related to the presence of $\widetilde \omega$  in \eqref{ds-qsol-mew-t} as opposed to \eqref{ds-qsol-mew-cpam-t}) in the leading-order asymptotics for subsequent values of $\xi$. 
On the contrary, the soliton wake arising at $\xi=\Vw$ has no effect on the leading-order asymptotics for $\xi>\Vw$. 
The numerical simulations of Figures \ref{numerics1} and \ref{numerics2} illustrate these remarks.
\end{remark}

\begin{remark}[\b{Multiple leading-order contributions from the poles}]
We find it quite remarkable that 
\textit{in the two mixed regimes a single pair of complex conjugate poles produces $O(1)$ contributions to the solution
at two different velocities}: 
the soliton velocity and the wake velocity,
as specified in Theorems~\ref{ds-twr-t} and~\ref{ds-ewr-t}.
(These predictions are validated by the numerical results in Figures~\ref{numerics1} and~\ref{numerics2}.)
To the best of our knowledge, this is the first time that such a phenomenon has been observed in the long-time asymptotic analysis of an integrable system, and is perhaps one of the main novelties in the results of the present work. 
Moreover, the numerical results in the bottom row of Figure~\ref{numerics2} suggest that \textit{the soliton-generated wake may comprise itself two different localized structures}. We emphasize however that, since these two structures propagate with the same velocity, in order to be able to differentiate between them one would have to compute the asymptotics by taking $x = \xi t + y$. 
Such a calculation is outside the scope of this work.
\end{remark}

\noindent
\textbf{Structure of the paper.}
In Section \ref{ds-rhp-s}, the solution of IVP \eqref{ds-fnls-ivp} for the focusing NLS equation is associated with the solution of a matrix Riemann-Hilbert problem via the inverse scattering transform. Furthermore, the four different long-time asymptotic patterns, namely the transmission, trap, trap/wake, and transmission/wake regimes, are motivated through the  behavior of the jump matrices of this Riemann-Hilbert problem. 
The transmission regime is analyzed in Section \ref{ds-esc-s}, resulting in the proof of Theorem \ref{ds-per-t}. The proof of Theorem \ref{ds-ptr-t} for the trap regime is provided in Section \ref{ds-trap-s}. The two mixed regimes are discussed in Section \ref{ds-mixed-s}, leading to the proofs of Theorems~\ref{ds-twr-t} and \ref{ds-ewr-t}. Finally, some concluding remarks are given in Section \ref{s:conclusions}.

\section{The Riemann-Hilbert Problem and Outline of the Asymptotic Analysis}
\label{ds-rhp-s}

The implementation of the inverse scattering transform method for IVP \eqref{ds-fnls-ivp}  begins with the integration of the Lax pair \eqref{ds-fnls-lp} for the  $2\times 2$ matrix-valued function $\Psi$ assuming as usual  that the solution $q$ of problem \eqref{ds-fnls-ivp} is given. This task is known as the direct problem.  Then,   $q$  is expressed in terms of a sectionally  meromorphic function $M$ which is defined via appropriate combinations of the two column vectors $\Psi_1$ and $\Psi_2$ of $\Psi$, and which satisfies a certain matrix Riemann-Hilbert problem. This portion of the analysis is known as the inverse problem.
Specifically, the discussion of the direct problem in Section~2 of~\cite{bm2017} motivates the following definition for the sectionally meromorphic matrix-valued function $M$:
\eee{\label{ds-mrh}
M(x, t, k)
=
\left\{
\def\arraystretch{2.1}
\begin{array}{ll}
\left( \dfrac{\Psi_{+1}(x, t, k)}{\bar a(k) d(k)}, \Psi_{-2}(x, t, k) \right) e^{-i\theta(\xi, k) t \sigma_3},
&k\in\mathbb C^+\setminus B^+,
\\
\left( \Psi_{-1}(x, t, k), \dfrac{\Psi_{+2}(x, t, k)}{a(k) d(k)} \right) e^{-i\theta(\xi, k) t \sigma_3}, 
&k\in\mathbb C^-\setminus B^-.
\end{array}
\right.
} 
In the above definition, we use the notation
$$
\mathbb C^{\pm} := \left\{k\in \mathbb C: \Im (k) \gtrless 0\right\},
\quad
B^+:= i\left[0, q_o\right], 
\quad
B^- := i\left[-q_o, 0\right]
$$
and  denote  by $\Psi_\pm$ the so-called Jost solutions, namely the simultaneous solutions of the Lax pair \eqref{ds-fnls-lp} with prescribed normalizations as $x\to\pm\infty$:
\eee{\label{ds-infcond}
\Psi_\pm(x, t, k)
=
\begin{pmatrix}
1 & {i(\lambda - k)}/{\bar q_\pm}
\\
{i(\lambda -k)}/{q_\pm} & 1
\end{pmatrix} e^{i\theta(\xi, k) t \sigma_3} \left[1+o(1) \right],
\quad
x\to \pm\infty.
}
(Recall that the quantities $\lambda$, $\theta$ and $\xi$ are defined by  \eqref{ds-lambda-def},  \eqref{ds-theta-def} and \eqref{ds-xi-def} respectively.)
Furthermore, we define the spectral function $a(k)$ along with  its Schwarz conjugate $\bar a(k)$ by
\eee{\label{ds-a-ab-def}
a(k) = \frac{\text{wr} \left[\Psi_{-1}(x, t, k), \Psi_{+2}(x, t, k) \right]}{d(k)}, 
\quad
\bar a(k) := \overline{a(\bar k)}
= 
\frac{\text{wr} \left[\Psi_{+1}(x, t, k), \Psi_{-2}(x, t, k) \right]}{d(k)},
} 
where ``$\text{wr}$'' denotes the Wronskian determinant and 
\eee{\label{ds-d-def}
d(k) := \frac{2\lambda(k)}{\lambda(k)+k}.
}
Importantly, the Wronskian determinants   appearing  in \eqref{ds-a-ab-def} are independent of $x$ and $t$, and hence the functions $a$ and $\bar a$  depend only on $k$.

\begin{figure}[t!]
\begin{center}
\includegraphics[height=4.5cm, width=5cm]{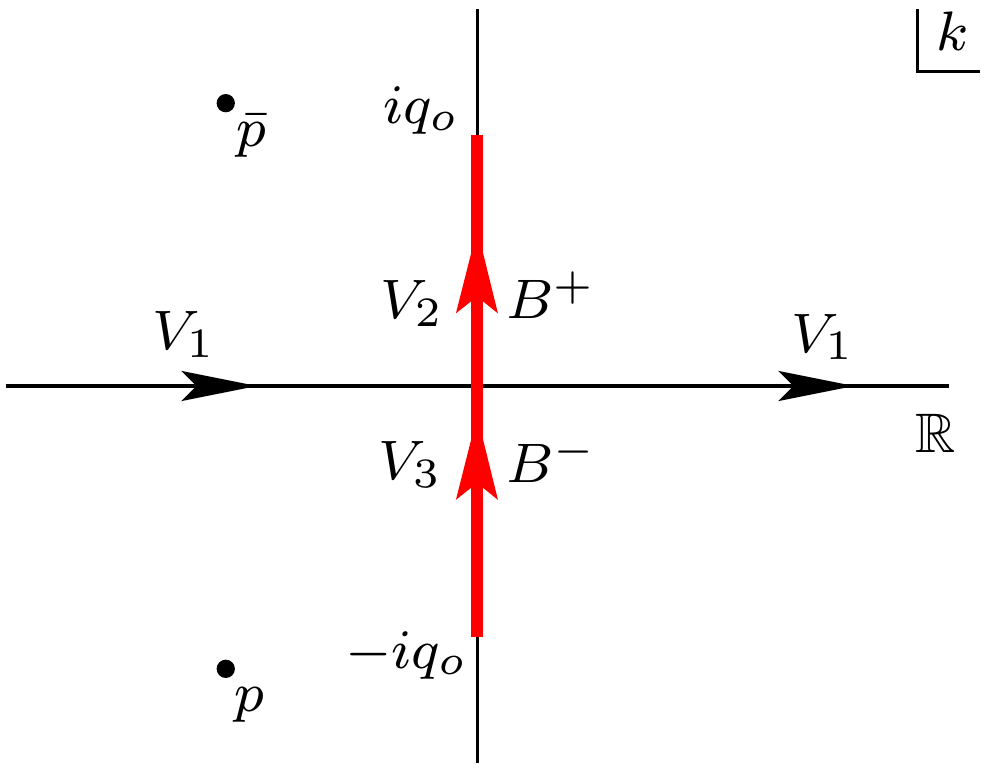}
\caption{The branch cut $B=B^+\cup B^-$,  the jumps $V_1$, $V_2$, $V_3$ of Riemann-Hilbert problem \eqref{ds-rhp-cpam}, and the conjugate pair of simple poles  $p, \bar p$.}
\label{ds-jumps-cpam-f}
\end{center}
\end{figure}

The definition \eqref{ds-mrh} of $M$, in combination with the analyticity properties of $\Psi_\pm$ (see \cite{bm2017} for more details), implies that the only sources of nonanalyticity of  $M$ are
\begin{enumerate}[label=(\roman*)]
\advance\itemsep 1mm
\item the \textit{continuous spectrum} 
\eee{
\Sigma := \mathbb R\cup B,
}
along which $M$ exhibits jump discontinuities, and
\item the possible zeros of the spectral function $a(k)$, which form the \textit{discrete spectrum} of the Riemann-Hilbert problem satisfied by $M$.
\end{enumerate} 
It was shown in \cite{bm2017} that, if there is no discrete spectrum, namely, if
\eee{\label{ds-cpam-as}
a(k)\neq 0 \quad \forall k \in \mathbb C^-\cup \Sigma,
} 
then the function $M(x, t, k)$ satisfies the following Riemann-Hilbert problem:
\sss{\label{ds-rhp-cpam}
\ddd{
M^+ &= M^- V_1,    &&k\in \mathbb R,
\\
M^+ &= M^- V_2, &&k\in B^+,
\\
M^+ &= M^- V_3, &&k\in B^-,
\\
M &= I +O\left(\tfrac 1k\right), \quad && k \to \infty,
}
}
where the jump matrices along the three contours $\mathbb R$, $B^+$, $B^-$ comprising the continuous spectrum $\Sigma$ are given by 
(see Figure~\ref{ds-jumps-cpam-f})
\sss{\label{ds-rhp-cpam-j}
\ddd{
V_1(x, t, k)
&=
\left(
\def\arraystretch{1.4}
\begin{array}{lr}
\dfrac{1}{d(k)}\left[1+r(k)\bar r(k)\right]
&
\bar r(k)e^{2i \theta(\xi, k) t} 
\\
r(k) e^{-2i \theta(\xi, k) t} 
&
d(k)
\end{array}
\right),
\\
V_2(x, t, k)
&=
\left(
\def\arraystretch{1.8}
\begin{array}{lr}
- \dfrac{\lambda(k)-k}{iq_-}\,  \bar r(k) \, e^{2i \theta(\xi, k) t}
&
 \dfrac{2\lambda(k)}{i\bar q_-} 
\\
 \dfrac{\bar q_-}{2i \lambda(k)}  \left[1+r(k)\bar r(k)\right] 
 &
 -\dfrac{\lambda(k)+k}{i\bar q_-} \, r(k) \,  e^{-2i \theta(\xi, k) t}
\end{array}
\right),
\\
V_3(x, t, k)
&=
\left(
\def\arraystretch{2}
\begin{array}{lr}
\dfrac{\lambda(k)+k}{iq_-}\, \bar r(k) \, e^{2i \theta(\xi, k) t}
& 
\dfrac{q_-}{2i\lambda(k)}\left[1+r(k)\bar r(k)\right]
\\
\dfrac{2\lambda(k)}{iq_-} 
&
\dfrac{\lambda(k)-k}{i\bar q_-}\, r(k)\,  e^{-2i \theta(\xi, k) t}
\end{array}
\right),
} 
}
with the reflection coefficient $r$ defined by
\eee{\label{ds-r-coef-def}
r(k) = -\frac{b(k)}{\bar a(k)}, \quad b(k):= \frac{\text{wr} \left[\Psi_{+1}(x, t, k), \Psi_{-1}(x, t, k) \right]}{d(k)}.
} 

Removing the assumption \eqref{ds-cpam-as}, i.e. allowing the spectral function $a(k)$ to vanish in $\mathbb C^-$, 
results in a Riemann-Hilbert problem  with a \textit{nonempty} discrete spectrum. 
In this work, we consider the simplest such scenario, according to which the initial data ${f}(x)$ of  IVP \eqref{ds-fnls-ivp} is such that $a(k)$ has a \textit{unique, simple zero} in $\mathbb C^-\setminus \Sigma$. That is, we assume that there exists a unique $p\in \mathbb C^-\setminus \Sigma$ such that $a(p)=0$ and, furthermore, $a'(p)\neq 0$. 
Correspondingly, the Schwarz conjugate $\bar a(k)$ of $a(k)$ possesses a \textit{unique, simple zero} $\bar p\in \mathbb C^+\setminus \Sigma$ and, by definition \eqref{ds-mrh}, $M$ is meromorphic in $\mathbb C\setminus \Sigma$ with two simple poles, at $k=p$ and at $k=\bar p$. 
Therefore,  in addition to the jumps $V_1, V_2, V_3$ along the continuous spectrum $\Sigma$, 
the Riemann-Hilbert problem for $M$  must be supplemented by suitable residue conditions at $p$ and $\bar p$. 
These can be computed as follows.

Since $a(p)=0$, by expression  \eqref{ds-a-ab-def}  we have that 
$
\text{wr} \left[\Psi_{-1}(x, t, p), \Psi_{+2}(x, t, p) \right] = 0
$
for all $x, t\in \mathbb R$.
In turn, since neither $\Psi_{-1}(x, t, p)$ nor $\Psi_{+2}(x, t, p)$ can be identically zero due to the normalization \eqref{ds-infcond},  
we infer that  there exists a constant $C_p\neq0$ such that
\sss{\eee{\label{ds-sym1}
\Psi_{+2}(x, t, p)
=
C_p \Psi_{-1}(x, t, p)
\quad  \forall x, t\in \mathbb R.
}
Similarly, evaluating \eqref{ds-a-ab-def} at $k=\bar p$  we obtain
\eee{\label{ds-sym2}
\Psi_{+1}(x, t, \bar p)
=
C_{\bar p} \Psi_{-2}(x, t, \bar p)
\quad  \forall x, t\in \mathbb R
}}
for some constant $C_{\bar p}  \neq 0$. 
Thus, 
\sss{\label{ds-cp-cpb-def}
\eee{\label{ds-cp-def}
\underset{k= p}{\text{Res}} \, \frac{\Psi_{+2}(x, t, k)}{a(k) d(k)} 
=
\frac{\Psi_{+2}(x, t, p)}{a'(p) d(p)}
=
c_p   \Psi_{-1}(x, t, p)
\quad 
\forall x, t \in\mathbb R,
\quad
c_p := \frac{C_p}{a'(p) d(p)}.
}
\eee{ \label{ds-cpb-def}
\underset{k=\bar p}{\text{Res}} \,
 \frac{\Psi_{+1}(x, t, k)}{\bar a(k) d(k)} 
=
\frac{\Psi_{+1}(x, t, \bar p)}{\bar a'(\bar p) d(\bar p)}
=
c_{\bar p}   \Psi_{-2}(x, t, \bar p)
\quad \forall x, t\in\mathbb R,
\quad
c_{\bar p} := \frac{C_{\bar p}}{\bar a'(\bar p) d(\bar p)}.
}
} 
Relations \eqref{ds-cp-def} and \eqref{ds-cpb-def} imply the following  residue conditions for  $M$:
\sss{\label{ds-res-cond}
\ddd{
&\underset{k=p}{\text{Res}} \, M(x, t, k) 
=
M(x, t, p)
\left(
\begin{array}{lr}
0 & c_p\, e^{2i\theta(\xi, p) t} 
\\
0 & 0
\end{array}
\right)
&&
\forall x, t \in\mathbb R,
\label{ds-res-p-cond}
\\
&\underset{k=\bar p}{\text{Res}} \, M(x, t, k) 
=
M(x, t, \bar p)
\left(
\begin{array}{lr}
0 & 0
\\
c_{\bar p}\, e^{-2i\theta(\xi, \bar p) t}  & 0
\end{array}
\right)
\quad
&&\forall x, t \in\mathbb R.
\label{ds-res-pb-cond}
}
}
The Riemann-Hilbert problem for the focusing NLS IVP \eqref{ds-fnls-ivp} in the presence of the discrete spectrum $\{p, \bar p\}$ 
comprises the empty-discrete-spectrum problem \eqref{ds-rhp-cpam} augmented with the residue conditions \eqref{ds-res-cond}. 
To ensure uniqueness of solutions of the above Riemann-Hilbert problem, one must also supplement it
with suitable growth conditions at the branch points \cite{bm2019}.

The $x$-part of the Lax pair \eqref{ds-fnls-lp} together with the definition \eqref{ds-mrh} and the asymptotic condition \eqref{ds-n-rhp-as-cond} yield the  solution of the IVP \eqref{ds-fnls-ivp} via the reconstruction formula
\eee{ \label{ds-q-recon-n}
q(x, t) = -2i\lim_{k \to \infty}  k M_{12}(x, t, k).
}

For the  purpose of  computing the long-time asymptotics, it is convenient to convert the residue conditions \eqref{ds-res-cond} into jump discontinuities. In particular, following \cite{m2008}, we  let $\p D_p^\ve $ and $\p D_{\bar p}^\ve$ be the positively oriented boundaries of the disks $D_p^\ve $ and $D_{\bar p}^\ve$ of radius $\ve$ centered at $p$ and $\bar p$ respectively, and define the function $N$ by
\eee{\label{ds-nrh}
N(x, t, k)
=
\begin{cases}
M(x, t, k) V_p(x, t, k), &k\in D_{p}^{\ve },
\\
M(x, t, k),
&k\in \mathbb C^-\setminus \left( B^-\cup \overline{D_{p}^{\ve }}\,\right),
\\
M(x, t, k) V_{\bar p}(x, t, k),
&k\in D_{\bar p}^{\ve },
\\
M(x, t, k),
&k\in \mathbb C^+\setminus \left( B^+\cup \overline{D_{\bar p}^{\ve }}\, \right),
\end{cases}
}
where the matrices $V_p$ and $V_{\bar p}$ are given by
\eee{\label{ds-vp-vpb-def}
V_p(x, t, k)
=
\left(
\def\arraystretch{1}
\begin{array}{lr}
1 & -\dfrac{c_p}{k-p}\, e^{2i\theta(\xi, p)t}
\\
0 & 1
\end{array}
\right),
\quad
V_{\bar p}(x, t, k)
=
\left(
\def\arraystretch{1}
\begin{array}{lr}
1 &0
\\
-\dfrac{c_{\bar p}}{k-\bar p}\, e^{-2i\theta(\xi, \bar p)t} & 1
\end{array}
\right).
}
Note that the residue conditions \eqref{ds-res-cond} imply that $N$ is analytic at $p$ and $\bar p$. 
Furthermore, the jumps of $N$ along the continuous spectrum $\Sigma$ are the same with those of $M$ since $N = M$ outside the disks $D_p^\ve$ and $D_{\bar p}^\ve$. 
Therefore,  $N(x, t, k)$ is analytic for $k\in \mathbb C\setminus \left( \Sigma\cup  \p D_p^\ve \cup \p D_{\bar p}^\ve\right)$ and satisfies the following Riemann-Hilbert problem:
\sss{\label{ds-n-rhp-intro}
\ddd{
&N^+ = N^- V_1,   &&k\in \mathbb R,
\\
&N^+ = N^-V_2, &&k\in B^+,
\\
&N^+ = N^-V_3, &&k\in B^-,
\\
&N^+ = N^-V_p,  &&k\in \p D_p^\ve,
\label{ds-n-rhp-intro-dp}
\\
&N^+ = N^-V_{\bar p}, &&k\in \p D_{\bar p}^\ve,
\label{ds-n-rhp-intro-dpb}
\\
&N = I +O\left(1/k\right),\qquad && k \to \infty,
\label{ds-n-rhp-as-cond}
}
}
with the jumps $V_1, V_2, V_3$ given by \eqref{ds-rhp-cpam-j} and the jumps $V_p, V_{\bar p}$ defined by \eqref{ds-vp-vpb-def}.
Note that the transformation \eqref{ds-nrh} does not affect the normalization as $k\to\infty$. Thus, the long-time asymptotic behavior of the solution $q$ of  the IVP \eqref{ds-fnls-ivp} for the focusing NLS equation can equivalently be obtained by determining the corresponding behavior of the solution $N$ of the Riemann-Hilbert problem \eqref{ds-n-rhp-intro}.

\vskip 3mm
\noindent
\textbf{Overview of the long-time asymptotic analysis.}
The time dependence of the jumps of Riemann-Hilbert problem \eqref{ds-n-rhp-intro} is dictated by the exponentials $e^{\pm i\theta t}$, which become highly oscillatory in the limit $t\to\infty$. Thus, a delicate analysis via the nonlinear steepest descent method of Deift and Zhou \cite{dz1993, dz1995} is required in order to extract the leading-order asymptotic contribution to the solution.
Like in the classical steepest descent method, the main idea behind the Deift-Zhou method is to deform the contours associated with the oscillatory jumps to appropriate regions of the complex $k$-plane where the exponentials $e^{\pm i \theta t}$ decay to zero as $t\to\infty$.
Hence, the first step in the asymptotic analysis of  problem \eqref{ds-n-rhp-intro}  consists of studying the sign structure of $\Re (i\theta)$ 
in the complex $k$-plane. Recall, however, that the controlling phase function $\theta$ depends parametrically on the similarity variable $\xi$. Thus, similarly to the use of the steepest descent method for  computing the long-time asymptotics of solutions of linear equations (e.g., see \cite{AS1981, Whitham1974}), the analysis begins by studying how the sign structure of $\Re(i\theta)$ changes as $\xi$ increases from $-\infty$ to $\infty$.

\begin{figure}[t!]
\begin{center}
\includegraphics[scale=1]{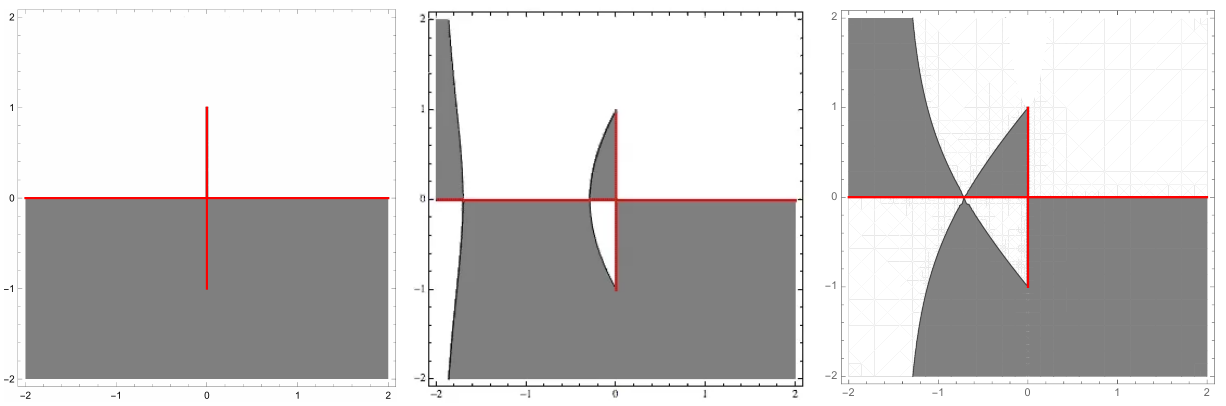}
\includegraphics[scale=1]{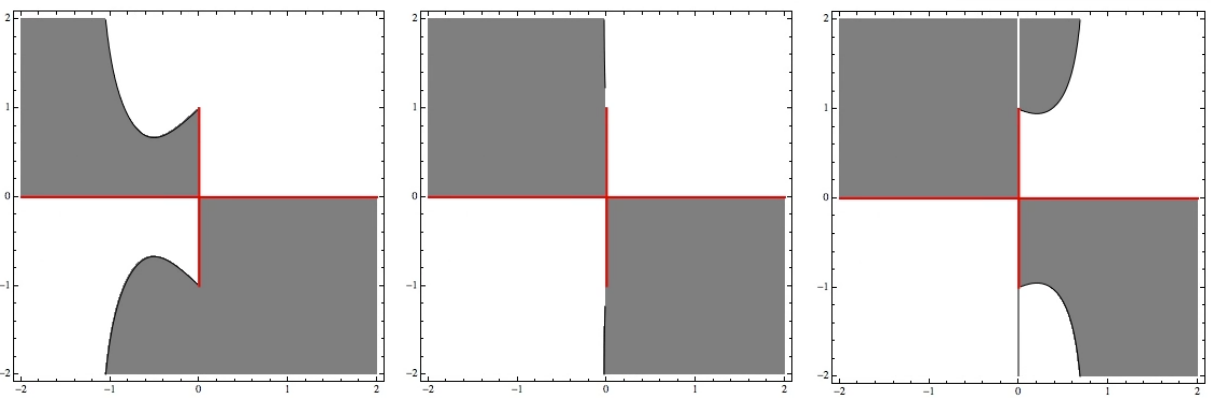}
\includegraphics[scale=1]{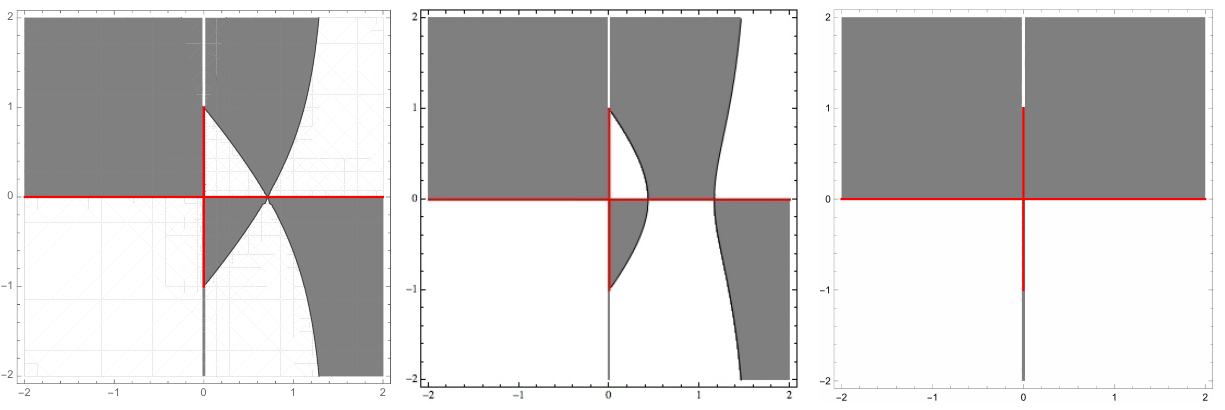}
\caption{
The sign of $\Re (i\theta)$ as $\xi$ increases from $-\infty$ to $+\infty$. %
\textit{Gray}: $\Re (i\theta)<0$; \textit{White}: $\Re (i\theta)>0$. 
\textit{Top row}: $\xi=-\infty$, $\xi \in (-\infty, \Vo)$, $\xi = \Vo$;
\textit{Middle row}: $\xi \in (\Vo, 0)$, $\xi = 0$, $\xi \in (0, -\Vo)$;
\textit{Bottom row}: $\xi = -\Vo$, $\xi = (-\Vo, \infty)$, and $\xi = +\infty$. 
\label{sign-structure-f}
}
\end{center}
\end{figure}

Let us first focus on the sign structure of $\Re (i\theta)$ for $\xi<0$ (i.e. $x<0$),  which is depicted in the first four frames of Figure~\ref{sign-structure-f}. Observe that, as $\xi$ increases from $-\infty$ to $0$, the sign of $\Re (i\theta)$ switches from negative to positive in the third quadrant and from positive to negative in the second quadrant, while it remains the same in the first and the fourth quadrant. More specifically, two  regions of positive sign emerge in the third quadrant: an unbounded region on the left of the point $k_1$, and a bounded region on the right of the point $k_2$ and on the left of the branch cut $B$, where
\eee{\label{ds-k1-k2-def}
k_1(\xi) := \frac18 \left(\xi - \sqrt{\xi^2 - \Vo^2}\right),
\quad
k_2(\xi) := \frac18 \left(\xi + \sqrt{\xi^2 - \Vo^2}\right)
}
are the two stationary points of $\theta$ with $\Vo$  defined by \eqref{ds-xistar-xip-def}.
 The two regions of positive sign grow continuously and remain disjoint until $\xi=\Vo$  (third frame in Figure~\ref{sign-structure-f}) where  $k_1(\Vo)=k_2(\Vo)=\frac{\Vo}{8}$. 
 Note that for $\xi\leqslant v_o$ the stationary points $k_1, k_2$ are real. 
 Subsequently, however, for $\Vo < \xi < 0$, the two stationary points become complex, and the two regions of positive sign merge to a \textit{single} region that eventually grows to occupy all of the third quadrant (fourth and fifth frames in Figure~\ref{sign-structure-f}). 
As a result, the value $\xi=\Vo$ is a bifurcation point in the analysis of the problem via the Deift-Zhou method. 

\begin{figure}[t!]
\begin{center}
\includegraphics[scale=0.325]{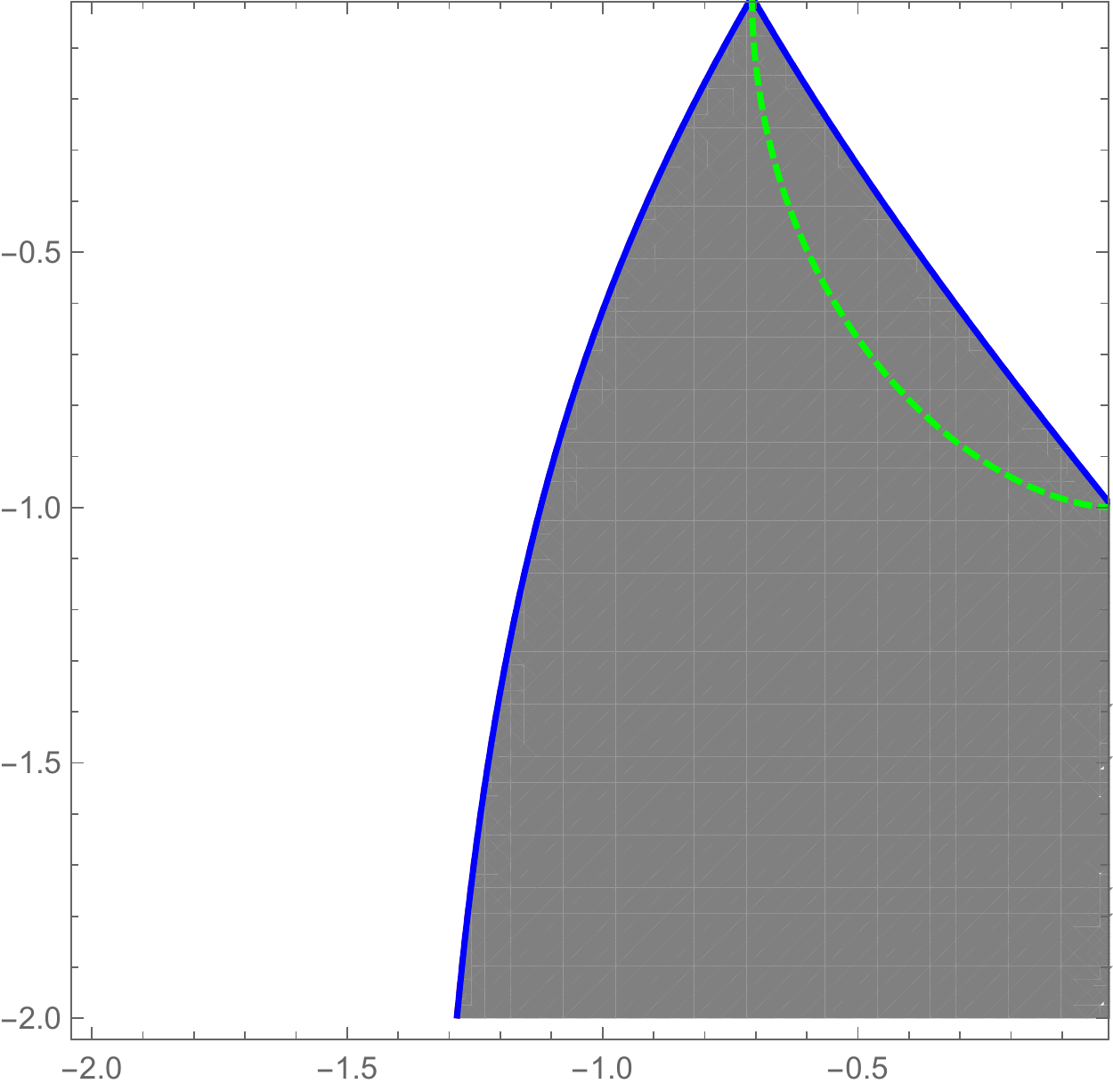}~%
\includegraphics[scale=0.325]{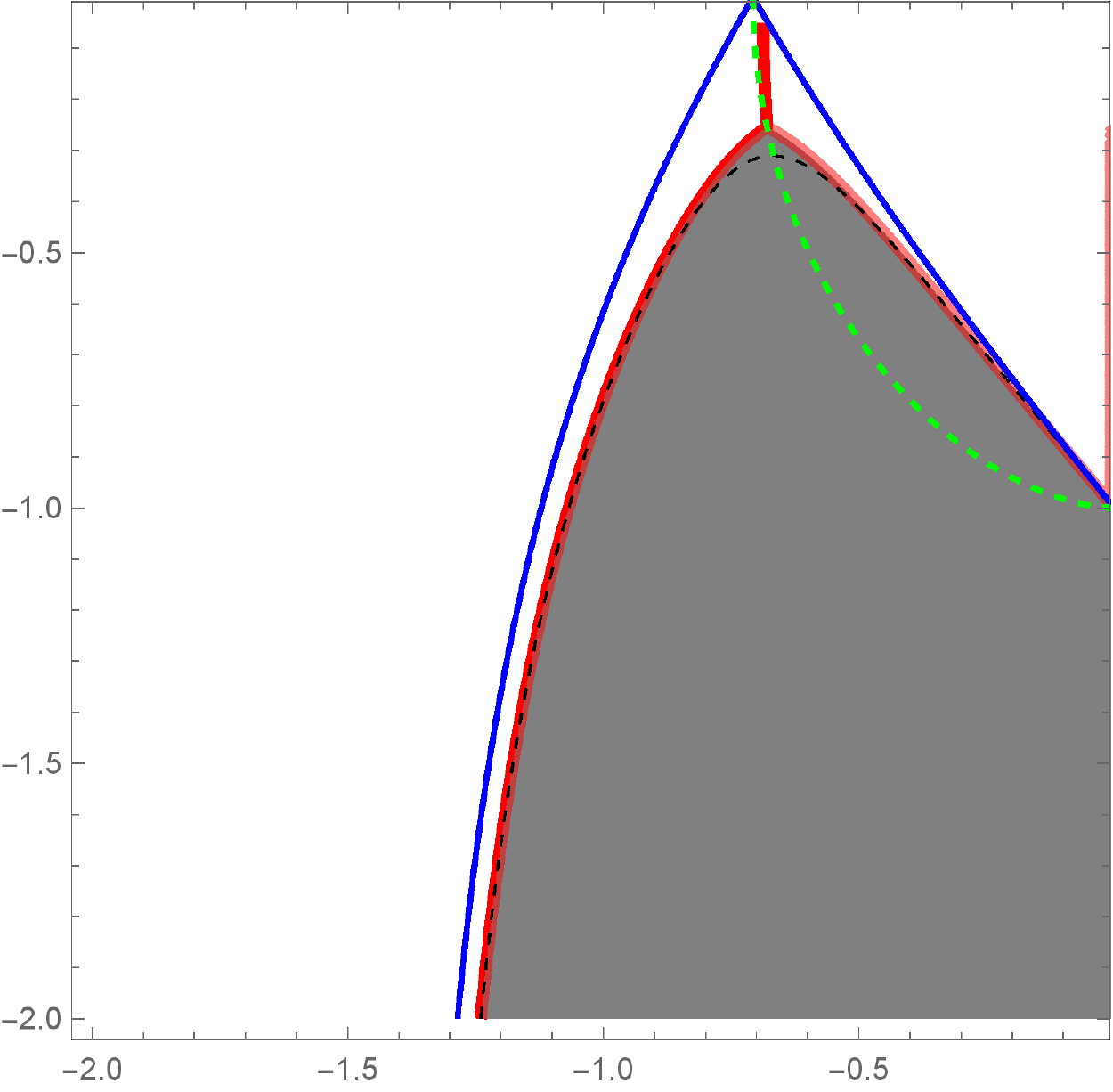}~%
\includegraphics[scale=0.325]{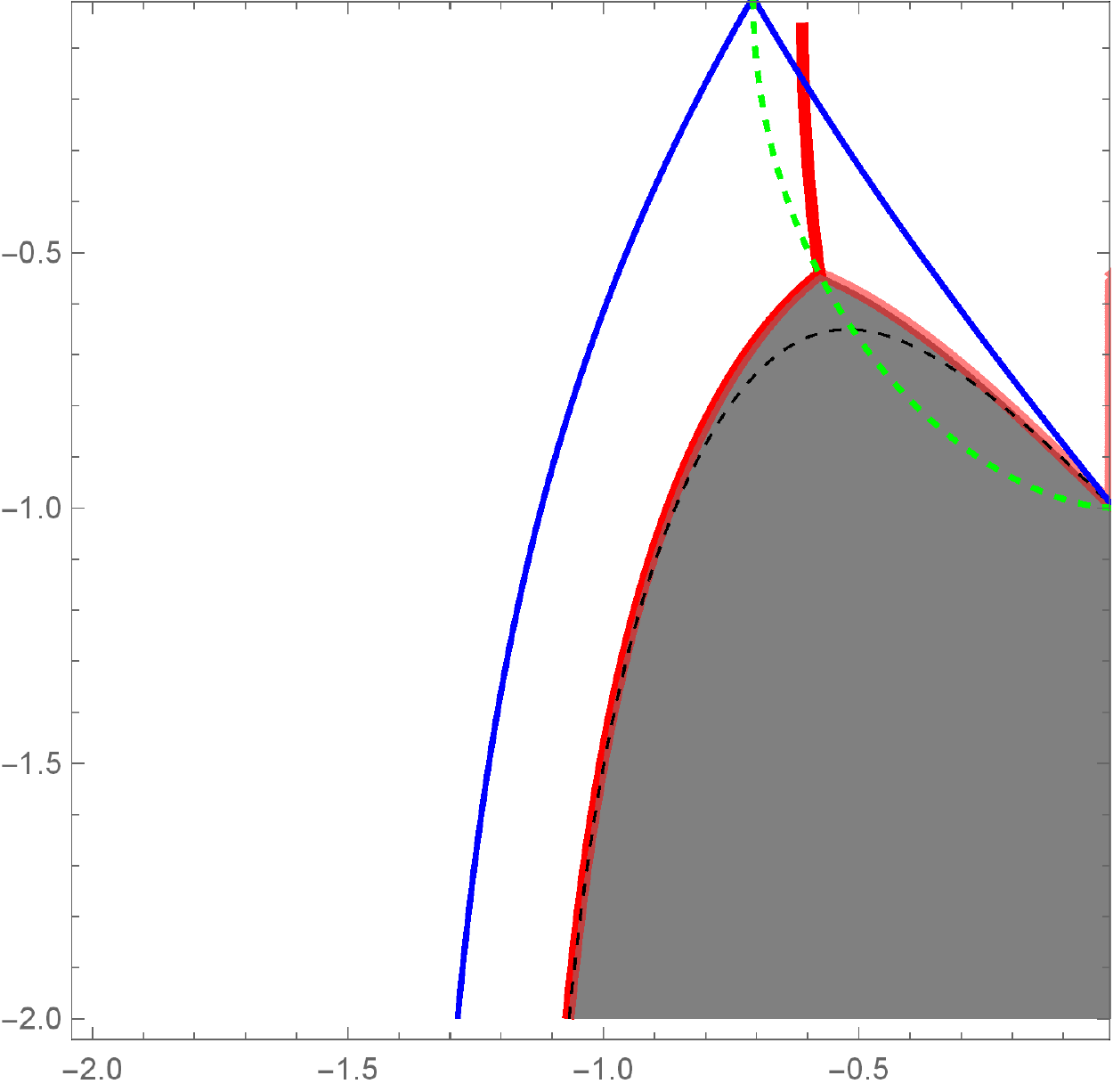}\\[1ex]
\includegraphics[scale=0.325]{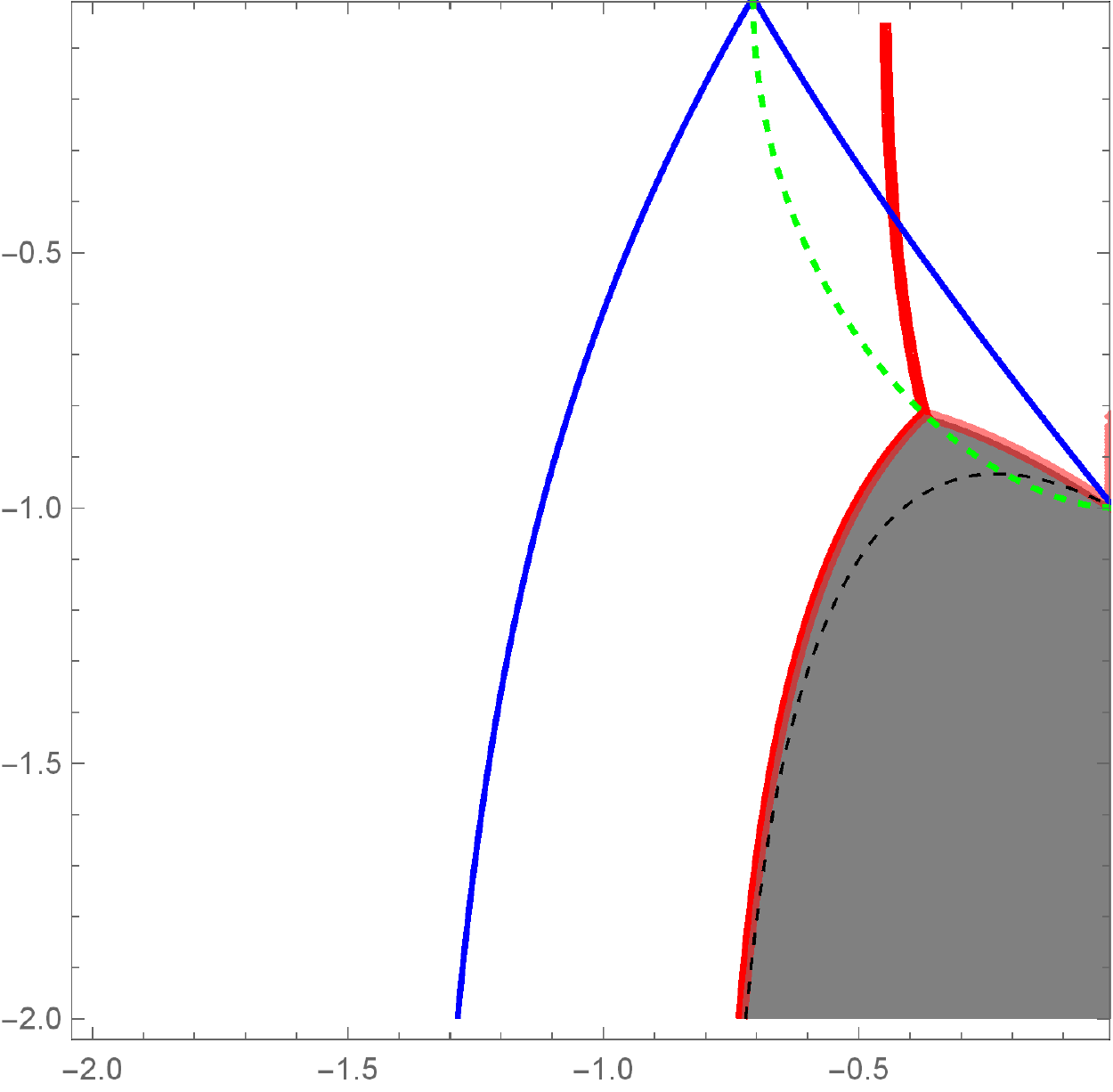}~%
\includegraphics[scale=0.325]{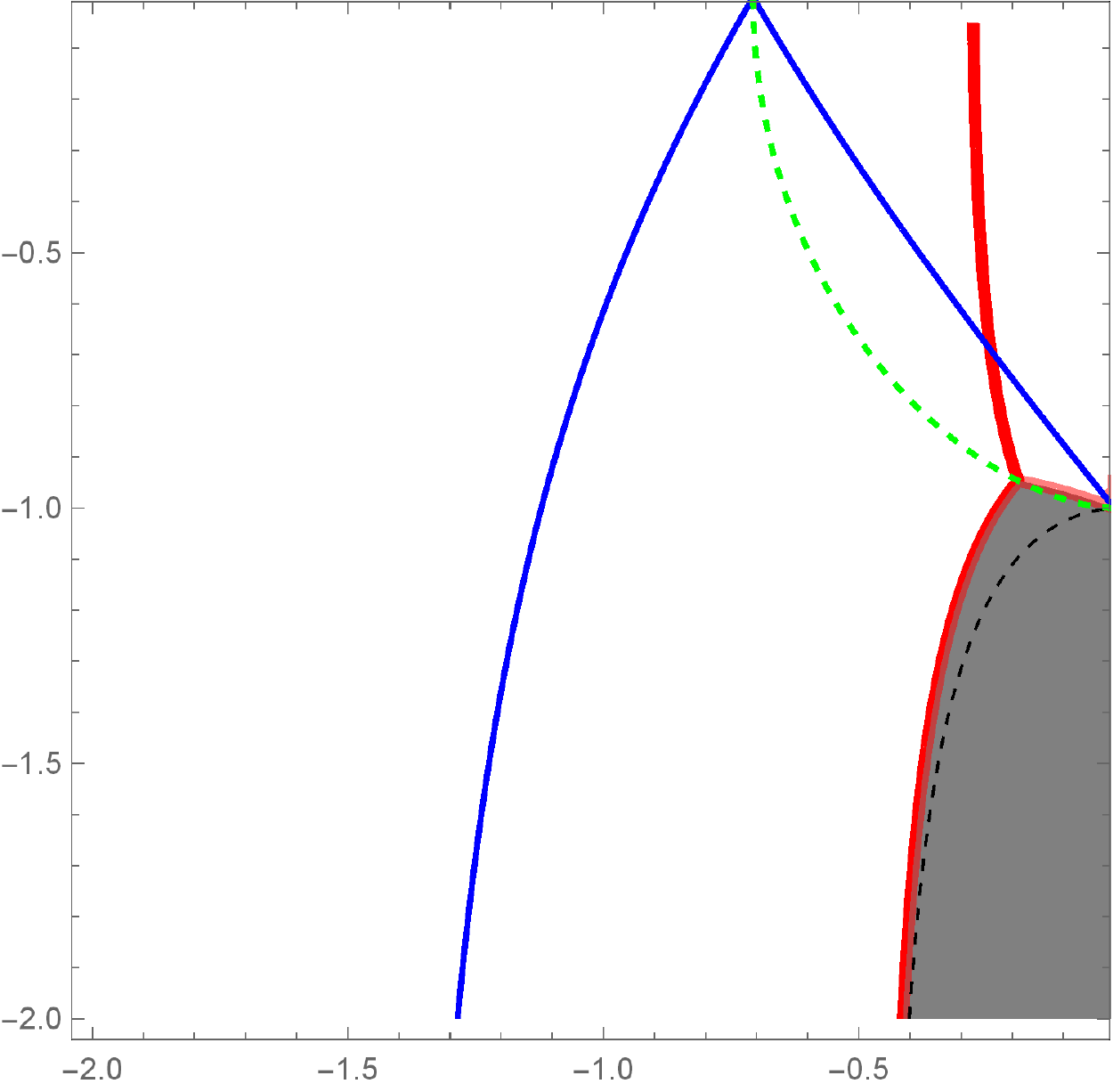}~%
\includegraphics[scale=0.325]{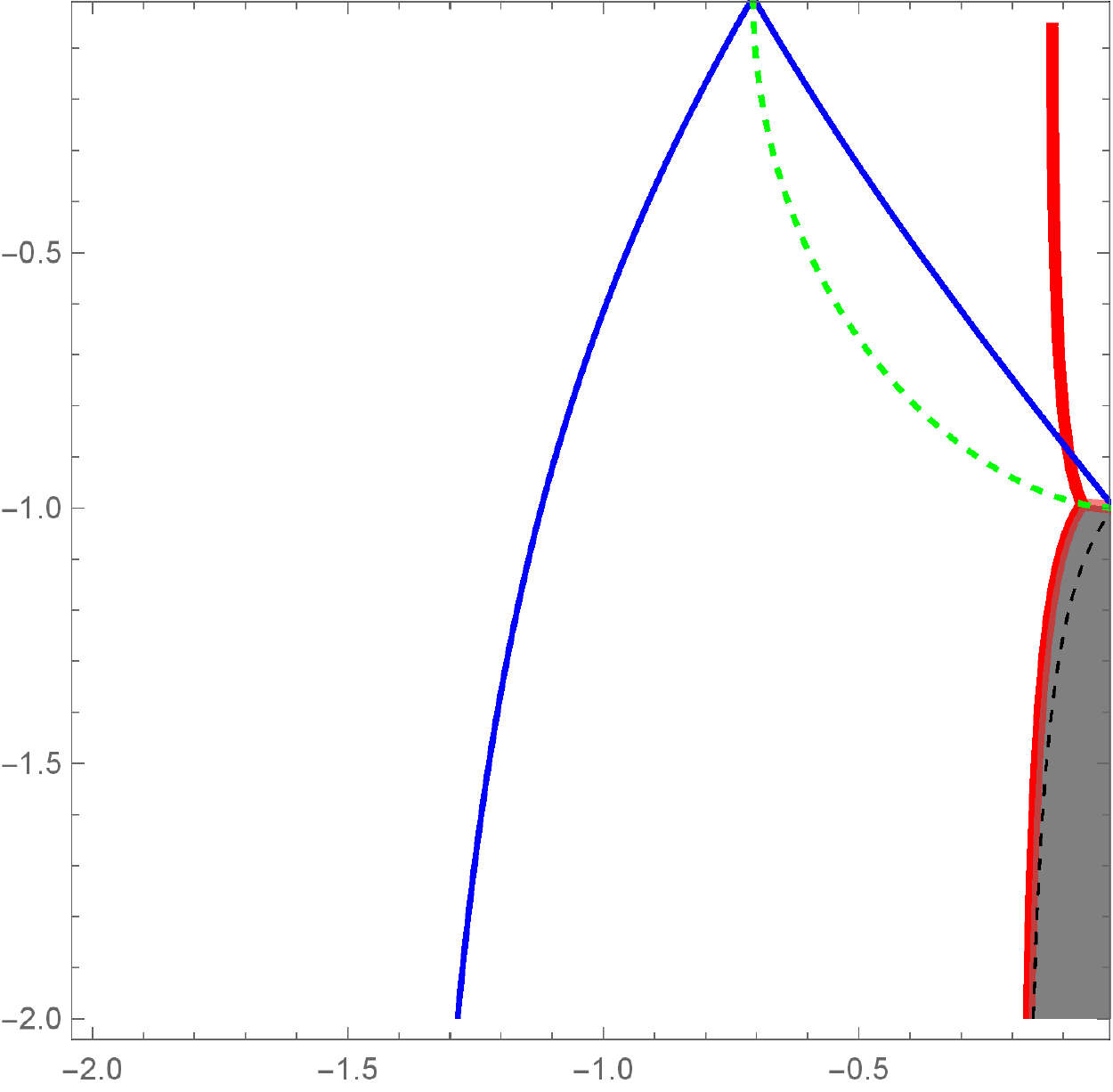}
\caption{
The sign of $\Re (ih)$ as $\xi$ increases gradually from $\Vo$ (top left plot) to just below $\xi=0$ (bottom right plot). %
\textit{Gray}: $\Re (ih)<0$; \textit{White}: $\Re (ih)>0$. 
Note $h\equiv \theta$ when $\xi = \Vo$.
\textit{Red curves}: The contours $\Im(h)=0$.
\textit{Green curve}: The path described by $\bar \alpha$ (defined by \eqref{ds-mod-eqs-i}), the intersection point of the contours $\Im(h)=0$.
\textit{Black curve}: The contour $\Im(\theta) = 0$. 
\textit{Blue curve}: The contour $\text{Im}\big[\theta(v_o, k)\big]=0$.
}
\label{f:h}
\end{center}
\end{figure}

More specifically, for $\xi<\Vo$,
the sign structure of $\text{Re}(i\theta)$ allows for two different factorizations of the jump $V_1$ along the real axis, 
both of which result in exponentially decaying contributions,
as will be explained in detail in Sections \ref{ds-esc-s} and \ref{ds-trap-s}. 
On the other hand, for $\Vo < \xi < 0$ it turns out that only one of the aforementioned factorizations can be employed. This results in an exponentially growing jump along a certain portion of the deformed jump contour, which is corrected by introducing a so-called $g$-function \cite{dvz1994, dvz1997} (see also Chapter 4 of \cite{kmm2003} in the context of semiclassical analysis).  
The corresponding transformation of the Riemann-Hilbert problem replaces the original
controlling phase function from $\theta$ to the Abelian integral $h$ defined by \eqref{ds-h-abel-i},
and is the reason why the asymptotics change dramatically as $\xi$ crosses $\Vo$.
The sign structure of $\Re(ih)$ in the third quadrant of the complex $k$-plane as $\xi$ decreases from $\Vo$ to 0
is shown in Figure~\ref{f:h}.

Of course, apart from the jumps $V_1$, $V_2$, $V_3$ along the continuous spectrum, the Riemann-Hilbert problem \eqref{ds-n-rhp-intro} also  involves the jumps $V_p, V_{\bar p}$ originating from the poles $p, \bar p$. 
Thus, another crucial value of $\xi$ now emerges, namely the value $\Vs$ for which $\text{Re}(i\theta)$ vanishes at $p$ and $\bar p$. Observe that $\Vs$ is the same  for $p$ and $\bar p$, since  
$\Re \left[i\theta(\xi, p)\right]=0 \Leftrightarrow \Re \left[i\theta(\xi, \bar p)\right]=0$ due to the symmetry $\theta(\xi, \bar k) = \overline{\theta(\xi, k)}$. Solving either of these equations, we obtain $\Vs$ in the explicit form \eqref{ds-xistar-xip-def}. 
Note that in the third quadrant, where $p$ lies, we have $\lambda_\re, \lambda_\im \leqslant 0$, thus $\Vs<0$.

In the range $(-\infty, \Vo)$, we shall see that the jumps $V_p, V_{\bar p}$ (equivalently, the poles $p, \bar p$) contribute to the leading-order asymptotics only when $\xi=\Vs$, provided that $p$ is such that $v_s \in (-\infty, \Vo)$. 
On the other hand, as explained above, in the range $(\Vo, 0)$  the phase function $\theta$ is replaced by the Abelian integral $h$ defined by \eqref{ds-h-abel-i}. Thus, the role of $\Vs$ is now played by the solutions of the equation
\eee{\label{ds-xiph-def}
\Re \left[ih(\xi, p)\right]=0.
}
The complicated form of $h$ does not allow us to solve equation \eqref{ds-xiph-def}  explicitly. It turns out, however, that, depending on the  location of $p$ inside the third quadrant, equation  \eqref{ds-xiph-def}  has either zero, one or two solutions in the interval $(\Vo, 0)$.
More specifically, as already noted in Section \ref{overview-s}, the third quadrant is divided into the four regions $D_1$, $D_2^+$, $D_2^-$, $D_3$ of Figure~\ref{ds-regions-f}, where for $\xi\in (\Vo, 0)$ equation \eqref{ds-xiph-def}  has no solutions in $D_1$, a unique solution $\Vd$ in $D_2^+$,  two solutions $\Vd<\Vw$ in $D_2^-$, and  a unique solution $\Vw$ in $D_3$. 
The mathematical description of the long-time asymptotic regimes that arise in these four regions is given in Theorems \ref{ds-per-t}-\ref{ds-ewr-t}. Before proceeding to the proofs of these results, we give a brief outline of the way in which the asymptotics unravels in each regime.
\vskip 3mm
\noindent 
$p\in D_1$:  \textit{The transmission regime.}
In this case, $\Vs<\Vo$ and, furthermore, equation \eqref{ds-xiph-def} has \textit{no} solution in the interval $(\Vo, 0)$ --- in fact, $\text{Re}(ih)(\xi, p)>0$ for all $\Vo<\xi < 0$. 
For $\xi<\Vs$, the jumps $V_p, V_{\bar p}$ decay exponentially and hence do not yield leading-order contributions. Thus, the dominant component of Riemann-Hilbert problem \eqref{ds-n-rhp-intro} in the limit $t\to\infty$ involves only the jumps along the continuous spectrum $\Sigma$, giving rise to the plane wave \eqref{ds-qsol-pw-t}. 
At $\xi=\Vs$, the jumps $V_p, V_{\bar p}$ switch from exponentially decaying to purely oscillatory. Consequently, they are now part of the dominant Riemann-Hilbert problem, generating the soliton \eqref{ds-esc-q-sol-lim-pw1-t}. Observe that this soliton propagates with velocity $\Vs$ and, since $|\Vs|>|\Vo|$, it eventually escapes to infinity outside the wedge $|\xi|<|\Vo|$  of Figure~\ref{ds-bifurc-f}.
For $\Vs<\xi<\Vo$, the jumps $V_p, V_{\bar p}$ grow exponentially. Nevertheless, it turns out that this growth can be converted into  decay  via an appropriate transformation. Hence, similarly to the range $\xi<\Vs$, the leading-order asymptotic behavior does not depend on $V_p, V_{\bar p}$ and is characterized by the plane wave \eqref{ds-qsol-pw-t}, but now with a \textit{phase shift} generated by the soliton that has arisen at $\xi=\Vs$.
Finally,  for $\Vo<\xi < 0$ the phase function switches from $\theta$ to $h$. Then, since $\text{Re}(ih)(\xi, p)>0$ for all $\Vo<\xi < 0$,  the jumps $V_p, V_{\bar p}$ do not contribute to the leading-order asymptotics. Hence,  no soliton is present in the range $\Vo<\xi < 0$ and the solution is asymptotically equal to the modulated elliptic wave \eqref{ds-qsol-mew-t} with the phase shift already generated by the soliton at $\xi=\Vs$  in the range  $\Vs<\xi<\Vo$.  
\vskip 3mm
\noindent
$p\in D_2^+$: \textit{The trap regime}.
In this case, $\Vs>\Vo$ and, in addition, equation \eqref{ds-xiph-def} has a \textit{unique} solution $\Vd$ in the interval $(\Vo, 0)$ --- in fact, it turns out that $\Vs<\Vd$.
Thus, for $\xi<\Vd$ the jumps $V_p, V_{\bar p}$ are not significant at leading order. In particular,  for $\xi<\Vo$ the leading-order asymptotics is given by the plane wave \eqref{ds-qsol-pw-t}, while for $\Vo < \xi < \Vd$ the solution is asymptotically equal to the modulated elliptic wave \eqref{ds-qsol-mew-cpam-t}. 
At $\xi=\Vd$, however, the jumps $V_p, V_{\bar p}$ become purely oscillatory and hence do contribute to the leading-order asymptotics, which is now given by the soliton \eqref{ds-trap-q-sol-lim-pw1-t}. Observe that, since $|\Vd|< |\Vo|$, this soliton is \textit{trapped forever} inside the  wedge $|\xi|<|\Vo|$ of Figure~\ref{ds-bifurc-f}. Furthermore, the fact that $|\Vd|< |\Vs|$ indicates that the soliton is \textit{delayed} by its interaction with the modulated elliptic wave. 
Finally, since $\Vd$ is the only solution of equation \eqref{ds-xiph-def}  in $(\Vo, 0)$, for $\Vd < \xi < 0$ the jumps $V_p, V_{\bar p}$ do not affect the leading-order asymptotics, which is now equal to the modulated elliptic wave \eqref{ds-qsol-mew-t}  with an additional \textit{phase shift}   generated by the soliton  at $\xi=\Vd$. 
\vskip 3mm
\noindent
$p\in D_2^-$: \textit{The trap}/\textit{wake regime.}
This case is  similar to the  trap regime apart from the fact that now equation \eqref{ds-xiph-def} has \textit{two} (as opposed to one) solutions in the interval $(\Vo, 0)$, namely $\Vd$ and $\Vw$ with $\Vd<\Vw$. 
Therefore, for $\xi < \Vw$ the asymptotics is the same with the one in the  trap regime, including the soliton that arises at $\xi=\Vd$. However, at $\xi=\Vw$ a new phenomenon emerges, namely the  soliton wake   \eqref{ds-trapwake-t}.
Importantly, contrary to the soliton (which induces a phase shift for $\xi>\Vd$), the soliton wake does \textit{not} affect  the leading-order asymptotics in the range $\Vw<\xi<0$.
\vskip 3mm
\noindent
$p\in D_3$: \textit{The transmission}/\textit{wake regime.}
This case is similar to the  transmission regime apart from the fact that equation \eqref{ds-xiph-def} now has a \textit{unique} solution $\Vw$ in the interval $(\Vo, 0)$ (as opposed to no solution). Thus, the leading-order asymptotics is the same with the one in the transmission regime except for $\xi=\Vw$, where the soliton wake \eqref{ds-trapwake-t} arises. Importantly, contrary to the soliton at $\xi=\Vs$ (which generates a phase shift for $\xi>\Vs$), the leading-order asymptotics for $\Vw<\xi<0$  are not affected by the soliton wake.

%
%
%
%
\section{The  Transmission Regime: Proof of Theorem \ref{ds-per-t}}
\label{ds-esc-s}

This regime arises when $p$ lies in the region $D_1$ of Figure \ref{ds-regions-f}, in which  case  $\Vs<\Vo$ and $\text{Re}(ih)(\xi, p)$ does not vanish in the interval $(\Vo,0)$. Thus, we split the interval $(-\infty, 0)$ into the following ranges: 
$\xi< \Vs$; $\xi=\Vs$; $\Vs<\xi< \Vo$; and $\Vo<\xi< 0$.

\subsection{The range $\xi<\Vs$: plane wave}
\label{ds-esc-pw1-ss}
\begin{figure}[t!]
\centering
\includegraphics[scale=1]{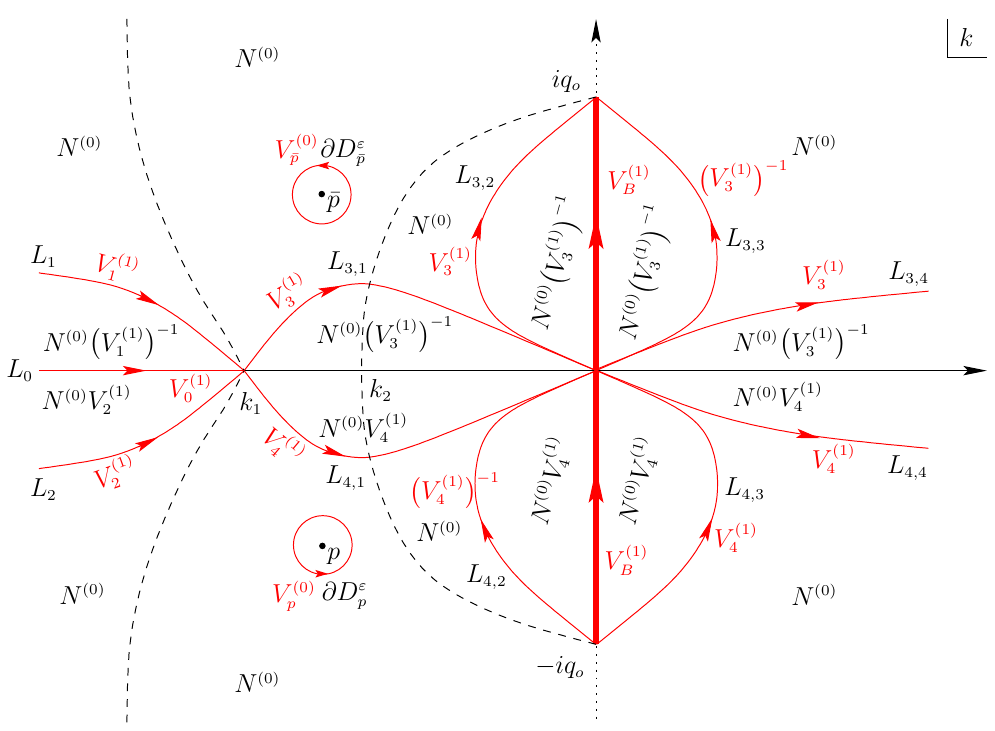}
\caption{Plane wave region in the  transmission regime: the first stage of the first deformation. 
A new function $N^{(1)}$ is defined in terms of the solution $N^{(0)}$ of Riemann-Hilbert problem \eqref{ds-n-rhp-intro} via different expressions in different regions of the complex $k$-plane. This allows us to eliminate the jump of \eqref{ds-n-rhp-intro} along $(k_1, \infty)$. The jumps $V_p = V_p^{(0)}$ and $V_{\bar p} = V_{\bar p}^{(0)}$ along the circles $\p D_p^\ve$ and $\p D_{\bar p}^\ve$ are not  affected at this stage.
\label{ds-esc-def1a-pw1-f}
}
\vspace{0mm}
\includegraphics[scale=1]{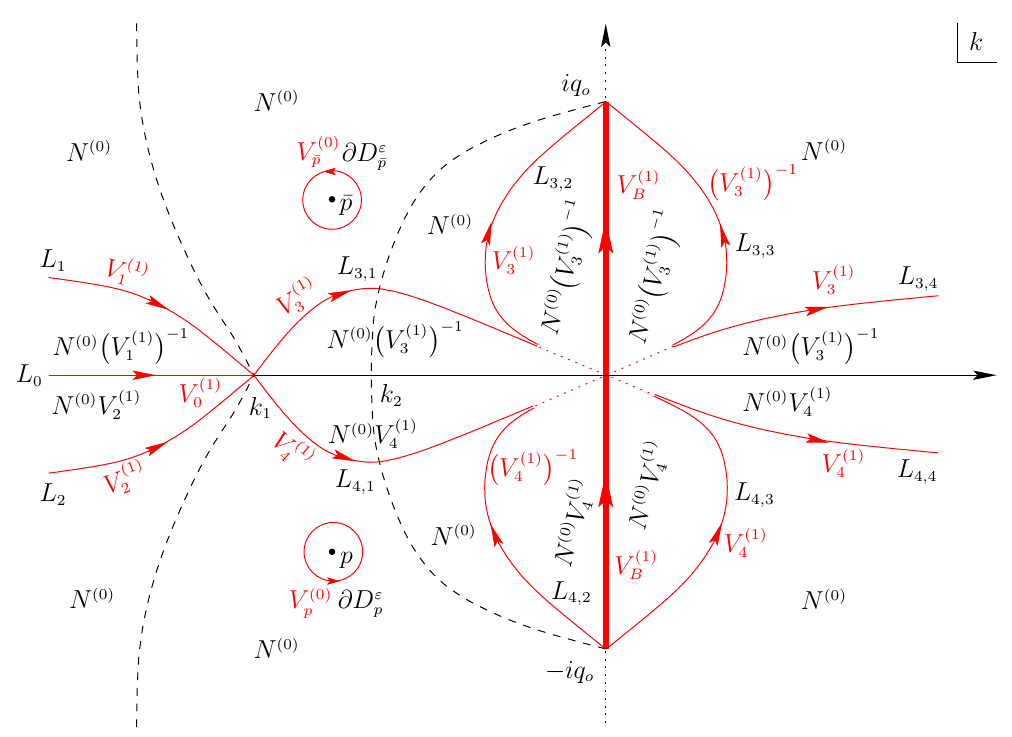}
\caption{Plane wave region in the  transmission regime: the second stage of the first deformation.
In the second quadrant, the function $N^{(1)}$ is defined in terms of $N^{(0)}$ by the same expression both below the contour $L_{3, 1}$ and to the right of the contour $L_{3, 2}$. Thus, $N^{(1)}$ does not have a jump along the overlapping portion between these two contours (dotted line), allowing one to lift them away from the origin. The same is true for the contour pairs $\{L_{3, 3}, L_{3,4}\}$, $\{L_{4,1}, L_{4, 2}\}$ and  $\{L_{4, 3}, L_{4, 4}\}$. The jumps along $\p D_p^\ve$ and $\p D_{\bar p}^\ve$ remain unchanged.
\label{ds-esc-def1c-pw1-f}
}
\end{figure}

\begin{figure}[t!]
\begin{center}
\includegraphics[scale=1]{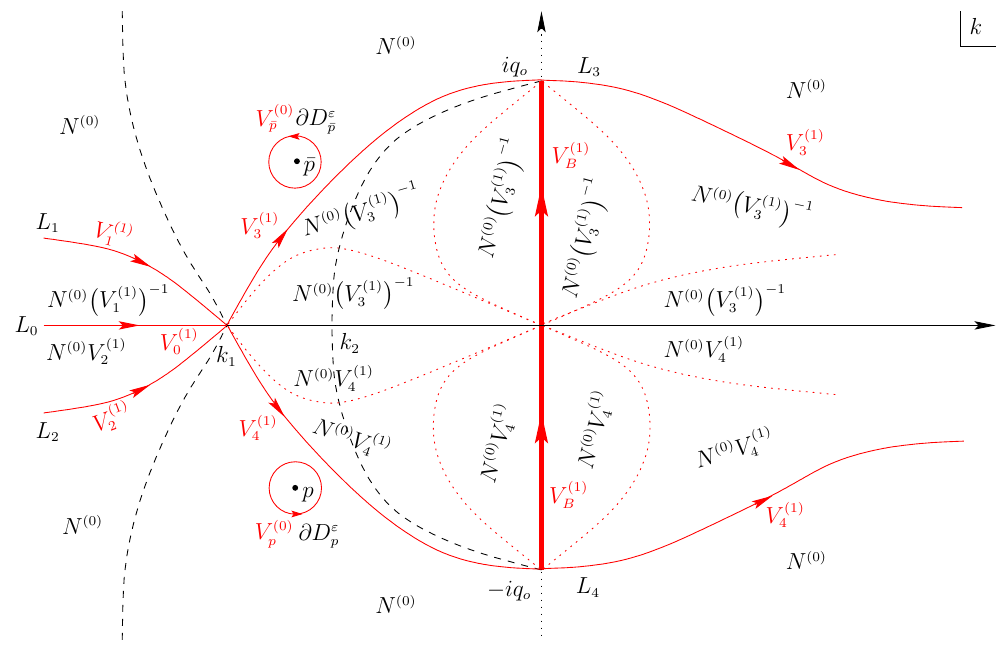}
\caption{Plane wave region in the  transmission regime:  the third stage of the first deformation.
Having lifted the jump contours away from the origin as shown in Figure \ref{ds-esc-def1c-pw1-f}, one can now adjust the definition of $N^{(1)}$ according to the present figure in order to move these jump contours outside the finite region defined by the branch cut $B$ and the dashed line through the stationary point $k_2$. This ensures that the relevant jumps occur along contours of  appropriate sign for $\text{Re}(i\theta)$. The jumps along $\p D_p^\ve$ and $\p D_{\bar p}^\ve$ are as before.
\label{ds-esc-def1d-pw1-f}
}
\vspace{7mm}
\includegraphics[scale=1]{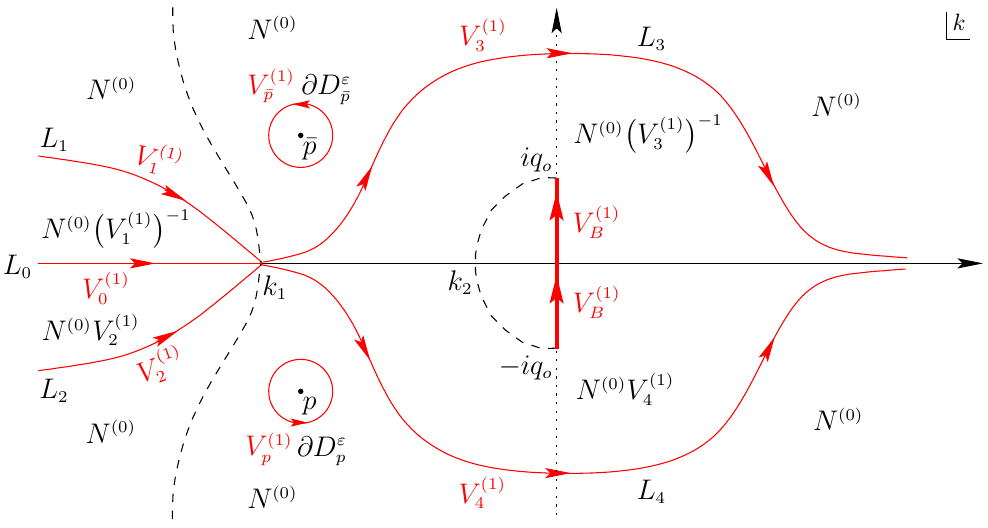}
\caption{Plane wave region in the  transmission regime: the fourth and final stage of the first deformation.
The jump contours $L_3$ and $L_4$ have been lifted away from the branch points $\pm iq_o$ similarly to \cite{bm2017}.
Overall, the jumps along $\p D_p^\ve$ and $\p D_{\bar p}^\ve$ have not changed in the transition from $N^{(0)}$ to $N^{(1)}$, i.e. $V_p^{(1)} = V_p^{(0)}$ and $V_{\bar p}^{(1)} = V_{\bar p}^{(0)}$. 
\label{ds-esc-def1b-pw1-f}
}
\end{center}
\end{figure}

In this range, we have $\Re (i\theta)(\xi, p)<0$ and $\Re (i\theta)(\xi, \bar p)>0$. Hence, the jumps $V_p$ and $V_{\bar p}$ given by \eqref{ds-vp-vpb-def} tend to the identity  as $t\to \infty$ and therefore are not expected to be part of the dominant component of Riemann-Hilbert problem \eqref{ds-n-rhp-intro} in the limit $t\to\infty$. 
Next, we shall show that this is indeed the case by performing several deformations of problem \eqref{ds-n-rhp-intro} in the spirit of the Deift-Zhou nonlinear steepest descent method. We emphasize that although some of these deformations are similar to those of the no-discrete-spectrum analysis of \cite{bm2017},  one now needs to carefully handle the jumps around the poles $p$, $\bar p$, which were not present in \cite{bm2017}.
\vskip 2mm
\noindent\textbf{First deformation.} 
This deformation is carried out in four stages. In each of these stages, a new function $N^{(1)}$ is defined in terms of the solution $N=N^{(0)}$ of Riemann-Hilbert problem \eqref{ds-n-rhp-intro}, as shown in Figures  \ref{ds-esc-def1a-pw1-f}-\ref{ds-esc-def1b-pw1-f}. Importantly, the jumps $V_p = V_p^{(0)}$ and $V_{\bar p} = V_{\bar p}^{(0)}$ are not affected by this deformation. In its final form,  the function $N^{(1)}$ is analytic in $\mathbb C\setminus \big(\bigcup_{j=0}^4 L_j \cup B\cup \p D_p^\ve \cup \p D_{\bar p}^\ve\big)$, satisfies the asymptotic condition
\eee{
N^{(1)} = I +O\left(\frac 1k\right), \quad  k \to \infty,
}
and possesses the following jump discontinuities along the contours $\bigcup_{j=0}^4 L_j \cup B\cup \p D_p^\ve \cup \p D_{\bar p}^\ve$, as shown in Figure \ref{ds-esc-def1b-pw1-f}:
\ddd{\label{ds-n1-jumps}
&V_B^{(1)}
=
V_B
=
\left(
\def\arraystretch{1.2}
\begin{array}{lr}
0 & \tfrac{q_-}{iq_o } \\
\tfrac{\bar q_-}{iq_o } & 0
\end{array}
\right), 
\quad 
V_0^{(1)}
=
\left(
\def\arraystretch{1.2}
\begin{array}{lr}
1+r\bar r & 0 \\
0 & \dfrac{1}{1+r\bar r}
\end{array}
\right),
\quad
V_1^{(1)}
=
\left(
\def\arraystretch{1.2}
\begin{array}{lr}
d^{-\frac 12} & \dfrac{d^{\frac 12}\, \bar r e^{2i \theta t } }{1+r\bar r} \\
0 & d^{\frac 12}
\end{array}
\right),
\nn\\
&V_2^{(1)}
=
\left(
\def\arraystretch{1.2}
\begin{array}{lr}
d^{-\frac 12} &  0 \\
\dfrac{d^{\frac 12} r e^{-2i \theta t } }{1+r\bar r} & d^{\frac 12}
\end{array}
\right),
\quad
V_3^{(1)}
=
\left(
\def\arraystretch{1.2}
\begin{array}{lr}
d^{-\frac 12} & 0 \\
d^{-\frac 12}r e^{-2i \theta t } & d^{\frac 12}
\end{array}
\right),
\quad
V_4^{(1)}
=
\left(
\def\arraystretch{1.2}
\begin{array}{lr}
d^{-\frac 12}
&
d^{-\frac 12}\bar r e^{2i \theta t } 
\\
0
&
d^{\frac 12}
\end{array}
\right),
\nn\\
&\hskip 5.3cm
V_p^{(1)} = V_p^{(0)}, \quad V_{\bar p}^{(1)} = V_{\bar p}^{(0)}.
}

\begin{figure}[t!]
\begin{center}
\includegraphics[scale=1]{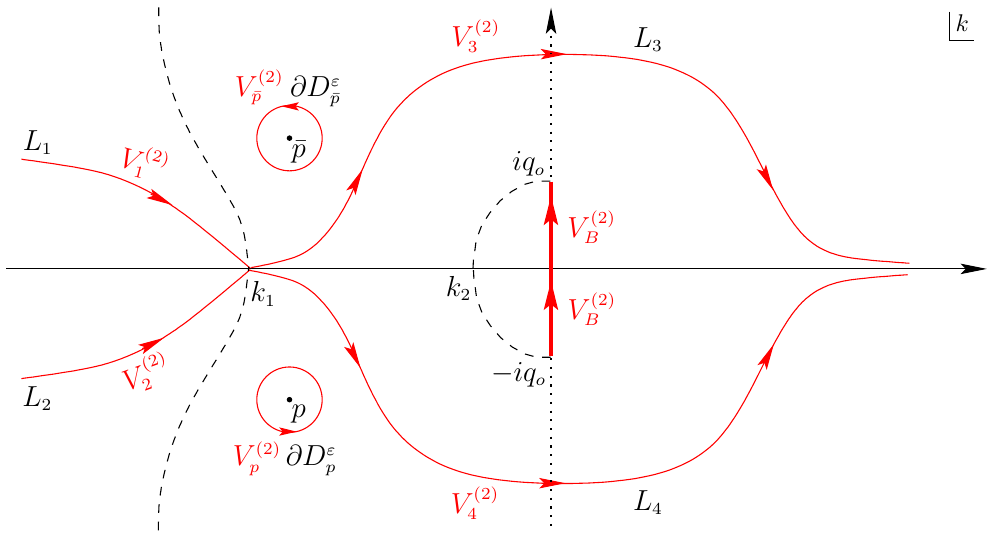}
\caption{Plane wave region in the  transmission regime: the second deformation.}
\label{ds-esc-def2-pw1-f}
\end{center}
\bigskip
\begin{center}
\includegraphics[scale=1]{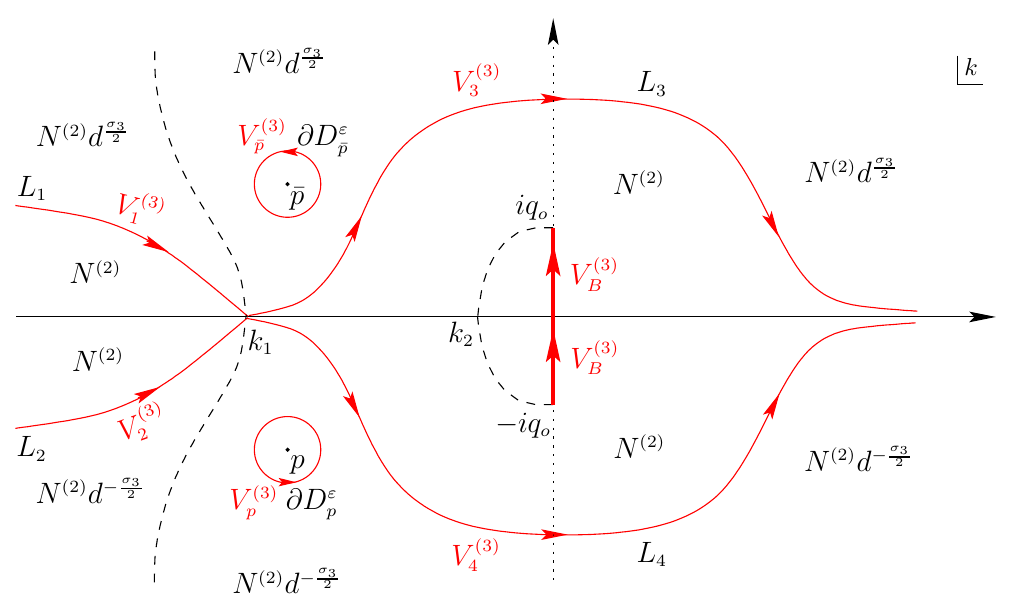}
\caption{Plane wave region in the  transmission regime: the third deformation.}
\label{ds-esc-def3-pw1-f}
\end{center}
\end{figure}

\vskip 3mm
\noindent\textbf{Second deformation.} 
The jump $V_0^{(1)}$  along the contour $L_0 := (-\infty, k_1)$ shown in Figure \ref{ds-esc-def1b-pw1-f} can be removed by means of the transformation
\eee{\label{ds-n2-def-pw1}
N^{(2)}(x, t, k) = N^{(1)}(x, t, k) \delta(\xi, k)^{-\sigma_3}, 
}
where the scalar function $\delta(\xi, k)$ is analytic in $\mathbb C\setminus (-\infty, k_1)$ and satisfies the Riemann-Hilbert problem
\sss{\label{ds-rhpdelta}
\ddd{
\delta^+ &= \delta^- \left(1+r \bar r\right),\quad &&k\in (-\infty, k_1),\label{ds-deljump}
\\
\delta &=1+O\left(\tfrac 1k \right), &&k \to \infty.
}
}
In fact, problem \eqref{ds-rhpdelta} can be solved explicitly via the Plemelj formulae to yield
\eee{\label{ds-delta-def}
\delta(\xi, k) = \exp\left\{\frac{1}{2i\pi}\int_{-\infty}^{k_1(\xi)} \frac{\ln\left[1+r(\nu)\bar r(\nu)\right]}{\nu-k} \, d\nu\right\},
\quad k\notin (-\infty, k_1).
}   
Through transformation \eqref{ds-n2-def-pw1}, the jumps of $N^{(1)}$  give rise to corresponding jumps for $N^{(2)}$. As shown in Figure \ref{ds-esc-def2-pw1-f}, these jumps occur along the contours $\bigcup_{j=1}^4 L_j \cup B$ and are given by
\ddd{\label{ds-n2-jumps}
&V_B^{(2)}
=
\left(
\def\arraystretch{1.2}
\begin{array}{lr}
0 & \tfrac{q_-}{iq_o} \delta^2 \\
\tfrac{\bar q_-}{iq_o} \delta^{-2} & 0
\end{array}
\right), 
\quad
V_1^{(2)}
=
\left(
\def\arraystretch{1.2}
\begin{array}{lr}
d^{-\frac 12} & \dfrac{d^{\frac 12}\, \bar r e^{2i \theta t } }{1+r\bar r} \delta^2 \\
0 & d^{\frac 12}
\end{array}
\right),
\quad
V_2^{(2)}
=
\left(
\def\arraystretch{1.2}
\begin{array}{lr}
d^{-\frac 12} &  0 \\
\dfrac{d^{\frac 12} r e^{-2i \theta t}}{1+r\bar r} \delta^{-2} & d^{\frac 12}
\end{array}
\right),
\nn\\
&\hskip 2cm
V_3^{(2)}
=
\left(
\def\arraystretch{1.2}
\begin{array}{lr}
d^{-\frac 12} & 0 \\
d^{-\frac 12}r e^{-2i \theta t} \delta^{-2} & d^{\frac 12}
\end{array}
\right),
\quad
V_4^{(2)}
=
\left(
\def\arraystretch{1.2}
\begin{array}{lr}
d^{-\frac 12}
&
d^{-\frac 12}\bar r e^{2i \theta t } \delta^2
\\
0
&
d^{\frac 12}
\end{array}
\right),
}
as well as along the disks $\p D_p^\ve$ and $\p D_{\bar p}^\ve$, where they read (modified for the first time)
\eee{
V_p^{(2)} 
= 
\left(
\def\arraystretch{1}
\begin{array}{lr}
1 & -\dfrac{c_p \, \delta^2(\xi, k)}{k-p}\,  e^{2i\theta(\xi, p)t}
\\
0 & 1
\end{array}
\right),
\quad
V_{\bar p}^{(2)} 
=\left(
\def\arraystretch{1}
\begin{array}{lr}
1 &0
\\
-\dfrac{c_{\bar p} \, \delta^{-2}(\xi, k)}{k-\bar p}\,  e^{-2i\theta(\xi, \bar p)t} & 1
\end{array}
\right).
}
Finally, the normalization condition \eqref{ds-delta-def} for $N^{(1)}$ is also satisfied by $N^{(2)}$.
 
\vskip 3mm
\noindent\textbf{Third deformation.}
The function $d(k)$ can be eliminated from the jump matrices along $\bigcup_{j=1}^4 L_j$ by introducing a new function $N^{(3)}$ defined in terms of $N^{(2)}$ according to Figure \ref{ds-esc-def3-pw1-f}. In particular, the jumps of $N^{(3)}$ along the contours $\bigcup_{j=1}^4 L_j \cup B$ are given by
\ddd{\label{ds-n3-jumps}
&V_B^{(3)}
=
\left(
\def\arraystretch{1.2}
\begin{array}{lr}
0 & \tfrac{q_-}{iq_o} \delta^2 \\
\tfrac{\bar q_-}{iq_o} \delta^{-2} & 0
\end{array}
\right), 
\quad
V_1^{(3)}
=
\left(
\def\arraystretch{1.2}
\begin{array}{lr}
1 & \dfrac{\bar r e^{2i \theta t } }{1+r\bar r} \delta^2 \\
0 & 1
\end{array}
\right),
\quad
V_2^{(3)}
=
\left(
\def\arraystretch{1.2}
\begin{array}{lr}
1 &  0 \\
\dfrac{r e^{-2i \theta t}}{1+r\bar r} \delta^{-2} & 1
\end{array}
\right),
\nn\\
&\hskip 2.7cm
V_3^{(3)}
=
\left(
\def\arraystretch{1.2}
\begin{array}{lr}
1 & 0 \\
r e^{-2i \theta t} \delta^{-2} & 1
\end{array}
\right),
\quad
V_4^{(3)}
=
\left(
\def\arraystretch{1.2}
\begin{array}{lr}
1
&
\bar r e^{2i \theta t} \delta^2
\\
0
&
1
\end{array}
\right).
}
Moreover, noting that $N^{(3)} = N^{(2)} d^{-\frac{\sigma_3}{2}}$ for  $k\in D_p^\ve$ and $N^{(3)} = N^{(2)} d^{\frac{\sigma_3}{2}}$ for $k\in D_{\bar p}^\ve$, we obtain
\eee{\label{ds-vp3-vpb3-def-pw1}
V_p^{(3)} 
=
\left(
\def\arraystretch{1}
\begin{array}{lr}
1 & -\dfrac{c_p \, \delta^2(\xi, k)\, d(k)}{k-p}\,  e^{2i\theta(\xi, p)t}
\\
0 & 1
\end{array}
\right),
\
V_{\bar p}^{(3)} 
=
\left(
\def\arraystretch{1}
\begin{array}{lr}
1 &0
\\
-\dfrac{c_{\bar p} \, \delta^{-2}(\xi, k) \, d(k)}{k-\bar p}\,  e^{-2i\theta(\xi, \bar p)t} & 1
\end{array}
\right).
}

\vskip 3mm
\noindent\textbf{Fourth deformation.}
Our final goal is to convert the jump along the branch cut $B$ into the \textit{constant} matrix $V_B$ given by \eqref{ds-n1-jumps}. This can be achieved by means of the global transformation 
\eee{\label{ds-esc-n4-def-pw1}
N^{(4)}(x, t, k) = N^{(3)}(x, t, k) e^{i g(\xi, k)\sigma_3},
}
where the function $g(\xi, k)$ is analytic in $\mathbb C\setminus B$ and satisfies the jump condition
\eee{\label{ds-esc-jump-g-pw1}
 e^{i(g^+ + g^-)} = \delta^2,\quad k\in B,
}
and the normalization condition
\eee{\label{ds-esc-g-rhp-pw1}
\frac{g}{\lambda}
=
O\left(\frac 1k\right), \quad k \to\infty.
}
Indeed, the jump condition \eqref{ds-esc-jump-g-pw1} implies that the jump of $N^{(4)}$ along $B$  is precisely $V_B$. Equations \eqref{ds-esc-jump-g-pw1} and \eqref{ds-esc-g-rhp-pw1} formulate a Riemann-Hilbert problem for $g$, which can be solved explicitly to yield 
\eee{ \label{ds-esc-g-def-pw1}
g(\xi, k)
=
\frac{\lambda(k)}{2i\pi^2}
\int_{\zeta\in B} \frac{1}{\lambda(\zeta)\left(\zeta-k\right)} \int_{-\infty}^{k_1(\xi)} \frac{\ln\left[1+r(\nu)\bar r(\nu)\right]}{\nu-\zeta}\, d\nu d\zeta, \quad k\notin B.
} 

Under transformation \eqref{ds-esc-n4-def-pw1}, the Riemann-Hilbert problem for $N^{(3)}$ turns into the following Riemann-Hilbert problem for  $N^{(4)}$:
\sss{\label{ds-esc-n4-rhp-pw1}
\ddd{
N^{(4)+} &= N^{(4)-} V_B, && k\in  B,
\\
N^{(4)+} &= N^{(4)-} V_j^{(4)},  && k\in L_j,\ j=1,2,3,4,
\\
N^{(4)+} &= N^{(4)-} V_p^{(4)},  && k\in \p D_p^\ve,
\\
N^{(4)+} &= N^{(4)-} V_{\bar p}^{(4)},  && k\in \p D_{\bar p}^\ve,
\\
N^{(4)} &=\left[I +O\left(\tfrac 1k\right)\right]e^{ig_\infty(\xi) \sigma_3},\quad && k \to \infty,
}
}
with the jump $V_B$  given by \eqref{ds-n1-jumps} and
\ddd{
&V_1^{(4)}
=
\left(
\def\arraystretch{1}
\begin{array}{lr}
1 & \dfrac{\bar r e^{2i(\theta t -g)} }{1+r\bar r}\,\delta^{2}
\\
0 & 1
\end{array}
\right),
\quad
V_2^{(4)}
=
\left(
\def\arraystretch{1}
\begin{array}{lr}
1 & 0 \\
\dfrac{r e^{-2i(\theta t -g)} }{1+r\bar r}\,\delta^{-2} & 1
\end{array}
\right),
\quad
V_3^{(4)}
=
\left(
\def\arraystretch{1}
\begin{array}{lr}
1
&
0
\\
r e^{-2i(\theta t -g)} \delta^{-2}
&
1
\end{array}
\right),
\nn\\
&V_4^{(4)}
=
\left(
\def\arraystretch{1}
\begin{array}{lr}
1
&
\bar r e^{2i(\theta t -g)}  \delta^{2}
\\
0
&
1
\end{array}
\right),
\quad
V_p^{(4)}
=
\left(
\def\arraystretch{1}
\begin{array}{lr}
1 & -\dfrac{c_p \, \delta^2(\xi, k)\, d(k)\, e^{-2ig(\xi, k)}}{k-p}\,  e^{2i\theta(\xi, p)t}
\\
0 & 1
\end{array}
\right),
\nn\\
&V_{\bar p}^{(4)} 
=
\left(
\def\arraystretch{1}
\begin{array}{lr}
1 &0
\\
-\dfrac{c_{\bar p} \, \delta^{-2}(\xi, k) \, d(k)\, e^{2ig(\xi, k)}}{k-\bar p}\,  e^{-2i\theta(\xi, \bar p)t} & 1
\end{array}
\right),
\label{ds-vpvpb-n4-pw1}
}
where the associated jump contours are shown in Figure \ref{ds-esc-def3-pw1-f} and  $g_\infty(\xi)$ is  the limit of $g(\xi, k)$ as $k\to \infty$, i.e.
\eee{ \label{ds-esc-ginf-def-pw1}
g_\infty(\xi) 
:=
\lim_{k\to\infty} g(\xi, k)
=
-\frac{1}{2i\pi^2}
\int_{\zeta\in B} \frac{1}{\lambda(\zeta)} \int_{-\infty}^{k_1(\xi)} \frac{\ln\left[1+r(\nu)\bar r(\nu)\right]}{\nu-\zeta}\, d\nu d\zeta.
} 
Importantly, expressing $g_\infty$ in terms of $\delta$ and using the symmetries $\overline{\lambda(k)} = \lambda(\bar k)$ and $\overline{\delta(\xi, k)} = \left[\delta(\xi, \bar k)\right]^{-1}$, we have that 
$\overline{g_\infty(\xi)}=g_\infty(\xi)$, i.e. that $g_\infty\in \mathbb R$. 

Observe that all the jumps of $N^{(4)}$ apart from $V_B$ tend to the identity exponentially fast in the limit $t\to\infty$. Hence, proceeding as in the appendix of \cite{bm2017}, we find that the contribution of these jumps is of order $O(t^{-1/2})$. Then, starting from the reconstruction formula \eqref{ds-q-recon-n} and applying the four successive deformations that lead to $N^{(4)}$, we eventually obtain
\eee{\label{ds-esc-q-recon-pw1-as}
q(x, t)
=
-2i \lim_{k\to\infty} \left[k N_{12}^\dom(x, t, k)\right] e^{ig_\infty(\xi)} + O\big(t^{-\frac 12}\big), \quad t\to\infty,
} 
where $N^\dom(x, t, k)$ satisfies the dominant component of Riemann-Hilbert problem \eqref{ds-esc-n4-rhp-pw1}, that is
\sss{\label{ds-esc-nmod-rhp-pw1}
\ddd{
N^{\dom+} &= N^{\dom-} V_B, && k\in  B,
\\
N^\dom &=\left[I +O\left(\tfrac 1k\right)\right]e^{ig_\infty(\xi) \sigma_3},\quad && k \to \infty.
}
}

The dominant problem \eqref{ds-esc-nmod-rhp-pw1} has been extracted from problem \eqref{ds-esc-n4-rhp-pw1} in a similar way with  problem \eqref{ds-rhp-mod-pw1} of Subsection \ref{ds-esc-pw1-lim-ss}. 
In fact, it is straightforward to verify that $N^\dom$ is given by the  explicit formula
\eee{\label{ds-esc-nmod-def-pw1}
N^\dom
=
\frac 12\, e^{ig_\infty(\xi) \sigma_3} 
\left(
\def\arraystretch{1.5}
\begin{array}{lr}
\Lambda(k)+\Lambda^{-1}(k)
&
-\dfrac{q_o}{\bar q_-}\left[\Lambda(k)-\Lambda^{-1}(k)\right]
\\
-\dfrac{q_o}{q_-} \left[ \Lambda(k)- \Lambda^{-1}(k)\right]
&
\Lambda(k)+\Lambda^{-1}(k) 
\end{array}
\right),
}
where
\eee{\label{ds-Lam-def}
\Lambda(k) := \left(\frac{k-iq_o}{k+iq_o}\right)^{\frac 14}.
}
Expressions \eqref{ds-esc-q-recon-pw1-as}  and \eqref{ds-esc-nmod-def-pw1}  yield the leading-order asymptotics \eqref{ds-qsol-pw-t}  in the range  $\xi<\Vs$ of the  transmission regime $p\in D_1$.
We note that, as expected from the fact that the discrete spectrum does not contribute at leading order for $\xi<\Vs$,  \eqref{ds-qsol-pw-t} is consistent with the result obtained for  $\xi<\Vo$ in the case of no discrete spectrum analyzed in \cite{bm2017}.

\subsection{The case $\xi = \Vs$: soliton on top of a plane wave}
\label{ds-esc-pw1-lim-ss}
The same four deformations that were performed for $\xi<\Vs$ yield once again Riemann-Hilbert problem \eqref{ds-esc-n4-rhp-pw1}. 
In particular, the jumps along $\p D_p^\ve$ and $\p D_{\bar p}^\ve$ read
\sss{\label{ds-vppb4-comb-pw1}
\ddd{
V_p^{(4)}
&=
\left(\def\arraystretch{1}\begin{array}{lr}1 & -\dfrac{c_p \, \delta^2(\Vs, k)\, d(k)\, e^{-2ig(\Vs, k)}}{k-p}\,  e^{2i\theta(\Vs, p)t}
\\
0 & 1
\end{array}
\right), 
\label{ds-vp-4-xipm}
\\
V_{\bar p}^{(4)} 
&=
\left(\def\arraystretch{1}\begin{array}{lr}1 &0
\\
-\dfrac{c_{\bar p} \, \delta^{-2}(\Vs, k) \, d(k)\, e^{2ig(\Vs, k)}}{k-\bar p}\,  e^{-2i\theta(\Vs, \bar p)t} & 1
\end{array}
\right).
\label{ds-vpb-4-xipm}
}
}
However, since $\Re (i\theta)(\Vs, p)=\Re (i\theta)(\Vs, \bar p)=0$, the time-dependent exponentials involved in the jumps \eqref{ds-vppb4-comb-pw1} are purely oscillatory (as opposed to decaying). That is, contrary to the range $\xi<\Vs$, the jumps $V_p^{(4)}$ and $V_{\bar p}^{(4)}$ no longer tend to the identity  as $t\to\infty$.
Hence, $V_p^{(4)}$ and $V_{\bar p}^{(4)}$ are now expected to contribute to the leading-order asymptotics (together, of course, with the jump $V_B$ along the branch cut $B$, which is constant) and, therefore, they must be included in the dominant component of problem \eqref{ds-esc-n4-rhp-pw1}.
Next, we extract the dominant component from the rest of the problem.

\vskip 3mm
\noindent
\textbf{Decomposition into dominant and error problems.}
Let $D_{k_1}^{\epsilon}$ be a disk centered at $k_1$ with radius $\epsilon$ sufficiently small so that $D_{k_1}^\epsilon  \cap \left(B\cup \overline{D_p^\ve \cup D_{\bar p}^\ve} \right) =\O$. Then, write the solution $N^{(4)}$ of problem \eqref{ds-esc-n4-rhp-pw1} in the form
\eee{\label{ds-n4-dec-pw1}
N^{(4)}
=
N^\err N^\asymp,
\quad 
N^\asymp
=
\begin{cases}
N^\dom, &k\in \mathbb C\setminus  
D_{k_1}^\epsilon,
\\
N^{k_1}, &k\in D_{k_1}^\epsilon,
\end{cases}
}
where the components $N^\dom$, $N^{k_1}$ and $N^\err$ are defined as follows:
\vskip 3mm
\begin{enumerate}[label=$\bullet$, leftmargin=4mm, rightmargin=0mm]
\advance\itemsep 3mm
\item
The function $N^\dom(\Vs t, t, k)$ is analytic in $\mathbb C\setminus \left(B\cup \p D_p^\ve \cup \p D_{\bar p}^\ve \right)$ and satisfies the Riemann-Hilbert problem
\sss{\label{ds-rhp-mod-pw1}
\ddd{
N^{\dom+} &= N^{\dom-} V_B,  && k\in B,
\label{ds-rhp-mod-B}
\\
N^{\dom+} &= N^{\dom-} V_p^{(4)}, && k\in \p D_p^\ve,
\label{ds-rhp-mod-p}
\\
N^{\dom+} &= N^{\dom-} V_{\bar p}^{(4)},  && k\in \p D_{\bar p}^\ve,
\label{ds-rhp-mod-pb}
\\
N^\dom &=  \left[I +O\left(\tfrac 1k\right)\right]e^{ig_\infty(\Vs) \sigma_3},\quad && k \to \infty,
\label{ds-rhp-mod-pw1-as}
}
}
with $V_B$ given by \eqref{ds-n1-jumps} and  $V_p^{(4)}$, $V_{\bar p}^{(4)}$ given by \eqref{ds-vppb4-comb-pw1}.
\item
The function $N^{k_1}(\Vs t, t, k)$ is analytic  in $D_{k_1}^{\epsilon }\setminus \bigcup_{j=1}^4 L_j$ with jumps
\eee{\label{ds-rhpP}
N^{k_1+}
=
N^{k_1-}
V_j^{(4)},
\quad
k\in \widehat L_j := L_j\cap D_{k_1}^{\epsilon },\ j=1, 2, 3, 4.
}
Note that nothing has been specified about $N^{k_1}$ outside the disk $D_{k_1}^\epsilon$.
\item
The function $N^\err(\Vs, t, k)$ is analytic in $\mathbb C
\setminus \big(\bigcup_{j=1}^4 \wc L_j \cup  \p D_{k_1}^\epsilon\big)$, where $\wc L_j := L_j\setminus \widehat L_j$, and satisfies the Riemann-Hilbert problem
\sss{\label{ds-rhpE}
\ddd{
N^{\err+} &= N^{\err-} 
V^\err,\quad && k\in {\textstyle \bigcup}_{j=1}^4 \wc L_j \cup  \p D_{k_1}^{\epsilon},
\label{ds-rhpEb}
\\
N^\err &= I +O\left(\tfrac 1k \right), && k \to \infty,
\label{ds-rhpE-as}
}
}
where
\eee{\label{ds-VEdef}
V^\err
=
\begin{cases}
N^\dom V_j^{(4)} (N^\dom)^{-1}, &k\in \wc L_j,
\\
N^{\asymp-}(V_D^\asymp)^{-1}(N^{\asymp-})^{-1}, &k\in \p D_{k_1}^\epsilon,
\end{cases}
}
and $V_D^\asymp$ is the yet unknown jump of $N^\asymp$ along the circle $\p D_{k_1}^\epsilon$.
\end{enumerate}

\vskip 3mm

Under the four successive deformations that lead to problem \eqref{ds-esc-n4-rhp-pw1}, the reconstruction formula \eqref{ds-q-recon-n} becomes
\eee{ \label{ds-esc-q-recon-pw1}
q(x, t)
=
-2i \lim_{k \to \infty}  k N_{12}^{(4)}(\Vs t, t, k) e^{ig_\infty(\Vs)},
} 
where we have also used the fact that $\delta, d \to 1$ as $k \to \infty$.
This formula combined with the decomposition \eqref{ds-n4-dec-pw1} and the asymptotic conditions \eqref{ds-rhp-mod-pw1-as} and \eqref{ds-rhpE-as} implies 
\eee{
q(x, t) = -2i\lim_{k \to\infty} k \left[ N_{12}^\dom(\Vs t, t, k) e^{ig_\infty(\Vs)}+ N_{12}^\err(\Vs t, t, k)\right].
}
The error problem \eqref{ds-rhpE} is precisely that of the plane wave region in \cite{bm2017}, since the jumps around $p$ and $\bar p$ are not part of this problem. Hence, as shown in \cite{bm2017},   $\lim_{k \to\infty} k  N_{12}^\err = O\big(t^{-1/2}\big)$. In turn, we obtain
\eee{\label{ds-q-recon-pw1-2}
q(x, t) = -2i\lim_{k \to\infty} k  N_{12}^\dom(\Vs t, t, k) e^{ig_\infty(\Vs)} + O\big(t^{-\frac 12}\big), \quad t\to\infty.
}
It remains to determine $N^\dom$, i.e. to solve the dominant Riemann-Hilbert problem  \eqref{ds-rhp-mod-pw1}.

\vskip 3mm
\noindent
\textbf{Solution of the dominant problem.}
We begin by converting the jumps  along the circles $\p D_p^\ve$ and $\p D_{\bar p}^\ve$ back to residue conditions at $p$ and $\bar p$. This is done by reverting transformation \eqref{ds-nrh}, i.e. by letting
\eee{\label{ds-esc-mmod-def-pw1}
M^\dom 
=
\left\{
\def\arraystretch{1.2}
\begin{array}{ll}
N^\dom  \big(V_p^{(4)}\big)^{-1}, &k\in D_p^{\ve },
\\
N^\dom, &k\in \mathbb C^-\setminus \left( B^-\cup \overline{D_{p}^{\ve}}\,\right),
\\
N^\dom \big(V_{\bar p}^{(4)}\big)^{-1}, &k\in D_{\bar p}^{\ve },
\\
N^\dom, &k\in \mathbb C^+\setminus \left( B^+\cup \overline{D_{\bar p}^{\ve}}\,\right).
\end{array}
\right.
}
Then, $M^\dom$ is the solution of the Riemann-Hilbert problem
\sss{\label{ds-rhp-mmod-pw1}
\ddd{
M^{\dom+} &= M^{\dom-} V_B,  && k\in B,
\label{ds-rhp-mmod-B-pw1}
\\
M^\dom &=  \left[I +O\left(\tfrac 1k\right)\right]e^{ig_\infty(\Vs) \sigma_3},\quad  && k \to \infty,
\\
\underset{k=p}{\text{Res}}\, M^\dom &= \left(0, \rho_p\,  M_1^\dom(p)\right),
\label{ds-rhp-mmod-p-pw1}
\\
\underset{k=\bar p}{\text{Res}}\,  M^\dom &= \left(\rho_{\bar p}\,  M_2^\dom(\bar p), 0\right),
\label{ds-rhp-mmod-pb-pw1}
}
}
where $M_1^\dom, M_2^\dom$ denote the two columns of $M^\dom$ and
\sss{\label{ds-rhop-def}
\ddd{
&\rho_p 
= 
c_p \delta^2(\Vs, p) d(p) e^{2i\left[\theta(\Vs, p)t-g(\Vs, p)\right]},
\\
&\rho_{\bar p} 
= 
c_{\bar p} \delta^{-2}(\Vs, \bar p) d(\bar p) e^{-2i\left[\theta(\Vs, \bar p)t-g(\Vs, \bar p)\right]}.
}
}
In fact, the expressions for $\rho_p$ and $\rho_{\bar p}$ can be simplified after noting that  the symmetry (see \cite{bm2017})
$$
\overline{\Psi_\pm(x,t,\bar k)} = -\sigma_*\Psi_\pm(x,t,k)\sigma_*,
\quad   
\sigma_*
:=
\begin{pmatrix}
0 & 1 \\ -1 & 0
\end{pmatrix},
$$
together with relations \eqref{ds-sym1} and \eqref{ds-sym2} imply $C_{\bar p} = -\overline{C_p}$. Then, recalling the Schwarz symmetries $\bar a'(\bar k) = \overline{a'(k)}$, $d(\bar k) = \overline{d(k)}$ and the definitions \eqref{ds-cp-def} and \eqref{ds-cpb-def} of $c_p$ and $c_{\bar p}$, we obtain
\eee{\label{ds-Cp-Cpb-sym}
c_{\bar p} = - \overline{c_p}.
}
Hence, noting in addition that $\theta(\xi, \bar k) = \overline{\theta(\xi, k)}$, $g(\xi, \bar k) = \overline{g(\xi, k)}$, $\delta(\xi, \bar k) = \overline{\delta^{-1}(\xi, k)}$, and since $\theta(\Vs, \bar p)\in \mathbb R$, we have
\ddd{\label{ds-R-def}
\rho_p 
= 
R_p \, e^{2i\theta(\Vs, p)t},
\quad
\rho_{\bar p} 
= 
-\overline{R_p} \, e^{-2i\theta(\Vs, p)t},
\quad
R_p:=C_p \frac{\delta^2(\Vs, p)  e^{-2i g(\Vs, p)}}{a'(p)},
}
which shows that $\rho_{\bar p} = -\overline{\rho_p}$.

We will solve problem \eqref{ds-rhp-mmod-pw1} by decomposing it into discrete and continuous spectrum components via the substitution
\eee{\label{ds-esc-mmodt-def-pw1-0}
M^\dom = \mathcal  M^\dom W.
}
Here,  $W$ is the solution of  the continuous spectrum component problem
\sss{\label{ds-esc-w-rhp-pw1}
\ddd{
W^+ &= W^- V_B, && k\in B,
\label{ds-esc-w-rhp-pw1-jump}
\\
W &= \left[I + O\left(\tfrac 1k \right)\right] e^{ig_\infty(\Vs) \sigma_3}, \quad && k\to \infty,
\label{ds-esc-w-asymp-pw1}
}
}
which is nothing but problem \eqref{ds-esc-nmod-rhp-pw1} evaluated at $\xi=\Vs$. Therefore, 
\eee{\label{ds-esc-w-def-pw1}
W
= 
\frac 12 \, e^{ig_\infty(\Vs) \sigma_3} 
\left(
\def\arraystretch{1.5}
\begin{array}{lr}
\Lambda(k)+\Lambda^{-1}(k) 
&
-\frac{q_o}{\bar q_-}\left[\Lambda(k)-\Lambda^{-1}(k)\right]
\\
-\frac{q_o}{q_-} \left[ \Lambda(k)- \Lambda^{-1}(k)\right]
&
\Lambda(k)+\Lambda^{-1}(k) 
\end{array}
\right)
}
with $\Lambda$ given by \eqref{ds-Lam-def}. 

Since $\det W \equiv 1\neq 0$,  we can rearrange \eqref{ds-esc-mmodt-def-pw1-0} to 
\eee{\label{ds-esc-mmodt-def-pw1}
\mathcal M^\dom  = M^\dom  W^{-1}
}
and hence deduce that  $\mathcal M^\dom$ does not have a jump along $B$, i.e. $\mathcal M^\dom$ is indeed the discrete spectrum component of $M^\dom$.
Moreover, since $M_1^\dom$ and $W$ are analytic at $p$, formula \eqref{ds-esc-mmodt-def-pw1}  and the residue condition \eqref{ds-rhp-mmod-p-pw1} imply
\sss{\label{ds-esc-mmodt-res-pw1}
\eee{
\underset{k=p}{\text{Res}}\, \mathcal M_1^\dom
=
-W_{21}(p)   \rho_p \, M_1^\dom(p).
}
Similarly, we find
\ddd{
\underset{k=p}{\text{Res}}\, \mathcal M_2^\dom
&=
W_{11}(p)   \rho_p \, M_1^\dom(p),
\\
\underset{k=\bar p}{\text{Res}}\, \mathcal M_1^\dom
&=
W_{22}(\bar p)   \rho_{\bar p} \, M_2^\dom(\bar p)
=
\overline{W_{11}(p)}   \rho_{\bar p} \, M_2^\dom(\bar p),
\\
\underset{k=\bar p}{\text{Res}}\, \mathcal M_2^\dom
&=
-W_{12}(\bar p)  \rho_{\bar p} \, M_2^\dom(\bar p)
=
\overline{W_{21}(p)}  \rho_{\bar p} \, M_2^\dom(\bar p),
}
}
where in the last two conditions we have also made use of the symmetries
\eee{\label{ds-w-syms}
W_{22}(\bar p)
=
\overline{W_{11}(p)},
\quad
W_{12}(\bar p)
=
-\overline{W_{21}(p)}.
}
Furthermore, \eqref{ds-esc-mmodt-def-pw1} in combination with the asymptotic conditions for $M^\dom$ and $W$ as $k\to\infty$ yield the following asymptotic condition for $\mathcal M^\dom$:
\ddd{\label{ds-esc-mmodt-asymp-pw1}
\mathcal M^\dom
&=
\left(\left[I+O\left(\tfrac 1k\right)\right] e^{ig_\infty(\Vs) \sigma_3}\right)
\left(\left[I+O\left(\tfrac 1k\right)\right] e^{ig_\infty(\Vs) \sigma_3}\right)^{-1}
\nn\\
&=
\left[I+O\left(\tfrac 1k\right)\right]  \left[I+O\left(\tfrac 1k\right)\right]^{-1} 
=
I+O\left(\tfrac 1k\right), \quad k\to\infty,
}
where we note that the $O(\frac 1k)$ term possibly involves in some form  the exponential $e^{ig_\infty(\Vs)\sigma_3}$.

In summary, $\mathcal M^\dom$ is analytic for $k\in\mathbb C\setminus\left\{ p, \bar p\right\}$, has simple poles at $p$ and $\bar p$ with associated residues satisfying \eqref{ds-esc-mmodt-res-pw1}, and satisfies the asymptotic condition \eqref{ds-esc-mmodt-asymp-pw1} as $k\to \infty$. Thus, Liouville's theorem implies
\eee{\label{ds-esc-mmodt-form-pw1}
\mathcal M^\dom
=
I + \frac{\underset{k=p}{\text{Res}}\, \mathcal M^\dom}{k-p} + \frac{\underset{k=\bar p}{\text{Res}}\, \mathcal M^\dom}{k-\bar p},
}
and hence it only remains to determine the residues of $\mathcal M^\dom$ at $p$ and $\bar p$. 
In fact, thanks to \eqref{ds-esc-mmodt-res-pw1} this amounts to computing the corresponding residues of $M^\dom$.
Combining \eqref{ds-esc-mmodt-def-pw1}, \eqref{ds-esc-mmodt-res-pw1}, \eqref{ds-w-syms} and \eqref{ds-esc-mmodt-form-pw1}, we find
\sss{\label{ds-esc-mmod-form-pw1}
\ddd{\label{ds-esc-mmod1-form-pw1}
M^\dom_1
&=
W_1 + W_{11} \left[
-
\frac{W_{21}(p) \rho_p \, M^\dom_1(p)}{k-p}
+
\frac{\overline{W_{11}(p)} \rho_{\bar p}\, M^\dom_2(\bar p)}{k-\bar p}
\right]
\nn\\
&\quad
+
W_{21} 
\left[
\frac{W_{11}(p) \rho_p\,  M^\dom_1(p)}{k-p}
+
\frac{\overline{W_{21}(p)} \rho_{\bar p}\, M^\dom_2(\bar p)}{k-\bar p}
 \right]
}
and
\ddd{\label{ds-esc-mmod2-form-pw1}
M^\dom_2
&=
W_2 + 
W_{12} \left[
-
\frac{W_{21}(p) \rho_p \, M^\dom_1(p)}{k-p}
+
\frac{\overline{W_{11}(p)} \rho_{\bar p}\, M^\dom_2(\bar p)}{k-\bar p}
\right]
\nn\\
&\quad
+
W_{22} \left[
\frac{W_{11}(p) \rho_p\,  M^\dom_1(p)}{k-p}
+
\frac{\overline{W_{21}(p)} \rho_{\bar p}\, M^\dom_2(\bar p)}{k-\bar p}
 \right].
}
}
Since $M^\dom_1$ is analytic at $p$, we can evaluate \eqref{ds-esc-mmod1-form-pw1} at $k=p$ to obtain
\sss{\label{ds-esc-mmod-comb}
\ddd{
M^\dom_1(p)
&=
W_1(p) 
+
\rho_p \left[-W_{11}'(p)W_{21}(p)+W_{11}(p)W_{21}'(p)\right] 
M^\dom_1(p)
\nn\\
&\quad
+ 
\rho_{\bar p}\, \frac{\left|W_{11}(p)\right|^2 + \left|W_{21}(p)\right|^2}{p-\bar p}  \, M^\dom_2(\bar p).
\label{ds-esc-mmod1p-pw1-eq}
}
Similarly, since $M^\dom_2$ is analytic at $\bar p$, evaluating \eqref{ds-esc-mmod2-form-pw1} at $k=\bar p$ and using  the symmetries \eqref{ds-w-syms} (which also apply for $W'$), we have
\ddd{
M^\dom_2(\bar p)
&=
W_2(\bar p) -  \rho_p\, \frac{\left|W_{11}(p)\right|^2+\left|W_{21}(p)\right|^2}{p-\bar p} \, M^\dom_1(p)
\nn\\
&\quad
+
\rho_{\bar p} \left[\overline{W_{11}'(p) W_{21}(p)
-
W_{21}'(p)W_{11}(p)}\right] M^\dom_2(\bar p).
\label{ds-esc-mmod2p-pw1-eq}
}
}
Equations \eqref{ds-esc-mmod-comb} form a system for $M^\dom_1(p)$ and $M_2^\dom(\bar p)$, which can be solved to yield
\sss{\label{ds-esc-mmod12-sol-pw1}
\ddd{
M^\dom_1(p)
&=
\frac{\mathcal B \rho_{\bar p} W_2(\bar p) + \left(1-\bar{\mathcal A} \rho_{\bar p}\right) W_1(p)}{\mathcal B^2 \rho_p \rho_{\bar p} + \left(1-\bar{\mathcal A} \rho_{\bar p}\right)\left(1+\mathcal A \rho_p\right)},
\\
M^\dom_2(\bar p)
&=
\frac{\left(1+\mathcal A\rho_p\right)W_2(\bar p) - \mathcal B \rho_p W_1(p)}{\mathcal B^2 \rho_p \rho_{\bar p} + \left(1-\bar{\mathcal A} \rho_{\bar p}\right)\left(1+\mathcal A \rho_p\right)},
}
}
where  
\eee{
\mathcal A
=
W_{11}'(p)W_{21}(p)-W_{11}(p)W_{21}'(p),
\quad
\mathcal B
=
\frac{\left|W_{11}(p)\right|^2 + \left|W_{21}(p)\right|^2}{p-\bar p}.
}
Actually, using formula \eqref{ds-esc-w-def-pw1}, we can simplify $\mathcal A$ and $\mathcal B$  to the constants
\eee{\label{ds-ab-sol-def}
\mathcal A= \frac{i\bar q_-}{2\left(p^2+q_o^2\right)},
\quad
\mathcal B=\frac{\left|p-iq_o\right|+\left|p+iq_o\right|}{2\left|p^2+q_o^2\right|^{\frac 12}\left(p-\bar p\right)}.
}
Expressions \eqref{ds-esc-mmod12-sol-pw1} combined with  \eqref{ds-esc-mmodt-res-pw1} yield the residues of $\mathcal M^\dom$ at $p$ and $\bar p$,  and hence $\mathcal M^\dom$ itself via formula \eqref{ds-esc-mmodt-form-pw1}. 

Having determined $\mathcal M^\dom$, we return to the reconstruction formula  \eqref{ds-q-recon-pw1-2} and note that transformations \eqref{ds-esc-mmod-def-pw1} and \eqref{ds-esc-mmodt-def-pw1} imply
\eee{\label{ds-q-recon-pw1-3}
q(x, t)
=
-2i\lim_{k \to\infty} k  \big(\mathcal M^\dom  W \big)_{12} \, e^{ig_\infty(\Vs)} + O\big(t^{-\frac 12}\big), \quad t\to\infty.
}
Furthermore, by the asymptotic conditions \eqref{ds-esc-w-asymp-pw1} and \eqref{ds-esc-mmodt-asymp-pw1} we have
\eee{\label{ds-esc-mcalmod-asymp-pw1}
W  =  e^{ig_\infty(\Vs)\sigma_3} + \frac{w}{k} + O\left(\frac{1}{k^2}\right),  \quad
\mathcal M^\dom  = I + \frac{\mu}{k} + O\left(\frac{1}{k^2}\right),   \quad k\to \infty,
}
where the matrix-valued functions $w$ and $\mu$ may depend on $x$ and $t$ but not on $k$.
Thus,
\eee{
\big(\mathcal M^\dom W\big)_{12}
=
\frac{w_{12} + \mu_{12}\, e^{-ig_\infty(\Vs)}}{k} + O\left(\frac{1}{k^2}\right), \quad  k\to \infty.
}
Hence, noting that $w_{12} = \frac{iq_-}{2} e^{ig_\infty(\Vs)}$ by 
formula \eqref{ds-esc-w-def-pw1}, we obtain
\eee{\label{ds-esc-q-sol-lim-pw1-0}
q(x, t)
=
e^{2ig_\infty(\Vs)} q_- 
- 
2i  \mu_{12} + O\big(t^{-\frac 12}\big), \quad t\to \infty.
}
Moreover, matching the second expansion in \eqref{ds-esc-mcalmod-asymp-pw1} with the large-$k$ expansion of \eqref{ds-esc-mmodt-form-pw1}, we infer
\eee{
\mu_{12} 
=
\Big(\underset{k=p}{\text{Res}} \, \mathcal M_2^\dom \Big)^{(1)} 
+ 
\Big(\underset{k=\bar p}{\text{Res}} \, \mathcal M_2^\dom\Big)^{(1)}.
}
Thus, using successively  \eqref{ds-esc-mmodt-res-pw1}, \eqref{ds-esc-mmod12-sol-pw1} and  \eqref{ds-w-syms}, we find
\eee{
\mu_{12} 
=
\frac{\left(1-\bar{\mathcal A} \rho_{\bar p}\right)\rho_p W_{11}(p)^2-\left(1+\mathcal A\rho_p\right)\rho_{\bar p}\overline{W_{21}(p)}^2-2\mathcal B \rho_p \rho_{\bar p} W_{11}(p)\overline{W_{21}(p)}}{\mathcal B^2 \rho_p \rho_{\bar p} + \left(1-\bar{\mathcal  A} \rho_{\bar p}\right)\left(1+\mathcal A \rho_p\right)}.
\label{ds-m12-pw1-temp1}
}
Substituting for $W$ via \eqref{ds-esc-w-def-pw1} and inserting the resulting expression in \eqref{ds-esc-q-sol-lim-pw1-0}, we conclude that
\ddd{\label{ds-esc-q-sol-lim-pw1}
q(x, t) 
&=
q_{\text{pw}}(\Vs)
-
\frac i2 \, e^{2ig_\infty(\Vs)} \left[\mathcal B^2 \rho_p \rho_{\bar p} + \left(1-\bar{\mathcal A} \rho_{\bar p}\right)\left(1+\mathcal A \rho_p\right)\right]^{-1}
\bigg\{
\left(1-\bar{\mathcal A} \rho_{\bar p}\right)\rho_p 
 \left[\Lambda(p)+  \Lambda^{-1}(p)\right]^2
\nn\\
&\quad
-
\left(1+\mathcal A\rho_p\right)\rho_{\bar p}\,
\frac{q_-}{\bar q_-} \left[\overline{\Lambda(p) - \Lambda^{-1}(p)}\right]^2  
+
2\mathcal B \rho_p \rho_{\bar p} \,  \frac{q_o}{\bar q_-} 
 \left[ \Lambda(p)+ \Lambda^{-1}(p)\right]
\left[\overline{\Lambda(p) - \Lambda^{-1}(p)}\right]  
\bigg\}
\nn\\
&\hskip 9cm
 + O\big(t^{-\frac 12}\big), \quad t\to \infty,
}
where $q_{\text{pw}}(\Vs)$ is the plane wave \eqref{ds-qpw-def} evaluated at $\xi=\Vs$, the quantities $\Lambda, \mathcal A, \mathcal B, \rho_p, \rho_{\bar p}$ are given by \eqref{ds-Lam-def}, \eqref{ds-ab-sol-def},  \eqref{ds-R-def} and the real constant $g_\infty(\Vs)$ is obtained by evaluating \eqref{ds-esc-ginf-def-pw1} at $\xi=\Vs$.
In fact, setting 
\eee{\label{ds-DE-def}
\Lambda_1 := \Lambda(p)+ \Lambda^{-1}(p),
\quad
\Lambda_2 := \overline{\Lambda(p)- \Lambda^{-1}(p)},
}
and substituting for $\rho_p, \rho_{\bar p}$ via  \eqref{ds-R-def} turns the leading-order asymptotics \eqref{ds-esc-q-sol-lim-pw1} into the form \eqref{ds-esc-q-sol-lim-pw1-t}-\eqref{ds-qstheta-def} given in Theorem \ref{ds-per-t}.

\subsection{The range $\Vs<\xi<\Vo$: plane wave with a phase shift}
\label{ds-esc-pw2-ss}
The analysis in this range is similar to the one for $\xi<\Vs$. Indeed, under the same series of deformations as in Subsection \ref{ds-esc-pw1-ss}, Riemann-Hilbert problem \eqref{ds-n-rhp-intro} can be transformed once again into Riemann-Hilbert problem \eqref{ds-esc-n4-rhp-pw1}. We note, in particular, that, since $p\in D_1$, for $\Vs<\xi<\Vo$ the point $p$ lies inside the \textit{unbounded} region of positive sign to the left of the stationary point $k_1$ (the unbounded  region in white inside the third quadrant of the second frame of Figure \ref{sign-structure-f}). Thus, all four stages of the first deformation for $\xi<\Vs$ can be repeated for $\Vs<\xi<\Vo$ in a way that leaves the jump along $\p D_p^\ve$ invariant. By symmetry, the same is true for the  jump along $\p D_{\bar p}^\ve$. Importantly, we shall see later that this is not the case for $p\in D_3$ (transmission/wake regime).

An important difference between the ranges $(-\infty, \Vs)$ and $(\Vs, \Vo)$, however, is that in the latter case the jumps $V_p^{(4)}$ and $V_{\bar p}^{(4)}$ defined by \eqref{ds-vpvpb-n4-pw1}   \textit{grow exponentially} as $t\to\infty$, since $\Re (i\theta)(\xi, p)>0$ and $\Re (i\theta)(\xi, \bar p)<0$  for all $\xi>\Vs$.
This is to be contrasted with the range $(-\infty, \Vs)$, where we recall that these jumps decayed exponentially to the identity and hence could be immediately neglected from the dominant Riemann-Hilbert problem.
Nevertheless, it turns out that the jumps along $\p D_p^\ve$ and $\p D_{\bar p}^\ve$ still do not contribute to the leading-order asymptotics. Along the lines of  \cite{dkkz1996}, this can be seen by applying  the following additional transformation  to problem \eqref{ds-esc-n4-rhp-pw1}:
\sss{\label{ds-esc-n4t-def-pw2}
\eee{
\widetilde N^{(4)} 
=
\begin{cases}
N^{(4)}  n^{\sigma_3}, & k\in \mathbb C\setminus \left(\, \overline{D_p^\ve} \cup  \overline{D_{\bar p}^\ve}\, \right),
\\
N^{(4)} J_p \, n^{\sigma_3}, & k\in D_p^\ve,
\\
N^{(4)} J_{\bar p} \, n^{\sigma_3}, & k\in D_{\bar p}^\ve,
\end{cases}
} 
where $n(k)$ is the piecewise-defined function
\eee{  
n(k)
=
\left\{
\def\arraystretch{1.2}
\begin{array}{ll}
\dfrac{k-\bar p}{k-p}, & k\in \mathbb C\setminus \left(\,\overline{D_p^\ve} \cup \overline{D_{\bar p}^\ve}\,\right),
\\
k-\bar p, & k\in D_p^\ve,
\\
\dfrac{1}{k-p}, & k\in D_{\bar p}^\ve,
\end{array}
\right.
} 
and the matrices $J_p(\xi, k)$ and $J_{\bar p}(\xi, k)$ are given by
\ddd{
J_p(\xi, k) 
&=
\left(
\def\arraystretch{2.4}
\arraycolsep=0pt
\begin{array}{lr}
\dfrac{1- \frac{n^2(p)}{n^2(k)}}{k-p} &  c_p \, d(k)\delta^2(\xi, k) e^{2i\left[\theta(\xi, p)t-g(\xi, k)\right]}
\\
- \dfrac{n^2(p)\, e^{-2i\left[\theta(\xi, p)t-g(\xi, k)\right]}}{c_p \, d(k)\delta^2(\xi, k) \left(k-\bar p\right)^2} &  k-p
\end{array}
\right), 
\\
J_{\bar p}(\xi, k)
&=
\left(
\arraycolsep=0pt
\def\arraystretch{2.8}
\begin{array}{lr}
k-\bar p &  -\dfrac{e^{2i\left[\theta(\xi, \bar p)t-g(\xi, k)\right]}}{n^2(\bar p) c_{\bar p}\, d(k) \delta^{-2}(\xi, k) \left(k-p\right)^2} 
\\
c_{\bar p}\, d(k)\delta^{-2}(\xi, k) e^{-2i\left[\theta(\xi, \bar p)t-g(\xi, k)\right]}  & \dfrac{1-\frac{n^2(k)}{n^2(\bar p)}}{k-\bar p}
\end{array}
\right).
}
}
Note importantly that $J_p$ is analytic in $D_p^\ve$ since the singularity of its $11$-element at $k=p$ is removable. Similarly,  $J_{\bar p}$ is analytic in $D_{\bar p}^\ve$. Therefore, $\widetilde N^{(4)}$ inherits the analyticity of $N^{(4)}$ and satisfies the following Riemann-Hilbert problem:
\sss{\label{ds-esc-n4t-rhp-pw2}
\ddd{
\widetilde N^{(4)+} &= \widetilde N^{(4)-}  \widetilde V_B^{(4)},  && k\in  B,
\\
\widetilde N^{(4)+} &= \widetilde N^{(4)-}  \widetilde V_j^{(4)}, &&k\in L_j,\ j=1,2,3,4,
\\
\widetilde N^{(4)+} 
&=
\widetilde N^{(4)-} 
\widetilde V_p^{(4)}, && k\in  \p D_p^\ve,
\\
\widetilde N^{(4)+} 
&=
\widetilde N^{(4)-} 
\widetilde V_{\bar p}^{(4)}, && k\in  \p D_{\bar p}^\ve,
\\
\widetilde N^{(4)} &=\left[I +O\left(\tfrac 1k\right)\right]e^{ig_\infty(\xi) \sigma_3},\quad && k \to \infty,
}
}
with $g_\infty$  defined by \eqref{ds-esc-ginf-def-pw1} and 
\ddd{\label{ds-esc-n4t-jumps-pw2}
&\hskip 3cm
\widetilde V_B^{(4)}
=
\left(
\def\arraystretch{1.2}
\begin{array}{lr}
0 & \frac{q_-}{iq_o} \, n^{-2} \\
\frac{\bar q_-}{iq_o}\, n^2 & 0
\end{array}
\right),
\nn\\
\widetilde V_1^{(4)}
&=
\left(
\def\arraystretch{1}
\begin{array}{lr}
1 & \dfrac{\bar r \,\delta^2 n^{-2} e^{-2ig}}{1+r\bar r} \, e^{2i\theta t}
\\
0 & 1
\end{array}
\right),
\quad
\widetilde V_2^{(4)}
=
\left(
\def\arraystretch{1}
\begin{array}{lr}
1 & 0 \\
\dfrac{r \delta^{-2} n^2 e^{2ig}}{1+r\bar r}\, e^{-2i\theta t} & 1
\end{array}
\right),
\nn\\
\widetilde V_3^{(4)}
&=
\left(
\def\arraystretch{1}
\begin{array}{lr}
1
&
0
\\
r  \delta^{-2} n^2 e^{2ig} e^{-2i\theta t}
&
1
\end{array}
\right),
\hskip 0.7cm
\widetilde V_4^{(4)}
=
\left(
\def\arraystretch{1}
\begin{array}{lr}
1
&
\bar r  \delta^2 n^{-2} e^{-2ig} e^{2i\theta t} 
\\
0
&
1
\end{array}
\right),
\nn\\
&\hskip 1.2cm
\widetilde V_p^{(4)} 
= 
\left(
\def\arraystretch{1}
\begin{array}{lr}
 1 & 0 \\ -\dfrac{n^2(p) \delta^{-2}(\xi, k) e^{2ig(\xi, k)}}{c_p \, d(k) \left(k-p\right)}\, e^{-2i\theta(\xi, p)t} & 1 
\end{array}
\right),
\nn\\
&\hskip 1.2cm
\widetilde V_{\bar p}^{(4)} 
= 
\left(
\def\arraystretch{1}
\begin{array}{lr}
1 & -\dfrac{n^{-2}(\bar p) \delta^2(\xi, k) e^{-2ig(\xi, k)}}{ c_{\bar p}\, d(k)\left(k-\bar p\right)} \, e^{2i \theta(\xi, \bar p)t}
\\
0 &1 
\end{array}
\right).
}

All the jumps of $\widetilde N^{(4)}$ with the exception of $\widetilde V_B^{(4)}$ tend to the identity exponentially fast as $t\to\infty$. 
Importantly, as a result of transformation  \eqref{ds-esc-n4t-def-pw2}, this includes the jumps $\widetilde V_p^{(4)}$ and $\widetilde V_{\bar p}^{(4)}$.
Hence, we anticipate that the leading-order contribution of problem \eqref{ds-esc-n4t-rhp-pw2} comes from the jump $\widetilde V_B^{(4)}$. 
As this jump depends on $k$ through the function $n$, prior to decomposing problem \eqref{ds-esc-n4t-rhp-pw2} into dominant and error components we employ yet one more transformation that converts   $\widetilde V_B^{(4)}$ into the constant jump $V_B$. 
Specifically, we let
\eee{\label{ds-esc-n5t-def-pw2}
\widetilde N^{(5)}(x, t, k) = \widetilde N^{(4)}(x, t, k) e^{i\widetilde g(k)\sigma_3},
}
where the function $\widetilde g(k)$ is analytic in $\mathbb C\setminus B$ and satisfies the Riemann-Hilbert problem
\sss{\label{ds-esc-jump-gt-pw2}
\ddd{
&e^{i(\widetilde g^+ + \widetilde g^-)} = n^{-2},\quad && k\in B,
\\
&\frac{\widetilde g}{\lambda}
=
O\left(\frac 1k\right), && k \to\infty.
}
}
The above problem can be solved explicitly via Plemelj's formulae to yield 
\eee{\label{ds-esc-gt-def-pw2}
\widetilde g(k)
=
-\frac{\lambda(k)}{\pi}
\int_{\zeta\in B} \frac{\ln\left(\frac{\zeta-\bar p}{\zeta-p}\right)}{\lambda(\zeta)\left(\zeta-k\right)} \, d\zeta, \quad k\notin B.
} 
Note that $\widetilde g$ does not depend on $\xi$. 
Combining problems \eqref{ds-esc-n4t-rhp-pw2} and \eqref{ds-esc-jump-gt-pw2}, we obtain the following Riemann-Hilbert problem for $\widetilde N^{(5)}$:
\sss{\label{ds-esc-n5t-rhp-pw2}
\ddd{
\widetilde N^{(5)+} &= \widetilde N^{(5)-}  V_B,  && k\in  B,
\\
\widetilde N^{(5)+} &= \widetilde N^{(5)-} \widetilde V_j^{(5)}, &&k\in L_j,\ j=1,2,3,4,
\\
\widetilde N^{(5)+} 
&=
\widetilde N^{(5)-} 
\widetilde V_p^{(5)}, && k\in  \p D_p^\ve,
\\
\widetilde N^{(5)+} 
&=
\widetilde N^{(5)-} 
\widetilde V_{\bar p}^{(5)}, && k\in  \p D_{\bar p}^\ve,
\\
\widetilde N^{(5)} &=\left[I +O\left(\tfrac 1k\right)\right]e^{i\left[g_\infty(\xi)+\widetilde g_\infty\right] \sigma_3},\quad && k \to \infty,
}
}
where $V_B$ is defined by \eqref{ds-n1-jumps}, the remaining jumps are given by
\ddd{\label{ds-esc-n5t-jumps-pw2}
&\widetilde V_1^{(5)}
=
\left(
\def\arraystretch{1}
\begin{array}{lr}
1 & \dfrac{\bar r \,\delta^2 n^{-2} e^{-2i\left(g+\widetilde g\right)}}{1+r\bar r} \, e^{2i\theta t}
\\
0 & 1
\end{array}
\right),
\quad
\widetilde V_2^{(5)}
=
\left(
\def\arraystretch{1}
\begin{array}{lr}
1 & 0 \\
\dfrac{r \delta^{-2} n^2 e^{2i\left(g+\widetilde g\right)}}{1+r\bar r}\, e^{-2i\theta t} & 1
\end{array}
\right),
\nn\\
&\widetilde V_3^{(5)}
=
\left(
\def\arraystretch{1}
\begin{array}{lr}
1
&
0
\\
r  \delta^{-2} n^2 e^{2i\left(g+\widetilde g\right)} e^{-2i\theta t}
&
1
\end{array}
\right),
\hskip 0.73cm
\widetilde V_4^{(5)}
=
\left(
\def\arraystretch{1}
\begin{array}{lr}
1
&
\bar r  \delta^2 n^{-2} e^{-2i\left(g+\widetilde g\right)} e^{2i\theta t} 
\\
0
&
1
\end{array}
\right),
\nn\\
&\hskip 1.7cm
\widetilde V_p^{(5)} 
= 
\left(
\def\arraystretch{1}
\begin{array}{lr}
 1 & 0 
 \\ 
 -\dfrac{n^2(p) \delta^{-2}(\xi, k) e^{2i\left[g(\xi, k)+\widetilde g(k)\right]}}{c_p \, d(k) \left(k-p\right)}\, e^{-2i\theta(\xi, p)t} & 1 
\end{array}
\right),
\nn\\
&\hskip 1.7cm
\widetilde V_{\bar p}^{(5)} 
= 
\left(
\def\arraystretch{1}
\begin{array}{lr}
1 & -\dfrac{n^{-2}(\bar p) \delta^2(\xi, k) e^{-2i\left[g(\xi, k)+\widetilde g(k)\right]}}{ c_{\bar p}\, d(k) \left(k-\bar p\right)} \, e^{2i \theta(\xi, \bar p)t}
\\
0 &1 
\end{array}
\right),
}
the real quantity $g_\infty(\xi)$ is defined by \eqref{ds-esc-ginf-def-pw1}, and the real constant $\widetilde g_\infty$ is the $O(1)$ term of the expansions of $\widetilde g(k)$ as $k\to\infty$, i.e.
\eee{\label{ds-esc-ginft-def-pw2}
\widetilde g_\infty 
=
\frac{1}{\pi} \int_{\zeta\in B} \frac{\ln\left(\frac{\zeta-\bar p}{\zeta-p}\right)}{\lambda(\zeta)}\, d\zeta.
}

At leading order, the jumps of problem \eqref{ds-esc-n5t-rhp-pw2} are the same with those of problem \eqref{ds-esc-n4-rhp-pw1}. Indeed, along $B$  the jump of both problems is equal to $V_B$,  while the remaining jumps in both cases tend to the identity as $t\to\infty$.  
Thus, at leading order, the only difference between the two problems  is the presence of the constant phase $\widetilde g_\infty$ in the normalization condition of problem \eqref{ds-esc-n5t-rhp-pw2}.
Therefore, noting that under transformations \eqref{ds-esc-n4t-def-pw2} and \eqref{ds-esc-n5t-def-pw2} the reconstruction formula \eqref{ds-esc-q-recon-pw1} becomes 
\eee{ \label{ds-esc-q-recon-pw2}
q(x, t)
=
-2i \lim_{k \to \infty} k \widetilde N_{12}^{(5)}(x, t, k) e^{i\left[g_\infty(\xi)+\widetilde g_\infty\right]},
} 
we conclude that the leading-order asymptotics  in the range $\Vs<\xi<\Vo$ is equal to the plane wave \eqref{ds-qsol-pw-t} up to a   \textit{constant phase shift} of $2\widetilde g_\infty$, i.e.
\eee{\label{ds-esc-qasym-pw2}
q(x, t) = e^{2i\left[g_\infty(\xi)+\widetilde g_\infty\right]} q_- + O\big(t^{-\frac 12}\big),\quad t\to\infty.
} 
This result shows that the byproduct of the interaction of the plane wave emerging for $\xi<\Vs$ with the soliton arising for $\xi=\Vs$ is the constant phase shift $2\widetilde g_\infty$  for $\Vs<\xi<\Vo$. In fact, switching to the uniformization variable $z(k) = k + \lambda(k)$, we can compute the integral \eqref{ds-esc-ginft-def-pw2}  via Cauchy's residue theorem and thereby obtain $\widetilde g_\infty$ in the explicit form 
\eee{
\widetilde g_\infty=2\text{arg}\left[p+\lambda(p)\right],
}
which corresponds to a phase shift of $4\text{arg}\left[p+\lambda(p)\right]$ for the plane wave \eqref{ds-esc-qasym-pw2}, in perfect agreement with the inverse scattering transform result of \cite{bk2014}. In turn, the leading-order asymptotics \eqref{ds-esc-qasym-pw2} assume the form \eqref{ds-esc-qasym-pw2-t} of Theorem \ref{ds-per-t}.

\subsection{The range $\Vo<\xi<0$: modulated elliptic wave}
\label{ds-esc-mew-ss}
\begin{figure}[t!]
\begin{center}
\includegraphics[scale=1]{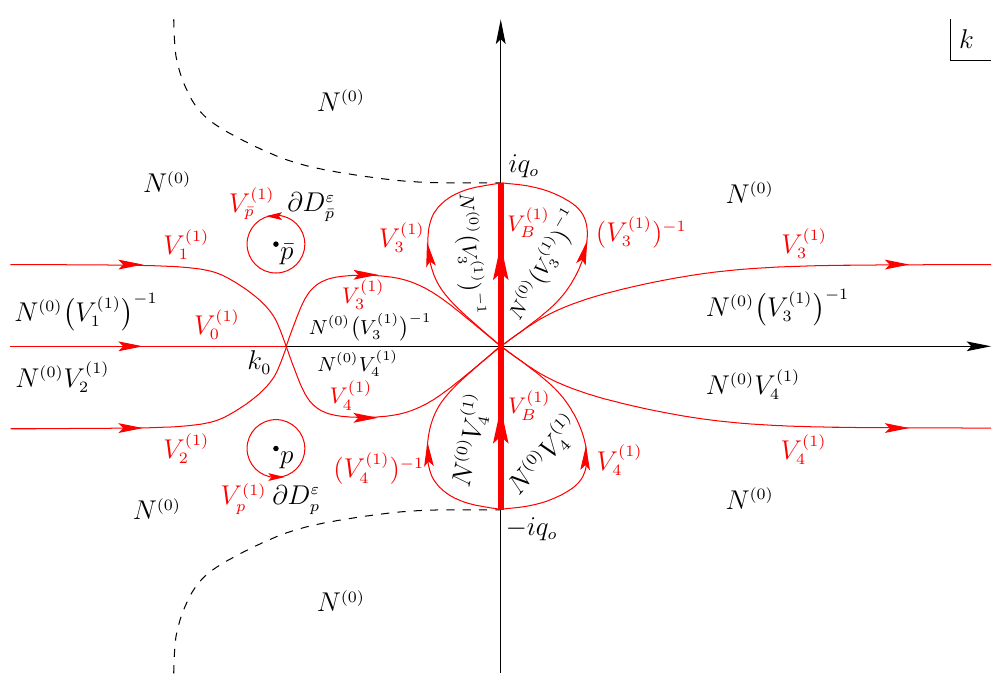}
\caption{Modulated elliptic wave region in the  transmission regime: the first stage of the first deformation. Since $p\in D_1$, it is possible to choose the disks $\overline{D_p^\ve}$ and $\overline{D_{\bar p}^\ve}$ to always lie in regions where $N^{(1)} = N^{(0)}$.}
\label{ds-esc-def1a-mew-f}
\end{center}
\end{figure}

In this range, the stationary points $k_1$ and $k_2$ of the phase function $\theta$ are complex (recall \eqref{ds-k1-k2-def}). 
This has a direct impact on the asymptotic analysis of Riemann-Hilbert problem \eqref{ds-n-rhp-intro}, since the deformations used for $\xi<\Vo$ (where $k_1$ and $k_2$ are real) are no longer effective. 
\vskip 3mm
\noindent
\textbf{First, second and third deformation.}
The first deformation consists of switching from the solution $N^{(0)}$ of problem \eqref{ds-n-rhp-intro} to the function $N^{(1)}$ defined in terms of $N^{(0)}$ by Figure \ref{ds-esc-def1a-mew-f}. This step is very similar to the first stage of the first deformation  in the plane wave region $\xi<\Vs$ (recall Figure \ref{ds-esc-def1a-pw1-f}) apart from the fact that the change of factorization of the jump along the real axis now takes place at the  point $k_o$, which is yet to be determined, instead of the stationary point $k_1$. 
The remaining three stages of the first deformation that were performed for $\xi<\Vs$ (recall Figures \ref{ds-esc-def1c-pw1-f}-\ref{ds-esc-def1b-pw1-f}) can also be carried out here, eventually allowing us to lift the jump contours $L_j$, $j=1,2,3,4,$ away from the origin as well as from the branch points $\pm i q_o$.
Importantly,  the fact that $p \in D_1$ allows us  to perform the first deformation  without modifying the jumps along $\p D_p^\ve$ and $\p D_{\bar p}^\ve$, since the various transformations can be adjusted so that the disks $\overline{D_p^\ve}$ and $\overline{D_{\bar p}^\ve}$  always lie in regions where $N^{(1)} = N^{(0)}$.

\begin{figure}[t!]
\begin{center}
\includegraphics[scale=1]{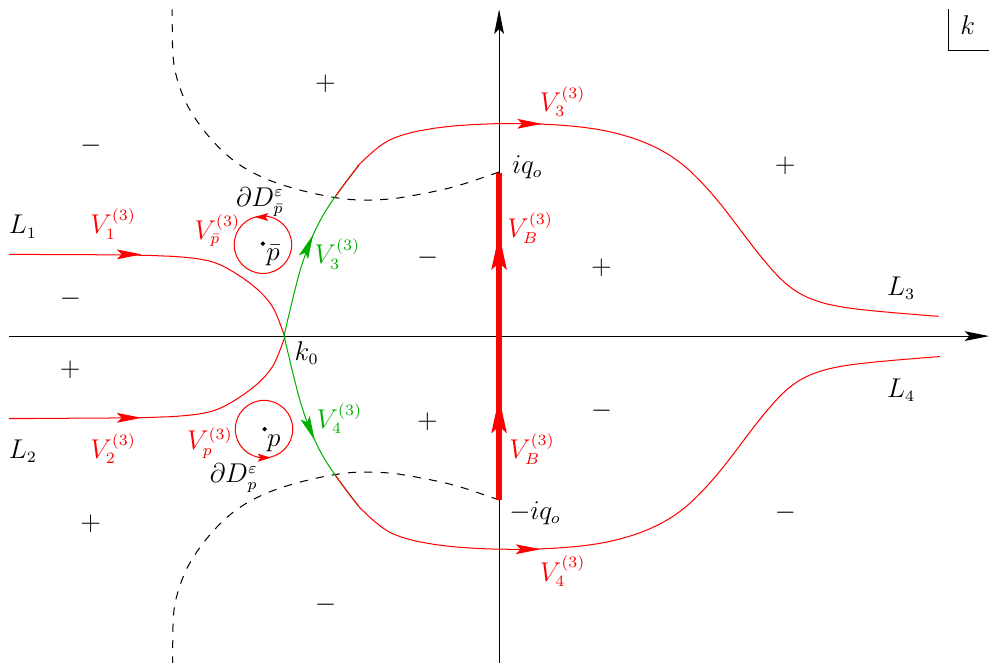}
\caption{Modulated elliptic wave region in the  transmission regime:  the jumps of $N^{(3)}$. As $t\to\infty$, the jumps $V_3^{(3)}$ and $V_4^{(3)}$  grow exponentially along the parts of the contours $L_3$ and $L_4$ connecting $k_o$ with the curve $\text{Re}(i\theta) = 0$ (dashed). Moreover, like in the first deformation (see Figure \ref{ds-esc-def1a-mew-f}), the deformed contours do not interfere with the disks $\overline{D_p^\ve}$ and $\overline{D_{\bar p}^\ve}$, leaving the corresponding jumps unaffected.}
\label{ds-esc-def3-mew-f}
\end{center}
\end{figure}

The second and the third deformation are identical to the corresponding ones in the plane wave region $\xi<\Vs$, leading to the function $N^{(3)}(x, t, k)$, which is analytic in $\mathbb C\setminus\big(\bigcup_{j=1}^4 L_j \cup  B \cup \p D_p^\ve \cup \p D_{\bar p}^\ve\big)$ and satisfies the Riemann-Hilbert problem
\sss{\label{ds-esc-n3t-rhp-mew}
\ddd{
 N^{(3)+} &=  N^{(3)-} V_B^{(3)},\quad && k\in B,
\\
 N^{(3)+} &=  N^{(3)-} V_j^{(3)}, && k\in L_j,\ j=1,2,3,4,
\\
 N^{(3)+} &=  N^{(3)-} V_p^{(3)}, && k\in  \p D_p^\ve,
\\
 N^{(3)+} &= N^{(3)-} V_{\bar p}^{(3)}, && k\in  \p D_{\bar p}^\ve,
\\
 N^{(3)} &=I +O\left(\tfrac 1k\right), && k \to \infty,
}
}
where the contours $L_j$ are depicted in Figure \ref{ds-esc-def3-mew-f} and the relevant jumps are given by \eqref{ds-n3-jumps} and \eqref{ds-vp3-vpb3-def-pw1} but with the function $\delta(\xi, k)$ now modified to
\eee{\label{ds-esc-del-def-mew}
\delta(\xi, k)
=
\exp
\left\{
\frac{1}{2i\pi}
\int_{-\infty}^{k_o(\xi)}
\frac{\ln\!\big[1+r(\nu)\bar r(\nu)\big]}{\nu -k}\, d\nu 
\right\},\quad k\notin(-\infty, k_o).
}

Importantly, we note that the jump $V_3^{(3)}$ \textit{grows exponentially} as $t\to\infty$ along the portion of the contour  $L_3$ colored in green in Figure \ref{ds-esc-def3-mew-f}, i.e. along the portion of $L_3$ that connects $k_o$ with the dashed curve $\text{Re}(i\theta)=0$ lying in the second quadrant of the complex $k$-plane.
The same is true for the jump $V_4^{(3)}$ and the green-colored portion of the contour $L_4$ that joins $k_o$ with the dashed curve $\text{Re}(i\theta)=0$  in the third quadrant of the complex $k$-plane.
This growth, which was not present for $\xi<\Vo$, can be removed with the help of appropriate factorizations and a \textit{time-dependent} version of transformation \eqref{ds-esc-n4-def-pw1}, as shown in the course of the following two deformations.

\vskip 3mm
\noindent
\textbf{Fourth deformation.}
The jumps $V_3^{(3)}$ and $V_4^{(3)}$ can be factorized in the form 
\eee{\label{ds-gr-facs}
V_3^{(3)} = V_5^{(4)} V_7^{(4)} V_5^{(4)},
\quad
V_4^{(3)} = V_6^{(4)} V_8^{(4)} V_6^{(4)},
}
where
\ddd{\label{ds-lense-jumps}
&V_5^{(4)}
=
\left(
\def\arraystretch{1.2}
\begin{array}{lr}
1 &\frac{\delta^{2}}{r}\, e^{2i \theta t } 
\\
0 &1
\end{array}
\right),
&&V_6^{(4)}
=
\left(
\def\arraystretch{1.2}
\begin{array}{lr}
1 & 0
\\
\frac{1}{ \bar r\delta^{2}}\, e^{-2i \theta t}  & 1
\end{array}
\right),
\nn\\
&V_7^{(4)}
=
\left(
\def\arraystretch{1.2}
\begin{array}{lr}
0 & -\frac{\delta^{2}}{r}\, e^{2i \theta t } 
\\
\frac{r}{\delta^{2}}\, e^{-2i \theta t } 
&
0
\end{array}
\right),
\quad
&&V_8^{(4)}
=
\left(
\def\arraystretch{1.2}
\begin{array}{lr}
0
&
\bar r\delta^{2} e^{2i \theta t } 
\\
-\frac{1}{\bar r \delta^{2}}\, e^{-2i \theta t } 
&
0
\end{array}
\right).
}
The advantage of the above factorization is that the matrices $V_5^{(4)}$ and $V_6^{(4)}$ each involve only one exponential and hence they have a definitive behavior as $t\to\infty$. In particular, in this limit $V_5^{(4)}$ and $V_6^{(4)}$  tend to the identity in regions of negative and positive sign of $\text{Re}(i\theta)$ respectively.
On the other hand, the matrices $V_7^{(4)}$ and $V_8^{(4)}$ still involve both of the exponentials $e^{\pm 2i\theta t}$ and so it is not possible to take their limit as $t\to\infty$. However, contrary to the original jumps $V_3^{(3)}$ and $V_4^{(3)}$, the matrices $V_7^{(4)}$ and $V_8^{(4)}$ are antidiagonal. This fact turns out to be crucial, as we will see in the fifth deformation below.

\begin{figure}[t]
\begin{center}
\includegraphics[scale=1]{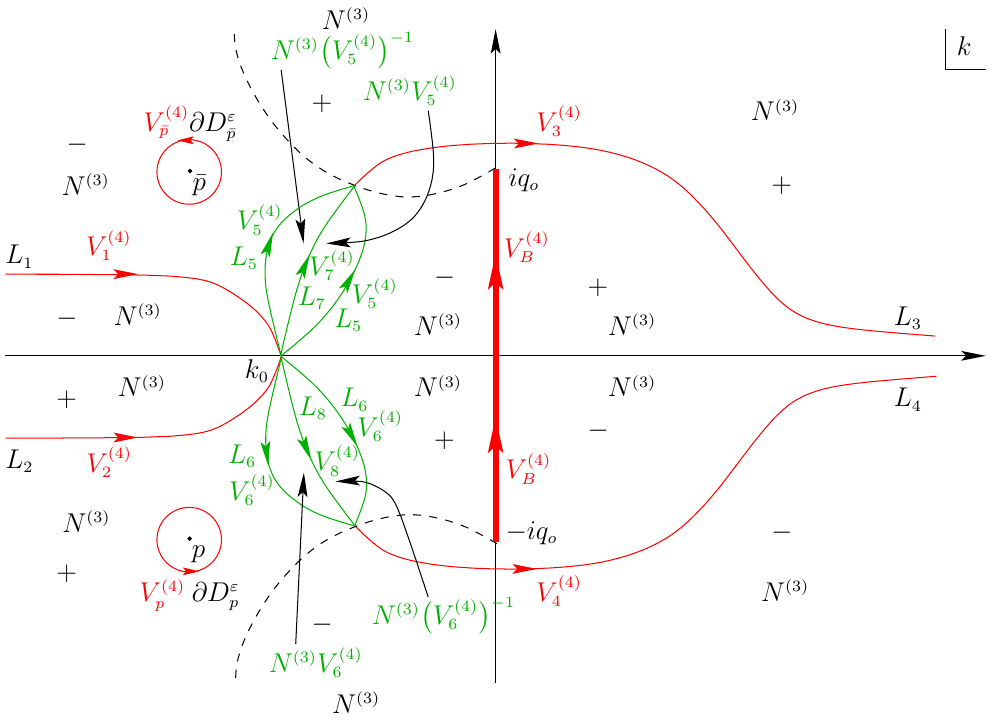}
\caption{Modulated elliptic wave region in the transmission regime:  the fourth deformation. The various jump contours have been chosen so as not to interfere with the disks 
$\overline{D_p^\ve}$ and $\overline{D_{\bar p}^\ve}$.}
\label{ds-esc-def4-mew-f}
\end{center}
\end{figure}

Using the factorization \eqref{ds-gr-facs}, we switch from $N^{(3)}$ to $N^{(4)}$ as shown in Figure \ref{ds-esc-def4-mew-f}. By this definition, $N^{(4)}(x, t, k)$ is analytic in $\mathbb C\setminus \big(\bigcup_{j=1}^4 L_j \cup  B \cup \p D_p^\ve \cup \p D_{\bar p}^\ve\big)$ and  satisfies the following Riemann-Hilbert problem:
\sss{\label{ds-esc-n3t-rhp-mew-2}
\ddd{
 N^{(4)+} &=  N^{(4)-} V_B^{(4)},\quad && k\in B,
\\
 N^{(4)+} &=  N^{(4)-} V_j^{(4)}, && k\in L_j,\ j=1, \ldots, 8,
\\
 N^{(4)+} &=  N^{(4)-} V_p^{(4)}, && k\in  \p D_p^\ve,
\\
 N^{(4)+} &= N^{(4)-} V_{\bar p}^{(4)}, && k\in  \p D_{\bar p}^\ve,
\\
 N^{(4)} &=I +O\left(\tfrac 1k\right), && k \to \infty,
}
}
where the contours $L_j$ are shown in Figure \ref{ds-esc-def4-mew-f}, the jumps $V_5^{(4)}, V_6^{(4)}, V_7^{(4)}, V_8^{(4)}$ are given by \eqref{ds-lense-jumps}, and the remaining jumps are as in problem \eqref{ds-esc-n3t-rhp-mew}.

The only growth surviving in problem \eqref{ds-esc-n3t-rhp-mew-2} after the fourth deformation is located in the $21$- and $12$-elements of the jumps $V_7^{(4)}$ and $V_8^{(4)}$ respectively. The fact that these jumps are antidiagonal allows us to eliminate this growth by employing the following transformation, which is essentially the mechanism leading to a modulated elliptic wave (as opposed to a plane wave).

\begin{figure}[t]
\begin{center}
\includegraphics[scale=1]{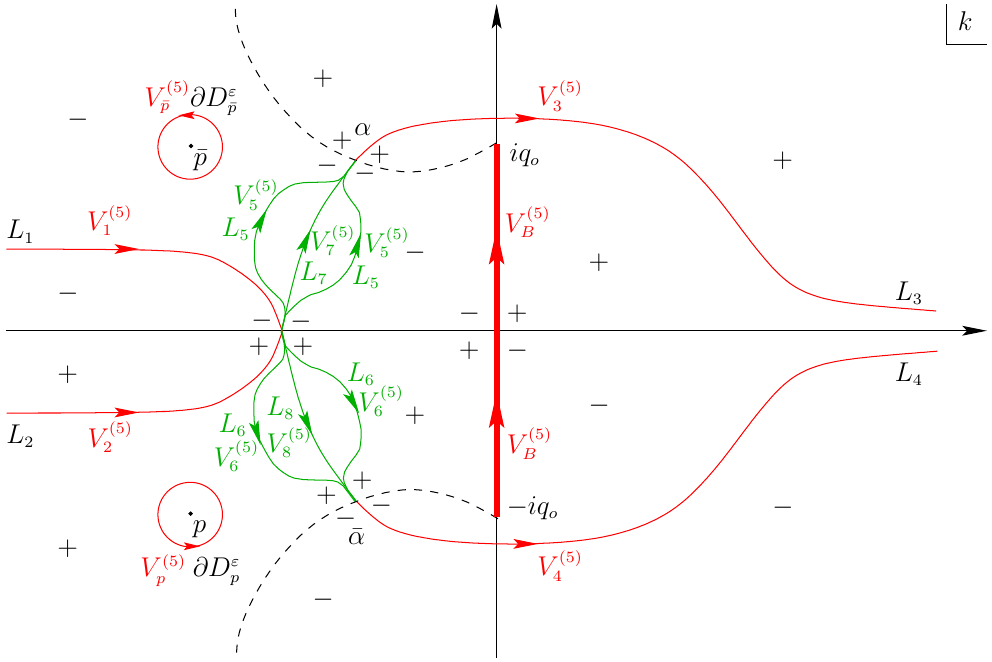}
\caption{Modulated elliptic wave region in the transmission regime:  the fifth deformation.}
\label{ds-esc-def5-mew-f}
\end{center}
\end{figure}

\vskip 3mm
\noindent
\textbf{Fifth deformation.}
Let the points $\alpha(\xi)$ and $\bar \alpha(\xi)$ be defined through the solution of the modulation equations \eqref{ds-mod-eqs-i}, which was shown in \cite{bm2017} to be unique for $\xi\in (\Vo, 0)$. In turn, let the point  $k_o(\xi)$ be given by \eqref{ds-ko-def-i}. Then, introduce the function $h(\xi, k)$ via the Abelian integral \eqref{ds-h-abel-i}. Note that $h$ involves the function $\gamma(\xi, k)$  defined by \eqref{ds-gamma-def-i}, which is made single-valued by taking  branch cuts along $B$ as well as along the curve
\eee{
\widetilde B:= L_7\cup (-L_8)
}
with the contours $L_7$ and $L_8$ depicted in Figure \ref{ds-esc-def5-mew-f}. We emphasize that  $\widetilde B$ must begin at $\bar \alpha$ and end at $\alpha$ by crossing the negative real axis at the point $k_o$  but is otherwise arbitrary for the moment. The function $\gamma$ is analytic for $k\in \mathbb C\setminus  B\cup \widetilde B$ and changes sign as $k$ crosses $B$ and $\widetilde B$. Via \eqref{ds-h-abel-i}, this induces analyticity of $h$ in $\mathbb C\setminus  B\cup \widetilde B$ as well as the following jump conditions along $B$ and $\widetilde B$:
\sss{\label{ds-rhp-h}
\ddd{
& h^+ + h^- = 0, && k\in B,\label{ds-jhL}
\\
&h^+ + h^-=\Omega, \quad &&k\in L_7\cup L_8,
\label{ds-jhLL}
}
}
where the real quantity $\Omega(\xi)$, which is independent of $k$, is defined by
\eee{\label{ds-Om-def}
\Omega(\xi) 
=
-4 \left(\int_{iq_o }^{\alpha}+\int_{-iq_o }^{\bar \alpha} \right)
\frac{\left(z-k_{o}\right)\left(z-\alpha\right)\left(z-\bar \alpha\right)}{ \gamma(z)}\, dz.
}
Moreover,  as shown in \cite{bm2017}, $\text{Re}(ih)$ has the same sign with $\text{Re}(i\theta)$ at infinity, near the origin, and near $\alpha$ and $\bar \alpha$.

The definition of $h$ and, more specifically, the jump conditions \eqref{ds-rhp-h} imply that the jumps of the function
\eee{\label{ds-m4r-0}
 N^{(5)}(x, t, k)
=
 N^{(4)}(x, t, k) e^{-i\left[h(\xi, k)-\theta(\xi,k)\right]t \sigma_3}
 }
along the contours $L_7$ and $L_8$ are bounded. Furthermore, those jumps that were bounded at the level of $N^{(4)}$ remain bounded at the level of  $N^{(5)}$ (see discussion below \eqref{ds-Ginf-mew}). 
Specifically, Riemann-Hilbert problem \eqref{ds-esc-n3t-rhp-mew-2} and transformation \eqref{ds-m4r-0} imply that  $N^{(5)}(x, t, k)$ is analytic for $k\in \mathbb C \setminus  \big(\bigcup_{j=1}^8 L_j \cup  B \cup \p D_p^\ve \cup \p D_{\bar p}^\ve\big)$ and satisfies the  jump conditions
\sss{\label{ds-n4t-rhp-mew}
\ddd{
 N^{(5)+} &=  N^{(5)-} V_B^{(5)},  && k\in B,
\\
 N^{(5)+} &=  N^{(5)-} V_j^{(5)}, && k\in L_j,\ j=1, \ldots, 8,
\\
 N^{(5)+} 
&=
 N^{(5)-} 
 V_p^{(5)}, && k\in  \p D_p^\ve,
\\
 N^{(5)+} 
&=
 N^{(5)-} 
 V_{\bar p}^{(5)}, && k\in  \p D_{\bar p}^\ve,
\\
 N^{(5)} &=
\left[I+O\left(\tfrac 1k\right)\right] e^{-iG_\infty(\xi) t  \sigma_3},\quad  && k \to \infty,
}
}
where the contours $L_j$ are  shown in Figure \ref{ds-esc-def5-mew-f} and 
\ddd{\label{ds-n5-jumps-mew}
&V_B^{(5)}
=
\left(
\def\arraystretch{1.2}
\begin{array}{lr}
0 & \tfrac{q_-}{iq_o} \delta^2  \\
\tfrac{\bar q_-}{iq_o} \delta^{-2}  & 0
\end{array}
\right), 
\quad
V_1^{(5)}
=
\left(
\def\arraystretch{1.2}
\begin{array}{lr}
1 & \frac{\bar r \delta^2}{1+r\bar r}  \, e^{2iht} \\
0 & 1
\end{array}
\right),
\quad
V_2^{(5)}
=
\left(
\def\arraystretch{1.2}
\begin{array}{lr}
1 &  0 \\
\frac{r\delta^{-2}}{1+r\bar r} e^{-2i h t} & 1
\end{array}
\right),
\nn\\
&V_3^{(5)}
=
\left(
\def\arraystretch{1.2}
\begin{array}{lr}
1 & 0 \\
r \delta^{-2} e^{-2i h t} & 1
\end{array}
\right),
\quad
V_4^{(5)}
=
\left(
\def\arraystretch{1.2}
\begin{array}{lr}
1
&
\bar r  \delta^2 e^{2i h t}
\\
0
&
1
\end{array}
\right),
\quad
V_5^{(5)}
=
\left(
\def\arraystretch{1.2}
\begin{array}{lr}
1 &\frac{\delta^{2}}{r}\, e^{2i h t } 
\\
0 &1
\end{array}
\right),
\nn\\
&V_6^{(4)}
=
\left(
\def\arraystretch{1.2}
\begin{array}{lr}
1 & 0
\\
\frac{1}{ \bar r\delta^{2}}\, e^{-2i h t}  & 1
\end{array}
\right),
\quad
V_7^{(5)}
=
\left(
\def\arraystretch{1.2}
\begin{array}{lr}
0 & -\frac{\delta^{2}}{r}\, e^{i \Omega t} 
\\
\frac{r}{\delta^{2}}\, e^{-i \Omega t} 
&
0
\end{array}
\right),
\quad
V_8^{(5)}
=
\left(
\def\arraystretch{1.2}
\begin{array}{lr}
0
&
\bar r\delta^{2} e^{i \Omega t} 
\\
-\frac{1}{\bar r \delta^{2}}\, e^{-i\Omega t} 
&
0
\end{array}
\right),
\nn\\
&V_p^{(5)} 
=
\left(
\def\arraystretch{1}
\begin{array}{lr}
1 & -\frac{c_p \, \delta^2(\xi, k)\, d(k)}{k-p}\,  e^{2i \left[h(\xi, k) + \theta(\xi, p) - \theta(\xi, k)\right]t}
\\
0 & 1
\end{array}
\right),
\nn\\
&V_{\bar p}^{(5)} 
=
\left(
\def\arraystretch{1}
\begin{array}{lr}
1 &0
\\
-\frac{c_{\bar p} \, \delta^{-2}(\xi, k) \, d(k)}{k-\bar p}\,  e^{-2i\left[h(\xi, k) + \theta(\xi, \bar p) - \theta(\xi, k)\right]t} & 1
\end{array}
\right),
}
with $\Omega(\xi)$ defined by \eqref{ds-Om-def} and the real quantity $G_\infty(\xi)$ given by
\eee{\label{ds-Ginf-mew}
G_\infty(\xi) 
=
-2 \left(\int_{iq_o }^\infty+\int_{-iq_o }^\infty \right)
\left[
\frac{\left(z-k_o\right)\left(z-\alpha\right)\left(z-\bar \alpha\right)}{\gamma(z)}
-
\left(z-\frac \xi 4\right)
\right] dz
-q_o^2.
}

\begin{remark}
The fourth deformation, which affects only the jumps along the green contours of Figure \ref{ds-esc-def3-mew-f} by opening those contours into the lenses comprising the contours $L_5$, $L_6$, $L_7$, $L_8$ of Figure \ref{ds-esc-def4-mew-f}, could also be performed \textit{after} the $g$-function deformation \eqref{ds-m4r-0}, leading  again to Riemann-Hilbert problem \eqref{ds-n4t-rhp-mew}.
Nonetheless, the order we have followed here has the advantage of revealing the basic form of the jumps along the contours of growth (see \eqref{ds-lense-jumps}) \textit{before} the introduction of the Abelian function $h$, thus allowing us to better motivate the desired properties that eventually lead to the definition of $h$. Of course, since for $\Vo < \xi < 0$ we switch the phase function from $\theta$ to $h$ via deformation \eqref{ds-m4r-0}, it should be emphasized that the contours $L_7$ and   $L_8$ of Figure~\ref{ds-esc-def5-mew-f} are not those of Figure~\ref{ds-esc-def4-mew-f}, but rather the contours connecting $\alpha$ and $\bar \alpha$ with $k_o$.
\end{remark}

The sign structure of $\text{Re}(ih)$ at infinity and near the origin together with the fact that, by definition, $h$ possesses precisely three critical points, namely $k_o$, $\alpha$, $\bar \alpha$, guarantee the existence of a neighborhood around the point $k_o$ where $\text{Re}(ih)<0$ in the second quadrant, $\text{Re}(ih)>0$ in the third quadrant, and $\text{Re}(ih)=0$ along $\mathbb R$ and the branch cut $\widetilde B$.\footnote{Indeed, if such a neighborhood did not exist then there would have to be one or more saddle points other than $k_o$ along the negative real axis, leading to a contradiction.} 
Thus, thanks to the sign structure of $\text{Re}(ih)$ near $\alpha$, it is always possible to have a contour from $k_o$ to $\alpha$ which lies to the right of the branch cut $\widetilde B$ and along which $\text{Re}(ih)<0$.\footnote{\label{footnote3}Indeed, the only way this could fail is if there were a  strip of  $\text{Re}(ih)>0$ connecting the branch cut $\widetilde B$ either with the real axis or with the branch cut $B$. The first scenario is not possible because if would create additional critical (saddle) points on the real axis. Moreover, the second scenario is also not realizable since, due to the continuity of $\text{Re}(ih)$ away from the cuts and the jump condition $\text{Re}(ih^+) = -\text{Re}(ih^-)$, the boundary of the strip away from $\widetilde B$ would have to be a zero-contour, i.e. a contour along which $\text{Re}(ih)=0$. But then, deforming $\widetilde B$ to this zero-contour we would get inconsistent jump conditions for $h$ along the part of the zero-contour that overlaps with $B$, since $\Omega\neq 0$ independently of the branch cuts $B$ and $\widetilde B$.} 
Similarly, thanks to the sign structure of $\text{Re}(ih)$ at infinity and near $k_o$ and $\alpha$, as well as to the fact that $h$ possesses precisely three critical points, we can always find a contour from $k_o$ to $\alpha$ which lies to the left of the branch cut $\widetilde B$ and along which $\text{Re}(ih)<0$.
Therefore, the contour $L_5$  of Figure \ref{ds-esc-def5-mew-f} can always be    chosen to satisfy $\text{Re}(ih)<0$. In turn, by symmetry, the contour $L_6$  can always be chosen to satisfy $\text{Re}(ih)>0$. Hence, the jumps $V_5^{(5)}$ and $V_6^{(5)}$ are guaranteed to decay exponentially to the identity  as $t\to\infty$. 
Analogous considerations ensure that $\text{Re}(ih)$ has the same sign structure with $\text{Re}(i\theta)$ along the contours $L_1$, $L_2$, $L_3$ and $L_4$ of Figure  \ref{ds-esc-def5-mew-f}. Thus, the jumps $V_1^{(5)}$, $V_2^{(5)}$, $V_3^{(5)}$ and $V_4^{(5)}$ inherit the behavior of their $N^{(4)}$-counterparts, i.e. they decay exponentially to the identity  as $t\to\infty$. 

Furthermore,  there exists at least one  zero-contour (i.e. a contour along which $\text{Re}(ih)=0$) connecting $\alpha$ and $\bar \alpha$ through $k_o$. This is because of the sign structure of $\text{Re}(ih)$ near $\alpha$ and $\bar \alpha$ as well as due to the jump condition \eqref{ds-jhLL} along $\widetilde B$, which implies that $\text{Re}(ih^+) = - \text{Re}(ih^-)$ since $\Omega\in\mathbb R$. Thus, either $\widetilde B$ is itself a zero-contour, or there exists a region of positive sign adjacent to $\widetilde B$ whose boundary will have to be a zero-contour due to the analyticity of $h$, the sign of $\text{Re}(ih)$ at infinity and near the origin, and the existence of precisely three critical points of $h$.
In the latter case, we can deform $\widetilde B$ to this zero-contour so that $\Re (ih)=0$ \textit{throughout} $\widetilde B$. 
In fact, any zero-contour connecting $\alpha$ and $\bar \alpha$ can only cross the negative real axis at $k_o$, since a zero-contour intersecting with the negative real axis at a point different than $k_o$ would require this point to a critical point, leading to a contradiction. 
Therefore, taking  into account once again the sign structure of $\text{Re}(ih)$ near $\alpha$ and $\bar \alpha$, we conclude that for $\xi\in (\Vo, 0)$ there exists a \textit{unique} zero-contour with endpoints   $\alpha$ and $\bar \alpha$ and  through the point $k_o$, namely the branch cut $\widetilde B$.

Furthermore, as $\xi$ increases from $\Vo$ to $0$ the branch cut $\widetilde B$ remains within the finite region enclosed by the trace of $\bar \alpha$  and $\alpha$ (the dashed green curve in Figure \ref{ds-regions-f} and its reflection through the real axis) and the branch cut $B$. Hence, for $p\in D_1$, as $\xi$ increases from $\Vo$ to $0$  the branch cut $\widetilde B$ remains to the right of the disks $\overline{D_p^\ve}$ and $\overline{D_{\bar p}^\ve}$ without interfering with $p$ and $\bar p$. The same is true for $p\in D_2^+$ since, as shown in Figure \ref{ds-regions-f}, this region lies by definition below the trace of $\bar \alpha$ (while $\widetilde B$ lies above that trace). On the other hand, if $p\in D_2^-\cup D_3$ then both $p$ and $\bar p$ are crossed by   $\widetilde B$ for some $\xi\in(\Vo, 0)$, this being the mechanism that generates the soliton wake   in the mixed regimes of Section \ref{ds-mixed-s}.

Next, recall that the transition from $\theta$ to $h$ in the jump matrices takes place at $\xi=\Vo$, where $h(\Vo, k)=\theta(\Vo, k)$ and $\alpha(\Vo)=\bar \alpha(\Vo)=k_o(\Vo)=\Vo/8$. Hence, at $\xi=\Vo$ the lower and upper dashed  curves of Figure \ref{ds-esc-def5-mew-f} are, respectively, the solid blue curve of Figure \ref{ds-regions-f} and its reflection through the real $k$-axis. 
A numerical investigation then shows that, as $\xi$ increases from $\Vo$ to $0$, the upper dashed curve remains convex and moves continuously upwards and to the right, eventually collapsing to the half-line $i[q_o, \infty)$. Analogously, the lower dashed curve remains concave and moves downwards and to the right, eventually collapsing to the half-line $i(-\infty, -q_o]$.
Hence, no points inside the regions $D_1$ and $D_3$ of Figure \ref{ds-regions-f} are crossed by the lower dashed curve of Figure \ref{ds-esc-def5-mew-f}  as $\xi$ increases from $\Vo$ to $0$. 
In particular, if $p\in D_1\cup D_3$, then the circle $\p D_p^\ve$, which is between the negative real $k$-axis the solid blue curve of Figure \ref{ds-regions-f} at $\xi=\Vo$, remains below the negative real $k$-axis and above the lower dashed curve of Figure \ref{ds-esc-def5-mew-f}  throughout the range $(\Vo, 0)$.
On the other hand, if $p\in D_2^+\cup D_2^-$, then it will be crossed at exactly one value of $\xi\in(\Vo, 0)$ by the lower dashed curve of Figure \ref{ds-esc-def5-mew-f}.
An analogous statement can be made for $\bar p$ by symmetry. 

The dashed  curves of Figure \ref{ds-esc-def5-mew-f} together with the branch cut $\widetilde B$ and the negative real axis make up the contours along which $\text{Re}(ih)=0$ on the left half of the complex $k$-plane. Hence, from the above-described behavior of these contours as $\xi$ increases from $\Vo$ to $0$, we conclude that
\begin{enumerate}[label=$\bullet$, leftmargin=4mm, rightmargin=0mm]
\advance\itemsep 1mm
\item If $p\in D_1$, then $\text{Re}(ih)(\xi, p)>0$ for all $\xi\in (\Vo, 0)$. Equivalently, the equation 
\eee{\label{ds-fnlsd-int-eq-xistar-0} 
\text{Re}(ih)(\xi, p) = 0
\Leftrightarrow
\text{Re}(ih)(\xi, \bar p) = 0
\Leftrightarrow
\int_{\bar p}^p dh(\xi, z) = 0
}
has no solution for $\xi\in (\Vo, 0)$.  Hence, in the  transmission regime $p\in D_1$ no soliton arises in the range $\xi\in (\Vo, 0)$.
\item If $p\in D_2^+$, then equation \eqref{ds-fnlsd-int-eq-xistar-0}  has a \textit{unique} solution $\Vd\in (\Vo, 0)$, which gives rise to a soliton (see Section \ref{ds-trap-s} for more details).
\item If $p\in D_2^-$, then equation \eqref{ds-fnlsd-int-eq-xistar-0} has two solutions in the interval $(\Vo, 0)$: one due to the crossing of $p$ by the lower dashed curve of \ref{ds-esc-def5-mew-f}, denoted by $\Vd$, and another one due to the crossing of $p$ by the branch cut $\widetilde B$, denoted by $\Vw$. Moreover, $\Vd<\Vw$ and the first solution corresponds to a soliton while the second one to a soliton wake (see Section \ref{ds-mixed-s} for more details).
\item Finally, if $p\in D_3$, then equation \eqref{ds-fnlsd-int-eq-xistar-0}  has a \textit{unique} solution $\Vw\in (\Vo, 0)$, which arises from the crossing of $p$ by the branch cut $\widetilde B$ and corresponds to a soliton wake (see Section \ref{ds-mixed-s} for more details).
\end{enumerate}

\begin{figure}[t!]
    \begin{center}
        \includegraphics[scale=0.335]{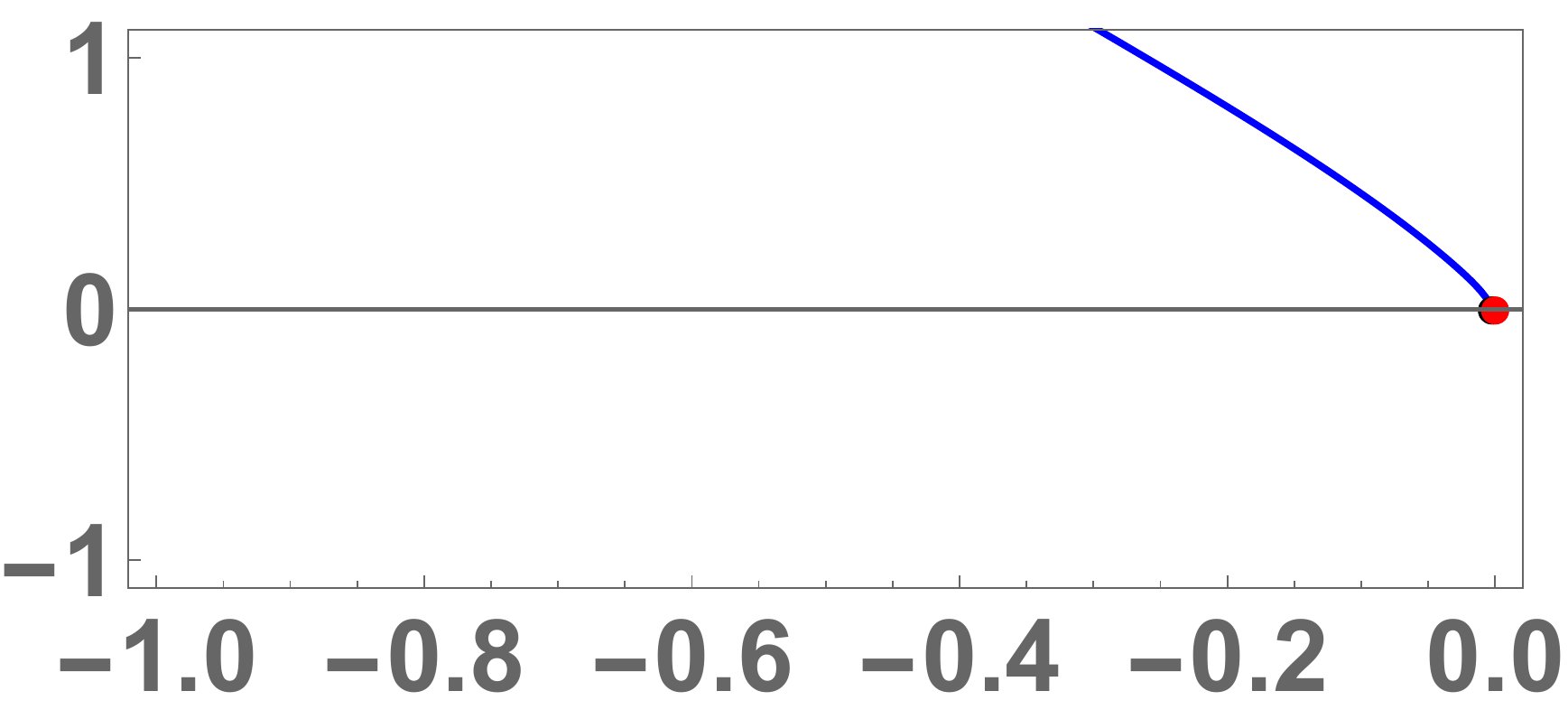}\qquad
        \includegraphics[scale=0.335]{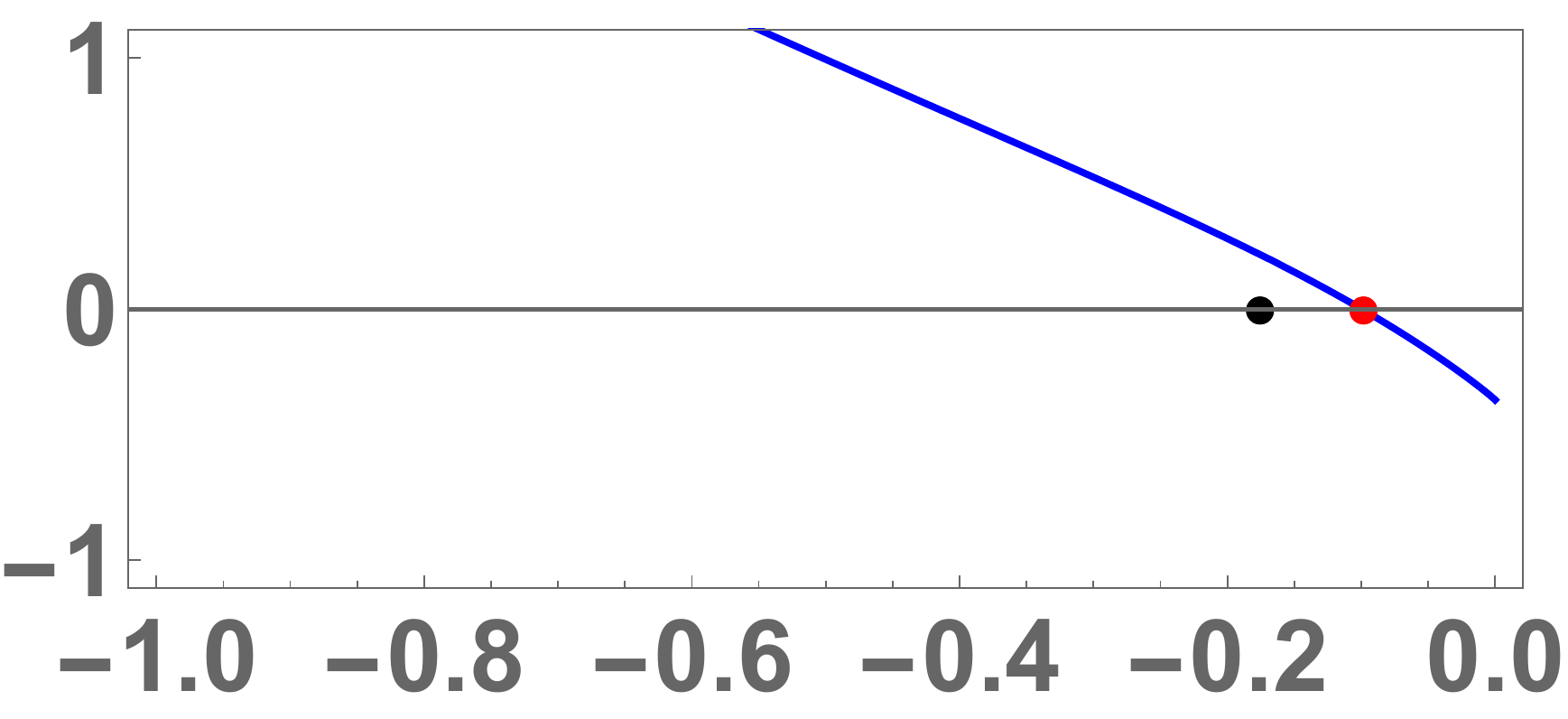}\\[2ex]
        \includegraphics[scale=0.335]{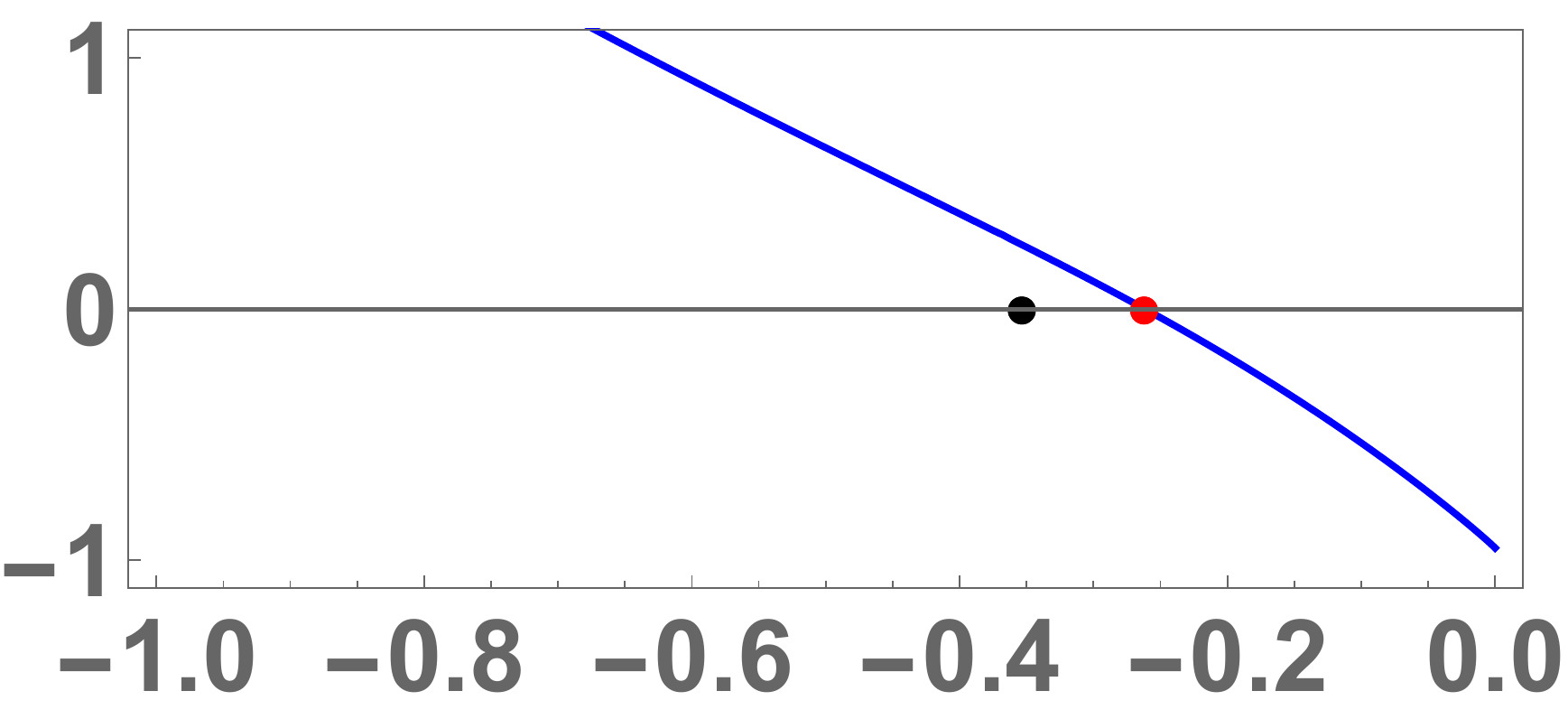}\qquad
        \includegraphics[scale=0.335]{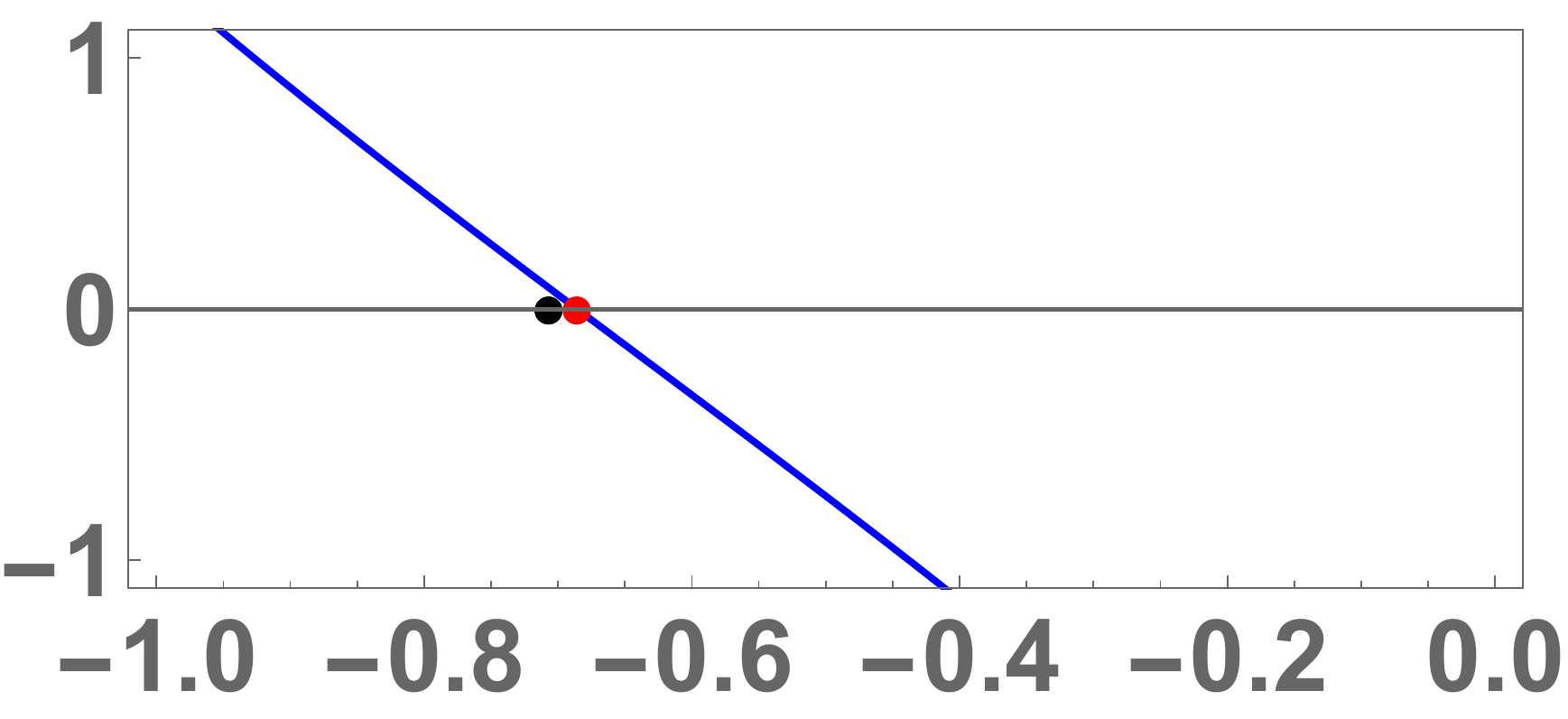}\\[2ex]
        \includegraphics[scale=0.335]{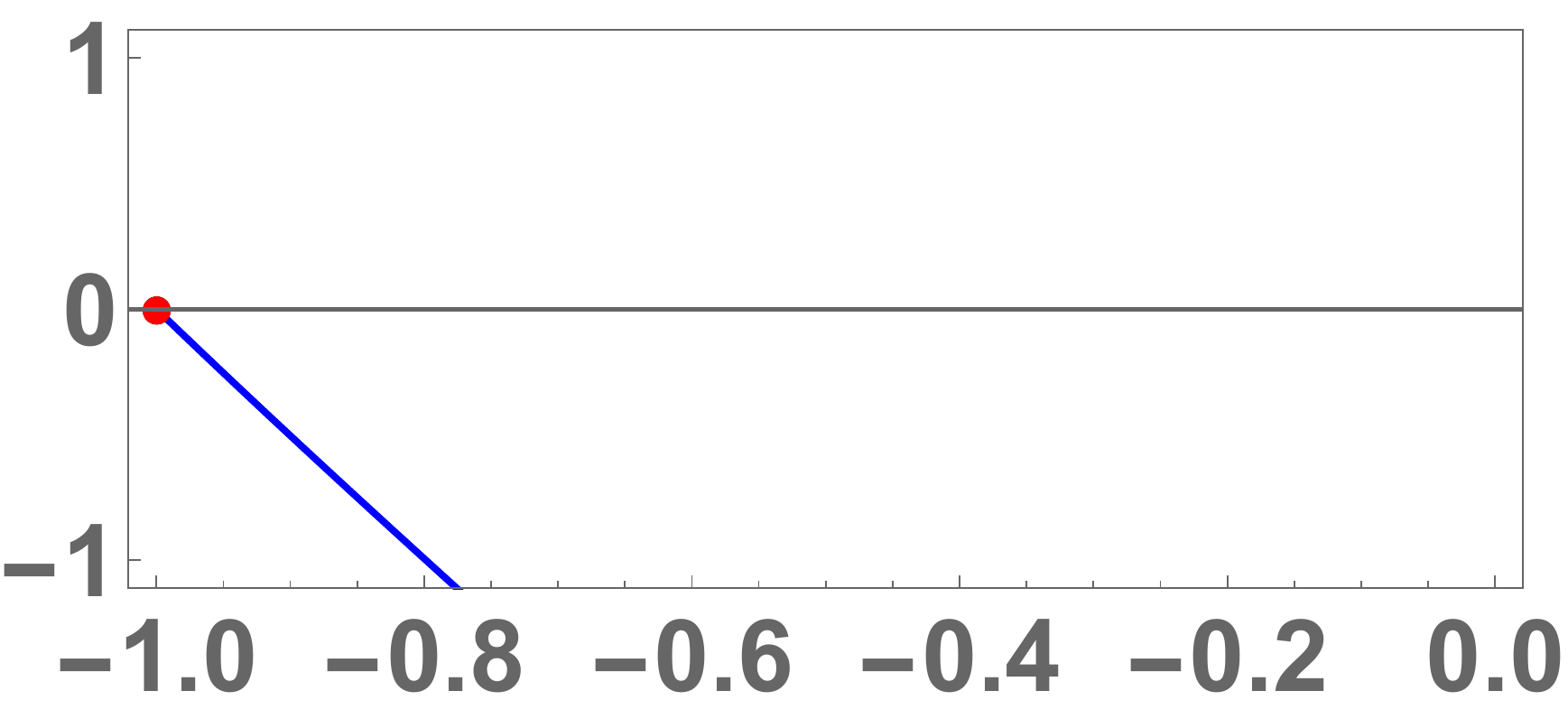}\qquad
        \includegraphics[scale=0.335]{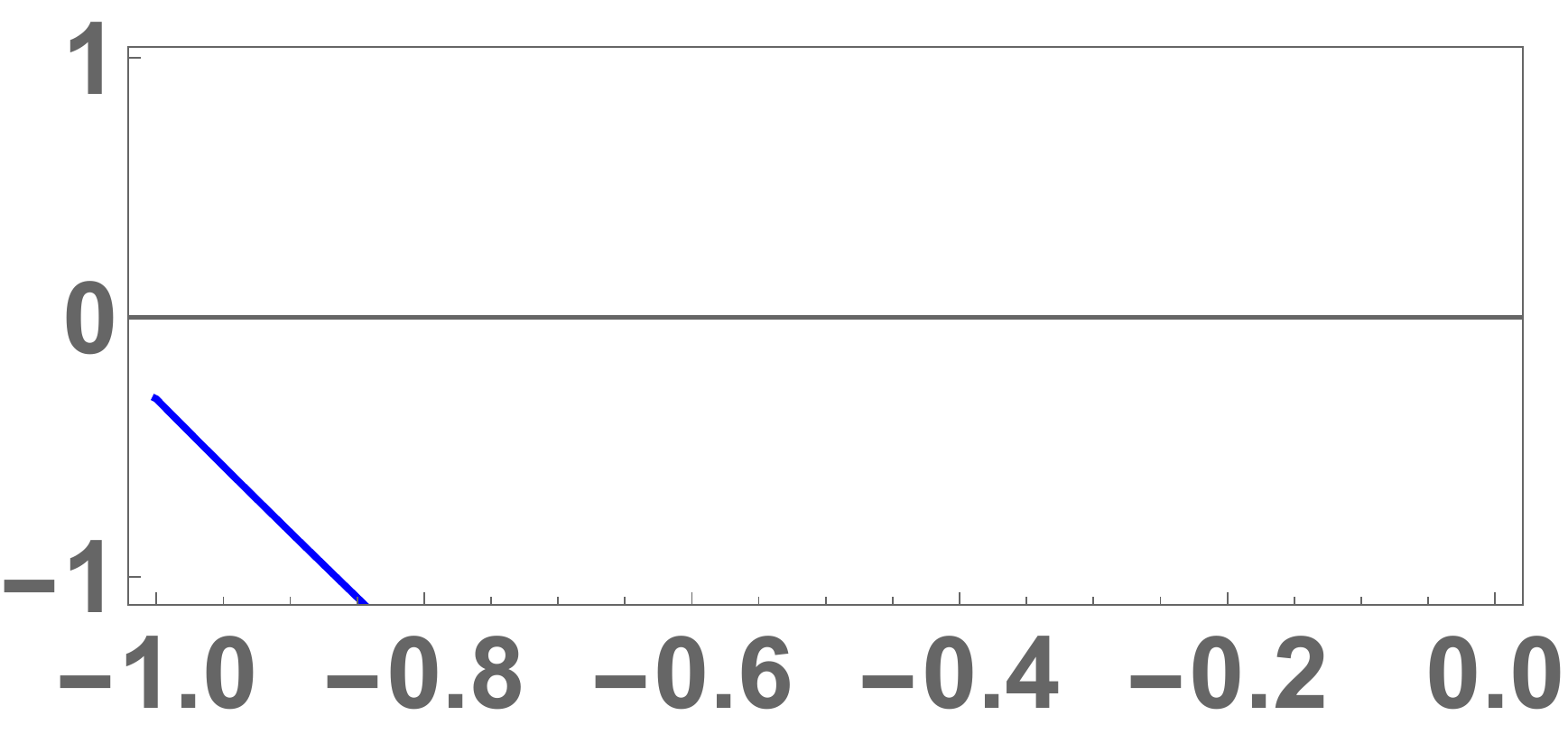}\\[2ex]
        \includegraphics[scale=0.335]{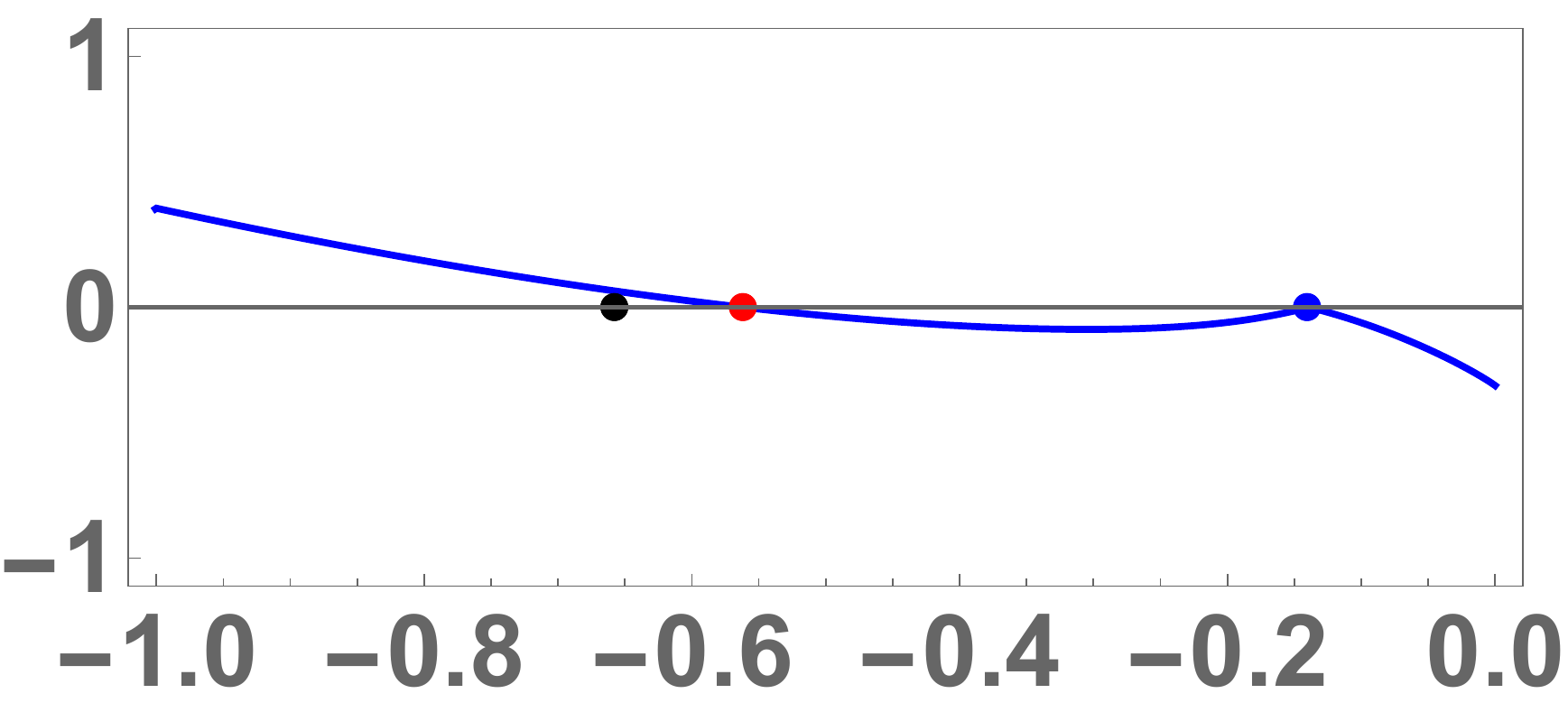}\qquad
        \includegraphics[scale=0.335]{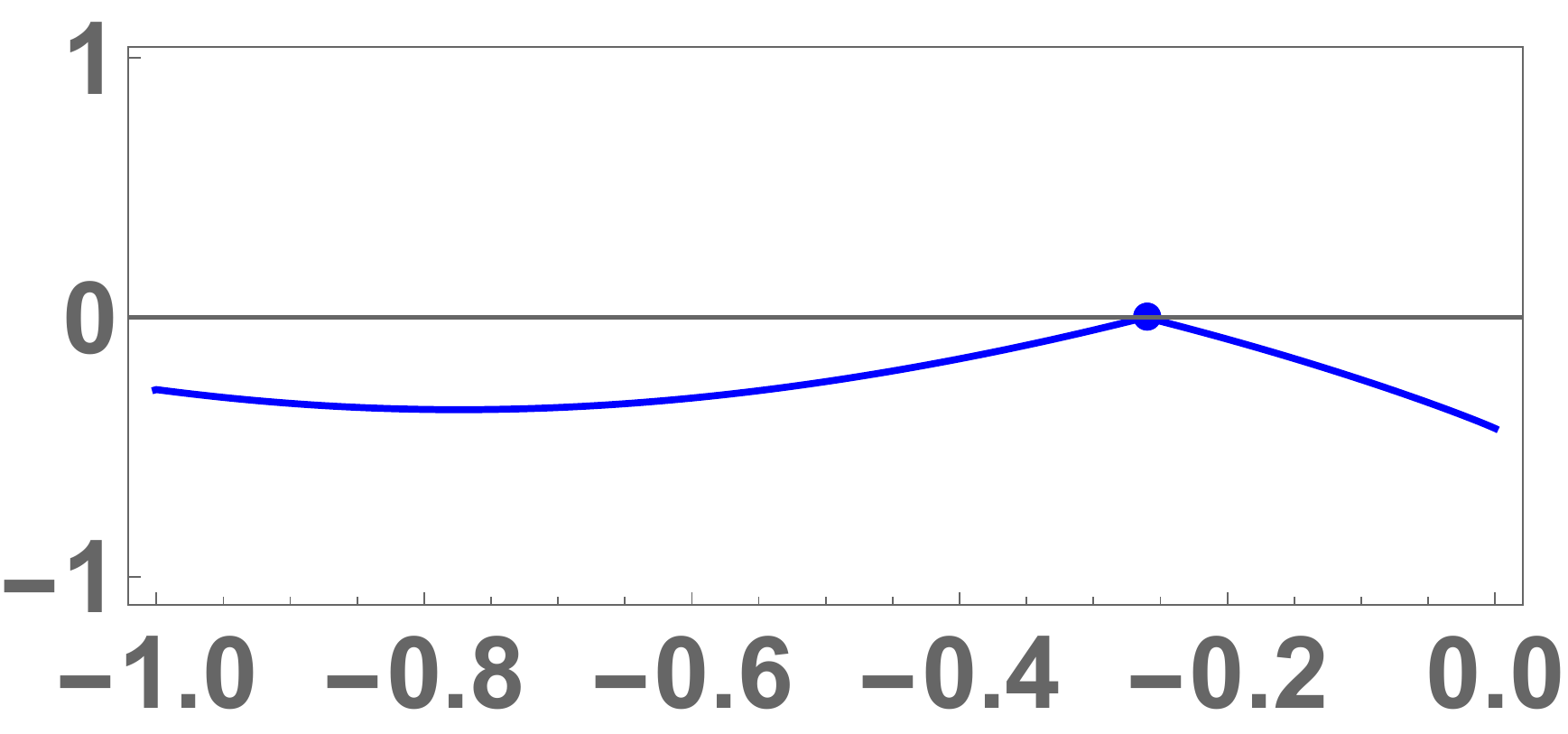}
        \caption{
Numerical evaluation  of the solutions of equation \eqref{ds-fnlsd-int-eq-xistar-0} for $q_o=1$ and the different choices of $p$ depicted by orange dots in Figure \ref{ds-regions-f}.  The horizontal axis corresponds to $\xi$,  the vertical axis to $\text{Re}(ih)(\xi, p)$, and the values on both axes are normalized by a factor of $|\Vo|$ so that the range $[-1, 0]$ of $\xi$ corresponds precisely to the modulated elliptic wave region $[\Vo, 0]$, where equation \eqref{ds-fnlsd-int-eq-xistar-0} is actually relevant.
The black dots denote the unperturbed soliton velocity $\Vs$ given by \eqref{ds-xistar-xip-def}, while the red and blue dots denote the solutions $\Vd$ and $\Vw$ of equation \eqref{ds-fnlsd-int-eq-xistar-0} corresponding to a soliton and a soliton wake respectively.
\textit{First row}:  $p = -1.1i\in D_2^+$  and $p = -0.082 - 1.1i\in D_2^+$.
\textit{Second row}: $p = -0.216-1.1i\in D_2^+$ and $p = -0.69 - 1.1i\in D_2^+$.
\textit{Third row}: $p = -1.146-1.1i\in D_2^+$ and $p = -1.24 - 1.1i\in D_1$.
\textit{Fourth row}: $p = -0.082-0.95i\in D_2^-$ and $p = -0.214 - 0.5i\in D_3$.
Note that in the left panel of the first and third rows $\Vs$ and $\Vd$ essentially coincide. In particular, the latter case, which concerns $D_2^+$, is in complete agreement with the numerical simulation shown in the center panel of Figure \ref{numerics2}. 
}
\label{h_root}
\end{center}
\end{figure}

The integral equation \eqref{ds-fnlsd-int-eq-xistar-0} can be expressed  in terms of the incomplete elliptic integrals of the first and second kind.
Importantly,  we note that when the poles coincide with the branch point $\pm iq_o$ equation \eqref{ds-fnlsd-int-eq-xistar-0} reduces to the  modulation equation \eqref{ds-mod-eq-bm-i}.
A numerical evaluation of the solutions of equation \eqref{ds-fnlsd-int-eq-xistar-0} for various choices of $p$ that cover all four possible regions $D_1$, $D_2^+$, $D_2^-$ and $D_3$ of the third quadrant is shown in Figure \ref{h_root}.

The only jumps of Riemann-Hilbert problem \eqref{ds-n4t-rhp-mew} that are not bounded as $t\to\infty$ are the ones along $\p D_p^\ve$ and $\p D_{\bar p}^\ve$, which grow exponentially since they are controlled by the sign of $\text{Re}(ih)(\xi, p)$.\footnote{\label{ds-p-foot}The fact that   $e^{2i\left[h(\xi, k) + \theta(\xi, p)-\theta(\xi, k)\right]t}$ is  controlled by $\text{Re}(ih)(\xi, p)$ can be seen by recalling that the jump $V_p^{(5)}$ along the circle $\p D_p^\ve$ originates from the residue condition at $p$. Eventually, whenever the dominant problem contains the contribution from $\p D_p^\ve$, the jump $V_p^{(5)}$ will be converted back to a (modified) residue condition at $p$. Thus, eventually we will end up evaluating the quantity $h(\xi, k) + \theta(\xi, p)-\theta(\xi, k)$ at $k=p$.}
Thus, similarly to Subsection \ref{ds-esc-pw2-ss}, we must employ the analogue of transformation \eqref{ds-esc-n4t-def-pw2} in order to convert this growth into decay. Before doing so, however, we apply the analogue of transformation \eqref{ds-esc-n4-def-pw1} in order to remove the $k$-dependence from the jumps along $B$ and $\widetilde B$.
\vskip 3mm
\noindent
\textbf{Sixth deformation.}
The jumps $V_B^{(5)}$, $V_7^{(5)}$ and $V_8^{(5)}$ can be made independent of $k$ by means of the transformation
\eee{\label{ds-esc-n5-def-mew} 
N^{(6)}(x, t, k) = N^{(5)}(x, t, k) e^{ig(\xi, k)\sigma_3},
}
where the function $g(\xi, k)$ is analytic in $\mathbb C\setminus (B\cup \widetilde B)$ and satisfies the following jump conditions:
\sss{\label{ds-pav1}
\ddd{
g^++g^-
&=
-i\ln\left(\delta^2\right), &&
k\in B,
\\
g^++g^-
&=
-i\ln\left(\tfrac{\delta^2}{r}\right) + \omega, \quad && k\in L_7,
\\
g^++g^-
&=
-i\ln\left( \delta^2 \bar r\right) + \omega, &&
k\in L_8,
}
}
with  $\delta(\xi, k)$ defined by \eqref{ds-esc-del-def-mew} and with the real quantity $\omega(\xi)$  defined by
\eee{
\label{ds-omr}
\omega(\xi)
=
i \,
\frac{
\displaystyle \int_{B} \frac{\ln \delta^2(\xi, \nu)}{  \gamma(\xi, \nu)}\, d\nu
+
\int_{\widetilde B^+} \frac{
\ln\left[\frac{\delta^2(\xi, \nu) }{r(\nu)}\right]}{\gamma(\xi, \nu)}\, d\nu 
+
\int_{\widetilde B^-} \frac{\ln\left[ \delta^2(\xi, \nu)  \, \bar r(\nu)\right]}{  \gamma(\xi, \nu)}\, d\nu 
}
{\displaystyle \int_{\widetilde B}  \frac{d\nu}{  \gamma(\xi, \nu)}},
}
where we have introduced the notation
\eee{
\widetilde B^\pm := \widetilde B \cap \mathbb C^\pm.
}
The solution of problem \eqref{ds-pav1} is obtained via the Plemelj formulae as
\eee{\label{ds-esc-g-def-mew}
g(\xi, k)
=
\frac{\gamma(\xi, k)}{2\pi}
\bigg[
\int_{B} \frac{\ln \delta^2(\xi, \nu)}{ \gamma(\xi, \nu )(\nu -k)}\, d\nu 
+
\int_{L_7} \frac{
\ln\left[\frac{\delta^2(\xi, \nu) }{r(\nu)}\right]+i\omega(\xi)}{ \gamma(\xi, \nu )(\nu -k)}\, d\nu 
-
\int_{L_8} \frac{\ln\left[ \delta^2(\xi, \nu)  \, \bar r(\nu)\right]+i\omega(\xi)}{\gamma(\xi, \nu )(\nu -k)}\, d\nu 
\bigg].
}
The presence of $\omega$ in the jump conditions \eqref{ds-pav1} ensures that $g(\xi, k)=O(1)$ as $k \to \infty$. Indeed, expression \eqref{ds-esc-g-def-mew} implies
\eee{\label{ds-ginfr}
g(\xi, k) = g_\infty(\xi) +O\left(\frac 1k\right), \quad k \to \infty,
}
where the real quantity $g_\infty(\xi)$ is defined by
\eee{\label{ds-esc-ginf-def-mew}
g_\infty(\xi)
=
-\frac{1}{2\pi}
\bigg[
\int_{B} \frac{\ln \delta^2(\xi, \nu)}{\gamma(\xi, \nu)}\,\nu d\nu
+
\int_{\widetilde B^+} \frac{
\ln\left[\frac{\delta^2(\xi, \nu) }{r(\nu)}\right]+i\omega(\xi)}{ \gamma(\xi, \nu)}\,\nu  d\nu 
+
\int_{\widetilde B^-} \frac{\ln\left[ \delta^2(\xi, \nu)  \, \bar r(\nu)\right]+i\omega(\xi)}{\gamma(\xi, \nu)}\,\nu  d\nu
\bigg]. 
}

In summary, the function $N^{(6)}$ defined by \eqref{ds-esc-n5-def-mew} is analytic in $\mathbb C \setminus  \big(\bigcup_{j=1}^6 L_j \cup  B \cup \widetilde B \cup \p D_p^\ve \cup \p D_{\bar p}^\ve\big)$ and satisfies the Riemann-Hilbert problem
\sss{\label{ds-esc-n6-rhp-mew}
\ddd{
N^{(6)+} &= N^{(6)-} V_B,  && k\in B,
\\
N^{(6)+} &= N^{(6)-}V_{\widetilde B}^{(6)},  && k\in \widetilde B,
\\
N^{(6)+} &= N^{(6)-}V_j^{(6)}, && k\in L_j,\ j=1, \ldots, 6,
\\
N^{(6)+}  
&=
N^{(6)-}
V_p^{(6)}, && k\in  \p D_p^\ve,
\\
N^{(6)+} 
&=
N^{(6)-} 
V_{\bar p}^{(6)}, && k\in  \p D_{\bar p}^\ve,
\\
N^{(6)}  &=\left[I +O\left(\tfrac 1k\right)\right]e^{i\left[g_\infty(\xi)-G_\infty(\xi) t\right]\sigma_3}, \quad && k \to \infty,
\label{ds-n5-asymp-mew}
}
}
where the jump along $B$ is given by \eqref{ds-n1-jumps}, the jump along $\widetilde B$ is equal to
\eee{\label{ds-esc-vb5-vbt5-def-mew}
V_{\widetilde B}^{(6)}
=
\left(
\def\arraystretch{1}
\begin{array}{lr}
0
&
-e^{i(\Omega  t - \omega)} 
\\
e^{-i(\Omega  t - \omega)} 
&
0
\end{array}
\right),
}
the jumps along  $\p D_p^\ve$ and $\p D_{\bar p}^\ve$ are given by
\sss{\label{ds-vp5-vpb5-def-mew}
\ddd{
V_p^{(6)} 
&=
\left(
\def\arraystretch{1}
\begin{array}{lr}
1 & -\frac{c_p \, \delta^2(\xi, k)\, d(k) e^{-2ig(\xi, k)}}{k-p}\,  e^{2i\left[h(\xi, k) + \theta(\xi, p) - \theta(\xi, k)\right]t}
\\
0 & 1
\end{array}
\right),
\\
V_{\bar p}^{(6)}  
&=
\left(
\def\arraystretch{1}
\begin{array}{lr}
1 &0
\\
-\frac{c_{\bar p} \, \delta^{-2}(\xi, k) \, d(k)e^{2ig(\xi, k)}}{k-\bar p}\,  e^{-2i\left[h(\xi, k) + \theta(\xi, \bar p) - \theta(\xi, k)\right]t} & 1
\end{array}
\right),
}
}
the jumps  along the  contours  $L_j$ of Figure \ref{ds-esc-def5-mew-f} are equal to
\sss{\label{ds-esc-n5-jumps-mew}
\ddd{
V_1^{(6)}
&=
\left(
\def\arraystretch{1}
\begin{array}{lr}
1
&
\frac{\bar r\,\delta^{2} e^{-2ig}}{1+r\bar r}\,e^{2i ht}
\\
0
&
1
\end{array}
\right), 
\quad
&&V_2^{(6)}
=
\left(
\def\arraystretch{1}
\begin{array}{lr}
1
&
0
\\
\frac{r \delta^{-2} e^{2ig}}{1+r\bar r}\,e^{-2iht}
&
1
\end{array}
\right), 
\\
V_3^{(6)}
&=
\left(
\def\arraystretch{1}
\begin{array}{lr}
1
&
0
\\
r \delta^{-2} e^{2ig}  e^{-2iht}
&
1
\end{array}
\right), 
\quad
&&V_4^{(6)}
=
\left(
\def\arraystretch{1}
\begin{array}{lr}
1
&
\bar r \delta^{2} e^{-2ig} e^{2iht}
\\
0
&
1
\end{array}
\right), 
\\
V_5^{(6)}
&=
\left(
\def\arraystretch{1}
\begin{array}{lr}
1
&
\frac{\delta^{2}e^{-2ig}}{r}\, e^{2iht} 
\\
0
&
1
\end{array}
\right),
\quad 
&&V_6^{(6)}
=
\left(
\def\arraystretch{1}
\begin{array}{lr}
1
&
0
\\
\frac{\delta^{-2}e^{2ig}}{ \bar r}\, e^{-2iht} 
&
1
\end{array}
\right),
}
}
and the real quantities $G_\infty(\xi)$ and $g_\infty(\xi)$ are given by \eqref{ds-Ginf-mew} and \eqref{ds-esc-ginf-def-mew} respectively.
\vskip 3mm
\noindent
\textbf{Converting growth into decay.}
The growing exponentials in the jumps $V_p^{(6)}$ and $V_{\bar p}^{(6)}$ can be converted into decaying ones via the analogue of transformation  \eqref{ds-esc-n4t-def-pw2}, i.e. by letting
\sss{\label{ds-esc-n5t-def-mew}
\eee{ 
\widetilde N^{(6)} 
=
\begin{cases}
N^{(6)}  n^{\sigma_3}, & k\in \mathbb C\setminus \left(\,\overline{D_p^\ve} \cup  \overline{D_{\bar p}^\ve}\,\right),
\\
N^{(6)} J_p \, n^{\sigma_3}, & k\in D_p^\ve,
\\
N^{(6)} J_{\bar p} \, n^{\sigma_3}, & k\in D_{\bar p}^\ve,
\end{cases}
} 
where
\ddd{
J_p(\xi, k) 
&=
\left(
\def\arraystretch{2.4}
\arraycolsep=0pt
\begin{array}{lr}
\dfrac{1- \frac{n^2(p)}{n^2(k)}}{k-p} &  c_p \, d(k)\delta^2(\xi, k) e^{2i\left[\left(h(\xi, k) + \theta(\xi, p)-\theta(\xi, k)\right)t-g\right]}
\\
- \dfrac{n^2(p)\, e^{-2i\left[\left(h(\xi, k) + \theta(\xi, p)-\theta(\xi, k)\right)t-g\right]}}{c_p \, d(k)\delta^2(\xi, k) \left(k-\bar p\right)^2} &  k-p
\end{array}
\right),  
\\
J_{\bar p}(\xi, k)
&=
\left(
\def\arraystretch{2.8}
\arraycolsep=0pt
\begin{array}{lr}
k-\bar p &  -\dfrac{e^{2i\left[\left(h(\xi, k) + \theta(\xi, \bar p)-\theta(\xi, k)\right)t-g\right]}}{n^2(\bar p) c_{\bar p}\, d(k) \delta^{-2}(\xi, k) \left(k-p\right)^2} 
\\
c_{\bar p}\, d(k)\delta^{-2}(\xi, k) e^{-2i\left[\left(h(\xi, k) + \theta(\xi, \bar p)-\theta(\xi, k)\right)t-g\right]}  & \dfrac{1-\frac{n^2(k)}{n^2(\bar p)}}{k-\bar p}
\end{array}
\right).
}
}
Then, $\widetilde N^{(6)}$ satisfies the Riemann-Hilbert problem
\sss{\label{ds-esc-n5t-rhp-mew}
\ddd{
\widetilde N^{(6)+} &= \widetilde N^{(6)-} \widetilde V_B^{(6)},  && k\in B,
\\
\widetilde N^{(6)+} &= \widetilde N^{(6)-} \widetilde V_{\widetilde B}^{(6)},  && k\in \widetilde B,
\\
\widetilde N^{(6)+} &= \widetilde N^{(6)-} \widetilde V_j^{(6)}, && k\in L_j,\ j=1, \ldots, 6,
\\
\widetilde N^{(6)+}  
&=
\widetilde N^{(6)-}  \widetilde V_p^{(6)}, && k\in  \p D_p^\ve,
\\
\widetilde N^{(6)+}  
&=
\widetilde N^{(6)-} \widetilde V_{\bar p}^{(6)}, && k\in  \p D_{\bar p}^\ve,
\\
\widetilde N^{(6)} &=\left[I +O\left(\tfrac 1k\right)\right]e^{i\left[g_\infty(\xi)-G_\infty(\xi) t\right]\sigma_3}, \quad && k \to \infty,
\label{ds-esc-n5t-asymp-mew}
}
}
where the jumps along $B$ and $\widetilde B$ are given by
\eee{\label{ds-esc-vb5t-vbt5t-def-mew}
\widetilde V_B^{(6)}
=
\left(
\def\arraystretch{1}
\begin{array}{lr}
0
&
\tfrac{q_-}{iq_o}   n^{-2}
\\
\tfrac{\bar q_-}{iq_o}  n^2
&
0
\end{array}
\right),
\quad
\widetilde V_{\widetilde B}^{(6)}
=
\left(
\def\arraystretch{1}
\begin{array}{lr}
0
&
-e^{i(\Omega  t - \omega)} n^{-2}
\\
e^{-i(\Omega  t - \omega)} n^2
&
0
\end{array}
\right),
}
the jumps along $\p D_p^\ve$ and $\p D_{\bar p}^\ve$ are equal to
\sss{\label{ds-vp5t-vpb5t-def-mew}
\ddd{
\widetilde V_p^{(6)} 
&=
\left(
\def\arraystretch{1}
\begin{array}{lr}
1 & 0
\\
-\dfrac{n^2(p)e^{2ig(\xi, k)}}{c_p \, \delta^2(\xi, k)\, d(k)(k-p)}\,  e^{-2i\left[h(\xi, k) + \theta(\xi, p) - \theta(\xi, k)\right]t}
 & 1
\end{array}
\right),
\\
\widetilde V_{\bar p}^{(6)}  
&=
\left(
\def\arraystretch{1}
\begin{array}{lr}
1 & -\dfrac{n^{-2}(\bar p) \, \delta^{2}(\xi, k) e^{-2ig(\xi, k)}}{c_{\bar p} \, d(k) (k-\bar p)}\,  e^{2i\left[h(\xi, k) + \theta(\xi, \bar p) - \theta(\xi, k)\right]t} 
\\
0 & 1
\end{array}
\right),
}
}
and the jumps   along the  contours  $L_j$ of Figure \ref{ds-esc-def5-mew-f} are given by
\sss{
\ddd{\label{ds-esc-n5t-jumps-mew}
\widetilde V_1^{(6)}
&=
\left(
\def\arraystretch{1}
\begin{array}{lr}
1
&
\dfrac{\bar r\,\delta^{2} e^{-2ig}}{1+r\bar r}\,e^{2i ht}  n^{-2}
\\
0
&
1
\end{array}
\right), 
\quad
&&\widetilde V_2^{(6)}
=
\left(
\def\arraystretch{1}
\begin{array}{lr}
1
&
0
\\
\frac{r \delta^{-2} e^{2ig}}{1+r\bar r}\,e^{-2iht} n^2
&
1
\end{array}
\right), 
\\
\widetilde V_3^{(6)}
&=
\left(
\def\arraystretch{1}
\begin{array}{lr}
1
&
0
\\
r \delta^{-2} e^{2ig}  e^{-2iht} n^2
&
1
\end{array}
\right), 
\quad
&&\widetilde V_4^{(6)}
=
\left(
\def\arraystretch{1}
\begin{array}{lr}
1
&
\bar r \delta^{2} e^{-2ig} e^{2iht} n^{-2}
\\
0
&
1
\end{array}
\right), 
\\
\widetilde V_5^{(6)}
&=
\left(
\def\arraystretch{1}
\begin{array}{lr}
1
&
\frac{\delta^{2}e^{-2ig}}{r}\, e^{2iht} n^{-2}
\\
0
&
1
\end{array}
\right),
\quad 
&&
\widetilde V_6^{(6)}
=
\left(
\def\arraystretch{1}
\begin{array}{lr}
1
&
0
\\
\frac{\delta^{-2}e^{2ig}}{ \bar r}\, e^{-2iht} n^2
&
1
\end{array}
\right).
}
}

Transformation \eqref{ds-esc-n5t-def-mew} has re-introduced $k$ in the jumps along $B$ and $\widetilde B$, which had been made $k$-independent via the sixth deformation. Thus, motivated by the plane wave region, where having a constant jump along $B$ allowed us to solve the dominant Riemann-Hilbert problem explicitly, we  next perform a final, seventh deformation in order to remove the $k$-dependence from the jumps $\widetilde V_B^{(6)}$ and $\widetilde V_{\widetilde B}^{(6)}$. 

\begin{remark}[\b{Order of deformations}]
In view of the above discussion, it becomes apparent that  the sixth deformation should have been postponed until  after transformation \eqref{ds-esc-n5t-def-mew}, since then the $k$-dependence from the jumps along $B$ and $\widetilde B$ would have to be removed only once instead of twice. However, the less efficient order of deformations that we have followed   has the advantage of revealing precisely which part of the overall phase of the modulated elliptic wave \eqref{ds-esc-qasym-mew-t} is generated by the soliton at $\xi=\Vs$ (namely, the constant $4 \textnormal{arg}\left[p+\lambda(p)\right]$ via the seventh deformation as shown in \eqref{ds-esc-ginft-def-mew}).
\end{remark}

\noindent
\textbf{Seventh deformation.}
Similarly to \eqref{ds-esc-n5-def-mew}, we eliminate the dependence on $k$ from the jumps $\widetilde V_B^{(6)}$ and $\widetilde V_{\widetilde B}^{(6)}$ by letting
\eee{\label{ds-esc-n6t-def-mew} 
\widetilde N^{(7)}(x, t, k) = \widetilde N^{(6)}(x, t, k) e^{i\widetilde g(\xi, k)\sigma_3},
}
where  the function $\widetilde g(\xi, k)$ is analytic in $\mathbb C\setminus (B\cup \widetilde B)$ and satisfies the jump conditions
\sss{\label{ds-esc-gt-jumps-mew}
\ddd{
\widetilde g^{\,+} + \widetilde g^{\,-}
&=
i\ln\left(n^2\right), &&
k\in B,
\\
\widetilde g^{\,+} + \widetilde g^{\,-}
&=
i\ln\left(n^2r\right) + \widetilde \omega, \quad && k\in L_7,
\\
\widetilde g^{\,+} + \widetilde g^{\,-}
&=
i\ln\left(\tfrac{n^2}{\bar r}\right) + \widetilde \omega, 
&&
k\in L_8,
}
}
with the real quantity $\widetilde \omega(\xi)$  given by
\eee{\label{ds-esc-omegat-def-mew}
\widetilde  \omega(\xi)
=
-i\, \frac{
\displaystyle
\int_{B} \frac{\ln n^2(\nu)}{  \gamma(\xi, \nu)}\, d\nu
+
\int_{L_7} \frac{
\ln\left[n^2(\nu) r(\nu)\right]}{  \gamma(\xi, \nu)}\, d\nu 
+
\int_{L_8} \frac{\ln\left[\frac{\bar r(\nu)}{n^2 (\nu)}\right]}{  \gamma(\xi, \nu)}\, d\nu 
}
{\displaystyle \int_{\widetilde B}  \frac{d\nu}{  \gamma(\xi, \nu)}}.
}
Similarly to \eqref{ds-esc-g-def-mew}, we have the explicit formula
\eee{\label{ds-esc-gt-def-mew}
\widetilde g(\xi, k)
=
-\frac{\gamma(\xi, k)}{2\pi}
\bigg\{
\int_{B} \frac{\ln n^2(\nu)}{ \gamma(\xi, \nu)(\nu -k)}\, d\nu 
+
\int_{L_7} \frac{
\ln\left[n^2(\nu)r(\nu)\right]-i\widetilde \omega(\xi)}{ \gamma(\xi, \nu)(\nu -k)}\, d\nu 
+
\int_{L_8} \frac{\ln\left[ \frac{ \bar r(\nu)}{n^2(\nu)}\right]+i\widetilde \omega(\xi)}{\gamma(\xi, \nu)(\nu -k)}\, d\nu 
\bigg\},
}
which implies
\eee{\label{ds-ginfrt}
\widetilde g(\xi, k) = \widetilde g_\infty(\xi) +O\left(\frac 1k\right), \quad k \to \infty,
}
with the real quantity $\widetilde g_\infty(\xi)$ given by
\eee{
\widetilde g_\infty(\xi)
=
\frac{1}{2\pi}
\bigg\{
\int_{B} \frac{\ln n^2(\nu)}{  \gamma(\xi, \nu)}\, \nu d\nu
+
\int_{L_7} \frac{
\ln\left[n^2(\nu)r(\nu)\right]
-i\widetilde \omega(\xi)}{\gamma(\xi, \nu)}\, \nu d\nu 
+
\int_{L_8} \frac{\ln\left[\frac{\bar r(\nu)}{n^2(\nu)}\right] +i\widetilde \omega(\xi)}{\gamma(\xi, \nu)}\, \nu d\nu 
\bigg\}.
}
It turns out that $\widetilde g_\infty$ is actually independent of $\xi$ and,  more precisely, 
\eee{\label{ds-esc-ginft-def-mew}
\widetilde g_\infty = 2\text{arg}\left[p+\lambda(p)\right]
}
like in the plane wave region. Eventually (see Remark \ref{ds-om-shift-r}), this implies that the effect of the soliton at $\xi=\Vo$ on the phase of the leading order asymptotics is the same both for  $\xi\in (\Vo, 0)$ and for $\xi\in(\Vs, \Vo)$.

The Riemann-Hilbert problem for $\widetilde N^{(6)}$ yields the following problem for $\widetilde N^{(7)}$:
\sss{\label{ds-esc-n6t-rhp-mew}
\ddd{
\widetilde N^{(7)+} &= \widetilde N^{(7)-} V_B,  && k\in B,
\\
\widetilde N^{(7)+} &= \widetilde N^{(7)-} \widetilde V_{\widetilde B}^{(7)},  && k\in \widetilde B,
\\
\widetilde N^{(7)+} &= \widetilde N^{(7)-} \widetilde V_j^{(7)}, && k\in L_j,\ j=1, \ldots, 6,
\\
\widetilde N^{(7)+}  
&=
\widetilde N^{(7)-}  \widetilde V_p^{(7)}, && k\in  \p D_p^\ve,
\\
\widetilde N^{(7)+}  
&=
\widetilde N^{(7)-}  \widetilde V_{\bar p}^{(7)}, && k\in  \p D_{\bar p}^\ve,
\\
\widetilde N^{(7)} &=\left[I +O\left(\tfrac 1k\right)\right]e^{i\left[g_\infty(\xi)+\widetilde g_\infty-G_\infty(\xi) t\right]\sigma_3}, \quad && k \to \infty,
\label{ds-esc-n6t-asymp-mew}
}
}
where the jump along $B$ is given by \eqref{ds-n1-jumps},  the jump along $\widetilde B$ is equal to 
\eee{\label{ds-esc-vb6t-vbt6t-def-mew}
\widetilde V_{\widetilde B}^{(7)}
=
\left(
\def\arraystretch{1}
\begin{array}{lr}
0
&
-e^{i(\Omega  t - \omega - \widetilde \omega)}  
\\
e^{-i(\Omega  t - \omega - \widetilde \omega)}  
&
0
\end{array}
\right)
}
with the real quantities $\Omega$, $\omega$ and $\widetilde \omega$ given by  \eqref{ds-Om-def}, \eqref{ds-omr} and \eqref{ds-esc-omegat-def-mew} respectively, the jumps along  $\p D_p^\ve$ and $\p D_{\bar p}^\ve$ are given by
\sss{\label{ds-vp6t-vpb6t-def-mew}
\ddd{
\widetilde V_p^{(7)} 
&=
\left(
\def\arraystretch{1}
\begin{array}{lr}
1 & 0
\\
-\dfrac{n^2(p)e^{2i\left[g(\xi, k)+\widetilde g(\xi, k)\right]}}{c_p \, \delta^2(\xi, k)\, d(k)(k-p)}\,  e^{-2i\left[h(\xi, k) + \theta(\xi, p) - \theta(\xi, k)\right]t}
 & 1
\end{array}
\right),
\\
\widetilde V_{\bar p}^{(7)}  
&=
\left(
\def\arraystretch{1}
\begin{array}{lr}
1 & -\dfrac{n^{-2}(\bar p) \, \delta^{2}(\xi, k) e^{-2i\left[g(\xi, k)+\widetilde g(\xi, k)\right]}}{c_{\bar p} \, d(k) (k-\bar p)}\,  e^{2i\left[h(\xi, k) + \theta(\xi, \bar p) - \theta(\xi, k)\right]t} 
\\
0 & 1
\end{array}
\right)
}
}
with the functions $d$, $n$, $\delta$, $g$ and $\widetilde g$ defined by \eqref{ds-d-def}, \eqref{ds-esc-n4t-def-pw2},  \eqref{ds-esc-del-def-mew}, \eqref{ds-esc-g-def-mew} and \eqref{ds-esc-gt-def-mew} respectively,
the jumps   along the  contours  $L_j$ of Figure \ref{ds-esc-def5-mew-f} are given by
\sss{\label{ds-esc-n6t-jumps-mew}
\ddd{
\widetilde V_1^{(7)}
&=
\left(
\def\arraystretch{1}
\begin{array}{lr}
1
&
\frac{\bar r\,\delta^{2} e^{-2i(g+\widetilde g)}}{1+r\bar r}\,e^{2i ht} n^{-2}
\\
0
&
1
\end{array}
\right), 
\quad
&&\widetilde V_2^{(7)}
=
\left(
\def\arraystretch{1}
\begin{array}{lr}
1
&
0
\\
\frac{r \delta^{-2} e^{2i(g+\widetilde g)}}{1+r\bar r}\,e^{-2iht} n^2
&
1
\end{array}
\right), 
\\
\widetilde V_3^{(7)}
&=
\left(
\def\arraystretch{1}
\begin{array}{lr}
1
&
0
\\
r \delta^{-2} e^{2i(g+\widetilde g)}  e^{-2iht} n^2
&
1
\end{array}
\right), 
\quad
&&\widetilde V_4^{(7)}
=
\left(
\def\arraystretch{1}
\begin{array}{lr}
1
&
\bar r \delta^{2} e^{-2i(g+\widetilde g)} e^{2iht} n^{-2}
\\
0
&
1
\end{array}
\right), 
\\
\widetilde V_5^{(7)}
&=
\left(
\def\arraystretch{1}
\begin{array}{lr}
1
&
\frac{\delta^{2}e^{-2i(g+\widetilde g)}}{r}\, e^{2iht} n^{-2}
\\
0
&
1
\end{array}
\right),
\quad 
&&
\widetilde V_6^{(7)}
=
\left(
\def\arraystretch{1}
\begin{array}{lr}
1
&
0
\\
\frac{\delta^{-2}e^{2i(g+\widetilde g)}}{ \bar r}\, e^{-2iht} n^2
&
1
\end{array}
\right),
}
}
the real quantities  $G_\infty$ and $g_\infty$  are defined by  \eqref{ds-Ginf-mew} and  \eqref{ds-esc-ginf-def-mew}, and the real constant $\widetilde g_\infty$ is given by \eqref{ds-esc-ginft-def-mew}.

\vskip 3mm
\noindent
\textbf{Decomposition into dominant and error problems.}
The jumps $\widetilde V_j^{(7)}$, $j=1, \ldots, 6$, do not contribute to the leading-order long-time asymptotics since they decay to the identity as $t\to\infty$ due to the sign structure of $\Re (ih)$ (see Figure \ref{ds-esc-def5-mew-f}). 
The same is true for the jumps $\widetilde V_p^{(7)}$ and $\widetilde V_{\bar p}^{(7)}$ since the exponentials involved in these jumps are controlled by the sign of $\text{Re}(ih)(\xi, p)$.$^{\ref{ds-p-foot}}$
Therefore,  the dominant component of Riemann-Hilbert problem \eqref{ds-esc-n6t-rhp-mew} must come from the jumps $V_B$ and $\widetilde V_{\widetilde B}^{(7)}$.
With these in mind, we decompose problem \eqref{ds-esc-n6t-rhp-mew} as follows. 

Let $D_{k_o}^{\epsilon}$, $D_{\alpha}^{\epsilon}$ and $D_{\bar \alpha}^{\epsilon}$ be disks of radius $\epsilon$ centered at $k_o$, $\alpha$ and $\bar \alpha$ respectively, where $\epsilon$ is sufficiently small so that these disks do not intersect with each other or with $B\cup \overline{D_p^\ve} \cup \overline{D_p^\ve}$. Then, write
\eee{\label{ds-n5t-decomp-mew}
\widetilde N^{(7)}
=
\widetilde N^\err \widetilde N^\asymp
}
where
\eee{
\widetilde N^\asymp
=
\begin{cases}
\widetilde N^\dom,  & k\in \mathbb C\setminus (D_{k_o}^{\epsilon}\cup D_{\alpha}^{\epsilon}\cup D_{\bar \alpha}^{\epsilon}),
\\
\widetilde N^D, & k\in D_{k_o}^{\epsilon}\cup D_{\alpha}^{\epsilon}\cup D_{\bar \alpha}^{\epsilon},
\end{cases}
}
and the functions $\widetilde N^\dom$, $\widetilde N^D$ and $\widetilde N^\err$ are defined as follows:
\vskip 3mm
\begin{enumerate}[label=$\bullet$, leftmargin=4mm, rightmargin=0mm]
\advance\itemsep 3mm
\item
$\widetilde N^\dom(x, t, k)$ is analytic in $\mathbb C\setminus(B \cup \widetilde B)$  and satisfies the Riemann-Hilbert problem 
\sss{\label{ds-esc-nmodt-rhp-mew}
\ddd{
\widetilde N^{\dom+} &= \widetilde N^{\dom-} V_B,  && k\in B,
\label{ds-rhpBbr}
\\
\widetilde N^{\dom+} &= \widetilde N^{\dom-}  \widetilde V_{\widetilde B}^{(7)}, && k\in \widetilde B,
\label{ds-rhpBcr}
\\
\widetilde N^\dom &=  
 \left[I+O\left(\tfrac 1k\right)\right]e^{i\left[g_\infty(\xi) +\widetilde g_\infty -G_\infty(\xi) t  \right]\sigma_3},\quad && k \to \infty.
 \label{ds-nmodt-asymp-mew}
}
}
\item
$\widetilde N^D(x, t, k)$ is analytic in $D_{k_o}^{\epsilon}\cup D_{\alpha}^{\epsilon}\cup D_{\bar \alpha}^{\epsilon}\setminus  \bigcup_{j=1}^{8}  L_j$ with jumps 
\eee{\label{ds-rhpPr} 
\widetilde N^{D+} = \widetilde N^{D-} \widetilde V_j^{(6)},
\quad
k\in \widehat L_j := L_j\cap \left(D_{k_o}^{\epsilon}\cup D_{\alpha}^{\epsilon}\cup D_{\bar \alpha}^{\epsilon}\right),
\ 
j=1, \ldots, 8,
}
as shown in  Figure \ref{ds-esc-ndt-jumps-mew-f}.
\begin{figure}[t]
\begin{center}
\includegraphics[width=0.238\textwidth]{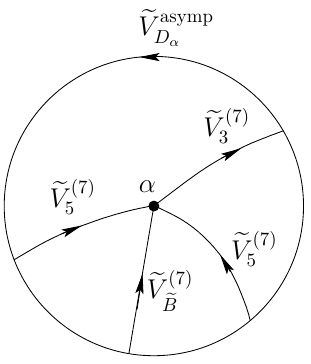}
\hskip 4mm
\includegraphics[width=0.241\textwidth]{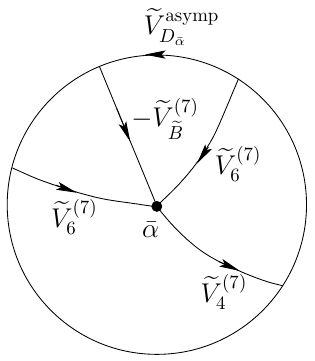}
\hskip 5mm
\includegraphics[width=0.233\textwidth]{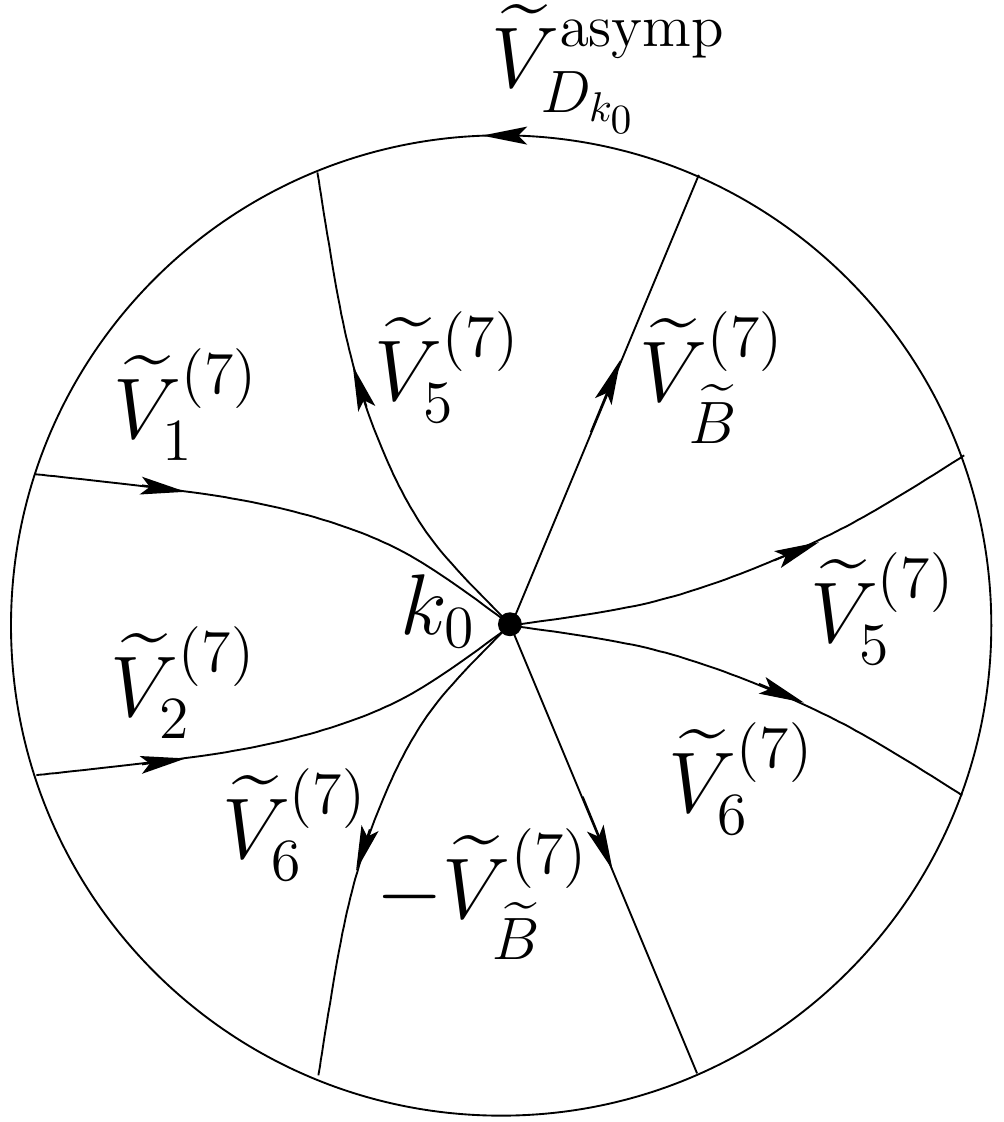}
\caption{Modulated elliptic wave   in the  transmission regime: The jumps of  $\widetilde N^D$  in the interior of and along the boundary of the disks $\overline{D_\alpha^\epsilon}$, $\overline{D_{\bar \alpha}^\epsilon}$ and $\overline{D_{k_o}^\epsilon}$. Note that, although the jumps $\widetilde V_{D_\alpha}^\asymp$, $\widetilde V_{D_{\bar \alpha}}^\asymp$, $\widetilde V_{D_{k_o}}^\asymp$ are unknown, they are equal to the identity up to $O(t^{-1/2})$ and hence do not affect the dominant problem.}
\label{ds-esc-ndt-jumps-mew-f}
\end{center}
\end{figure}

\item
$\widetilde N^\err(x, t, k)$ is analytic in $\mathbb C \setminus (\bigcup_{j=1}^{6} \wc L_j \cup \p D_{k_o}^{\epsilon}\cup \p D_{\alpha}^{\epsilon}\cup \p D_{\bar \alpha}^{\epsilon}\cup \p D_p^\ve \cup \p D_{\bar p}^\ve\big)$ with $\wc L_j := L_j\setminus ( \overline{D_{k_o}^{\epsilon}}\cup \overline{D_{\alpha}^{\epsilon}}\cup \overline{D_{\bar \alpha}^{\epsilon}})$ and satisfies the Riemann-Hilbert problem (see Figure \ref{ds-esc-nerrt-jumps-mew-f})
\sss{\label{ds-rhpEr}
\ddd{
\widetilde N^{\err+} &= \widetilde N^{\err-} \widetilde V^\err, \quad
&&  k\in {\textstyle \bigcup}_{j=1}^{6} \wc L_j \cup \p D_{k_o}^{\epsilon}\cup \p D_{\alpha}^{\epsilon}\cup \p D_{\bar \alpha}^{\epsilon}\cup \p D_p^\ve \cup \p D_{\bar p}^\ve,
\label{ds-rhpEbr}
\\
\widetilde N^\err  &= I +O\left(\tfrac 1k\right), && k \to \infty, \label{ds-net-asymp-mew}
}
}
where
\eee{\label{ds-VEdefr}
\widetilde V^\err
=
\begin{cases}
\widetilde N^\dom \widetilde V_j^{(7)} (\widetilde N^\dom)^{-1}, &k\in \wc L_j,
\\
\widetilde N^\dom \widetilde V_p^{(7)} (\widetilde N^\dom)^{-1}, &k\in \p D_p^\ve,
\\
\widetilde N^\dom \widetilde V_{\bar p}^{(7)} (\widetilde N^\dom)^{-1}, &k\in \p D_{\bar p}^\ve,
\\
\widetilde N^{\asymp-}(\widetilde V_D^\asymp)^{-1}(\widetilde N^{\asymp-})^{-1}, &k\in \p D_{k_o}^{\epsilon}\cup \p D_{\alpha}^{\epsilon}\cup \p D_{\bar \alpha}^{\epsilon},
\end{cases}
}
and 
\eee{\label{ds-VDasymp}
\widetilde V_D^\asymp
=
\begin{cases}
\widetilde V_{D_\alpha}^\asymp, &k\in \p D_{\alpha}^\epsilon,
\\
\widetilde V_{D_{\bar \alpha}}^\asymp, &k\in \p D_{\bar \alpha}^\epsilon,
\\
\widetilde V_{D_{k_o}}^\asymp, &k\in \p D_{k_o}^\epsilon.
\end{cases}
}
Importantly, despite the fact that the jump $\widetilde V_D^\asymp$  is unknown, in \cite{bm2017} it was estimated to be equal to the identity up to $O(t^{-1/2})$  and hence it does not affect the leading-order asymptotics.
\begin{figure}[t]
\begin{center}
\includegraphics[scale=1]{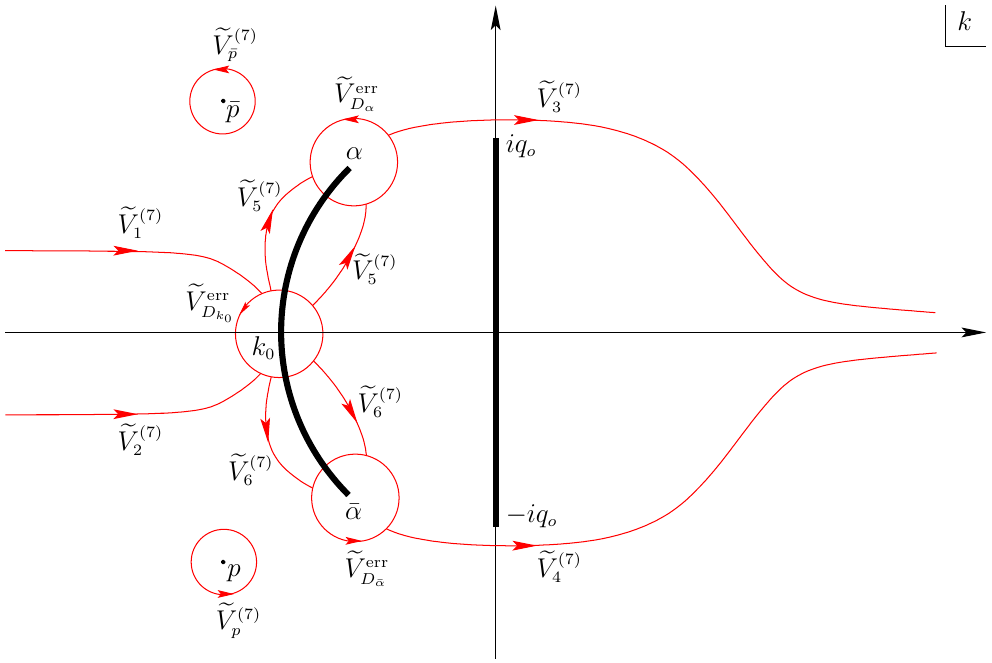}
\caption{Modulated elliptic wave in the  transmission regime: The jumps of $\widetilde N^\err$.}
\label{ds-esc-nerrt-jumps-mew-f}
\end{center}
\end{figure}
\end{enumerate}

\begin{figure}[t]
\begin{center}
\includegraphics[scale=.5]{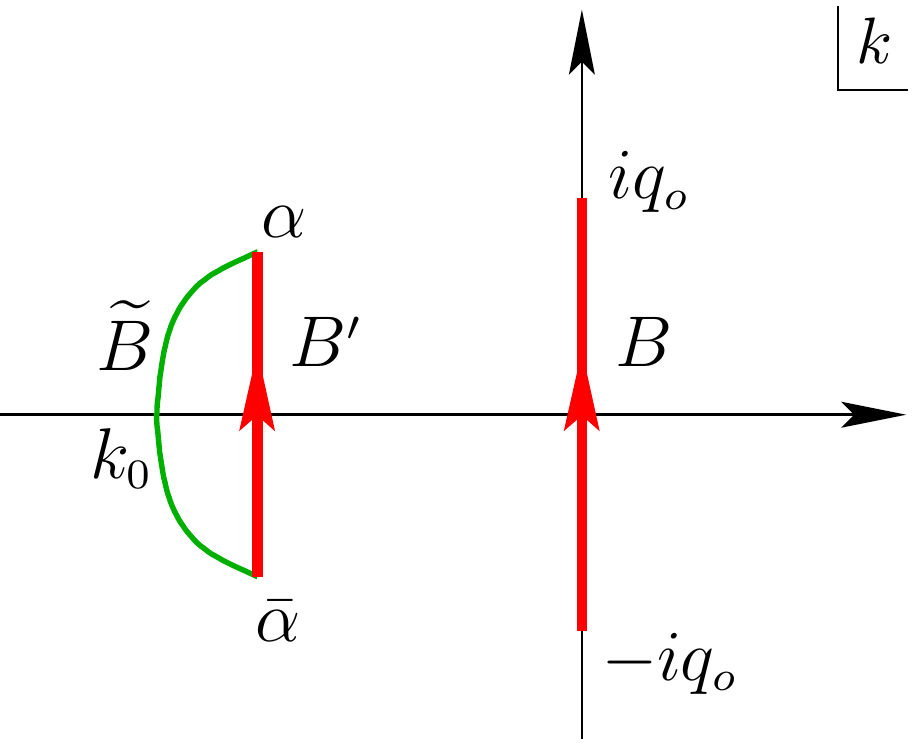}
\caption{Modulated elliptic wave in the  transmission regime: 
The jump contours of $\widetilde N^\dom$. The fact that the jump $\widetilde V_{\widetilde B}^{(7)}$ is independent of $k$ has allowed us to deform $\widetilde B$ to the straight line segment  $B'$  from $\bar\alpha$ to $\alpha$. 
Note that the original and deformed versions of $\widetilde N^\dom$ agree outside the finite region $\mathcal D$ enclosed by $\widetilde B$ and $B'$.}
\label{ds-esc-nmodt-jumps-mew-f}
\end{center}
\end{figure}

\begin{figure}[t]
\begin{center}
\includegraphics[scale=0.425]{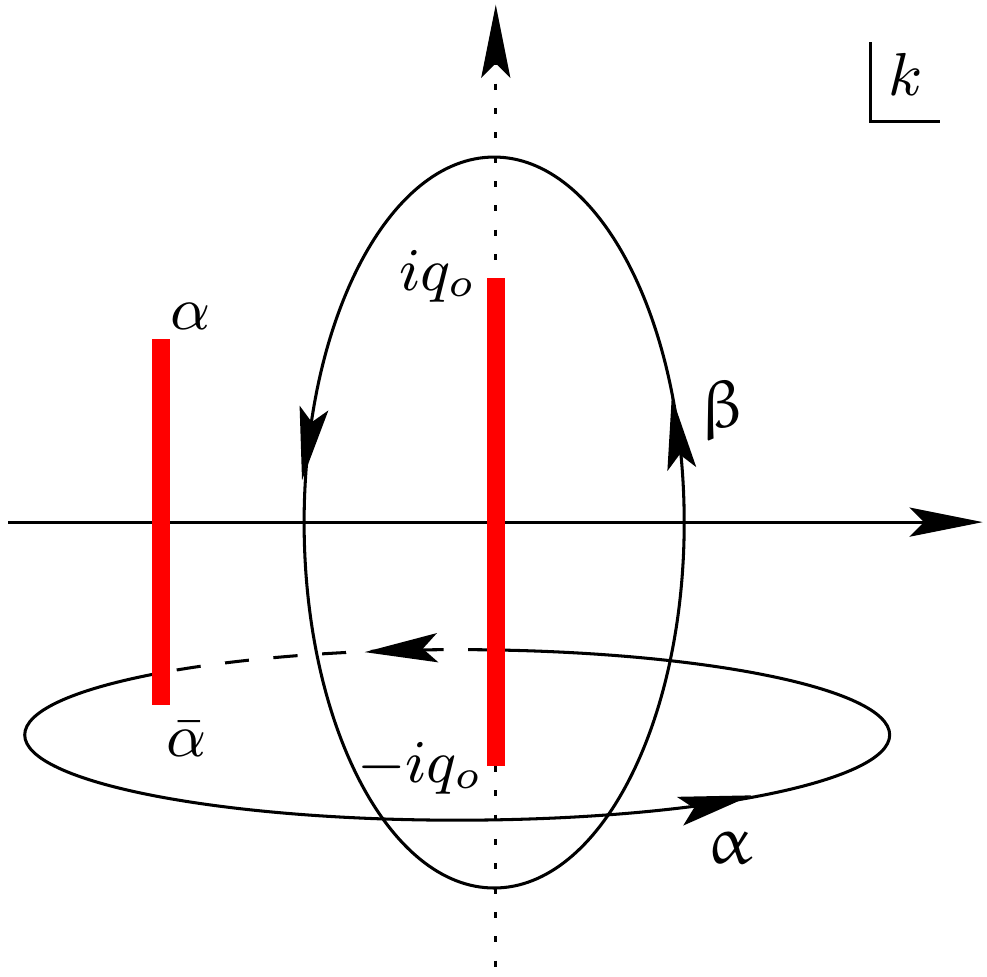}
\caption{Modulated elliptic wave  in the  transmission regime: The basis $\{\upalpha, \upbeta\}$ of cycles for the genus-1 Riemann surface associated with the function $\Gamma(\xi, k)$.
The cycle  $\upbeta$  is a closed, anti-clockwise contour that encircles the branch cut $B$ while lying on the first sheet  of the Riemann surface. The cycle $\upalpha$ consists of an anti-clockwise contour that begins from the left of the branch cut $B'$ and approaches the branch cut $B$ from the right while lying on the first sheet, and  then returns to $B'$ via the second sheet (dashed portion).
}
\label{ds-abcycles}
\end{center}
\end{figure}

\vskip 3mm
\noindent
\textbf{Solution of the dominant problem.}
We begin by noting that, at the level of the dominant problem \eqref{ds-esc-nmodt-rhp-mew}, since the jump $\widetilde V_{\widetilde B}^{(7)}$ is independent of $k$, the jump contour $\widetilde B$ can be deformed  to the straight line segment $B'$ from $\bar \alpha$ to $\alpha$ so that the jump contours of problem \eqref{ds-esc-nmodt-rhp-mew} are as shown in Figure \ref{ds-esc-nmodt-jumps-mew-f}.
Problem \eqref{ds-esc-nmodt-rhp-mew} was solved in \cite{bm2017} in the case of $\widetilde\omega=\widetilde g_\infty=0$.  Adapting that analysis  to account for the presence of $\widetilde\omega$ and $\widetilde g_\infty$, we obtain
\eee{
\label{ds-esc-ntmod-sol-mew}
\widetilde N^\dom(x, t, k)
=
e^{i\left[g_\infty(\xi)+\widetilde g_\infty -G_\infty(\xi) t  \right]\sigma_3}
 \widetilde{\mathcal N}^{-1}(\xi, \infty, c) \, \widetilde{\mathcal N}(\xi, k,  c)
}
with 
\eee{\label{ds-ncaldef}
\hspace*{-2mm}
\widetilde{\mathcal N}(\xi, k, c)
=
\frac 12
\def\arraystretch{1.5}
\left(
\begin{array}{lr}
\left[\eta(\xi, k)+\eta^{-1}(\xi, k)\right]\@{\widetilde N}_1(\xi, k, c)
& 
i\left[\eta(\xi, k)-\eta^{-1}(\xi, k)\right]\@{\widetilde N}_2(\xi, k, c)
\\
-i\left[\eta(\xi, k)-\eta^{-1}(\xi, k)\right]\@{\widetilde N}_1(\xi, k, -c)
&
\left[\eta(\xi, k)+\eta^{-1}(\xi, k)\right]\@{\widetilde N}_2(\xi, k, -c)
\end{array}
\right)
}
and
\eee{
\label{ds-ncalasym}
\widetilde{\mathcal N}(\xi, \infty, c)
:=
 \lim_{k \to\infty}
\widetilde{\mathcal N}(\xi, k, c),
}
where the function $\eta$ with branch cuts along $B$ and $B'$ (see Figure \ref{ds-esc-nmodt-jumps-mew-f}) is defined by
\eee{\label{ds-pdef} 
\eta(\xi, k)
:=
\left[
\frac
{\left(k-iq_o \right)\left(k-\alpha\right)}
{\left(k+iq_o \right)\left(k-\bar \alpha \right)}
\right]^{\frac 14}
}
and where $\@{\widetilde N}_1$ and $\@{\widetilde N}_2$ denote the first and second component of the vector-valued function 
\eee{
\label{ds-mcaldef-til}
\@{\widetilde N} (\xi, k, c) 
:=
\left(
\frac{
 \Theta\big(-\frac{\Omega t}{2\pi}+\frac{\omega+\widetilde \omega}{2\pi}+\frac{i\ln\left(\frac{\bar q_-}{iq_o}\right)}{2\pi}+\upnu(k)+c \big)}
 {  \sqrt{\frac{iq_o }{\bar q_-}} \ \Theta\left(\upnu(k)+c \right)},
\frac{\Theta\big(-\frac{\Omega t}{2\pi}+\frac{\omega+\widetilde \omega}{2\pi}+\frac{i\ln\left(\frac{\bar q_-}{iq_o}\right)}{2\pi}-\upnu(k)+c \big)}
{ \sqrt{\tfrac{\bar q_-}{iq_o }}\   \Theta\left(-\upnu(k)+c \right)}
\right).
}
In the above definition, the dependence of $\Omega$, $\omega$, $\widetilde \omega$, $\upnu$ and $c$ on $\xi$ has been suppressed for convenience. Moreover, $\Theta(k)=\Theta(\xi, k)$ is the following variant of the third Jacobi theta function:
\eee{\label{ds-Theta-def}
\Theta(\xi, k)
=
\theta_3(\pi k, e^{i\pi \tau(\xi)}),
\quad
\theta_3(z, \varrho)
:=
\sum_{\ell\in\mathbb Z} e^{2i \ell  z} \varrho^{\ell^2},
}
with Riemann period   
\eee{\label{ds-tau-a}
\tau(\xi) 
:= \left(\oint_{\upbeta} \frac{dk}{ \Gamma(\xi, k)}\right)^{-1} \oint_{\upalpha}   \frac{dk}{\Gamma(\xi, k)}
=
\frac{iK\big(\sqrt{1-m^2}\,\big)}{K\left(m\right)},
}
where the function $\Gamma$ is defined in terms of the function $\gamma$ of  \eqref{ds-gamma-def-i} by
\eee{\label{ds-gamma0-def}
\Gamma(\xi, k) = \left\{\begin{array}{ll} \gamma(\xi, k), & k\in \mathbb C \setminus \overline{\mathcal D}, \\   -\gamma(\xi, k), & k\in \mathcal D,\end{array}\right.
}
where $\mathcal D$ denotes the finite region enclosed by $\widetilde B$ and $B'$ (see Figure \ref{ds-esc-nmodt-jumps-mew-f}).
This definition implies that $\Gamma$ has branch cuts along $B$ and $B'$, i.e. the branch cut $\widetilde B$ of $\gamma$ has been deformed to $B'$ in the case of $\Gamma$.
The cycles $\{\upalpha, \upbeta\}$ of the genus-1 Riemann surface associated with   $\Gamma$   are depicted in Figure \ref{ds-abcycles}.
Furthermore, the Abelian map $\upnu$ in the arguments of the $\Theta$-functions in \eqref{ds-mcaldef-til} is defined by  
\eee{\label{ds-vdefr0}
\upnu(k)
=
\upnu(\xi, k) 
=
\left(\oint_{\upbeta} \frac{d\nu}{ \Gamma(\xi, \nu)}\right)^{-1}
 \int_{iq_o }^k  \frac{d\nu}{\Gamma(\xi, \nu)},
}
and, finally,
\eee{\label{ds-c-choice} 
c
=
c(\xi)
:=
\upnu\!\left(\frac{q_o \alpha_\re}{q_o +\alpha_\im}\right)+\frac 12\left(1+\tau \right).
}

\begin{remark}[\b{Analyticity of $\widetilde{\mathcal N}(\xi, k, c)$}]\label{ds-c-rem}
The definition \eqref{ds-c-choice} of $c$ ensures that the only possible singularity of $\widetilde{\mathcal N}(\xi, k, c)$ on the first sheet of the Riemann surface may occur at $k=\frac{q_o \alpha_\re}{q_o +\alpha_\im}$. This is because $\Theta(-\upnu(k)+c)$ vanishes whenever $-\upnu(k)+c=\frac 12(1+\tau)+\mathbb Z + \tau \mathbb Z$ and $\upnu(k)$ is injective on each sheet of the Riemann surface as an Abelian map.
Furthermore, this singularity is actually removable since it is the unique (finite) zero of $\eta-\eta^{-1}$ on the first sheet of the Riemann surface. Hence, all four entries of $\widetilde{\mathcal N}(\xi, k, c)$  are analytic away from the branch cuts $B$ and $B'$.
\end{remark}

\begin{remark}[\b{Invertibility of $\widetilde{\mathcal N}(\xi, \infty, c)$}]
Since $\lim_{k\to\infty}\eta(\xi, k)= 1$  and $\Theta(k) = \Theta(-k)$, 
letting 
\eee{\label{ds-vdefr}
\upnu_\infty=\upnu_\infty(\xi) := \lim_{k\to\infty} \upnu(\xi, k)
}
we have
\eee{
\det \widetilde{\mathcal N}(\xi, \infty, c)
=
\frac{\Theta\big(-\frac{\Omega t}{2\pi}+\frac{\omega+\widetilde \omega}{2\pi}+\frac{i\ln\left(\frac{\bar q_-}{iq_o}\right)}{2\pi}+\upnu_\infty+c \big)\Theta\big(-\frac{\Omega t}{2\pi}+\frac{\omega+\widetilde \omega}{2\pi}+\frac{i\ln\left(\frac{\bar q_-}{iq_o}\right)}{2\pi}-\upnu_\infty-c \big)}{\Theta^2\big(\upnu_\infty+c \big)}.
\nn
}
The denominator of this expression is always nonzero thanks to the choice of $c$ (see Remark \ref{ds-c-rem}). Moreover, noting that
$
\Delta:=\left\{-\Omega t+\omega+\widetilde \omega+i\ln\left[\bar q_-/(iq_o)\right]\right\}/2\pi\in\mathbb R
$
we observe that  subtracting or adding $\Delta$ to  $\upnu_\infty+c$ does not affect the imaginary part of the argument of the $\Theta$-functions in the numerator of $\det \widetilde{\mathcal N}(\xi, \infty, c)$. 
Thus, recalling that the zeros of $\Theta(k)$ are located at $k = \frac 12 \left(1+\tau\right) + \mathbb Z + \tau \mathbb Z$ and noting that $\tau$ is purely imaginary, we deduce that  $\det \widetilde{\mathcal N}(\xi, \infty, c)$ is nonzero and hence $\widetilde{\mathcal N}(\xi, \infty, c)$ is invertible, as required by \eqref{ds-esc-ntmod-sol-mew}.
\end{remark}

Starting from the reconstruction formula  \eqref{ds-q-recon-n}  and applying the successive deformations that lead from $N=N^{(0)}$ to $\widetilde N^{(7)}$ while keeping in mind that $n, d, \delta\to 1$ as $k\to \infty$, we obtain
\eee{\label{ds-q-recon-mew}
q(x, t)
=
-2i \lim_{k \to \infty} k  \widetilde N^{(7)}_{12}(x, t, k)  e^{i\left[g_\infty(\xi)+\widetilde g_\infty -G_\infty(\xi) t  \right]}.
}
Furthermore, according to the decomposition \eqref{ds-n5t-decomp-mew}, for large $k$  we have
\eee{
\widetilde N^{(7)}_{12} = \widetilde N_{11}^\err \widetilde N_{12}^\dom + \widetilde N_{12}^\err \widetilde N_{22}^\dom.
}
Hence, using also the asymptotic conditions \eqref{ds-nmodt-asymp-mew} and \eqref{ds-net-asymp-mew}, we find
\eee{\label{ds-q-recon-mew-2}
q(x, t)
=
-2i \lim_{k \to \infty} k \widetilde N_{12}^\dom(x, t, k) e^{i\left[g_\infty(\xi) + \widetilde g_\infty -G_\infty(\xi) t  \right]} 
 -2i \lim_{k \to \infty}  k \widetilde N_{12}^\err(x, t, k).
}
All of the jumps of $\widetilde N^\err$, including those along $\p D_p^\ve$ and $\p D_{\bar p}^\ve$,  tend to the identity exponentially fast as $t\to\infty$. Hence, the second term in \eqref{ds-q-recon-mew-2} is of lower order. In fact, similarly to \cite{bm2017} (see also \cite{bv2007}) we have
\eee{\label{ds-net-est-mew}
\lim_{k \to \infty}  k \widetilde N_{12}^\err(x, t, k) 
=
O\big(t^{-\frac 12}\big),
\quad t\to \infty.
}
Moreover, since the original and deformed versions of $\widetilde N^\dom$ agree outside the finite region $\mathcal D$ enclosed by $\widetilde B$ and $B'$ (see Figure \ref{ds-esc-nmodt-jumps-mew-f}) and hence in the limit $k\to \infty$,  formula \eqref{ds-esc-ntmod-sol-mew} together with the expansion
\eee{\label{ds-pasym}
\eta(\xi, k) = 1-\frac{i(q_{o} +\alpha_\im)}{2k}+O\left(\frac{1}{k^2}\right),\quad k \to \infty,
}
imply
\eee{\label{ds-qsol-mbr}
\lim_{k\to\infty}  k \widetilde N_{12}^\dom(x, t, k) 
=
\frac 12 \left(q_o +\alpha_\im\right)
\frac{\@{\widetilde N}_2(\xi, \infty, c)}{\@{\widetilde N}_1(\xi, \infty, c)}
\,
e^{i\left[g_\infty(\xi) + \widetilde g_\infty-G_\infty(\xi) t \right]}.
}
Therefore, at leading order the reconstruction formula \eqref{ds-q-recon-mew-2} yields the modulated elliptic wave
\begin{multline}\label{ds-esc-qsol-mew}
q(x, t)
=
\frac{q_o \left(q_o +\alpha_\im\right) }{\bar q_-}
\frac{\Theta\big(\!-\frac{\Omega   t}{2\pi}+\frac{\omega+\widetilde \omega}{2\pi}+\frac{i\ln\left(\frac{\bar q_-}{iq_o}\right)}{2\pi}-\upnu_\infty +c \big)
\,
\Theta\left(\upnu_\infty +c \right)}
{\Theta\big(\!-\frac{\Omega   t}{2\pi}+\frac{\omega+\widetilde \omega}{2\pi}+\frac{i\ln\left(\frac{\bar q_-}{iq_o}\right)}{2\pi}+\upnu_\infty +c \big) 
\,
\Theta\left(-\upnu_\infty +c \right)}
\,
e^{2i\left[g_\infty(\xi) + \widetilde g_\infty - G_\infty(\xi)  t  \right]}
\\
+
O\big(t^{-\frac 12}\big),
\quad t\to \infty,
\end{multline}
where the real quantities
$\alpha_\im$,  
$\Omega$,
$\omega$, 
$\widetilde \omega$,
$G_\infty$,
$g_\infty$,
depend only on $\xi$ and 
are given respectively by 
\eqref{ds-mod-eqs-i},
\eqref{ds-Om-def},
\eqref{ds-omr},
\eqref{ds-esc-omegat-def-mew},  
\eqref{ds-Ginf-mew},
\eqref{ds-esc-ginf-def-mew}, 
the real constant $\widetilde g_\infty$ is given by \eqref{ds-esc-ginft-def-mew}, and the  quantities $c(\xi)$ and $\upnu_\infty(\xi)$ are defined by \eqref{ds-c-choice} and \eqref{ds-vdefr} respectively.
In fact, as shown in \cite{bm2017},  $\Omega$ can be expressed as
\eee{\label{ds-Omega-exp}
\Omega(\xi)
=
\frac{\pi \left|\alpha+iq_o\right|}{K(m)}  \left(\xi-2\alpha_\re\right)
}
with $K(m)$ being the complete elliptic integral of the first kind with elliptic modulus $m$ obtained via the modulation equations \eqref{ds-mod-eqs-i}.
Then, performing some straightforward manipulations of the relevant theta functions, we can write   \eqref{ds-esc-qsol-mew}    in the more explicit form \eqref{ds-esc-qasym-mew-t}-\eqref{ds-qsol-mew-t} with
\eee{\label{ds-Xdef}
X_o = X_o(\xi) := \frac{1}{2\pi}\left[\omega(\xi) - i\ln\Big(\frac{q_-}{q_o }\Big)\right] + \frac14.
}

The proof of Theorem \ref{ds-per-t} for the leading-order asymptotics in the transmission regime $p\in D_1$ is complete.

\begin{remark}[\b{Phase and position shifts}]
\label{ds-om-shift-r}
Setting $\widetilde \omega=\widetilde g_\infty=0$ in \eqref{ds-esc-qasym-mew-t}-\eqref{ds-qsol-mew-t} gives  \eqref{ds-qsol-mew-cpam-t}, which is precisely the modulated elliptic wave of \cite{bm2017}. That is, the effect of the soliton arising at $\xi=\Vs$ on the leading-order asymptotics for $\xi \in (\Vo, 0)$ is the constant phase shift $2\widetilde g_\infty=4\text{arg}\left[p+\lambda(p)\right]$ as well as a position shift related to the presence of the quantity $\widetilde \omega$.
\end{remark}

%
%
%
%
\section{The  Trap Regime: Proof of Theorem \ref{ds-ptr-t}}
\label{ds-trap-s}

This regime arises for $p$  inside the region $D_2^+$ of Figure \ref{ds-regions-f}. In that case, we have  $\Vs>\Vo$ i.e. $\text{Re}(i\theta)(\xi, p)<0$ throughout the interval $(-\infty, \Vo)$ and hence no soliton arises there. Furthermore, as already noted in the context of the fifth deformation of Subsection \ref{ds-esc-mew-ss}, for $p\in D_2^+$ the equation $\text{Re}(ih)(\xi, p)=0$ has a \textit{unique} solution $\Vd$ in the interval $(\Vo,0)$ (in fact, it turns out that $\Vd>\Vs$).
Thus, we split the range $(-\infty, 0)$ into the   subintervals $\xi< \Vo$; $\Vo<\xi< \Vd$; $\xi=\Vd$; and $\Vd<\xi< 0$.

\subsection{The range $\xi<\Vo$: plane wave}
\label{ds-trap-ss}
No soliton arises in this range since the asymptotics  is dictated by the phase function $\theta$ and the fact that $\Vo<\Vs$ means that $\text{Re}(i\theta)(\xi, p)<0$ throughout $(-\infty, \Vo)$.
Therefore, the analysis required is the same with the one carried out for $\xi<\Vs$ in the transmission regime (see Subsection \ref{ds-esc-pw1-ss}) and the leading-order asymptotics is given by  the plane wave \eqref{ds-qsol-pw-t}.

\subsection{The range $\Vo<\xi<\Vd$: modulated elliptic wave}
\label{ds-trap-mew1-ss}
\begin{figure}[t]
\begin{center}
\includegraphics[scale=1]{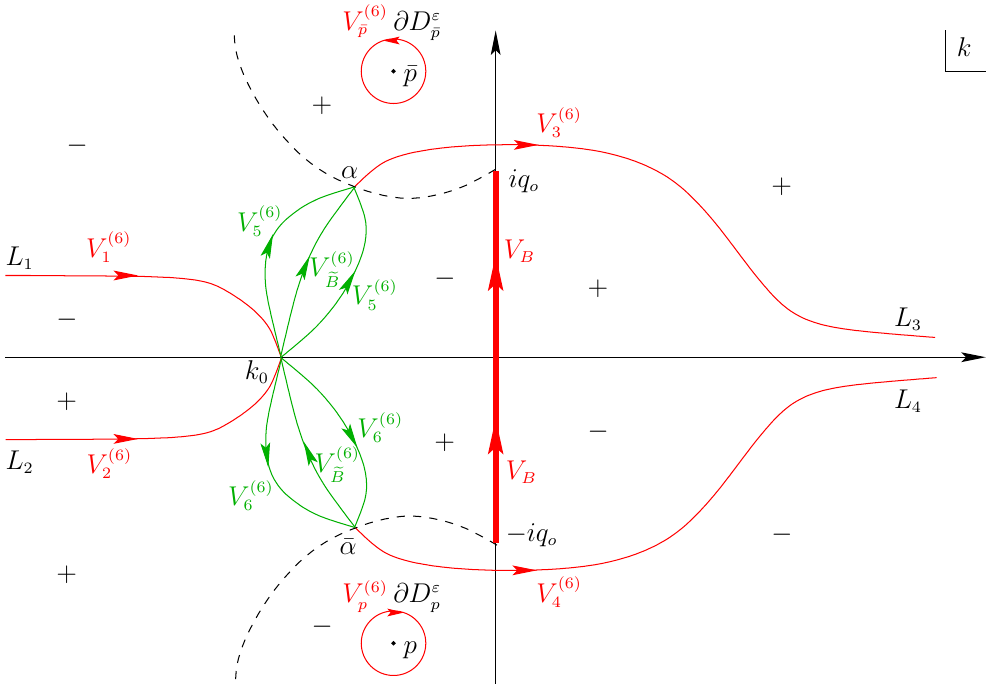}
\caption{Modulated elliptic wave in the trap regime:  the jumps of  Riemann-Hilbert problem \eqref{ds-esc-n6-rhp-mew}.
Contrary to the corresponding region in the  transmission regime, the jumps along $\p D_p^\ve$ and $\p D_{\bar p}^\ve$ tend to the identity as $t\to\infty$.
}
\label{ds-trap-def4-mew-f}
\end{center}
\end{figure}

This range is very similar to the range $(\Vo ,  0)$  of the  transmission regime. In particular, applying the first six deformations of Subsection \ref{ds-esc-mew-ss} we arrive again at Riemann-Hilbert problem \eqref{ds-esc-n6-rhp-mew} for the function $N^{(6)}$.
For $p\in D_2^+$, however, $\Re (ih)(\xi, p)<0$ as opposed to $\Re (ih)(\xi, p)>0$ (see Figure~\ref{ds-trap-def4-mew-f} and the  relevant discussion in Subsection \ref{ds-esc-mew-ss}). Therefore, the jumps along $\p D_p^\ve$ and $\p D_{\bar p}^\ve$  now tend to the identity exponentially fast as $t\to\infty$, allowing us to proceed to the decomposition of problem \eqref{ds-esc-n6-rhp-mew} into dominant and error components directly,  without the need for  transformations \eqref{ds-esc-n5t-def-mew} and \eqref{ds-esc-n6t-def-mew}.

Indeed, performing the analogue of decomposition \eqref{ds-n5t-decomp-mew} and proceeding as in Subsection \ref{ds-esc-mew-ss}, we find  
\eee{\label{ds-q-recon-mew-3}
q(x, t)
=
-2i \lim_{k \to \infty} k N_{12}^\dom(x, t, k) e^{i\left[g_\infty(\xi) -G_\infty(\xi) t  \right]} 
+
O\big(t^{-\frac 12}\big), 
}
where $N^\dom$ denotes the solution of the dominant component of Riemann-Hilbert problem \eqref{ds-esc-n6-rhp-mew} in the case $p\in D_2^+$. Specifically, as expected from the discussion above,  $N^\dom$ satisfies problem \eqref{ds-esc-nmodt-rhp-mew} with $\widetilde \omega = \widetilde g_\infty = 0$, i.e. $N^\dom$ is analytic in $\mathbb C\setminus(B \cup \widetilde B)$ with
\sss{\label{ds-trap-nmod-rhp-mew1}
\ddd{
N^{\dom+} &= N^{\dom-} V_B,  && k\in B,
\\
N^{\dom+} &= N^{\dom-} V_{\widetilde B}^{(6)}, && k\in \widetilde B,
\\
N^\dom &=  
 \left[I+O\left(\tfrac 1k\right)\right]e^{i\left[g_\infty(\xi)  -G_\infty(\xi) t  \right]\sigma_3},\quad && k \to \infty.
}
}
Problem \eqref{ds-trap-nmod-rhp-mew1} arises  in the case of empty discrete spectrum analyzed in \cite{bm2017}. 
Actually, thanks to the fact that the jump $V_{\widetilde B}^{(6)}$ is independent of $k$, the jump contour $\widetilde B$ in problem \eqref{ds-trap-nmod-rhp-mew1}  can be deformed  to the straight line segment $B'$ from $\bar \alpha$ to $\alpha$ (see Figure \ref{ds-esc-nmodt-jumps-mew-f}).
Then, following \cite{bm2017}, we obtain the solution of this deformed problem as
\eee{
\label{ds-trap-nmod-sol-mew1}
N^\dom(x, t, k)
=
e^{i\left[g_\infty(\xi) - G_\infty(\xi) t  \right]\sigma_3}
\mathcal N^{-1}(\xi, \infty, c) \,  \mathcal N(\xi, k, c),
}
where
\eee{\label{ds-trap-ncal-def-mew1}
\mathcal N(\xi, k, c)
=
\frac 12
\def\arraystretch{1.5}
\left(
\begin{array}{lr}
\left[\eta(\xi, k)+\eta^{-1}(\xi, k)\right]\@N_1(\xi, k, c)
& 
i\left[\eta(\xi, k)-\eta^{-1}(\xi, k)\right]\@N_2(\xi, k, c)
\\
-i\left[\eta(\xi, k)-\eta^{-1}(\xi, k)\right]\@N_1(\xi, k, -c)
&
\left[\eta(\xi, k)+\eta^{-1}(\xi, k)\right]\@N_2(\xi, k, -c)
\end{array}
\right)
}
and 
\eee{
\mathcal N(\xi, \infty, c) := \lim_{k\to\infty} \mathcal N(\xi, k, c)
}
with the function $\eta$ defined by \eqref{ds-pdef} and with  $\@N_1$ and $\@N_2$ denoting the first and second column of the vector-valued function
\eee{
\label{ds-mcaldef}
\@N (\xi, k, c) 
=
\left(
\frac{
 \Theta\big(-\frac{\Omega t}{2\pi}+\frac{\omega}{2\pi}+\frac{i\ln\left(\frac{\bar q_-}{iq_o}\right)}{2\pi}+\upnu(k)+c \big)}
 {  \sqrt{\frac{iq_o }{\bar q_-}}\ \Theta\left(\upnu(k)+c \right)},
\frac{\Theta\big(-\frac{\Omega t}{2\pi}+\frac{\omega}{2\pi}+\frac{i\ln\left(\frac{\bar q_-}{iq_o}\right)}{2\pi}-\upnu(k)+c \big)}
{ \sqrt{\tfrac{\bar q_-}{iq_o }}  \ \Theta\left(-\upnu(k)+c \right)}
\right),
}
where $\Omega$, $\omega$, $\upnu$  and $c$ are given by  \eqref{ds-Omega-exp}, \eqref{ds-omr}, \eqref{ds-vdefr0} and \eqref{ds-c-choice} respectively. 
We note that formula \eqref{ds-trap-nmod-sol-mew1} is consistent with formula \eqref{ds-esc-ntmod-sol-mew} after  setting $\widetilde \omega = \widetilde g_\infty = 0$.

Recall that the original and deformed versions of $N^\dom$ agree outside the finite region $\mathcal D$ enclosed by $\widetilde B$ and $B'$ (see Figure \ref{ds-esc-nmodt-jumps-mew-f}) and hence in the limit $k\to \infty$. Thus, inserting the solution \eqref{ds-trap-nmod-sol-mew1} in the reconstruction formula \eqref{ds-q-recon-mew-3} and utilizing the explicit form \eqref{ds-Omega-exp} of $\Omega$ together with the theta functions manipulations performed in \cite{bm2017}, we obtain the leading-order asymptotics \eqref{ds-trap-qasym-mew-t}-\eqref{ds-qsol-mew-cpam-t}.

\subsection{The case $\xi = \Vd$: soliton on top of a modulated elliptic wave}
\label{ds-trap-mew1-lim-ss}
Recall that in the trap regime currently under consideration  the value $\Vd$ is the unique solution of equation \eqref{ds-xih-def} in the interval $(\Vo, 0)$. That is, for $p\in D_2^+$ the quantities $\Re (ih)(\xi, p)$ and $\Re (ih)(\xi, \bar p)$ vanish inside $(\Vo, 0)$ only at $\xi=\Vd$. 
In turn, the jumps $V_p^{(6)}$ and $V_{\bar p}^{(6)}$ given by \eqref{ds-vp5-vpb5-def-mew} become part of the dominant component of Riemann-Hilbert problem \eqref{ds-esc-n6-rhp-mew} only for $\xi=\Vd$.
Indeed, as noted earlier$^{\ref{ds-p-foot}}$ and will be confirmed below,  whenever these jumps are part of the dominant problem they are eventually converted to residue conditions at $p$ and $\bar p$. Thus, the relevant exponentials reduce to $e^{\pm 2i h(\xi, p)t}$, which for $p\in D_2^+$ are purely oscillatory (as opposed to growing or decaying) only for $\xi = \Vd$.
On the other hand,  thanks to the global sign structure of $\Re (ih)$ (see Figure \ref{ds-trap-def4-mew-f}) the jumps $V_j^{(6)}$ tend to the identity exponentially fast as $t\to\infty$, like in the range $(\Vo, \Vd)$.  

\begin{figure}[t!]
\begin{center}
\includegraphics[width=0.238\textwidth]{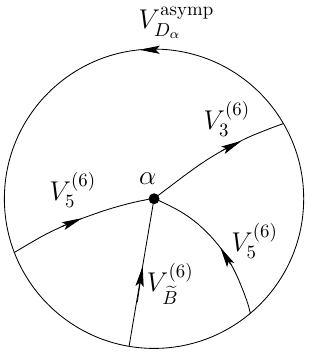}
\hskip 4mm
\includegraphics[width=0.241\textwidth]{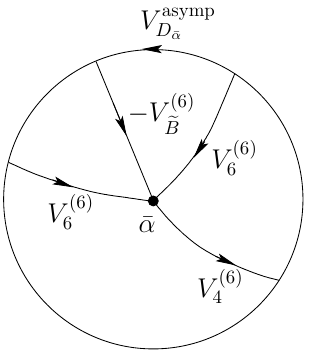}
\hskip 5mm
\includegraphics[width=0.233\textwidth]{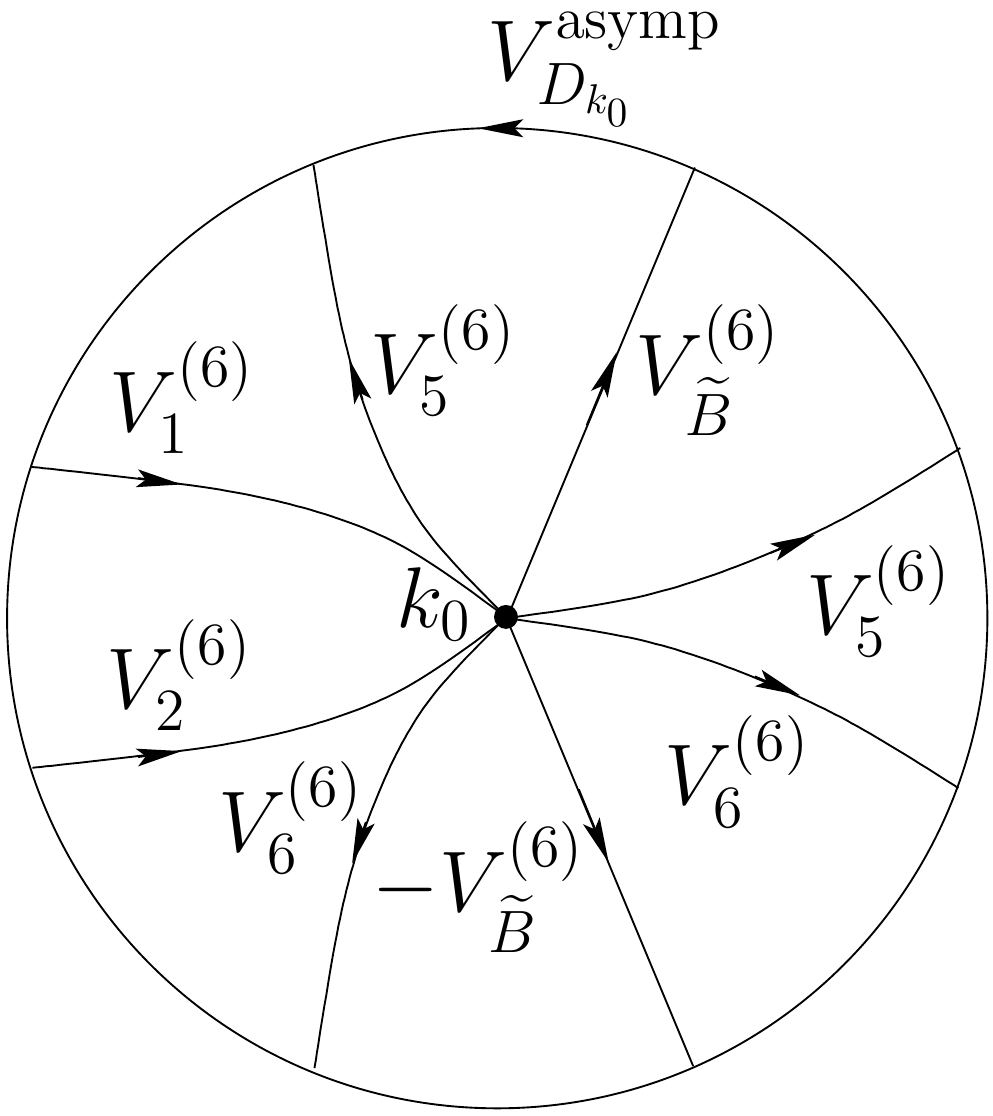}
\caption{Modulated elliptic wave  in the  trap regime: The jumps of  $N^D$  in the interior of and along the boundary of the disks $\overline{D_\alpha^\epsilon}$, $\overline{D_{\bar \alpha}^\epsilon}$ and $\overline{D_{k_o}^\epsilon}$. Although the jumps $V_{D_\alpha}^\asymp$, $V_{D_{\bar \alpha}}^\asymp$, $V_{D_{k_o}}^\asymp$ are unknown, they are equal to the identity up to $O(t^{-1/2})$ and hence do not affect the dominant problem.}
\label{ds-trap-nd-jumps-mew1-f}
\end{center}
\end{figure}

Following the above remarks, for $\xi = \Vd$ we write the solution of problem \eqref{ds-esc-n6-rhp-mew} as 
\eee{\label{ds-n5-decomp-mew1-lim}
N^{(6)}
=
N^\err N^\asymp,
}
where, for disks $D_{k_o}^{\epsilon}$, $D_{\alpha}^{\epsilon}$, $D_{\bar \alpha}^{\epsilon}$  of radius $\epsilon$ centered at $k_o$, $\alpha$, $\bar \alpha$ and such that they do not intersect with each other or with $B\cup \overline{D_p^\ve} \cup \overline{D_p^\ve}$, we let 
\eee{
N^\asymp
=
\begin{cases}
N^\dom,  & k\in \mathbb C\setminus (D_{k_o}^{\epsilon}\cup D_{\alpha}^{\epsilon}\cup D_{\bar \alpha}^{\epsilon}),
\\
N^D,  & k\in D_{k_o}^{\epsilon}\cup D_{\alpha}^{\epsilon}\cup D_{\bar \alpha}^{\epsilon},
\end{cases}
}
and define the functions $N^\dom$, $N^D$ and $N^\err$  as follows:
\vskip 3mm
\begin{enumerate}[label=$\bullet$, leftmargin=4mm, rightmargin=0mm]
\advance\itemsep 3mm
\item
$N^\dom(\Vd t, t, k)$ is analytic in $\mathbb C\setminus(B \cup \widetilde B \cup \p D_p^\ve \cup \p D_{\bar p}^\ve)$  and satisfies the Riemann-Hilbert problem 
\sss{\label{ds-trap-nmod-rhp-mew1-lim}
\ddd{
N^{\dom+} &= N^{\dom-} V_B,  && k\in B,
\label{ds-trap-nmod-jumpb-mew1-lim}
\\
N^{\dom+} &= N^{\dom-} V_{\widetilde B}^{(6)}, && k\in \widetilde B,
\label{ds-trap-nmod-jumpbt-mew1-lim}
\\
N^{\dom+} &= N^{\dom-} V_p^{(6)}, && k\in \p D_p^\ve,
\\
N^{\dom+} &= N^{\dom-} V_{\bar p}^{(6)}, && k\in \p D_{\bar p}^\ve,
\\
N^\dom &=  
 \left[I+O\left(\tfrac 1k\right)\right]e^{i\left[g_\infty(\Vd) -G_\infty(\Vd) t  \right]\sigma_3},\quad && k \to \infty.
 \label{ds-trap-nmod-asymp-mew1-lim}
}
}
\item
$N^D(\Vd t, t, k)$ is analytic in $D_{k_o}^{\epsilon}\cup D_{\alpha}^{\epsilon}\cup D_{\bar \alpha}^{\epsilon}\setminus  \bigcup_{j=1}^{8}  L_j$ with jumps 
\eee{
N^{D+} = N^{D-} V_j^{(6)},
\quad
k\in \widehat L_j := L_j\cap \left(D_{k_o}^{\epsilon}\cup D_{\alpha}^{\epsilon}\cup D_{\bar \alpha}^{\epsilon}\right),
\ 
j=1, \ldots, 8,
}
as shown in  Figure \ref{ds-trap-nd-jumps-mew1-f}.
\item
$N^\err(\Vd t, t, k)$ is analytic in $\mathbb C
\setminus \big(\bigcup_{j=1}^{6} \wc L_j \cup \p D_{k_o}^{\epsilon}\cup \p D_{\alpha}^{\epsilon}\cup \p D_{\bar \alpha}^{\epsilon}\big)$ with $\wc L_j := L_j\setminus ( \overline{D_{k_o}^{\epsilon}}\cup \overline{D_{\alpha}^{\epsilon}}\cup \overline{D_{\bar \alpha}^{\epsilon}})$ and satisfies the Riemann-Hilbert problem 
\sss{\label{ds-nerr-trap-sol}
\ddd{
N^{\err+} &= N^{\err-} \, V^\err, \quad
&&  k\in {\textstyle \bigcup}_{j=1}^{6} \wc L_j \cup \p D_{k_o}^{\epsilon}\cup \p D_{\alpha}^{\epsilon}\cup \p D_{\bar \alpha}^{\epsilon},
\\
N^\err &= I +O\left(\tfrac 1k\right), && k \to \infty, \label{ds-trap-ne-asymp-mew1-lim}
}
}
with
\eee{
V^\err
=
\begin{cases}
N^\dom V_j^{(6)} (N^\dom)^{-1}, &k\in \wc L_j,
\\
N^{\asymp-}(V_D^\asymp)^{-1}(N^{\asymp-})^{-1}, &k\in \p D_{k_o}^{\epsilon}\cup \p D_{\alpha}^{\epsilon}\cup \p D_{\bar \alpha}^{\epsilon},
\end{cases}
}
and
\eee{
V_D^\asymp
=
\begin{cases}
V_{D_\alpha}^\asymp, &k\in \p D_{\alpha}^\epsilon,
\\
V_{D_{\bar \alpha}}^\asymp, &k\in \p D_{\bar \alpha}^\epsilon,
\\
V_{D_{k_o}}^\asymp, &k\in \p D_{k_o}^\epsilon.
\end{cases}
}
\end{enumerate}

Although $V_D^\asymp$ is unknown, in  \cite{bm2017} it was shown that the contribution of the error problem \eqref{ds-nerr-trap-sol} to the leading-order asymptotics is of $O\big(t^{-1/2}\big)$. Therefore, starting from the reconstruction formula \eqref{ds-q-recon-n} and applying the six  deformations that lead to $N^{(6)}$,  we find
\eee{\label{ds-trap-q-recon-asymp-mew1}
q(x, t)
=
-2i \lim_{k \to \infty} k   N_{12}^\dom(\Vd t, t, k) e^{i\left[g_\infty(\Vd) - G_\infty(\Vd) t  \right]} + O\big(t^{-\frac 12}\big),
\quad t\to \infty.
}
It remains to determine $N^\dom$.
 
\vskip 3mm
\noindent
\textbf{Solution of the dominant problem.} 
We now determine the solution $N^\dom$ of problem \eqref{ds-trap-nmod-rhp-mew1-lim}  via the same series of steps followed in Subsection \ref{ds-esc-pw1-lim-ss}.
First, we convert the jumps $V_p^{(6)}$ and $V_{\bar p}^{(6)}$ along the circles $\p D_p^\ve$ and $\p D_{\bar p}^\ve$  to  residue conditions at $p$ and $\bar p$ via the transformation
\eee{\label{ds-trap-mmod-def-mew1}
M^\dom 
=
\begin{cases}
N^\dom  \big(V_p^{(6)}\big)^{-1}, &k\in D_p^{\ve },
\\
N^\dom, &k\in \mathbb C^-\setminus \big(B^-\cup \widetilde B^- \cup \overline{D_{p}^{\ve}}\,\big),
\\
N^\dom \big(V_{\bar p}^{(6)}\big)^{-1}, &k\in D_{\bar p}^{\ve },
\\
N^\dom, &k\in \mathbb C^+\setminus \big(B^+\cup \widetilde B^+ \cup \overline{D_{\bar p}^{\ve}}\, \big).
\end{cases}
}
Indeed, the function $M^\dom$ maintains the jumps of $N^\dom$ along $B$ and $\widetilde B$ but has simple poles at $p$ and $\bar p$ instead of jumps along $\p D_p^\ve$ and $\p D_{\bar p}^\ve$. Specifically, $M^\dom$ satisfies the Riemann-Hilbert problem
\sss{\label{ds-trap-rhp-mmod-mew1}
\ddd{
&M^{\dom+} = M^{\dom-} V_B,  && k\in B,
\\
&M^{\dom+} = M^{\dom-} V_{\widetilde B}^{(6)},  && k\in \widetilde B,
\\
&M^\dom =  \left[I +O\left(\tfrac 1k\right)\right]e^{i\left[g_\infty(\Vd)-G_\infty(\Vd)t\right] \sigma_3},\quad  && k \to \infty,
\\
& \underset{k=p}{\text{Res}}\, M^\dom = \left(0, \rho_p\,  M_1^\dom(p)\right),
\label{ds-trap-rhp-mmod-p-mew1}
\\
& \underset{k=\bar p}{\text{Res}}\,  M^\dom = \left(\rho_{\bar p}\,  M_2^\dom(\bar p), 0\right),
\label{ds-trap-rhp-mmod-pb-mew1}
}
}
where $M_1^\dom, M_2^\dom$ denote the first and second column of the matrix $M^\dom$ and
\sss{\label{ds-rho-def}
\ddd{
&\rho_p = c_p \delta^2(\Vd, p) d(p) e^{2i\left[h(\Vd, p)t-g(\Vd, p)\right]},
\\
&\rho_{\bar p} = c_{\bar p} \delta^{-2}(\Vd, \bar p) d(\bar p) e^{-2i\left[h(\Vd, \bar p)t-g(\Vd, \bar p)\right]}.
}
}
In fact, similarly to Subsection \ref{ds-esc-pw1-lim-ss}, we can use the definitions \eqref{ds-cp-def} and \eqref{ds-cpb-def} of $c_p$ and $c_{\bar p}$ together with the symmetries $C_{\bar p} = - \overline{C_p}$, $d(\bar k) = \overline{d(k)}$, $\bar a'(\bar k) = \overline{a'(k)}$, $\delta(\xi, \bar k) = \overline{\delta^{-1}(\xi, k)}$,  $g(\xi, \bar k) = \overline{g(\xi, k)}$, $h(\xi, \bar k) = \overline{h(\xi, k)}$ and the fact that $h(\Vd, \bar p)\in \mathbb R$ to write expressions \eqref{ds-rho-def} as
\ddd{\label{ds-R-def-trap}
\rho_p 
= 
R_p \, e^{2ih(\Vd, p)t},
\quad
\rho_{\bar p} 
= 
-\overline{R_p} \, e^{-2ih(\Vd, p)t},
\quad
R_p:=C_p \frac{\delta^2(\Vd, p)  e^{-2i g(\Vd, p)}}{a'(p)}.
}
Note that the writing \eqref{ds-R-def-trap} reveals that $\rho_{\bar p} = -\overline{\rho_p}$.

In order to solve problem \eqref{ds-trap-rhp-mmod-mew1}, it is convenient to let
\eee{\label{ds-trap-mmodt-def-mew1}
M^\dom = \mathcal M^\dom W
}
with $W$ being the solution of the continuous spectrum component problem
\sss{\label{ds-trap-nmod-rhp-mew1-w}
\ddd{
W^+ &= W^- V_B,  && k\in B,
\\
W^+ &= W^- V_{\widetilde B}^{(6)}, && k\in \widetilde B,
\\
W &=  
 \left[I+O\left(\tfrac 1k\right)\right]e^{i\left[g_\infty(\Vd)  -G_\infty(\Vd) t  \right]\sigma_3},\quad && k \to \infty.
 \label{ds-w-as-mew}
}
}
Observe that problem \eqref{ds-trap-nmod-rhp-mew1-w} is simply problem \eqref{ds-trap-nmod-rhp-mew1} evaluated at $\xi=\Vd$. Therefore, as for problem \eqref{ds-trap-nmod-rhp-mew1}, its jump contour $\widetilde B$ in problem \eqref{ds-trap-nmod-rhp-mew1-w} can be deformed to the straight line segment $B'$ connecting $\bar \alpha$ with $\alpha$ (see Figure \ref{ds-esc-nmodt-jumps-mew-f}), and the solution of this deformed problem is given by formula \eqref{ds-trap-nmod-sol-mew1} as 
\eee{\label{ds-trap-W-mew1}
W
=
e^{i\left[g_\infty(\Vd) - G_\infty(\Vd) t  \right]\sigma_3}
\mathcal N^{-1}(\Vd, \infty, c) \,  \mathcal N(\Vd, k, c)
}
with $\mathcal N$ defined by \eqref{ds-trap-ncal-def-mew1}.
Note further that $\det W$ inherits the analyticity of $W$ away from $B$ and $\widetilde B$, while \eqref{ds-trap-nmod-rhp-mew1-w} implies that
\ddd{
\det W^+ &= \det W^- \det V_B,  && k\in B,
\nn\\
\det W^+ &= \det W^- \det V_{\widetilde B}^{(6)}, \quad &&  k\in \widetilde B.
\nn
}
Therefore, since $\det V_B \equiv  \det V_{\widetilde B}^{(6)} \equiv 1$, we deduce that $\det W$ does not have jumps along $B$ and $\widetilde B$, i.e. $\det W$ is entire in $k$. Moreover, the asymptotic condition \eqref{ds-w-as-mew} implies that $\lim_{k\to \infty} \det W =1$. Thus, we conclude via Liouville's theorem that $\det W = 1$ for all $k\in\mathbb C$.

Combining \eqref{ds-trap-mmodt-def-mew1} and \eqref{ds-trap-nmod-rhp-mew1-w}, we find that the discrete component $\mathcal M^\dom$ of $M^\dom$  is analytic in $\mathbb C\setminus\left\{ p, \bar p\right\}$ and has simple poles at $p$ and $\bar p$ with the following residues:
\sss{\label{ds-trap-mmodt-res-mew1}
\ddd{
\underset{k=p}{\text{Res}}\, \mathcal M_1^\dom
&=
-W_{21}(p)  \rho_p \, M_1^\dom(p),
\\
\underset{k=p}{\text{Res}}\, \mathcal M_2^\dom
&=
W_{11}(p)   \rho_p \, M_1^\dom(p),
\\
\underset{k=\bar p}{\text{Res}}\, \mathcal M_1^\dom
&=
W_{22}(\bar p)   \rho_{\bar p} \, M_2^\dom(\bar p),
\\
\underset{k=\bar p}{\text{Res}}\, \mathcal M_2^\dom
&=
-W_{12}(\bar p)   \rho_{\bar p} \, M_2^\dom(\bar p).
}
}
Moreover, $\mathcal M^\dom$ satisfies the asymptotic condition 
\eee{\label{ds-trap-mmod-asymp-mew1}
\mathcal M^\dom
=
I+O\left(\tfrac 1k\right), \quad k\to\infty.
}
Then, arguing as in Subsection \ref{ds-esc-pw1-lim-ss}, we deduce that
\eee{\label{ds-trap-mmodt-form-mew1}
\mathcal M^\dom
=
I + \frac{\underset{k=p}{\text{Res}}\, \mathcal M^\dom}{k-p} + \frac{\underset{k=\bar p}{\text{Res}}\, \mathcal M^\dom}{k-\bar p}.
}
Thus, in order to determine $\mathcal M^\dom$ it suffices to determine its two residues at $p$ and $\bar p$.

From \eqref{ds-trap-mmodt-def-mew1} we have
\sss{
\ddd{
M_1^\dom
&=
W_1
+
W_{11} 
\left[
-
\frac{W_{21}(p)\,  \rho_p \, M_1^\dom(p)}{k-p} 
+ 
\frac{W_{22}(\bar p)\,  \rho_{\bar p} \, M_2^\dom(\bar p)}{k-\bar p} 
\right]
\nn\\
&\quad
+W_{21} 
\left[
\frac{W_{11}(p)\,  \rho_p \, M_1^\dom(p)}{k-p} 
- 
\frac{W_{12}(\bar p)\,  \rho_{\bar p} \, M_2^\dom(\bar p)}{k-\bar p} 
\right]
}
and
\ddd{
M_2^\dom
&=
W_2
+
W_{12} 
\left[
-
\frac{W_{21}(p)\,  \rho_p \, M_1^\dom(p)}{k-p} 
+ 
\frac{W_{22}(\bar p)\,  \rho_{\bar p} \, M_2^\dom(\bar p)}{k-\bar p} 
\right]
\nn\\
&\quad
+W_{22} 
\left[
\frac{W_{11}(p)\,  \rho_p \, M_1^\dom(p)}{k-p} 
- 
\frac{W_{12}(\bar p)\,  \rho_{\bar p} \, M_2^\dom(\bar p)}{k-\bar p} 
\right].
}
}
Evaluating the first of the above equations at $k=p$ and the second one at $k=\bar p$ (recall that $M_1^\dom$ and $M_2^\dom$ are analytic at $p$ and $\bar p$ respectively), we obtain the system
\sss{\label{ds-trap-mmod12-sys-mew1-0}
\ddd{
M_1^\dom(p)
&=
W_1(p) 
+ 
\rho_{\bar p}  \, \frac{W_{11}(p)W_{22}(\bar p) - W_{21}(p) W_{12}(\bar p)}{p-\bar p} \, M_2^\dom(\bar p)
\nn\\
&\quad
+
\rho_p \left[W_{21}'(p) W_{11}(p)-W_{11}'(p) W_{21}(p)\right]  M_1^\dom(p),
\\
M_2^\dom(\bar p)
&=
W_2(\bar p)
+
\rho_p \, \frac{W_{12}(\bar p) W_{21}(p) - W_{22}(\bar p) W_{11}(p)}{p-\bar p} \, M_1^\dom(p)
\nn\\
&\quad
+  
\rho_{\bar p}  \left[W_{12}'(\bar p)W_{22}(\bar p) - W_{22}'(\bar p) W_{12}(\bar p)\right] M_2^\dom(\bar p),
}
}
which can be solved to yield
\sss{\label{ds-trap-mmod12-sys-sol-mew1}
\ddd{
M_1^\dom(p)
&=
\frac{- \mathcal B\rho_{\bar p} W_2(\bar p)+\left(1+\mathcal C\rho_{\bar p}\right) W_1(p)}{\mathcal B^2 \rho_p\rho_{\bar p} + \left(1+ \mathcal C\rho_{\bar p}\right)\left(1+\mathcal A\rho_p\right)},
\\
M_2^\dom(\bar p)
&=
 \frac{\mathcal B\rho_p W_1(p)+\left(1+\mathcal A\rho_p\right) W_2(\bar p)}{\mathcal B^2 \rho_p\rho_{\bar p} + \left(1 + \mathcal C\rho_{\bar p}\right)\left(1+\mathcal A\rho_p\right)},
}
}
where
\sss{\label{ds-trap-abc-def}
\ddd{
\mathcal A &= W_{11}'(p)W_{21}(p)-W_{11}(p)W_{21}'(p),
\\
\mathcal B &= \frac{W_{21}(p)W_{12}(\bar p)-W_{11}(p)W_{22}(\bar p)}{p-\bar p},
\\
\mathcal C &= W_{22}'(\bar p)W_{12}(\bar p)-W_{12}'(\bar p)W_{22}(\bar p).
}
}
Expressions \eqref{ds-trap-mmod12-sys-sol-mew1} determine $\mathcal M^\dom$ through  \eqref{ds-trap-mmodt-form-mew1} and the residue relations \eqref{ds-trap-mmodt-res-mew1}. 

Having computed $\mathcal M^\dom$, we return to the reconstruction formula \eqref{ds-trap-q-recon-asymp-mew1} which upon \eqref{ds-trap-mmod-def-mew1} and \eqref{ds-trap-mmodt-def-mew1} reads
\eee{\label{ds-trap-q-recon-sol}
q(x, t) 
=
-2i\lim_{k \to\infty} k  \big(\mathcal M^\dom W\big)_{12} e^{i\left[g_\infty(\Vd)-G_\infty(\Vd)t\right]} + O\big(t^{-\frac 12}\big), \quad t\to\infty.
}
Now, the asymptotic conditions \eqref{ds-w-as-mew} and \eqref{ds-trap-mmod-asymp-mew1} imply 
\eee{
W  =  e^{i\left[g_\infty(\Vd)-G_\infty(\Vd)t\right]\sigma_3} + \frac{w}{k} + O\left(\frac{1}{k^2}\right),  \quad
\mathcal M^\dom  = I + \frac{\mu}{k} + O\left(\frac{1}{k^2}\right),   \quad k\to \infty,
\nn
}
where the matrix-valued functions $w$ and $\mu$ may depend on $x$ and $t$ but not on $k$. Therefore, 
\eee{
\big(\mathcal M^\dom W\big)_{12}
=
\frac{w_{12} + \mu_{12}\, e^{-i\left[g_\infty(\Vd)-G_\infty(\Vd)t\right]}}{k} + O\left(\frac{1}{k^2}\right), \quad  k\to \infty.
\nn
}
Substituting for $w_{12}$ via  \eqref{ds-trap-W-mew1} (note that the original and deformed versions of $W$ agree outside the finite region $\mathcal D$ enclosed by $\widetilde B$ and $B'$ in Figure \ref{ds-esc-nmodt-jumps-mew-f}  and hence as $k\to \infty$) and for $\mu_{12}$ via  \eqref{ds-trap-mmodt-res-mew1} and \eqref{ds-trap-mmodt-form-mew1} yields the leading-order asymptotics \eqref{ds-trap-q-sol-lim-pw1-t}-\eqref{ds-qsh-def} via the reconstruction formula \eqref{ds-trap-q-recon-sol}.

\begin{remark}[\b{Dependence on $g_\infty$ and $G_\infty$}]
\label{ds-overall-phase-r}
Formula \eqref{ds-trap-W-mew1} implies that $W_{11}$ and $W_{12}$ depend on $g_\infty$ and $G_\infty$ through the exponential $e^{i\left[g_\infty(\Vd) -  G_\infty(\Vd) t  \right]}$ while $W_{21}$ and $W_{22}$ instead contain the exponential $e^{-i\left[g_\infty(\Vd) -  G_\infty(\Vd) t  \right]}$. 
Hence, the quantities $\mathcal A, \mathcal B, \mathcal C$ defined by \eqref{ds-trap-abc-def} are independent of $g_\infty$ and $G_\infty$, and the overall dependence of the leading-order asymptotics \eqref{ds-trap-q-sol-lim-pw1-t} on $g_\infty$ and $G_\infty$ comes through a factor of $e^{2i\left[g_\infty(\Vd) -  G_\infty(\Vd) t  \right]}$.
\end{remark}

\subsection{The range $\Vd< \xi <0$: modulated elliptic wave with a phase shift}
\label{ds-trap-mew2-ss}
This range can be handled identically to the range $\Vo<\xi<0$ of the transmission regime that was analyzed in Subsection \ref{ds-esc-mew-ss}. Consequently, the leading-order asymptotics is characterized once again by \eqref{ds-esc-qsol-mew} as the modulated elliptic wave \eqref{ds-qsol-mew-t}  with a phase shift of $4\textnormal{arg}\left[p+\lambda(p)\right]$.

\vspace*{3mm}

The proof of Theorem \ref{ds-ptr-t} for the leading-order asymptotics in the trap regime $p\in D_2^+$ is complete.

%
%
%
%
\section{The Mixed Regimes: Proof of Theorems \ref{ds-twr-t} and \ref{ds-ewr-t}}
\label{ds-mixed-s}

In Sections \ref{ds-esc-s} and \ref{ds-trap-s}, we showed that the scenarios $p\in D_1$ and $p\in D_2^+$ give rise to \textit{pure} asymptotic regimes, namely a  transmission regime (Theorem \ref{ds-per-t}) and a  trap regime (Theorem \ref{ds-ptr-t}) respectively.
We now proceed to the analysis of the remaining two regions of Figure \ref{ds-regions-f}, namely  $D_2^-$ and $D_3$. We shall show that these regions correspond to \textit{mixed} asymptotic regimes, specifically a trap/wake regime (Theorem \ref{ds-twr-t}) and a   transmission/wake regime (Theorem \ref{ds-ewr-t})  respectively.

\subsection{The  trap/wake regime}
\label{ds-trap-wake-ss}

\begin{figure}[t!]
\begin{center}
\includegraphics[scale=1]{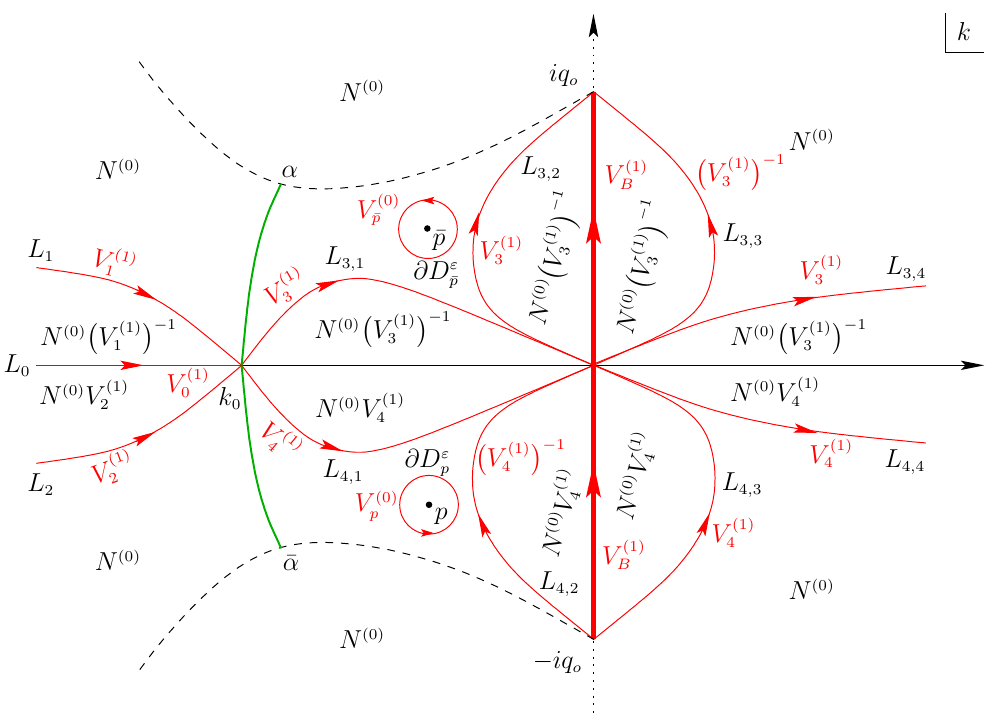}
\caption{Modulated elliptic wave in the trap/wake ($p\in D_2^-$) and transmission/wake ($p\in D_3$) regimes in  the ranges $\Vd<\xi<\Vw$ (for $p\in D_2^-$) and $v_o<\xi<\Vw$ (for $p\in D_3$): the initial stage of the first deformation. The jumps along $\p D_p^\ve$ and $\p D_{\bar p}^\ve$ are not affected.}
\label{ds-wake-pw2-f}
\end{center}
\end{figure}

Recall that for $p\in D_2^-$ we have $\Vo<\Vs<0$. Furthermore, as noted in Subsection \ref{ds-esc-mew-ss}, the integral equation \eqref{ds-fnlsd-int-eq-xistar-0} possesses exactly two solutions in the interval $(\Vo, 0)$: $\Vd$, which corresponds to the crossing of the pole $p$ by the dashed black curve in the third quadrant of Figure \ref{ds-wake-pw2-f} (for $\xi<\Vd$, the pole lies below this curve), and $\Vw>\Vd$, which corresponds to the crossing of $p$ by the branch cut $\widetilde B$  (green contour connecting $\alpha$ and $\bar \alpha$ in Figure \ref{ds-wake-pw2-f}). 
Note that the latter crossing can happen only if $p$ lies on the right (as opposed to the left) of $\widetilde B$ immediately after $\xi=\Vd$, and this is the way one distinguishes the trap/wake regime $p\in D_2^-$ from the trap regime $p\in D_2^+$.

For $\xi<\Vo$, the deformations performed in the trap regime  can be repeated  to lead once again to Riemann-Hilbert problem \eqref{ds-esc-n4-rhp-pw1}.
Furthermore, like in the trap regime, for $\xi<\Vo$ the dominant component of this problem only involves the jump along the branch cut $B$ since  the jumps along $\p D_p^\ve$ and $\p D_{\bar p}^\ve$ tend to the identity exponentially fast as $t\to \infty$ due to the fact that $\text{Re}(i\theta)(\xi, p)<0$ and $\text{Re}(i\theta)(\xi, \bar p)>0$ throughout the interval $(-\infty, \Vo)$.
Thus, the leading-order asymptotics for $p\in D_2^-$  and $\xi\in (-\infty, \Vo)$ is the same with the one of the trap regime, i.e. it is described by the plane wave \eqref{ds-qsol-pw-t}.

For $\Vo<\xi<\Vd$, the phase function switches from $\theta$ to $h$ via transformation \eqref{ds-m4r-0} and we eventually arrive at Riemann-Hilbert problem \eqref{ds-esc-n6-rhp-mew}. Moreover, we still have $\text{Re}(ih)(\xi, p)<0$ and $\text{Re}(ih)(\xi, \bar p)>0$, thus the jumps along $\p D_p^\ve$ and $\p D_{\bar p}^\ve$ still do not contribute to the leading-order asymptotics, which is  described by the modulated elliptic wave \eqref{ds-qsol-mew-cpam-t}.

At $\xi=\Vd$, we have $\text{Re}(ih)(\Vd, p)=\text{Re}(ih)(\Vd, \bar p)=0$. Thus, as explained in Subsection \ref{ds-trap-mew1-lim-ss}, the jumps along $\p D_p^\ve$ and $\p D_{\bar p}^\ve$ now contribute to the leading-order asymptotics, which is described  by \eqref{ds-trap-q-sol-lim-pw1-t} as the soliton \eqref{ds-qsh-def} on top of the modulated elliptic wave \eqref{ds-qsol-mew-cpam-t} evaluated at $\Vd$.

\begin{figure}[t!]
\begin{center}
\includegraphics[scale=1]{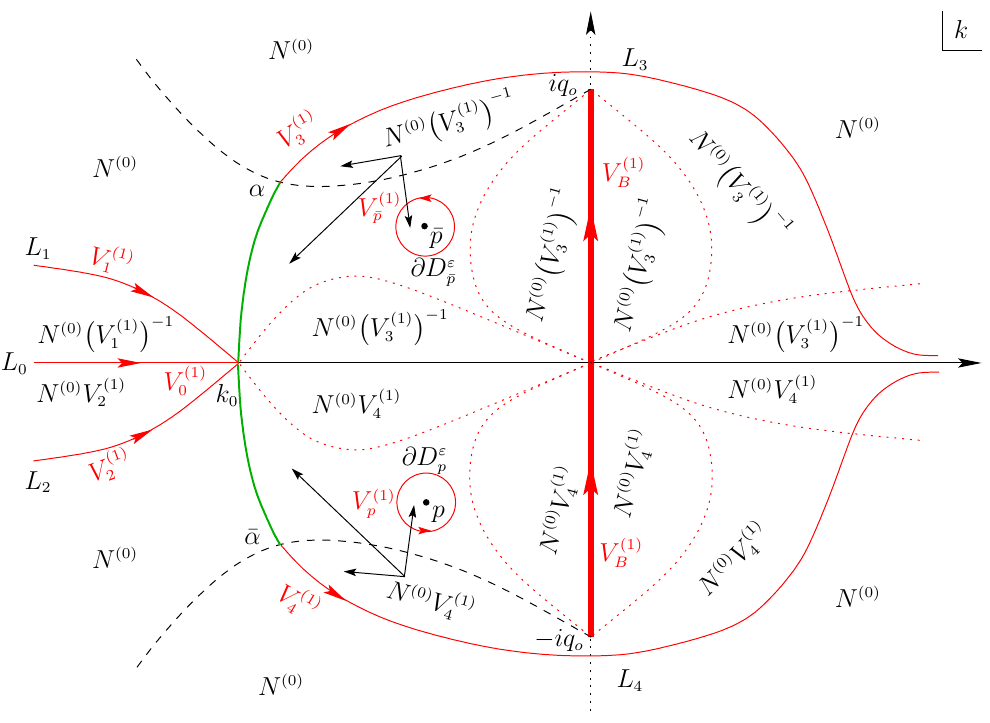}
\caption{Modulated elliptic wave in the trap/wake ($p\in D_2^-$) and transmission/wake ($p\in D_3$) regimes in  the ranges $\Vd<\xi<\Vw$ (for $p\in D_2^-$) and $v_o<\xi<\Vw$ (for $p\in D_3$): the final stage of the first deformation. The jumps along $\p D_p^\ve$ and $\p D_{\bar p}^\ve$ have now changed.}
\label{ds-wake-pw2-2-f}
\end{center}
\end{figure}

\begin{figure}[ht!]
\begin{center}
\includegraphics[scale=1]{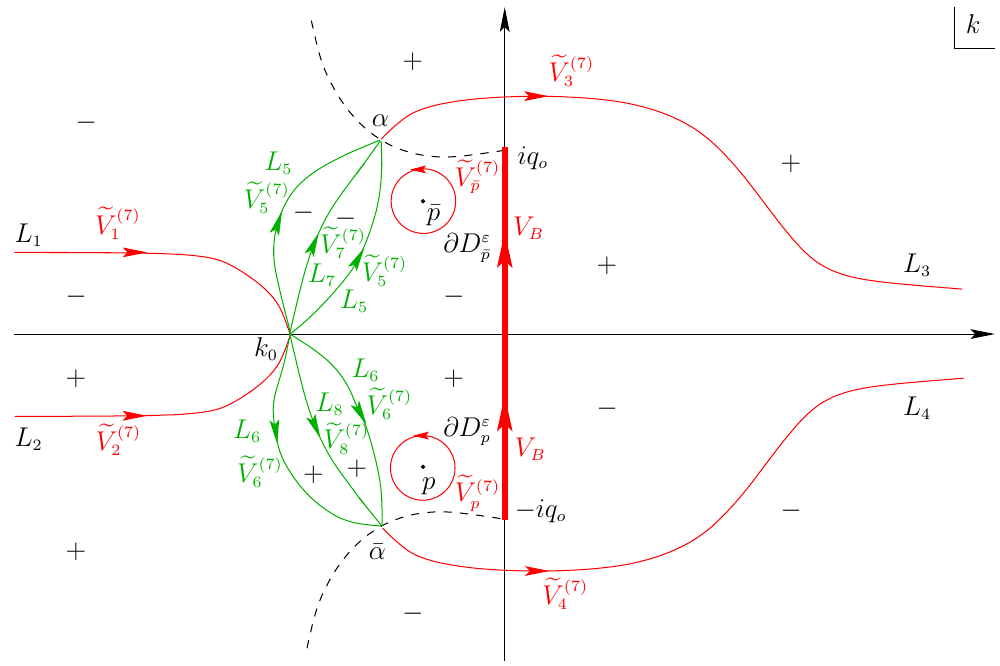}
\caption{Modulated elliptic wave  in the trap/wake  ($p\in D_2^-$) and transmission/wake regimes ($p\in D_3$): the jumps of the Riemann-Hilbert problem \eqref{ds-esc-n6t-rhp-mew}, the final problem in the ranges $\xi \in (\Vd, \Vw)\cup (\Vw, 0)$ (for $p\in D_2^-$) and $\xi \in (\Vo, \Vw)\cup (\Vw, 0)$ (for $p\in D_3$).}
\label{ds-trap-wake-def4-mew-f}
\end{center}
\end{figure}

The range $\Vd<\xi<\Vw$, which is not present in the trap regime since $\Vw$ does not arise there, requires a modification of the first four deformations. Specifically, while the first stage of the first deformation remains the same (compare Figure \ref{ds-wake-pw2-f}  with Figure \ref{ds-esc-def1a-mew-f}), the poles $p$ and $\bar p$ now lie on the \textit{right} of the branch cut $\widetilde B$. We emphasize that this is the defining difference between the trap/wake regime $p\in D_2^-$ and the trap regime $p\in D_2^+$, since in the latter case the poles are always on the \textit{left} of $\widetilde B$ for $\xi>\Vd$ (see also Figure \ref{f:h}).
Thus, for $\Vd<\xi<\Vw$ in the trap/wake regime, in order to lift the jump along $[k_o, 0]$ away from the real axis and onto  $\widetilde B$, the remaining stages of the first  deformation are adjusted from those of the trap regime to the factorization shown in Figure~\ref{ds-wake-pw2-2-f}. 
Then, applying the second deformation \eqref{ds-n2-def-pw1} with $\delta$ given by \eqref{ds-esc-del-def-mew} and the third deformation as shown in Figure \ref{ds-esc-def3-pw1-f} but with the disks $D_p^\ve$ and $D_{\bar p}^\ve$ now lying \textit{between} the contours $L_3$ and $L_4$, we obtain the analogue of  Figure \ref{ds-esc-def3-mew-f}, the only difference now being that the poles $p$  and $\bar p$ lie on the \textit{right} of $\widetilde B$. 
Subsequently, proceeding as in Subsection \ref{ds-esc-mew-ss}, we eventually arrive at the deformed Riemann-Hilbert problem of Figure \ref{ds-trap-wake-def4-mew-f}, which can be handled in the same way with Riemann-Hilbert problem \eqref{ds-esc-n6t-rhp-mew}. Indeed, since $\text{Re}(ih)(\xi, p)>0$ and $\text{Re}(ih)(\xi, \bar p)<0$, the jumps along $\p D_p^\ve$ and $\p D_{\bar p}^\ve$ are not significant at leading order and the corresponding  asymptotics is given by \eqref{ds-esc-qasym-mew-t} as the modulated elliptic wave \eqref{ds-qsol-mew-t}  with a phase shift of $4\textnormal{arg}\left[p+\lambda(p)\right]$.

For $\Vw<\xi<0$, the poles $p$ and $\bar p$ lie on the left of $\widetilde B$ and, therefore, the analysis is identical to the one for the trap regime in the range $\Vd<\xi<0$, leading once again to the asymptotics  \eqref{ds-esc-qasym-mew-t}.

It remains to analyze the case $\xi=\Vw$, which corresponds to the crossing of the poles $p$ and $\bar p$ by the branch cut $\widetilde B$ of $h$  (green contour connecting $\alpha$ and $\bar \alpha$ in Figure \ref{ds-wake-pw2-f}) as the latter  sweeps the region to its right en route to collapsing onto $B$ as $\xi \to 0^-$. Since $\widetilde B$ is also a zero-contour of $\text{Re}(ih)$, this is the mechanism giving rise to the second solution of equation \eqref{ds-fnlsd-int-eq-xistar-0}  when $p\in D_2^-$, since $\text{Re}(ih)(\xi, k)$ vanishes along $\widetilde B$ for all $\xi$ and hence $\text{Re}(ih)(\Vw, p)=0$.
As previously emphasized, the crossing of $p$ and $\bar p$ by $\widetilde B$ occurs for $p\in D_2^-$ but \textit{not} for $p\in D_2^+$, since in the latter case the poles are always on the left of $\widetilde B$ in the range $\Vd<\xi<0$.

\vskip 3mm

\begin{figure}[t!]
\begin{center}
\includegraphics[scale=1]{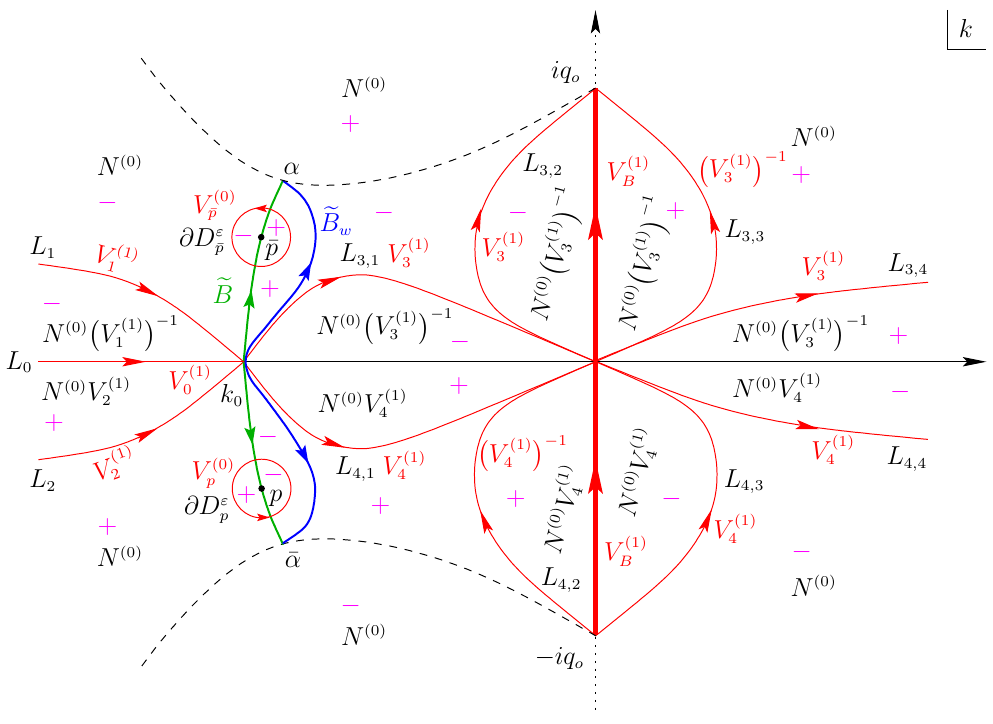}
\caption{Modulated elliptic wave in the trap/wake ($p\in D_2^-$) and transmission/wake ($p\in D_3$) regimes at $\xi=\Vw$: the sign structure of $\text{Re}(ih_w)$ and the initial stage of the first deformation. Note that the jumps along $\p D_p^\ve$ and $\p D_{\bar p}^\ve$ are not affected. Furthermore, the branch cut $\widetilde B_w$ of the function $h_w$ (blue contour) lies on the right of the zero-contour $\widetilde B$ of $\text{Re}(ih_w)$   (green contour) as well as on the right of the circles $\p D_p^\ve$ and $\p D_{\bar p}^\ve$. The finite region enclosed by $\widetilde B$ and $\widetilde B_w$ is denoted by $\mathcal R$.}
\label{ds-wake-vw-f}
\end{center}
\end{figure}

\noindent \textbf{The case $\xi=\Vw$: soliton wake.} 
For this value of $\xi$, the poles  $p$ and $\bar p$ lie on $\widetilde B$ (depicted in green in Figure \ref{ds-wake-pw2-f}), which is both a zero-contour  for $\text{Re}(ih)$ and a branch cut for $h$ (along with the branch cut $B=i[-q_o, q_o]$). 
For this reason, in view of the fifth deformation of Riemann-Hilbert problem \eqref{ds-n-rhp-intro} (see \eqref{ds-m4r-0}), it is convenient to switch from $h$ to a function $h_w$ which does not have a branch cut along $\widetilde B$ but which is such that $\text{Re}(ih_w)$ still vanishes along $\widetilde B$. 
More specifically, we define
\eee{\label{ds-h_w-def}
h_w(k) = \left\{\begin{array}{ll} 
h(\Vw, k), & k\in \mathbb C\setminus \overline{\mathcal R}, 
\\ 
\Omega (\Vw) - h(\Vw, k), & k\in  \mathcal R, 
\end{array}\right.
}
where the real constant $\Omega(\Vw)$ is given by \eqref{ds-Omega-exp} and $\mathcal R$ is the finite region enclosed by $\widetilde B$ and the contour $\widetilde B_w$ shown in blue in Figure \ref{ds-wake-vw-f}.
It is straightforward to see that the function $h_w$ (i) has branch cuts along $B$ and  $\widetilde B_w$, and (ii)  is continuous along $\widetilde B$.
Indeed, recall (see \eqref{ds-jhLL}) that along $\widetilde B$ we have $h^+ + h^- =  \Omega(\Vw)$. 
Hence,  according to the definition \eqref{ds-h_w-def} of $h_w$, along $\widetilde B$ we have $h_w^+  =\left(\Omega(\Vw)-h\right)^+ = \Omega(\Vw)-h^+ = h^-=h_w^-$, i.e.  $h_w$ is continuous along $\widetilde B$. 
On the other hand, along the contour $\widetilde B_w$ shown in blue in Figure \ref{ds-wake-vw-f}  we have  $h_w^+  = \left(\Omega(\Vw)-h\right)^+ = \Omega(\Vw) - h = \Omega(\Vw) - h_w^-$ (having used the fact that $h$ is continuous along $\widetilde B_w$). Hence, $h_w$ is discontinuous along $\widetilde B_w$ with $h_w^+ + h_w^- = \Omega(\Vw)$. Furthermore, $h_w$ is discontinuous along $B$ since it is equal to $h$ on both sides of $B$. Therefore, $h_w$ has branch cuts along $B$ and $\widetilde B_w$ and is continuous along $\widetilde B$.

\begin{figure}[t!]
\begin{center}
\includegraphics[scale=1]{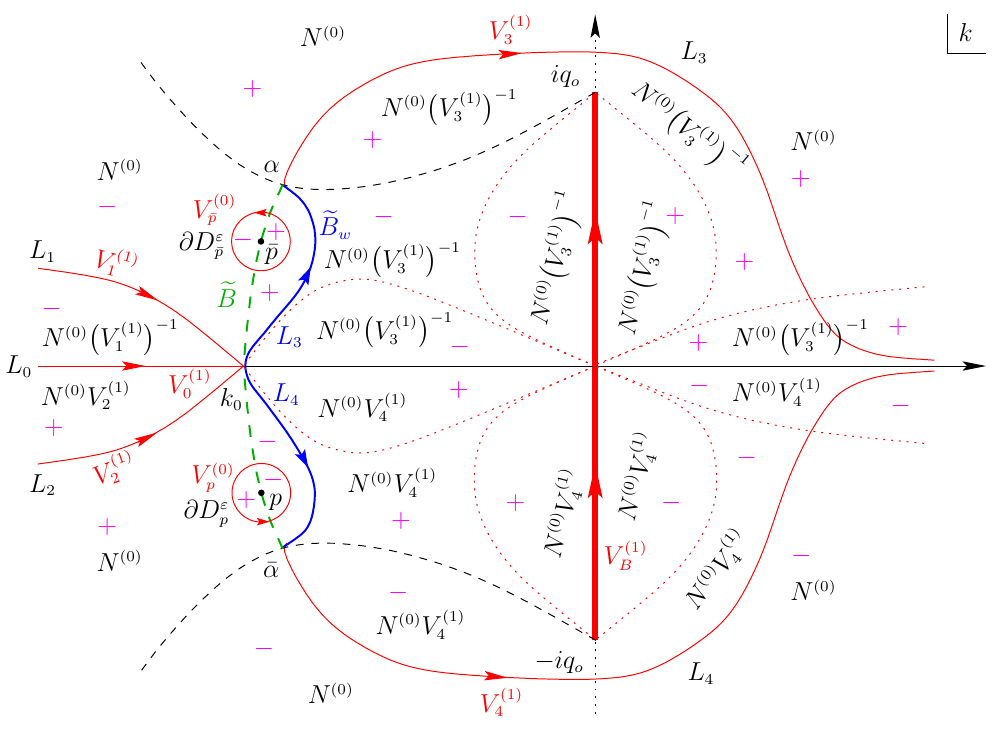}
\caption{Modulated elliptic wave in the trap/wake ($p\in D_2^-$) and transmission/wake ($p\in D_3$) regimes at $\xi=\Vw$: the final stage of the first deformation. The jumps along $\p D_p^\ve$ and $\p D_{\bar p}^\ve$ are not affected.}
\label{ds-wake-vw-2-f}
\end{center}
\end{figure}

The sign structure of $\text{Re}(ih_w)$ is shown in Figure \ref{ds-wake-vw-f}. With this in mind, we perform the first deformation according to Figure \ref{ds-wake-vw-f} and then  deform the contours $L_{3,1}$, $L_{3, 2}$ and $L_{4, 1}$, $L_{4,2}$ to the contours $L_3$ and $L_4$ of Figure \ref{ds-wake-vw-2-f}, which depicts the final stage of the first deformation. Note that $L_3$ consists of the upper half of the branch cut $\widetilde B_w$ as well as of the red contour starting from $\alpha$ and curving around $iq_o$ and down towards the positive real axis. Similarly, $L_4$ consists of both the lower half of $\widetilde B_w$ and the red curve emanating from $\bar \alpha$ and directed upwards towards the positive real axis.

\begin{figure}[t!]
\begin{center}
\includegraphics[scale=1]{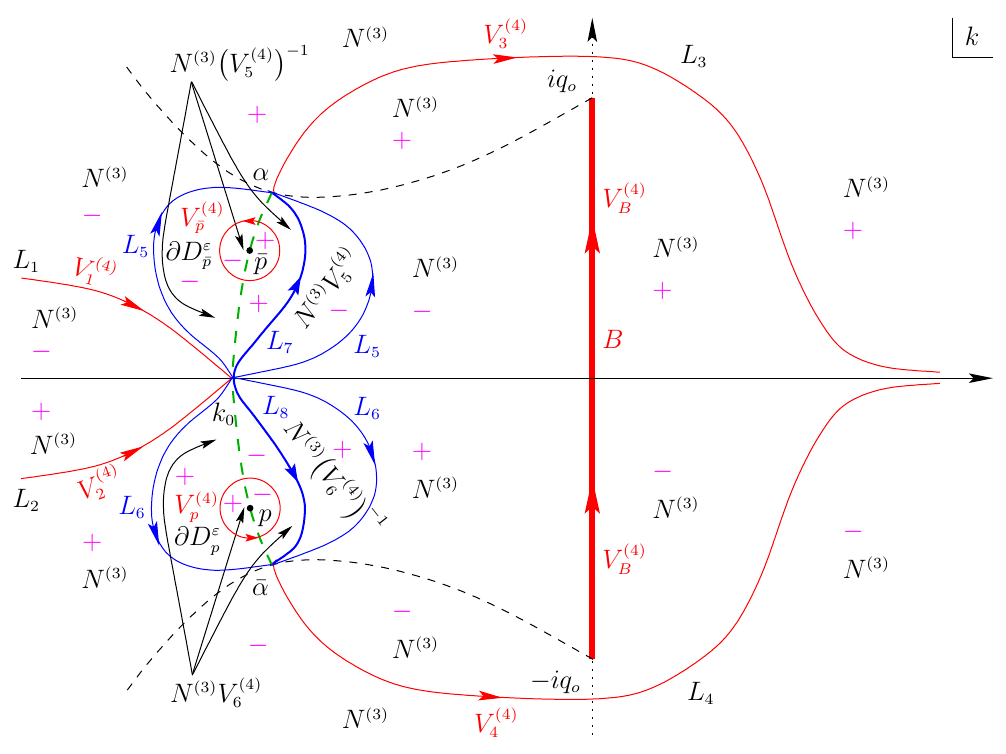}
\caption{Modulated elliptic wave in the trap/wake ($p\in D_2^-$) and transmission/wake ($p\in D_3$) regimes at $\xi=\Vw$: the fourth deformation. We recall that the dashed green contour is the branch cut $\widetilde B$ of $h$, along which $\text{Re}(ih)=\text{Re}(ih_w)=0$, while the function $h_w$ has branch cuts along $B=i[-q_o, q_o]$ and $\widetilde B_w = L_7 \cup (-L_8)$.}
\label{ds-wake-vw-3-f}
\end{center}
\end{figure}

The second and third deformations are identical to those performed in the trap regime, leading to Riemann-Hilbert problem \eqref{ds-esc-n3t-rhp-mew} but with the contours $L_j$, $j=1, 2, 3, 4$ as shown in Figure \ref{ds-wake-vw-2-f}.
For the fourth deformation, we use the factorizations \eqref{ds-gr-facs} to ``open up the lenses'' off the portions of the contours $L_3$ and $L_4$ that lie along $\widetilde B_w$  as shown in Figure \ref{ds-wake-vw-3-f}, where the jumps $V_j^{(4)}$, $j=5,6,7,8$ are given by \eqref{ds-lense-jumps}.
This figure provides the analogue of Figure \ref{ds-esc-def4-mew-f}, where the role of the contour $\widetilde B$ (which is a branch cut for $h$ and a zero-contour for $\text{Re}(ih)$ and, by definition~\eqref{ds-h_w-def}, for $\text{Re}(ih_w)$) is now held by the contour $\widetilde B_w$ (which is a branch cut for $h_w$). 
The difference between Figures \ref{ds-esc-def4-mew-f} and \ref{ds-wake-vw-3-f} is that in the latter case the disks $\overline{D_p^\ve}$ and $\overline{D_{\bar p}^\ve}$ lie between the contours $L_6$, $L_8$ and $L_5$, $L_7$ respectively. This is necessary in order for the contours $L_5$ and $L_6$ to lie in regions where the associated jumps $V_5^{(4)}$ and $V_6^{(4)}$ decay to the identity as $t\to \infty$.
Hence, in the fourth deformation shown in Figure \ref{ds-wake-vw-3-f}, $N^{(3)}$ changes to $N^{(3)} V_6^{(4)}$ in $D_p^\ve$ and to $N^{(3)} \big(V_5^{(4)}\big)^{-1}$ in $D_{\bar p}^\ve$ (in Figure \ref{ds-esc-def4-mew-f}, $N^{(3)}$ remains invariant inside the two disks). 
Consequently, $N^{(4)}$ satisfies the Riemann-Hilbert problem \eqref{ds-esc-n3t-rhp-mew-2} but with the jumps along $\p D_p^\ve$ and $\p D_{\bar p}^\ve$ now given by
\eee{
V_p^{(4)} = \big(V_6^{(4)}\big)^{-1} V_p^{(3)} V_6^{(4)}, 
\quad
V_{\bar p}^{(4)} =  V_5^{(4)} V_{\bar p}^{(3)} \big(V_5^{(4)}\big)^{-1},
}
where we recall that
\eee{
V_p^{(3)} = 
\left(
\def\arraystretch{1}
\begin{array}{lr}
1 & -\dfrac{c_p \, \delta^2(\Vw, k)\, d(k)}{k-p}\,  e^{2i\theta(\Vw, p)t}
\\
0 & 1
\end{array}
\right),
\quad
V_{\bar p}^{(3)}
=
\left(
\def\arraystretch{1}
\begin{array}{lr}
1 &0
\\
-\tfrac{c_{\bar p} \, \delta^{-2}(\Vw, k) \, d(k)}{k-\bar p}\,  e^{-2i\theta(\Vw, \bar p)t} & 1
\end{array}
\right)
}
and
\eee{
V_5^{(4)}
=
\left(
\def\arraystretch{1.2}
\begin{array}{lr}
1 &\tfrac{\delta^{2}}{r}\, e^{2i \theta t } 
\\
0 &1
\end{array}
\right),
\quad
V_6^{(4)}
=
\left(
\def\arraystretch{1.2}
\begin{array}{lr}
1 & 0
\\
\frac{1}{ \bar r\delta^{2}}\, e^{-2i \theta t}  & 1
\end{array}
\right).
}

Next, we switch from $N^{(4)}$ to $N^{(5)}$ via the analogue of transformation \eqref{ds-m4r-0}, now involving $h_w$  defined by \eqref{ds-h_w-def} instead of $h$:
\eee{\label{ds-wake-def5}
 N^{(5)}(\Vw t, t, k)
=
 N^{(4)}(\Vw t, t, k)  e^{-i\left[h_w(k)-\theta(\Vw,k)\right]t \sigma_3}.
 }
This transformation results in the analogue of Riemann-Hilbert problem \eqref{ds-n4t-rhp-mew}, where all relevant jumps are given by \eqref{ds-n5-jumps-mew} with $h$ replaced by $h_w$  except for the jumps along $\p D_p^\ve$ and $\p D_{\bar p}^\ve$, which are equal to
\sss{
\ddd{
V_p^{(5)} &= e^{i\left(h_w-\theta\right)t \sigma_3} \big(V_6^{(4)}\big)^{-1} V_p^{(3)} V_6^{(4)} e^{-i\left(h_w-\theta\right)t \sigma_3}, 
\\
V_{\bar p}^{(5)} &=  e^{i\left(h_w-\theta\right)t \sigma_3} V_5^{(4)} V_{\bar p}^{(3)} \big(V_5^{(4)}\big)^{-1} e^{-i\left(h_w-\theta\right)t \sigma_3}.
}
}
Finally, we perform the sixth deformation similarly to \eqref{ds-esc-n5-def-mew}, i.e. we let 
\eee{\label{ds-wake-def6}
N^{(6)}(\Vw t, t, k) = N^{(5)}(\Vw t, t, k) e^{ig_w(k) \sigma_3},
}
where the function $g_w$, which is the analogue of the function $g$ involved in \eqref{ds-esc-n5-def-mew}, is analytic in $\mathbb C\setminus (B\cup \widetilde B_w)$ and satisfies the following jump conditions:
\sss{\label{ds-wake-g_w-rhp}
\ddd{
g_w^+ + g_w^- &= -i\ln\left(\delta^2\right), && k\in B,
\\
g_w^+ + g_w^- &= -i\ln\left(\tfrac{\delta^2}{r}\right) + \omega_w, \quad && k\in L_7,
\\
g_w^+ + g_w^- &= -i\ln\left(\delta^2 \bar r\right) + \omega_w, && k\in L_8,
}
}
with the contours $B$, $L_7$, $L_8$ as in Figure \ref{ds-wake-vw-3-f}, the function $\delta(\Vw, k)$ given by \eqref{ds-esc-del-def-mew} and the real constant $\omega_w$  defined by
\eee{\label{ds-wake-omega'-def}
\omega_w
= 
 i \, 
\frac{
\displaystyle \int_{B} \frac{\ln \delta^2(\Vw, \nu)}{\gamma_w(\nu)}\, d\nu
+
\int_{\widetilde B_w^+} \frac{
\ln\left[\frac{\delta^2(\Vw, \nu) }{r(\nu)}\right]}{\gamma_w(\nu)}\, d\nu 
+
\int_{\widetilde B_w^-} \frac{\ln\left[ \delta^2(\Vw, \nu)  \, \bar r(\nu)\right]}{  \gamma_w(\nu)}\, d\nu 
}
{\displaystyle \int_{\widetilde B_w}  \frac{d\nu}{\gamma_w(\nu)}},
}
where 
\eee{
\widetilde B_w^\pm := \widetilde B_w \cap \mathbb C^\pm
}
and the function $\gamma_w$ is defined in terms of the function $\gamma$ (see \eqref{ds-gamma-def-i}) by
\eee{\label{ds-gamma'-def}
\gamma_w(k) = \left\{\begin{array}{ll} \gamma(\Vw, k), & k\in \mathbb C\setminus \overline{\mathcal R}, \\   -\gamma(\Vw, k), & k\in  \mathcal R, \end{array}\right.
}
where we recall that $\mathcal R$ is the finite region enclosed by $\widetilde B$ and $\widetilde B_w$ (see Figure \ref{ds-wake-vw-f}).
Recalling further that $\gamma$ is analytic in $\mathbb C\setminus  B\cup \widetilde B$ and changes sign as $k$ crosses $B$ and $\widetilde B$, we deduce that $\gamma_w$ has branch cuts along $B$ and $\widetilde B_w$, across which it changes sign, but is continuous as $k$ crosses $\widetilde B$.
Dividing the jumps of problem \eqref{ds-wake-g_w-rhp} by $\gamma_w$ and using Plemelj's formulae, we obtain 
\eee{\label{ds-wake-g_w}
g_w(k) = \frac{\gamma_w(k)}{2\pi}
\bigg[
\int_{B} \frac{\ln \delta^2(\Vw, \nu)}{ \gamma_w(\nu )(\nu -k)}\, d\nu 
+
\int_{L_7} \frac{
\ln\left[\frac{\delta^2(\Vw, \nu) }{r(\nu)}\right]+i\omega_w}{ \gamma_w(\nu )(\nu -k)}\, d\nu 
-
\int_{L_8} \frac{\ln\left[ \delta^2(\Vw, \nu)  \, \bar r(\nu)\right]+i\omega_w}{\gamma_w(\nu)(\nu -k)}\, d\nu 
\bigg].
}
We note that the presence of $\omega_w$ in the above formula ensures that $g_w(k) = g_{w, \infty} + O\left(\frac 1k\right)$ as $k\to \infty$ with the real constant $g_{w, \infty}$ given by
\eee{\label{ds-ginf'-def}
g_{w, \infty} 
= 
-\frac{1}{2\pi}
\bigg[
\int_{B} \frac{\ln \delta^2(\Vw, \nu)}{\gamma_w(\nu)}\,\nu d\nu
+
\int_{\widetilde B_w^+} \frac{
\ln\left[\frac{\delta^2(\Vw, \nu) }{r(\nu)}\right]+i\omega_w}{ \gamma_w(\nu)}\,\nu  d\nu 
+
\int_{\widetilde B_w^-} \frac{\ln\left[ \delta^2(\Vw, \nu)  \, \bar r(\nu)\right]+i\omega_w}{\gamma_w(\nu)}\,\nu  d\nu
\bigg]. 
}
Transformation \eqref{ds-wake-def6} results in the analogue of Riemann-Hilbert problem \eqref{ds-esc-n6-rhp-mew}, i.e.
\sss{\label{ds-wake-vw-def6}
\ddd{
N^{(6)+} &= N^{(6)-} V_B,  && k\in B,
\\
N^{(6)+} &= N^{(6)-}V_{\widetilde B_w}^{(6)},  && k\in \widetilde B_w,
\\
N^{(6)+} &= N^{(6)-}V_j^{(6)}, && k\in L_j,\ j=1, \ldots, 6,
\\
N^{(6)+}  
&=
N^{(6)-}
V_p^{(6)}, && k\in  \p D_p^\ve,
\\
N^{(6)+} 
&=
N^{(6)-} 
V_{\bar p}^{(6)}, && k\in  \p D_{\bar p}^\ve,
\\
N^{(6)}  &=\left[I +O\left(\tfrac 1k\right)\right]e^{i\left[g_{w, \infty}-G_\infty(\Vw) t\right]\sigma_3}, \quad && k \to \infty,
\label{ds-wake-vw-def6-asymp}
}
}
where the jump along $B$ is given by \eqref{ds-n1-jumps}, the jump along $\widetilde B_w$ is equal to
\eee{\label{ds-wake-vbtilde'}
V_{\widetilde B_w}^{(6)}
=
\left(
\def\arraystretch{1}
\begin{array}{lr}
0
&
-e^{i(\Omega(\Vw)  t - \omega_w)} 
\\
e^{-i(\Omega(\Vw)  t - \omega_w)} 
&
0
\end{array}
\right)
}
with the real constants $\Omega(\Vw)$ and $\omega_w$ given by \eqref{ds-Omega-exp} and \eqref{ds-wake-omega'-def} respectively, 
the jumps  along the  contours  $L_j$ of Figure \ref{ds-wake-vw-3-f} are equal to
\sss{\label{ds-wake-rem-jumps}
\ddd{
V_1^{(6)}
&=
\left(
\def\arraystretch{1}
\begin{array}{lr}
1 & \frac{\bar r\,\delta^{2} e^{-2ig_w}}{1+r\bar r}\,e^{2i h_wt}
\\
0 & 1
\end{array}
\right), 
\quad
&&V_2^{(6)}
=
\left(
\def\arraystretch{1}
\begin{array}{lr}
1 & 0
\\
\frac{r \delta^{-2} e^{2ig_w}}{1+r\bar r}\,e^{-2ih_wt} & 1
\end{array}
\right), 
\\
V_3^{(6)}
&=
\left(
\def\arraystretch{1}
\begin{array}{lr}
1 & 0
\\
r \delta^{-2} e^{2ig_w}  e^{-2ih_wt} & 1
\end{array}
\right), 
\quad
&&V_4^{(6)}
=
\left(
\def\arraystretch{1}
\begin{array}{lr}
1 & \bar r \delta^{2} e^{-2ig_w} e^{2ih_wt}
\\
0 & 1
\end{array}
\right), 
\\
V_5^{(6)}
&=
\left(
\def\arraystretch{1}
\begin{array}{lr}
1 & \frac{\delta^{2}e^{-2ig_w}}{r}\, e^{2ih_wt} 
\\
0 & 1
\end{array}
\right),
\quad 
&&V_6^{(6)}
=
\left(
\def\arraystretch{1}
\begin{array}{lr}
1 & 0
\\
\frac{\delta^{-2}e^{2ig_w}}{ \bar r}\, e^{-2ih_wt}  & 1
\end{array}
\right),
}
}
the jumps along  $\p D_p^\ve$ and $\p D_{\bar p}^\ve$ are given by
\sss{\label{ds-wake-vp-vpb-6}
\ddd{
V_p^{(6)} 
&=
\left(
\arraycolsep=1.4pt
\def\arraystretch{1.5}
\begin{array}{lr}
1 - \frac{c_p d(k)}{\bar r(k) \left(k-p\right)} \, e^{-2i\left[\theta(\Vw, k) - \theta(\Vw, p)\right]t} & -\frac{c_p \delta^2(\Vw, k) d(k) e^{-2ig_w(k)}}{k-p}\,  e^{2i\left[h_w(k) - \theta(\Vw, k) + \theta(\Vw, p)\right]t}
\\
\frac{c_p d(k) e^{2ig_w(k)}}{\bar r^2(k) \delta^2(\Vw, k) \left(k-p\right)} \, e^{-2i\left[h_w(k)+\theta(\Vw, k) - \theta(\Vw, p)\right]t} & 1 + \frac{c_p d(k)}{\bar r(k) \left(k-p\right)} \, e^{-2i\left[\theta(\Vw, k) - \theta(\Vw, p)\right]t}
\end{array}
\right),
\\
V_{\bar p}^{(6)}  
&=
\left(
\arraycolsep=-2pt
\def\arraystretch{1}
\begin{array}{lr}
1 - \frac{c_{\bar p} d(k)}{r(k) \left(k-\bar p\right)} \, e^{2i\left[\theta(\Vw, k) - \theta(\Vw, \bar p)\right]t}
& 
\frac{c_{\bar p} \delta^2(\Vw, k) d(k) e^{-2ig_w(k)}}{r^2(k) \left(k-\bar p\right)} \, e^{2i\left[h_w(k) + \theta(\Vw, k) - \theta(\Vw, \bar p)\right]t}
\\
-\frac{c_{\bar p} \delta^{-2}(\Vw, k) d(k) e^{2ig_w(k)}}{k-\bar p}\,  e^{-2i\left[h_w(k) - \theta(\Vw, k) + \theta(\Vw, \bar p)\right]t} 
& 1 + \frac{c_{\bar p} d(k)}{r(k) \left(k-\bar p\right)} \, e^{2i\left[\theta(\Vw, k) - \theta(\Vw, \bar p)\right] t}
\end{array}
\right),
}
}
and the real constants $G_\infty(\Vw)$ and $g_{w, \infty}$ are given by \eqref{ds-Ginf-mew} and \eqref{ds-ginf'-def} respectively.

The sign structure of $\text{Re}(ih_w)$  shown in Figure \ref{ds-wake-vw-3-f} indicates that the leading-order contribution to the solution of problem \eqref{ds-wake-vw-def6} in the limit $t\to \infty$ comes from the jumps along $B$, $\widetilde B_w$, $\p D_p^\ve$ and $\p D_{\bar p}^\ve$. Indeed, observe that while the jumps along the contours $L_j$, $j=1, \ldots, 6$, decay to the identity exponentially fast as $t\to \infty$, those along $B$ and $\widetilde B_w$ are purely oscillatory. Furthermore, noting that the jumps along $\p D_p^\ve$ and $\p D_{\bar p}^\ve$ will eventually be transformed to residue conditions at $p$ and $\bar p$ respectively, we see that the contributions of these jumps are also purely oscillatory as $t\to \infty$ since $\text{Re}(ih_w)(p) = \text{Re}(ih_w)(\bar p) = 0$ (recall that $p$ and $\bar p$ lie on the dashed green contour $\widetilde B$ of Figure \ref{ds-wake-vw-3-f}, along which $\text{Re}(ih_w)$ vanishes). 
This analysis motivates a decomposition of $N^{(6)}$ entirely analogous to \eqref{ds-n5-decomp-mew1-lim} and eventually leads to the asymptotic formula \eqref{ds-trap-q-recon-asymp-mew1}, i.e.
\eee{\label{ds-wake-q-asymp}
q(x, t)
=
-2i \lim_{k \to \infty} k   N_{12}^\dom(\Vw t, t, k) e^{i\left[g_{w, \infty} - G_\infty(\Vw) t  \right]} + O\big(t^{-\frac 12}\big),
\quad t\to \infty,
}
where $N^\dom$ is the solution of the dominant component of Riemann-Hilbert problem \eqref{ds-wake-vw-def6}, i.e.
\sss{\label{ds-wake-ndom-rhp}
\ddd{
N^{\dom+} &= N^{\dom-} V_B,  && k\in B,
\\
N^{\dom+} &= N^{\dom-} V_{\widetilde B_w}^{(6)}, && k\in \widetilde B_w,
\\
N^{\dom+} &= N^{\dom-} V_p^{(6)}, && k\in \p D_p^\ve,
\\
N^{\dom+} &= N^{\dom-} V_{\bar p}^{(6)}, && k\in \p D_{\bar p}^\ve,
\\
N^\dom &=  
 \left[I+O\left(\tfrac 1k\right)\right]e^{i\left[g_{w, \infty} -G_\infty(\Vw) t  \right]\sigma_3},\quad && k \to \infty.
}
}

The transformation
\eee{\label{ds-wake-ndom-mdom-trans}
M^\dom 
=
\begin{cases}
N^\dom  \big(V_p^{(6)}\big)^{-1}, &k\in D_p^{\ve },
\\
N^\dom, &k\in \mathbb C^-\setminus \big( B^-\cup \widetilde B_w^- \cup \overline{D_{p}^{\ve}}\,\big),
\\
N^\dom \big(V_{\bar p}^{(6)}\big)^{-1}, &k\in D_{\bar p}^{\ve },
\\
N^\dom, &k\in \mathbb C^+\setminus \big(B^+\cup \widetilde B_w^+ \cup  \overline{D_{\bar p}^{\ve}}\,\big),
\end{cases}
}
which is the analogue of transformation \eqref{ds-trap-mmod-def-mew1}, allows us to turn the jumps of $N^\dom$ along $D_p^\ve$ and $D_{\bar p}^\ve$ into residue conditions for $M^\dom$ at $p$ and $\bar p$. 
In particular, note that $V_p^{(6)}$ and $V_{\bar p}^{(6)}$ are meromorphic inside the disks $D_p^\ve$ and $D_{\bar p}^\ve$, their only singularities being simple poles at $p$ and $\bar p$ respectively. 
Furthermore, since $a(p) = \bar a(\bar p) = 0$, it follows from the definition \eqref{ds-r-coef-def} of the reflection coefficient $r(k)$ that $\frac{1}{\bar r(p)} = \frac{1}{r(\bar p)} = 0$. Thus, the singularity at $k=p$ is removable from all elements of the matrix $V_p^{(6)}$ except for the $12$-element. Similarly, the singularity at $k=\bar p$ is removable from all elements of the matrix $V_{\bar p}^{(6)}$ except for the $21$-element. Therefore, employing transformation \eqref{ds-wake-ndom-mdom-trans}, we convert problem \eqref{ds-wake-ndom-rhp} for $N^\dom$ to the following problem for $M^\dom$:
\sss{\label{ds-wake-mdom-rhp}
\ddd{
&M^{\dom+} = M^{\dom-} V_B,  && k\in B,
\\
&M^{\dom+} = M^{\dom-} \widetilde V_{\widetilde B_w}^{(6)},  && k\in \widetilde B_w,
\\
&M^\dom =  \left[I +O\left(\tfrac 1k\right)\right] e^{i\left[g_{w, \infty} -G_\infty(\Vw) t\right]\sigma_3},\quad  && k \to \infty,
\\
& \underset{k=p}{\text{Res}}\, M^\dom = \left(0, {\rho_p}_w \,  M_2^\dom(p)\right),
\\
& \underset{k=\bar p}{\text{Res}}\, M^\dom = \left({\rho_{\bar p}}_w\,   M_1^\dom(\bar p), 0\right),
}
}
where $M_1^\dom, M_2^\dom$ denote the two columns of $M^\dom$ and
\eee{\label{ds-rrb-til}
{\rho_p}_w 
=  c_p \delta^2(\Vw, p) d(p) \, e^{2i\left[h_w(p)t-g_w(p)\right]},
\quad
{\rho_{\bar p}}_w = c_{\bar p} \delta^{-2}(\Vw, \bar p) d(\bar p) \, e^{-2i\left[h_w(\bar p)t -g_w(\bar p)\right]},
}
which similarly to \eqref{ds-R-def-trap} can be expressed in the form
\ddd{\label{ds-R-def-wake}
{\rho_p}_w 
= 
{R_p}_w \, e^{2ih_w(p)t},
\quad
{\rho_{\bar p}}_w 
= 
-\overline{{R_p}_w} \, e^{-2ih(p)t},
\quad
{R_p}_w:=C_p \frac{\delta^2(\Vw, p)  e^{-2i g_w(p)}}{a'(p)},
}
revealing that ${\rho_{\bar p}}_w = -\overline{{\rho_p}_w}$.

Similarly to the previous sections, we solve problem \eqref{ds-wake-mdom-rhp} by employing the factorization
\eee{\label{ds-wake-factor}
M^\dom = \mathcal M^\dom W_w,
}
where $W_w$ is the solution of the continuous spectrum component of problem \eqref{ds-wake-mdom-rhp}, i.e.
\sss{\label{ds-wake-ww-rhp}
\ddd{
W_w^+ &= W_w^- V_B,  && k\in B,
\\
W_w^+ &= W_w^- V_{\widetilde B_w}^{(6)}, && k\in \widetilde B_w,
\\
W_w &=  
 \left[I+O\left(\tfrac 1k\right)\right]e^{i\left[g_{w, \infty}  -G_\infty(\Vw) t  \right]\sigma_3},\quad && k \to \infty,
}
}
and $\mathcal M^\dom$ solves the discrete spectrum component  of problem \eqref{ds-wake-mdom-rhp}, i.e.   $\mathcal M^\dom$ is analytic in $\mathbb C\setminus\left\{ p, \bar p\right\}$ and satisfies the residue conditions
\sss{\label{ds-wake-mcal-rhp}
\ddd{
\underset{k=p}{\text{Res}}\, \mathcal M_1^\dom
&=
-{W_w}_{21}(p)  {\rho_p}_w \, M_1^\dom(p),
\quad
&&\underset{k=p}{\text{Res}}\, \mathcal M_2^\dom
=
{W_w}_{11}(p)   {\rho_p}_w \, M_1^\dom(p),
\\
\underset{k=\bar p}{\text{Res}}\, \mathcal M_1^\dom
&=
{W_w}_{22}(\bar p)   {\rho_{\bar p}}_w \, M_2^\dom(\bar p),
&&\underset{k=\bar p}{\text{Res}}\, \mathcal M_2^\dom
=
-{W_w}_{12}(\bar p)   {\rho_{\bar p}}_w \, M_2^\dom(\bar p),
}
and the asymptotic condition 
\eee{\label{ds-wake-mcal-asymp}
\mathcal M^\dom
=
I+O\left(\tfrac 1k\right), \quad k\to\infty.
}
}
\indent Problem \eqref{ds-wake-ww-rhp} is entirely analogous to problem \eqref{ds-trap-nmod-rhp-mew1-w}. In fact, similarly to problem \eqref{ds-trap-nmod-rhp-mew1-w}, since the jump $V_{\widetilde B_w}^{(6)}$ is independent of $k$, the jump contour $\widetilde B_w$ in problem \eqref{ds-wake-ww-rhp} can be deformed to the straight line segment $B'$ connecting $\bar \alpha$ to $\alpha$ (see Figure \ref{ds-esc-nmodt-jumps-mew-f}). The solution of this deformed problem is then given by the analogue of formula \eqref{ds-trap-W-mew1}, i.e.
\eee{\label{ds-wake-ww-sol}
W_w
=
e^{i\left[g_{w, \infty} - G_\infty(\Vw) t  \right]\sigma_3}
\mathcal N_w^{-1}(\infty, c) \,  \mathcal N_w(k, c),
}
where $\mathcal N_w$ is defined similarly to \eqref{ds-trap-ncal-def-mew1} by
\eee{\label{ds-wake-Nwcal-def}
\mathcal N_w(k, c)
=
\frac 12
\def\arraystretch{1.5}
\left(
\begin{array}{lr}
\left[\eta(k)+\eta^{-1}(k)\right]{\@N_w}_1(k, c)
& 
i\left[\eta(k)-\eta^{-1}(k)\right]{\@N_w}_2(k, c)
\\
-i\left[\eta(k)-\eta^{-1}(k)\right]{\@N_w}_1(k, -c)
&
\left[\eta(k)+\eta^{-1}(k)\right]{\@N_w}_2(k, -c)
\end{array}
\right)
}
and 
\eee{
\mathcal N_w(\infty, c) := \lim_{k\to\infty} \mathcal N_w(k, c)
}
with $\eta$ defined by  \eqref{ds-pdef} and with  ${\@N_w}_1$ and ${\@N_w}_2$ denoting the first and second column of the vector-valued function
\eee{
\label{ds-Nw-def}
\@N_w (k, c) 
=
\left(
\frac{
 \Theta\big(-\frac{\Omega(\Vw) t}{2\pi}+\frac{\omega_w}{2\pi}+\frac{i\ln\left(\frac{\bar q_-}{iq_o}\right)}{2\pi}+\upnu(k)+c  \big)}
 {  \sqrt{\frac{iq_o }{\bar q_-}}\ \Theta\left(\upnu(k)+c  \right)},
\frac{\Theta\big(-\frac{\Omega(\Vw) t}{2\pi}+\frac{\omega_w}{2\pi}+\frac{i\ln\left(\frac{\bar q_-}{iq_o}\right)}{2\pi}-\upnu(k)+c  \big)}
{ \sqrt{\tfrac{\bar q_-}{iq_o }}  \ \Theta\left(-\upnu(k)+c  \right)}
\right),
}
where  $c=c(\Vw)$ and $\upnu(k) = \upnu(\Vw, k)$ are given by formulae \eqref{ds-vdefr0} and \eqref{ds-c-choice} evaluated at $\xi=\Vw$.

Furthermore, arguing as in Subsection \ref{ds-esc-pw1-lim-ss}, we infer that the solution of problem \eqref{ds-wake-mcal-rhp} takes the form
\eee{\label{ds-wake-mcaldom-form}
\mathcal M^\dom
=
I + \frac{\underset{k=p}{\text{Res}}\, \mathcal M^\dom}{k-p} + \frac{\underset{k=\bar p}{\text{Res}}\, \mathcal M^\dom}{k-\bar p}.
}
In addition, we compute
\sss{\label{ds-wake-mp}
\ddd{
M_1^\dom(p)
&=
\frac{-\mathcal B_w {\rho_{\bar p}}_w {W_w}_2(\bar p)+\left(1+\mathcal C_w {\rho_{\bar p}}_w\right) {W_w}_1(p)}{\mathcal B_w^2 {\rho_p}_w {\rho_{\bar p}}_w + \left(1+ \mathcal C_w {\rho_{\bar p}}_w\right)\left(1+\mathcal A_w {\rho_p}_w\right)},
\\
M_2^\dom(\bar p)
&=
 \frac{\mathcal B_w {\rho_p}_w {W_w}_1(p)+\left(1+\mathcal A_w {\rho_p}_w \right) {W_w}_2(\bar p)}{\mathcal B_w^2 {\rho_p}_w {\rho_{\bar p}}_w + \left(1 + \mathcal C_w {\rho_{\bar p}}_w\right)\left(1+\mathcal A_w {\rho_p}_w \right)},
}
}
where
\sss{\label{ds-abc-w-def}
\ddd{
\mathcal A_w &= {W_w}_{11}'(p){W_w}_{21}(p)-{W_w}_{11}(p){W_w}_{21}'(p),
\\
\mathcal B_w &= \frac{{W_w}_{21}(p){W_w}_{12}(\bar p)-{W_w}_{11}(p){W_w}_{22}(\bar p)}{p-\bar p},
\\
\mathcal C_w &= {W_w}_{22}'(\bar p){W_w}_{12}(\bar p)-{W_w}_{12}'(\bar p){W_w}_{22}(\bar p).
}
}
Combining expressions \eqref{ds-wake-mcal-rhp} and \eqref{ds-wake-mp}, we obtain $\mathcal M^\dom$ through the representation \eqref{ds-wake-mcaldom-form}.
In turn, proceeding as in Subsection \ref{ds-trap-mew1-lim-ss} we obtain the leading-order asymptotics  for the focusing NLS IVP \eqref{ds-fnls-ivp} at $\xi=\Vw$ in the form \eqref{ds-trapwake-t}-\eqref{ds-trapwake-qw-def}.
Finally, similarly to Remark \ref{ds-overall-phase-r}, we note that the overall dependence of the asymptotic solution \eqref{ds-trapwake-t} on $g_{w, \infty}$ and $G_\infty(\Vw)$   is expressed by a factor of $e^{2i\left[g_{w, \infty}  -  G_\infty(\Vw) t  \right]}$.

The proof of Theorem \ref{ds-twr-t} for the leading-order asymptotics in the trap/wake regime $p\in D_2^-$ is complete.

\begin{remark}[\b{Soliton versus soliton wake}] 
\label{ds-phase-shift-r}
We recall that the soliton arising at $\xi=\Vd$ induces a phase shift in the asymptotics for $\Vd<\xi<0$. This is because in the transition from $\Vd^-$ to $\Vd^+$ the quantity $\text{Re}(ih)$ switches sign from negative (Figure \ref{ds-trap-def4-mew-f}) to positive  (Figure \ref{ds-trap-wake-def4-mew-f}) along $\p D_p^\ve$. Hence, in the latter case the additional transformation \eqref{ds-esc-n5t-def-mew} must be employed in order to convert growth into decay in the jump along $\p D_p^\ve$.
On the other hand, the soliton wake arising at $\xi=\Vw$ does not cause a phase shift in the asymptotics for $\Vw<\xi<0$. To see this, recall that the wake is created at $\xi=\Vw$ because at that value of $\xi$ the contour $\widetilde B$, along which $\text{Re}(ih)$ vanishes for all $\xi$, crosses the pole $p$. Hence, $\text{Re}(ih)(\Vw, p) = \text{Re}(ih_w)(p)=0$ and the jump along $\p D_p^\ve$ contributes to the leading-order asymptotics. However,  the quantity $\text{Re}(ih)$ is positive along $\p D_p^\ve$ \textit{both} right before \textit{and} right after the crossing with $\widetilde B$ (see Figure \ref{ds-trap-wake-def4-mew-f}). Consequently, the jump along $\p D_p^\ve$ remains bounded in the transition from $\Vw^-$ to $\Vw^+$ (recall that transformation \eqref{ds-esc-n5t-def-mew} has already been applied for $\xi>\Vd$) and hence no further transformations are required for $\Vw<\xi<0$.
\end{remark}

\begin{remark}[\b{$h$ versus $h_w$}]
Although the case $\xi=\Vw$ was analyzed by switching from the phase function $h$ (used for all $\xi\neq \Vw$) to the phase $h_w$ defined by \eqref{ds-h_w-def}, it would have still been possible to obtain the asymptotic result \eqref{ds-trapwake-t}-\eqref{ds-trapwake-qw-def} by adhering to $h$. However, the fact that for $\xi=\Vw$ the poles $p$ and $\bar p$ lie along the branch cut $\widetilde B$ of $h$ would have made the analysis significantly more complicated. In particular, even the very first step of the analysis, namely, transformation \eqref{ds-nrh} which converts the residue conditions at $p$ and $\bar p$ to jumps along the circles $\p D_p^\ve$ and $\p D_{\bar p}^\ve$, would have resulted in additional jump conditions inside the disks $D_p^\ve$ and $D_{\bar p}^\ve$ due to fact that for $\xi=\Vw$ these disks are crossed by the branch cut $\widetilde B$ of $h$. Switching from $h$ to $h_w$, whose branch cut $\widetilde B_w$ does not intersect with $\overline{D_p^\ve}$ and $\overline{D_{\bar p}^\ve}$, significantly simplifies the analysis for the case $\xi = \Vw$. Of course, for all $\xi\neq \Vw$ the poles are away from $\widetilde B$ and hence the switch from $h$ to $h_w$ is  not necessary.
\end{remark}

\subsection{The  transmission/wake regime}
\label{ds-esc-wake-ss}

Recall that for $p\in D_3$ we have $\Vs<\Vo$ and, furthermore, the integral equation \eqref{ds-fnlsd-int-eq-xistar-0} possesses a unique solution $\Vw$ in the interval $(\Vo, 0)$, which corresponds to the crossing of  the pole $p$ by the branch cut $\widetilde B$.

\begin{figure}[t!]
\begin{center}
\includegraphics[scale=1]{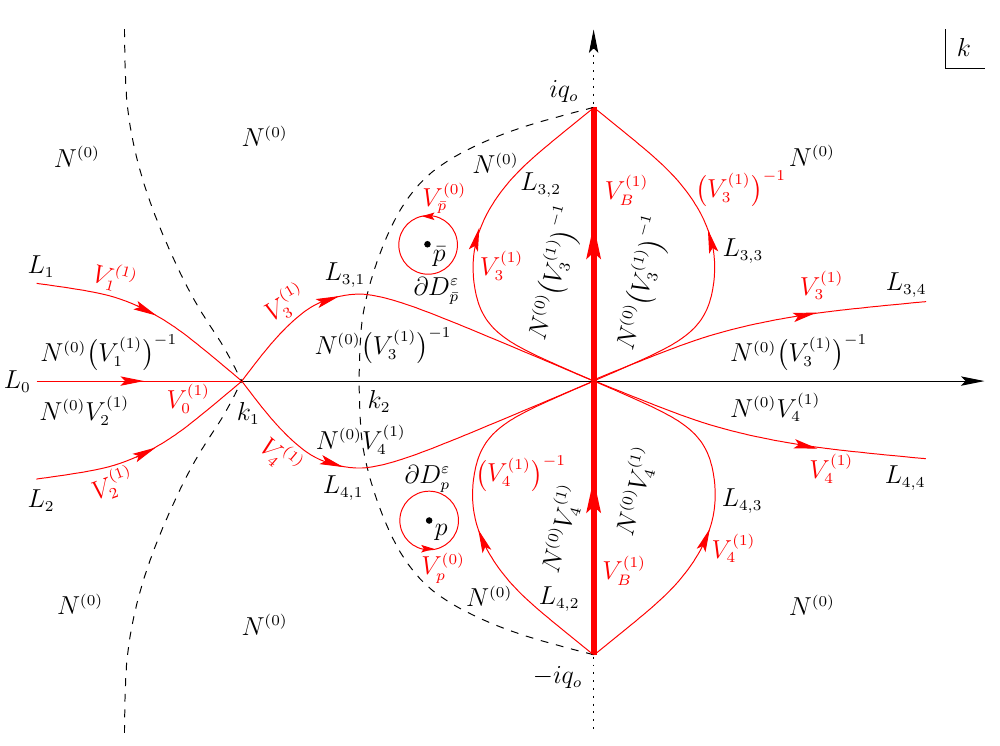}
\caption{Plane wave in the transmission/wake regime ($p\in D_3$) for $\Vs<\xi<\Vo$: the first stage of the first deformation. The jumps along $\p D_p^\ve$ and $\p D_{\bar p}^\ve$ are unaffected.}
\label{ds-esc-scen2-pw2-f}
\end{center}
\end{figure}

\begin{figure}[t!]
\begin{center}
\includegraphics[scale=1]{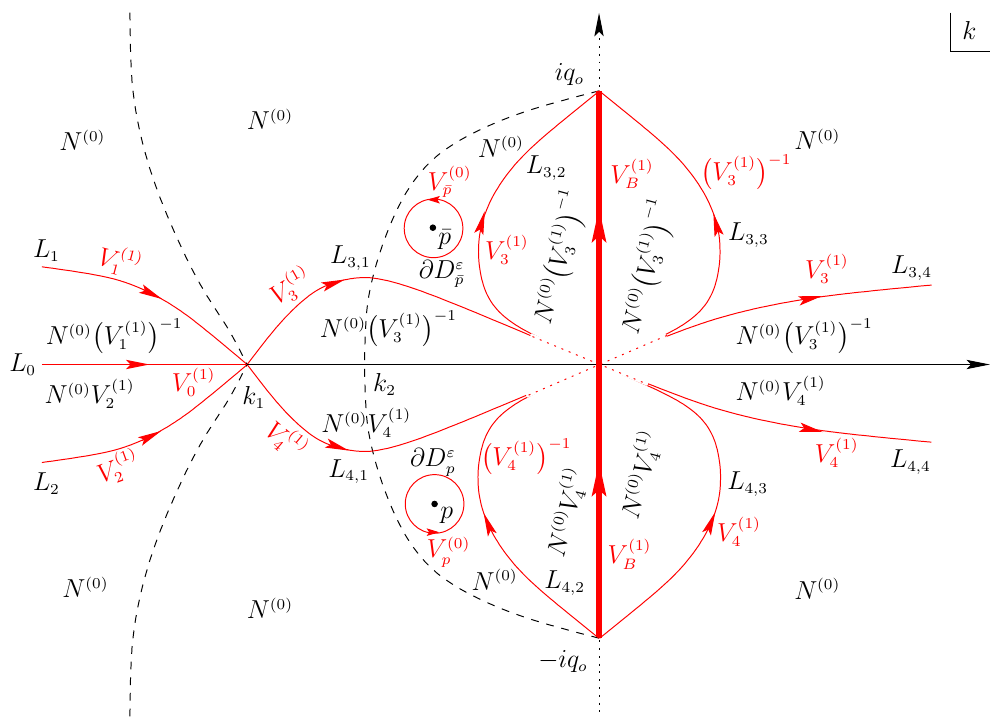}
\caption{Plane wave in the transmission/wake regime ($p\in D_3$) for $\Vs<\xi<\Vo$: the second stage of the first deformation. The overlapping portions of the contours $L_{3,1}$ and $L_{3,2}$, as well as of the contours $L_{4,1}$ and $L_{4,2}$, have been removed since by definition $N^{(1)}$ does not have a jump there. Hence, the contours $L_{3,1}$, $L_{3,2}$, $L_{4,1}$, $L_{4,2}$ have been lifted away from the origin. The jumps along $\p D_p^\ve$ and $\p D_{\bar p}^\ve$ remain unchanged.}
\label{ds-esc-scen2-2-pw2-f}
\end{center}
\end{figure}

\begin{figure}[t!]
\begin{center}
\includegraphics[scale=1.1]{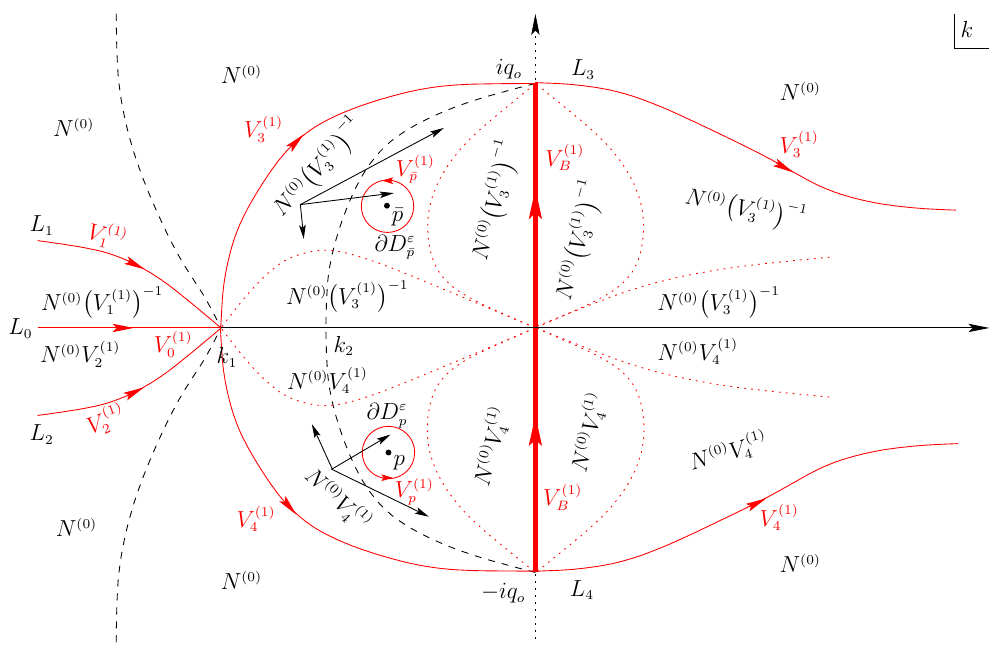}
\caption{Plane wave   in the transmission/wake regime ($p\in D_3$) for $\Vs<\xi<\Vo$: the third stage of the first deformation. The jumps along $\p D_p^\ve$ and $\p D_{\bar p}^\ve$ have now changed.}
\label{ds-esc-scen2-3-pw2-f}
\end{center}
\end{figure}

\begin{figure}[t!]
\begin{center}
\includegraphics[scale=1.1]{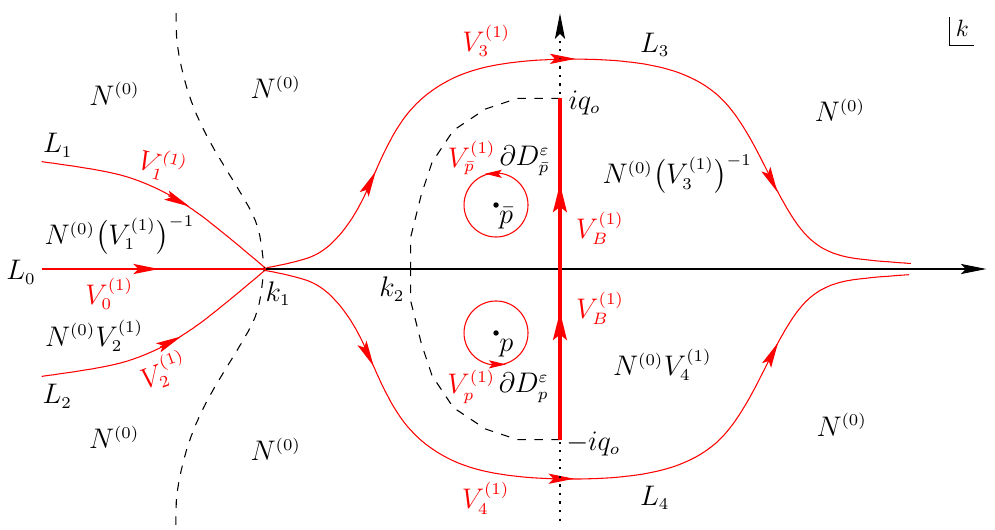}
\caption{Plane wave  in the transmission/wake regime ($p\in D_3$) for $\Vs<\xi<\Vo$: the fourth and final stage of the first deformation. 
The jump contours $L_3$ and $L_4$ have been lifted away from the branch points $\pm iq_o$ similarly to \cite{bm2017}.
}
\label{ds-esc-def1b-pw2-f}
\end{center}
\end{figure}

For $\xi<\Vo$, performing the deformations of Subsection \ref{ds-esc-pw1-ss} of the transmission regime, we obtain Riemann-Hilbert problem \eqref{ds-esc-n4-rhp-pw1}. Then, like in the transmission regime, since $\text{Re}(i\theta)(\xi, p)<0$ throughout the interval $(-\infty, \Vs)$ the leading-order asymptotics is described by the plane wave \eqref{ds-qsol-pw-t}. 
At $\xi=\Vs$, we have $\text{Re}(i\theta)(\Vs, p)=0$. Hence, Riemann-Hilbert problem \eqref{ds-esc-n4-rhp-pw1} can be analyzed like in Subsection \ref{ds-esc-pw1-lim-ss} to yield the asymptotics in the form \eqref{ds-esc-q-sol-lim-pw1-t} as the soliton \eqref{ds-qstheta-def} on top of the plane wave \eqref{ds-qpw-def} evaluated at $\Vs$.
For $\Vs<\xi<\Vo$, the fact that $p\in D_3$ means that   $p_\re>k_2$, as opposed to $p_\re<k_1$ of the case $p\in D_1$. That is, if $p\in D_3$ then $p$ is crossed by the portion of the curve $\text{Re}(i\theta)=0$ that goes through $\pm iq_o$ and $k_2$ (as opposed to the one going through $k_1$). Hence, after the crossing $p$ lies inside the finite region enclosed by the curve $\text{Re}(i\theta)=0$ and the branch cut $B$ (see Figure \ref{ds-esc-scen2-pw2-f}) as opposed to the unbounded region on the left of $\text{Re}(i\theta)=0$ and $k_1$. 
For this reason, the analysis of Subsection \ref{ds-esc-pw2-ss} for $p\in D_1$ now needs to be  modified as described below.
\vskip 3mm
\noindent
\textbf{First deformation} (Figures \ref{ds-esc-scen2-pw2-f}-\ref{ds-esc-def1b-pw2-f})\textbf{.}
Choose the contours $L_{4,1}, L_{4,2}$ and $L_{3, 1}, L_{3, 2}$ so that they do not intersect with the disks $\overline{D_p^\ve}$ and $ \overline{D_{\bar p}^\ve}$. Then, as shown in Figures \ref{ds-esc-scen2-pw2-f}-\ref{ds-esc-scen2-3-pw2-f},  in order to deform $L_4$ outside the bounded region of ``wrong'' (i.e. positive) sign, we eventually need to set $N^{(1)} = N^{(0)} V_4^{(1)}$ inside the disk $D_p^\ve$, as opposed to the regime $p\in D_1$ in which the fact that  $p_\re<k_1$ allows us to always have $N^{(1)} = N^{(0)}$ in $D_p^\ve$ (see Figures \ref{ds-esc-def1a-pw1-f}-\ref{ds-esc-def1b-pw1-f}). The situation is analogous for the disk $D_{\bar p}^\ve$. 
Therefore, for $p\in D_3$ the  Riemann-Hilbert problem for $N^{(1)}$ reads
\sss{\label{ds-esc-n1-rhp-pw2}
\ddd{
N^{(1)+} &= N^{(1)-} V_B^{(1)}, \quad && k\in  B,
\\
N^{(1)+} &= N^{(1)-} V_j^{(1)},  && k\in L_j,\ j=1,2,3,4,
\\
N^{(1)+} &= N^{(1)-} V_p^{(1)},  && k\in \p D_p^\ve,
\\
N^{(1)+} &= N^{(1)-} V_{\bar p}^{(1)},  && k\in \p D_{\bar p}^\ve,
\\
N^{(1)} &= I +O\left(\tfrac 1k\right), && k \to \infty,
}
}
where the matrices $V_B^{(1)}$ and $V_j^{(1)}$, $j=1,2,3,4,$ are defined as in \eqref{ds-n1-jumps} but the matrices $V_p^{(1)}$ and $V_{\bar p}^{(1)}$ are  given instead by
\ddd{\label{ds-esc-vp1-vpb1-def-pw2}
V_p^{(1)} 
&=
\big(V_4^{(1)}\big)^{-1} V_p^{(0)} V_4^{(1)}
=
\left(
\def\arraystretch{1}
\begin{array}{lr}
1 & -\dfrac{c_p \, d(k)}{k-p}\,  e^{2i\theta(\xi, p)t}
\\
0 & 1
\end{array}
\right),
\\
V_{\bar p}^{(1)} 
&=
V_3^{(1)}  V_{\bar p}^{(0)}  \big(V_3^{(1)}\big)^{-1}
=
\left(
\def\arraystretch{1}
\begin{array}{lr}
1 &0
\\
-\dfrac{c_{\bar p} \, d(k)}{k-\bar p}\,  e^{-2i\theta(\xi, \bar p)t} & 1
\end{array}
\right).
}

\vskip 3mm
\noindent
\textbf{Second deformation.}
This deformation is identical to \eqref{ds-n2-def-pw1} of Subsection \ref{ds-esc-pw1-ss} and results in the Riemann-Hilbert problem
\sss{\label{ds-esc-n2-rhp-pw2}
\ddd{
N^{(2)+} &= N^{(2)-} V_B^{(2)}, \quad && k\in  B,
\\
N^{(2)+} &= N^{(2)-} V_j^{(2)},  && k\in L_j,\ j=1,2,3,4,
\\
N^{(2)+} &= N^{(2)-} V_p^{(2)},  && k\in \p D_p^\ve,
\\
N^{(2)+} &= N^{(2)-} V_{\bar p}^{(2)},  && k\in \p D_{\bar p}^\ve,
\\
N^{(2)} &= I +O\left(\tfrac 1k\right), && k \to \infty,
}
}
with $V_B^{(2)}$ and $V_j^{(2)}$, $j=1,2,3,4$, as in \eqref{ds-n2-jumps} but with
\sss{\label{ds-esc-vp2-vpb2-def-pw2}
\ddd{
V_p^{(2)} 
&=
\delta^{\sigma_3} V_p^{(1)} \delta^{-\sigma_3}
=
\left(
\def\arraystretch{1}
\begin{array}{lr}
1 & -\dfrac{c_p \, \delta^2(\xi, k) d(k)}{k-p}\,  e^{2i\theta(\xi, p)t}
\\
0 & 1
\end{array}
\right),
\\
V_{\bar p}^{(2)} 
&=
\delta^{\sigma_3} V_{\bar p}^{(0)} \delta^{-\sigma_3}
=
\left(
\def\arraystretch{1}
\begin{array}{lr}
1 &0
\\
-\dfrac{c_{\bar p}\, \delta^{-2}(\xi, k)d(k)}{k-\bar p}\,  e^{-2i\theta(\xi, \bar p)t} & 1
\end{array}
\right).
}
}

\begin{figure}[t]
\begin{center}
\includegraphics[scale=1]{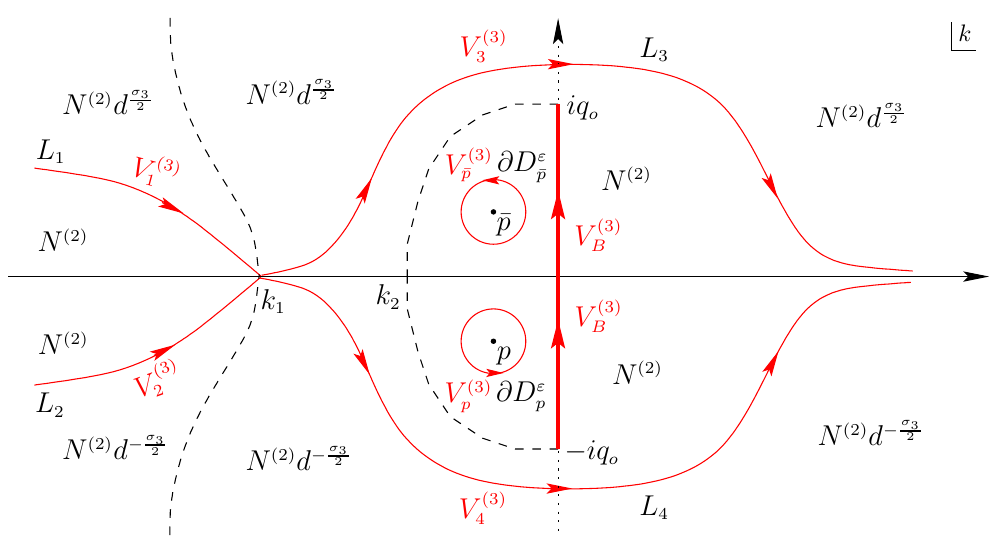}
\caption{Plane wave   in the transmission/wake regime ($p\in D_3$) for $\Vs<\xi<\Vo$: the third deformation.}
\label{ds-esc-def3-pw2-f}
\end{center}
\end{figure}

\vskip 3mm
\noindent
\textbf{Third deformation} (Figure \ref{ds-esc-def3-pw2-f})\textbf{.}
This deformation is different than the one of  Figure \ref{ds-esc-def3-pw1-f} in that the disks $\overline{D_p^\ve}$ and $\overline{D_{\bar p}^\ve}$ now lie between the contours $L_4$ and $L_3$ and hence inside these disks we have $N^{(3)} = N^{(2)}$. Thus, the jumps \eqref{ds-esc-vp2-vpb2-def-pw2} now remain invariant while the remaining jumps of problem \eqref{ds-esc-n2-rhp-pw2} are modified as in the second deformation \eqref{ds-n2-def-pw1} but  with $d^\frac 12$ now holding the role of $\delta$. Eventually, we find that $N^{(3)}$ satisfies the same problem as in Subsection \ref{ds-esc-pw1-ss}.

\vskip 3mm
\noindent
\textbf{Fourth deformation.}
This is identical to \eqref{ds-esc-n4-def-pw1}, leading to Riemann-Hilbert problem \eqref{ds-esc-n4-rhp-pw1}.

\vskip 3mm

In summary, in the range $\Vs<\xi<\Vo$ the original Riemann-Hilbert problem \eqref{ds-n-rhp-intro} can be deformed to Riemann-Hilbert problem \eqref{ds-esc-n4-rhp-pw1} both for $p\in D_1$ and for $p\in D_3$. Thus, performing the analysis of Subsection \ref{ds-esc-pw2-ss}, we obtain once again the phase-shifted plane wave \eqref{ds-esc-qasym-pw2-t}.

For $\Vo<\xi<0$, the phase function changes from $\theta$ to $h$. Recall that for $p\in D_1$ the asymptotics is given by the phase-shifted modulated elliptic wave \eqref{ds-esc-qasym-mew-t} throughout the range $(\Vo, 0)$. Now, however, $p\in D_3$ and hence, as noted in the relevant discussion of Subsection \ref{ds-esc-mew-ss}, there is a value $\Vw\in (\Vo, 0)$ for which the contour $\widetilde B$ crosses the pole $p$ en route to collapsing onto $B$. Indeed, this value is the unique solution of equation \eqref{ds-fnlsd-int-eq-xistar-0} in the interval $(\Vo, 0)$. 
Thus, similarly to the case $p\in D_2^-$, a soliton wake arises from the dominant component of Riemann-Hilbert problem \eqref{ds-esc-n6t-rhp-mew} at $\xi=\Vw$. In fact, the dominant problem is precisely that of  Subsection \ref{ds-trap-wake-ss}.
Therefore, at $\xi=\Vw$ the leading-order asymptotics is characterized by  \eqref{ds-trapwake-t} as the soliton wake \eqref{ds-trapwake-qw-def} on top of the modulated elliptic wave $q_{\text{mew}, w}(t)$. 
Finally, for $\Vo < \xi < \Vw$ and $\Vw<\xi<0$ the jumps along $\p D_p^\ve$ and $\p D_{\bar p}^\ve$ are not part of the dominant problem and, like in the transmission regime, the leading-order asymptotics is given by the phase-shifted modulated elliptic wave \eqref{ds-esc-qasym-mew-t}. 
We note that, like in the trap/wake regime (see Remark \ref{ds-phase-shift-r}), the soliton at $\Vs$  induces a phase shift  of $4\arg\left[p+\lambda(p)\right]$ in the asymptotics but no phase shift is generated by the soliton wake at $\Vw$.

The proof of Theorem \ref{ds-ewr-t} for the leading-order asymptotics in the transmission/wake regime $p\in D_3$ is complete.

\section{Conclusions}
\label{s:conclusions}

In summary, we have characterized the interactions between solitons and localized disturbances in focusing media governed by the NLS equation.  We reiterate that the main points of novelty of the results are on one hand the existence of a trapping regime, in which the velocity of the soliton differs from that of the case without radiation, and on the other hand the existence of mixed transmission/wake and trap/wake regimes, in which a single discrete eigenvalue gives rise to $O(1)$ contributions at two different velocities in the long-time asymptotics.  

The applicability of the deformations used in the present work requires that one can extend the reflection coefficient into the complex plane.  As in \cite{bm2017}, this can be done as long as the potential decays to background sufficiently rapidly (according to \eqref{ds-fnls-ic-space}) so as to ensure the existence of a Bargmann strip of analyticity.

Regarding the wake formulae, we note that the expression for the solution at the wake coincides formally with that of a soliton on top of an elliptic background.  On the other hand, crossing the wake does not result in an additional phase shift for the solution, while crossing a soliton does (cf. Remark \ref{ds-shifts-r}).
We also emphasize that the asymptotic expressions for the solution are not uniform with respect to $\xi$, as one can see by taking the limit of the various expressions as  $\xi \to \Vs$, $\xi\to\Vd$ and $\xi\to \Vw$.

We should mention that there exists previous literature on the interaction of solitons and the radiation on a nonzero background for integrable systems by the Riemann-Hilbert approach. See, for example, \cite{krugerteschl} for the case of the Toda lattice and \cite{andreievteschl} for the case of the KdV equation.  To the best of our knowledge, however, none of the previous cases studied in the literature give rise to the phenomena presented here, in which a localized disturbance results in a change of the soliton velocity and/or the production of a wake.

We believe that the asymptotic formulae giving rise to a soliton on top of an elliptic wave should be a limiting reduction of 3-phase solutions of the focusing NLS equation.  In this regard, we should mention that, in \cite{BBEIM}, the authors consider elliptic solutions of the focusing NLS equation as well as solutions corresponding to a nonlinear superposition of a soliton and hyperelliptic solutions.  (We also note that reductions of 2-phase solutions of the focusing NLS equations, which give rise to solitons on a constant background, were studied in \cite{bertolagiavedoni}.)  The authors of \cite{BBEIM} show that, in the genus-1 reduction, their solution reduces to the cnoidal wave solution of focusing NLS, namely (4.5.1) in \cite{BBEIM}.  Importantly, however, the dn solution is not the most general periodic solution of the focusing NLS equation (e.g., see \cite{kamchatnov, k1990, deconincksegal}).  More precisely, the dn solution is just one of the special cases corresponding to a trivial phase.  It is also the case that the modulate elliptic waves arising in the long-time asymptotics are not simply dn solutions.  Therefore, it is doubtful that the formulae for the soliton on top of an elliptic wave in our work reduce to those in \cite{BBEIM}. 

The asymptotic expressions in our work remain valid in the limit $\xi \to 0$.  As shown in \cite{blm2016}, in this limit $m\to1$ and the solution reduces to the well-known sech-shaped soliton solution of the focusing NLS equation.  However, some details of the derivation are different in this case and hence we omit the details for brevity.  It is also the case that the asymptotic formulae remain valid in the limit  $p_\re\to 0$, i.e. when the discrete eigenvalue lies on the imaginary axis.  In this case, the velocity $\Vs$ is zero (cf. \eqref{ds-xistar-xip-def}).  Thus, one does not see the soliton in the plane wave region.  Moreover, when $\xi=0$, the points $\alpha$,  $\bar \alpha$ and $k_o$ in the definition \eqref{ds-h-abel-i} of $h$ collapse to $iq_o$, $-iq_o$ and $0$, respectively, and hence $h(0, k) = -2\left(k^2+q_o^2\right)$.  Thus, for $p_\re=0$ we have $\text{Im}(h)(0, p)=0$, i.e. for $p_\re=0$ the imaginary parts of  both $\theta$ and $h$ are  zero at $\xi=0$.  Therefore, in this case the velocity $\Vd$ of the trapped soliton coincides with the unperturbed velocity $\Vs$.

Another interesting special case is that of a pole $p$ lying on the branch cut $i[-q_o, q_o]$, which gives rise to an Akhmediev breather.  There are four different considerations: 
(i)~Akhmediev breathers are periodic in $x$, and are therefore outside the class of initial conditions for which the inverse scattering transform formalism of \cite{bk2014,bm2017} applies (namely, constant nonzero boundary conditions);  
(ii)~Neglecting the direct problem in the inverse scattering transform, one could still consider the Riemann-Hilbert problem with a pole in the branch cut and ask what happens then. Nonetheless, even the simple formulation of a Riemann-Hilbert problem with a pole in the branch cut requires some care;  
(iii)~Akhmediev breathers are homoclinic in $t$ (i.e. they decay to the background as $t\to\pm\infty$) and hence they do not appear in the long-time asymptotics.  However, one can still see the result of their presence in the phase difference of the background before/after the breather;  
(iv)~Indeed, one can look at the case of a pole $p$ on the branch cut as a limit of the case of $p \in D_3$. Therefore, the case of an Akhmediev breather can be viewed as a limit of the trap/wake scenario.  Then, the analysis shows that one will see a wake located at $\xi = 0$ (for the same reasons as those outlined in the previous paragraphs).

As usual, the inverse scattering transform is formulated under the assumption of existence and uniqueness of solutions. The well-posedness of IVP  \eqref{ds-fnls-ivp-u} for short times with initial conditions in suitable Sobolev spaces was recently proved in \cite{m2017}.  The question of global well-posedness for initial conditions in the space \eqref{ds-fnls-ic-space} is still open.  In general, the issue of existence and uniqueness of solutions of the Riemann-Hilbert problems associated with the inverse scattering transform for integrable nonlinear partial differential equations is a nontrivial one \cite{bdt1988,Zhou,trogdonolver}.  Therefore, since this issue is peripheral to the main thrust of this work, we do not consider it here.
However, we note that the asymptotic results provide an explicit solution of the Riemann-Hilbert problem (and hence of the IVP) modulo the solution of the error Riemann-Hilbert problem, which is a small-norm problem and, therefore, is expected to have a unique solution.  On the other hand, whether the corresponding solution of the NLS equation belongs to the same function space \eqref{ds-fnls-ic-space} remains an interesting open question.  We note that proving well-posedness of an IVP in a given function space through the inverse scattering transform is in general a nontrivial problem.

\vspace*{3mm}
\noindent
\textbf{Acknowledgements.}
We would like to thank M. Ablowitz, G. El and P. Miller for many insightful discussions. This work was partially supported by the National Science Foundation under grant number DMS-1614623. Finally, we are indebted to the anonymous referees whose constructive criticism led to the improvement of our paper.

%
%
%
%
\setcounter{section}{0}
\renewcommand{\thesection}{\Alph{section}}
\section{Appendix: Soliton Solutions}

The pure one-soliton solution of the focusing NLS IVP \eqref{ds-fnls-ivp} can be derived by solving Riemann-Hilbert problem \eqref{ds-rhp-cpam} together with the residue conditions~\eqref{ds-res-cond} in the case of a zero reflection coefficient, i.e. by solving the problem
\sss{\label{ds-rhp-soliton1}
    \ddd{
        &M^+(x, t, k) = M^-(x, t, k)  V_{1, 0}(k),    &&k\in \mathbb R,
        \\
        &M^+(x, t, k)  = M^-(x, t, k)  V_{2, 0}(k), &&k\in B^+,
        \\
        &M^+(x, t, k)  = M^-(x, t, k)  V_{3, 0}(k), &&k\in B^-,
        \\
        &M(x, t, k)  = I +O\left(\tfrac 1k\right), \quad && k \to \infty,
        \\
        &\underset{k=p}{\text{Res}} \, M(x, t, k) 
        =
        \left(0, c_p  \, e^{2i\vartheta(x, t, p)} M_1(x, t, p)\right), 
       && x, t \in\mathbb R,
        \\
        &\underset{k=\bar p}{\text{Res}} \, M(x, t, k) 
        =
        \left(c_{\bar p}  \, e^{-2i\vartheta(x, t, \bar p)} M_2(x, t, \bar p), 0 \right),
        \quad && x, t \in\mathbb R,
    }
}
where the relevant jump matrices are given by 
\eee{
        V_{1, 0}(k)
        =
        \begin{pmatrix}
            \dfrac{1}{d(k)}
            &
            0 
            \\
            0 
            &
            d(k)
        \end{pmatrix},\ \
        V_{2, 0}(k)
        =
        \begin{pmatrix}
            0
            &
            \dfrac{2\lambda(k)}{i\bar q_-} 
            \\
            \dfrac{\bar q_-}{2i \lambda(k)}
            &
            0
        \end{pmatrix},\ \
        V_{3, 0}(k)
        =
        \begin{pmatrix}
            0
            & 
            \dfrac{q_-}{2i\lambda(k)}
            \\
            \dfrac{2\lambda(k)}{iq_-} 
            &
            0
        \end{pmatrix}
} 
with the functions $\lambda$ and $d$ defined by \eqref{ds-lambda-def} and \eqref{ds-d-def}, and where the phase function $\vartheta$ is defined as
\eee{\label{ds-vtheta-def}
\vartheta(x, t, k) :=  \lambda(k)  \left(x-2k t\right)
}
with $M_1$ and $M_2$ denoting the first and second column of $M$ respectively.

The jump $V_{1, 0}$ of problem \eqref{ds-rhp-soliton1} along $\mathbb R$ can be eliminated via the transformation
\eee{\label{ds-mm1}
M^{(1)}
= 
\left\{
\def\arraystretch{2.4}
\begin{array}{ll}
M 
\left(
\def\arraystretch{1}
\begin{array}{cc}
            d^{\frac 12}
            &
            0 
            \\
            0 
            &
            d^{-\frac 12}
\end{array}
\right) 
 & k\in \mathbb C^+\setminus B^+,
\\
M 
\left(
\def\arraystretch{1}
\begin{array}{cc}
            d^{-\frac 12}
            &
            0 
            \\
            0 
            &
            d^{\frac 12}
\end{array}
\right) & k\in \mathbb C^-\setminus B^-,
\end{array}
\right.
}
which implies the following Riemann-Hilbert problem for $M^{(1)}$:
\sss{\label{e:soliton-M1}
    \ddd{
        &M^{(1)+} = M^{(1)-} V_B,    &&k\in B,
        \\
        &M^{(1)} = I +O\left(\tfrac 1k\right), \quad && k \to \infty,
        \\
        &\underset{k=p}{\text{Res}} \, M^{(1)}(x, t, k) 
        =
        \left(0, \varrho_{p}(x, t) M_1^{(1)}(x, t, p)\right),
        &&  x, t \in\mathbb R, \label{e:soliton-M1residue-a}
        \\
        &\underset{k=\bar p}{\text{Res}} \, M^{(1)}(x, t, k) 
        =
        \left(-\overline{\varrho_p(x, t)} M_2^{(1)}(x, t, \bar p), 0 \right),
        \quad
        && x, t \in\mathbb R, \label{e:soliton-M1residue-b}
    }
}
where the jump matrix $V_B$ is defined in \eqref{ds-n1-jumps} and, recalling  the definitions \eqref{ds-cp-cpb-def} and the symmetry \eqref{ds-Cp-Cpb-sym}, we have introduced the quantity
\ddd{\label{e:soliton-rhop}
    \varrho_{p}(x, t) := \frac{C_{p}}{a'(p)} e^{2i\vartheta(x, t, p)}.
}

Problem \eqref{e:soliton-M1} can be solved by using the factorization
\eee{\label{e:soliton-M1M2W}
M^{(1)} = M^{(2)} W,
}
where $W$ is the solution of the continuous spectrum component of problem \eqref{e:soliton-M1}, i.e.
\sss{
    \ddd{
        W^+ &= W^- V_B,    &&k\in B,
        \\
        W &= I +O\left(\tfrac 1k\right), \quad && k \to \infty,
    }
}
and, similarly to problem \eqref{ds-esc-w-rhp-pw1}, is given by the explicit formula
\eee{\label{e:soliton-W}
    W
    =
    \frac 12 
    \begin{pmatrix}
        \Lambda(k)+\Lambda^{-1}(k)
        &
        -\frac{q_o}{\bar q_-}\left[\Lambda(k)-\Lambda^{-1}(k)\right]
        \\
        -\frac{q_o}{q_-} \left[ \Lambda(k)- \Lambda^{-1}(k)\right]
        &
        \Lambda(k)+\Lambda^{-1}(k) 
    \end{pmatrix}
}
with $\Lambda(k)$ defined by \eqref{ds-Lam-def}.

In turn, $M^{(2)}$ is the solution of the discrete spectrum component of problem \eqref{e:soliton-M1}, i.e. $M^{(2)}$ is analytic for all $k\in \mathbb C$ apart from the poles $p$ and $\bar p$, where it satisfies the residue conditions
\sss{\label{e:soliton-M2residue}
    \ddd{
        &\underset{k=p}{\text{Res}}\, M_1^{(2)}
        =
        -W_{21}(p) c_p \, d(p)  M_1^{(1)}(p), \quad &&
        \underset{k=p}{\text{Res}}\, M_2^{(2)}
        =
        W_{11}(p)   c_p \, d(p)  M_1^{(1)}(p),
        \\
        &\underset{k=\bar p}{\text{Res}}\, M_1^{(2)}
        =
        W_{22}(\bar p)   c_{\bar p} \, d(\bar p)  M_2^{(1)}(\bar p),
        &&
        \underset{k=\bar p}{\text{Res}}\, M_2^{(2)}
        =
        -W_{12}(\bar p)  c_{\bar p} \, d(\bar p)  M_2^{(1)}(\bar p).
    }
}
Furthermore, $M^{(2)}$ satisfies the asymptotic condition
\eee{
M^{(2)} = I +O\left(\frac 1k\right), \quad  k \to \infty.
}
Therefore, similarly to Subsection \ref{ds-esc-pw1-lim-ss}, we infer that $M^{(2)}$ is of the form
\eee{\label{soliton-M3}
    M^{(2)}
    =
    I + \frac{\underset{k=p}{\text{Res}}\, M^{(2)}}{k-p} + \frac{\underset{k=\bar p}{\text{Res}}\, M^{(2)}}{k-\bar p}.
}
Expressions \eqref{e:soliton-M1M2W}, \eqref{e:soliton-W}, \eqref{e:soliton-M2residue} and \eqref{soliton-M3} yield
\sss{
\ddd{
&M^\dom_1(p)
        =
        \frac{-\mathcal B \overline{\varrho_p} W_2(\bar p) + \left(1+\overline{\mathcal A \varrho_p}\right) W_1(p)}{\left(1+\overline{\mathcal A \varrho_p}\right)\left(1+\mathcal A \varrho_p\right)-\mathcal B^2 |\varrho_p|^2},
\\
&M^\dom_2(\bar p)
        =
        \frac{\left(1+\mathcal A\varrho_p\right)W_2(\bar p) - \mathcal B \varrho_p W_1(p)}{\left(1+\overline{\mathcal A \varrho_p}\right)\left(1+\mathcal A \varrho_p\right)-\mathcal B^2 |\varrho_p|^2}
}
}
with the constants $\mathcal A$, $\mathcal B$ given by \eqref{ds-ab-sol-def}. Hence, in view of \eqref{soliton-M3} and \eqref{e:soliton-M2residue}, the function $M^{(2)}$ has been determined. 

Then, reverting the transformations  \eqref{e:soliton-M1M2W} and \eqref{ds-mm1} we obtain the solution $M$ of problem \eqref{ds-rhp-soliton1} which, combined with the reconstruction formula \eqref{ds-q-recon-n}, yields the pure one-soliton solution of IVP~\eqref{ds-fnls-ivp} for the focusing NLS equation  in the form
\ddd{\label{e:soliton}
    q(x, t) 
    &=
    q_- 
    -
    \frac i2 \left\{\left[1+\overline{\mathcal A \varrho_p(x, t)}\right] \left[1+\mathcal A \varrho_p(x, t)\right]-\mathcal B^2 \left|\varrho_p(x, t)\right|^2\right\}^{-1}
 \nn\\
 &\ 
   \cdot \bigg\{\!
    \left[1+\overline{\mathcal A \varrho_p(x, t)}\right]\varrho_p(x, t)
    \left[\Lambda(p)+  \Lambda^{-1}(p)\right]^2
    +
    \left[1+\mathcal A\varrho_p(x, t)\right]\overline{\varrho_p}(x, t)\,
    \frac{q_-}{\bar q_-} \left[\overline{\Lambda(p) - \Lambda^{-1}(p)}\right]^2  
    \nn\\
    &\hskip 1cm
    -
    2\mathcal B \left|\varrho_p(x, t)\right|^2 \,  \frac{q_o}{\bar q_-} 
    \left[ \Lambda(p)+ \Lambda^{-1}(p)\right]
    \left[\overline{\Lambda(p) - \Lambda^{-1}(p)}\right]  
    \bigg\}.
}
Actually, letting
\eee{\label{chi-psi-def}
\chi(x,t) := -2\text{Im}\big[\vartheta(x, t, p)\big] + \ln\left|\frac{C_p}{a'(p)}\right|,
\quad
\psi(x,t) := 2\text{Re}\big[\vartheta(x, t, p)\big] + \arg\left(\frac{C_p}{a'(p)}\right)
}
allows us to express the quantity $\varrho_p$ defined by \eqref{e:soliton-rhop} as 
\eee{
\varrho_p(x,t) = e^{\chi(x,t)+i \psi(x,t)}.
}
In turn,  formula \eqref{e:soliton} takes the more compact form
\eee{\label{e:soliton1}
    q(x, t) 
    =
    q_-  + 
    \frac{e^{\chi(x, t)}\left( \bar{\mathcal A} \Lambda_1^2\bar q_- + \mathcal A \Lambda_2^2q_-  - 2 \mathcal B\Lambda_1 \Lambda_2 q_o\right) + e^{i\psi(x, t)} \Lambda_1^2\bar q_-  +  e^{-i\psi(x, t)} \Lambda_2^2q_-}
    {4i\bar{q}_-\left[\sqrt{|\mathcal A|^2-\mathcal B^2}\cosh\!\big(\chi(x, t) + \ln\sqrt{|\mathcal A|^2-\mathcal B^2}\, \big) + \text{Re}\left(\mathcal A e^{i\psi(x, t)}\right)\right]}
}
with the constants $\Lambda_1$ and $\Lambda_2$ given by \eqref{ds-DE-def}. 
It now becomes evident that the pure one-soliton is localized along the line $\chi(x, t)+\ln\sqrt{|\mathcal A|^2-\mathcal B^2} = 0$ which in view of \eqref{chi-psi-def} is equivalent to
\eee{\label{e:soliton-location}
\text{Im}\big[\vartheta(x, t, p)\big] = \frac 12\left( \ln\left|\frac{C_p}{a'(p)}\right|+ \  \ln\sqrt{|\mathcal A|^2-\mathcal B^2} \right).
}

Expressing the pure one-soliton in the non-standard form \eqref{e:soliton1} allows us to compare it against the leading-order asymptotics \eqref{ds-esc-q-sol-lim-pw1-t}, since both expressions involve two portions: the background (first term)  and a traveling wave part (second term). To perform this comparison, we calculate the long-time asymptotics of \eqref{e:soliton1}.
The pure one-soliton propagates along the line specified by equation \eqref{e:soliton-location}. Noting that $\vartheta(x, t, k) = \theta(\xi, k) t$ (cf. definitions \eqref{ds-vtheta-def} and \eqref{ds-theta-def}), we infer that a necessary condition for  equation \eqref{e:soliton-location} to hold in the limit $t\to \infty$  is that $\Im[\theta(\xi,p)] = 0$.
This last equation, however, amounts to $\xi=\Vs$ (recall  \eqref{ds-xistar-xip-def}-\eqref{vs_def}). Thus, we consider three cases: $\xi < \Vs$ (left of soliton); $\xi > \Vs$ (right of soliton); $\xi = \Vs$.

If  $\xi < \Vs$ then $\Im[\theta(\xi,p)] > 0$ (recall Figure \ref{sign-structure-f}). Hence, as $t\to\infty$ we have $\chi(x,t) \to -\infty$ and, in turn, $e^{\chi(x,t)} \to 0$ and $\cosh\!\big(\chi(x,t) + \ln\sqrt{|\mathcal A|^2-\mathcal B^2}\, \big) \to \infty$. Therefore, for $\xi<\Vs$ we obtain
\eee{
q(x,t) = q_- + o(1),\quad t\to\infty,
}
in agreement with the asymptotics \eqref{ds-qsol-pw-t} apart from the real constant phase $g_\infty(\xi)$, which originates from the radiation of Riemann-Hilbert problem \eqref{ds-rhp-cpam}, \eqref{ds-res-cond} and which vanishes once the reflection coefficient in this problem is set to zero.

If $\xi > \Vs$ then $\Im[\theta(\xi,p)] < 0$. Hence, as $t\to\infty$ we have $\chi(x,t) \to\infty$ and, therefore, $e^{\chi(x,t)} \to \infty$ and $\cosh\!\big(\chi(x,t) + \ln\sqrt{|\mathcal A|^2-\mathcal B^2}\,\big) \to \infty$. Thus, expressing $\cosh$ in exponential form we obtain
\eee{\label{ds-qsol-right-as}
    q(x,t) = q_-  + 
    \frac{\bar{\mathcal A} \Lambda_1^2\bar q_- + \mathcal A \Lambda_2^2q_-  - 2 \mathcal B\Lambda_1 \Lambda_2 q_o}
    {2i\bar{q}_-\left(|\mathcal A|^2-\mathcal B^2\right)} + o(1),\quad t\to\infty.
}
Through algebraic manipulations, it can be shown that the leading-order term of \eqref{ds-qsol-right-as} is equal to $q_+$, which is consistent with the fact that propagation along speeds $\xi>\Vs$ always remains to the right of the soliton at $\xi=\Vs$.

Finally, if $\xi = \Vs$ then $\Im[\theta(\xi,p)] = 0$ and hence $\Im[\vartheta(x, t, p)] = 0$, which implies
$$
\chi(\Vs t, t) = \ln\left|\mathcal R_p\right|,
\quad
\psi(\Vs t, t) =  2\theta(\Vs, p) t + \arg\left(\mathcal R_p\right),
\quad
\mathcal R_p := \frac{C_p}{a'(p)}.
$$
Then,  \eqref{e:soliton1} becomes
\ddd{\label{e:soliton1-vs}
&q(\Vs t, t)  
\\
&= q_-  + \frac{\left|\mathcal R_p\right|\left(\bar{\mathcal A} \Lambda_1^2\bar q_- + \mathcal A \Lambda_2^2q_-  - 2 \mathcal B\Lambda_1 \Lambda_2 q_o\right) +  e^{i\left[2 \theta(\Vs,p) t + \arg(\mathcal R_p)\right]} \Lambda_1^2\bar q_- +  e^{-i\left[2 \theta(\Vs,p) t + \arg(\mathcal R_p)\right]} \Lambda_2^2q_-}
{4i\bar q_-\left\{\sqrt{|\mathcal A|^2-\mathcal B^2}\cosh\big[\ln\big(|\mathcal R_p| \sqrt{|\mathcal A|^2-\mathcal B^2}\,\big)\big] + \Re \big(\mathcal A e^{i\left[2 \theta(\Vs,p) t + \arg(\mathcal R_p)\right]}\big)\right\}}.
\nn
}
The exact one-soliton solution \eqref{e:soliton1-vs} is the same with the leading-order asymptotics \eqref{ds-esc-q-sol-lim-pw1-t} except for three points:
(i) the background, which is  $q_-$ in \eqref{e:soliton1-vs} and   $q_-e^{2ig_\infty(\Vs)}$ in \eqref{ds-esc-q-sol-lim-pw1-t}; 
(ii)~an overall phase of $e^{2ig_\infty(\Vs)}$, which is present in \eqref{ds-esc-q-sol-lim-pw1-t} but not in \eqref{e:soliton1-vs};
(iii) the quantity $\mathcal R_p$ in \eqref{e:soliton1-vs}, which is replaced by $R_p$ in \eqref{ds-esc-q-sol-lim-pw1-t}, where
$R_p = \mathcal R_p\delta^2(\Vs,p)e^{-2ig(\Vs,p)}$.
However, this variation is exclusively due to the presence of radiation in  IVP \eqref{ds-fnls-ivp}. Indeed, setting the reflection coefficient equal to zero in the definitions \eqref{ds-delta-def}, \eqref{ds-esc-g-def-pw1} and \eqref{ds-esc-ginf-def-pw1} of $\delta$, $g$ and $g_\infty$ yields $\delta(\Vs, p) = g(\Vs, p) = g_\infty(\Vs)=0$, i.e. in the absence of radiation \eqref{ds-esc-q-sol-lim-pw1-t} would be identical to \eqref{e:soliton1-vs}.

%
%
%
\let\reftitle=\textit
\def\doibase{https://doi.org/}

\vspace{7mm}
\noindent
Gino Biondini \hspace{3.2cm} Sitai Li  \hspace{4.42cm}  Dionyssios Mantzavinos
\\
Department of Mathematics \hspace{0.8cm} Department of Mathematics \hfill Department of Mathematics
\\
State University of New York \hspace{0.6cm} University of Michigan \hspace{1.85cm}  University of Kansas
\\
Buffalo, NY 14260 \hspace{2.46cm} Ann Arbor, MI 48109 \hspace{2.05cm} Lawrence, KS 66045 \\
E-mail: \textit{biondini@buffalo.edu} \hspace{0.82cm} E-mail: \textit{sitaili@umich.edu} \hspace{1.44cm}  E-mail: \textit{mantzavinos@ku.edu} 

\end{document}